\documentclass[review]{article}

\usepackage[text={6.75in,9.5in},centering,letterpaper]{geometry}
\usepackage{amssymb,amsmath}
\usepackage{graphicx}
\usepackage{subcaption}

\newtheorem{theorem}{Theorem}[section]
\newtheorem{lemma}[theorem]{Lemma}

\newtheorem{corollary}[theorem]{Corollary}
\newtheorem{remark}[theorem]{Remark}

\catcode`\@=11
\@addtoreset{equation}{section}
\renewcommand{\theequation}{\thesection.\arabic{equation}}

\DeclareMathOperator{\sech}{sech}

\title{Multi-front dynamics in spatially inhomogeneous Allen-Cahn equations\thanks{\today}}

\author{Robbin Bastiaansen\thanks{Mathematical Institute, Utrecht University, 3508 TC, Utrecht, The Netherlands \& Department of Physics, Institute for Marine and Atmospheric Research Utrecht, Utrecht University,
Utrecht, The Netherlands ({\tt r.bastiaansen@uu.nl})}
\and Arjen Doelman\thanks{Mathematical Institute, Leiden University, 2300 RA, Leiden, Netherlands ({\tt doelman@math.leidenuniv.nl})}
\and Tasso J. Kaper\thanks{Department of Mathematics and Statistics, Boston University, Boston, MA 02215, USA
({\tt tasso@math.bu.edu})} }

\begin{document}

\maketitle

\begin{abstract}
Recent studies of biological, chemical, and physical pattern-forming systems have started to go beyond the classic `near onset' and `far from equilibrium' theories for homogeneous systems to include the effects of spatial heterogeneities. In this article, we build a conceptual understanding of the impact of spatial heterogeneities on the pattern dynamics of reaction-diffusion models.
We consider the simplest setting of an explicit, scalar, bi-stable Allen-Cahn equation driven by a general small-amplitude spatially-heterogeneous term $\varepsilon F(U,U_x,x)$. In the first part, we perform an analysis of the existence and stability of stationary one-, two- and $N$-front patterns for general spatial heterogeneity $F(U,U_x,x)$. In addition, we explicitly determine the $N$-th order system of ODEs that governs the evolution of the front positions of general $N$-front patterns to leading order. In the second part, we focus on a particular class of spatial heterogeneities where $F(U,U_x,x) = H'(x) U_x + H''(x) U$ with $H$ either spatially periodic or localised. For spatially periodic heterogeneities, we show that the fronts of a multi-front pattern will get `pinned' if the distances between successive fronts are sufficiently large, {\it i.e.}, the multi-front pattern is attracted to a nearby stable stationary multi-front pattern. For localised heterogeneities, we determine all stationary $N$-front patterns, and show that these are unstable for $N > 1$. We find instead slowly evolving `trains' of $N$-fronts that collectively travel to $\pm \infty$, either with slowly decreasing or increasing speeds.
\end{abstract}
\vspace{0.2cm}
\\ \\
{\bf Highlights}
\begin{itemize}
\item Spatial heterogeneities can prevent front  annihilation and coarsening in the Allen-Cahn equation
\item Spatial heterogeneities can create rich new multi-front patterns in the Allen-Cahn equation
\item Spatially periodic heterogeneities lead to front pinning in multi-front patterns
\item Multi-fronts form in clusters through intersections of invariant manifolds, and one per cluster is stable
\item Localized spatial forcing creates dynamically-evolving multi-front patterns
\end{itemize}

\section{Introduction}
\begin{figure}
\centering
\begin{subfigure}[t]{0.32\textwidth}
\centering
\includegraphics[width=\textwidth]{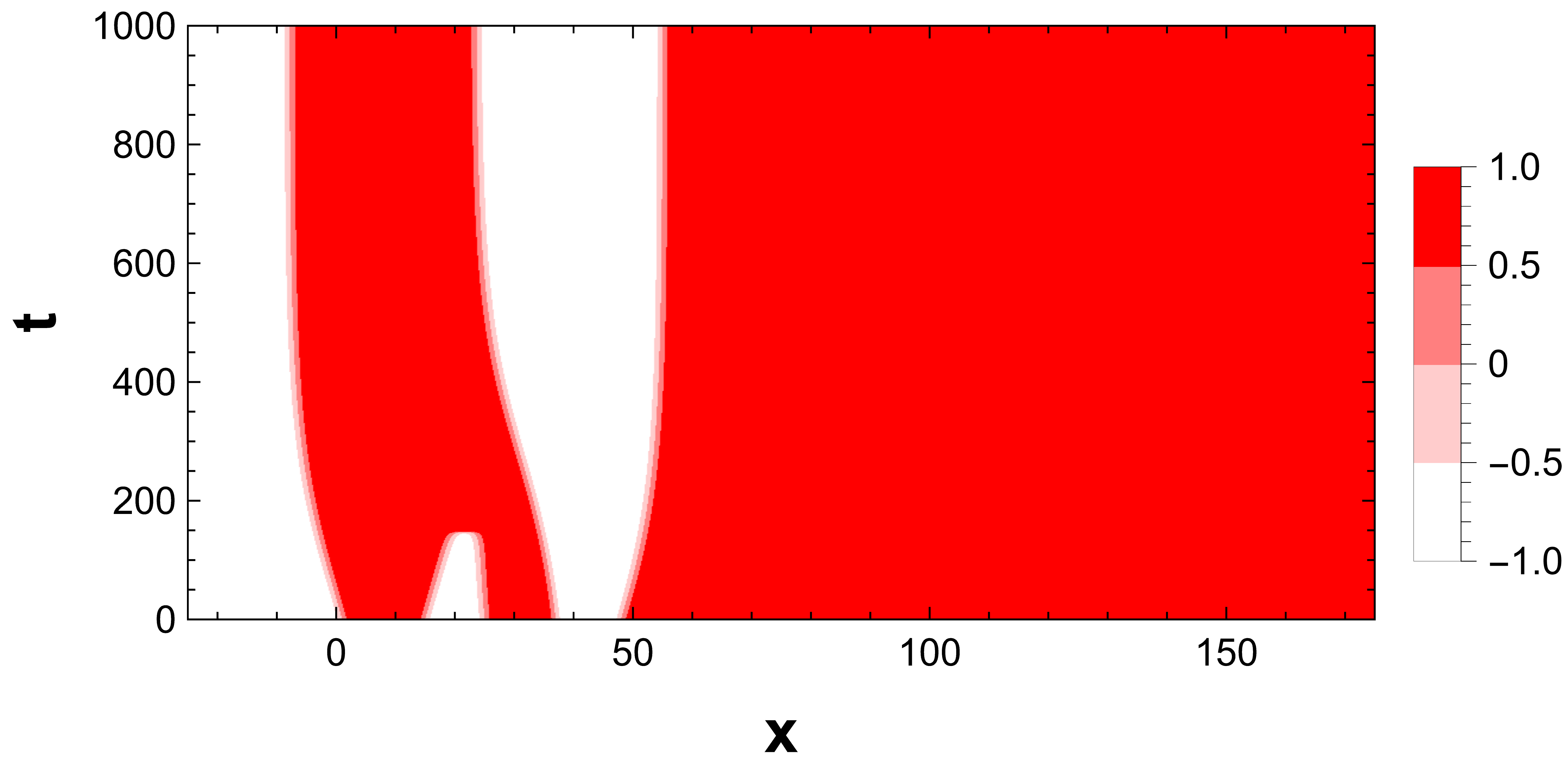}
\end{subfigure}
~
\begin{subfigure}[t]{0.32\textwidth}
\centering
\includegraphics[width = \textwidth]{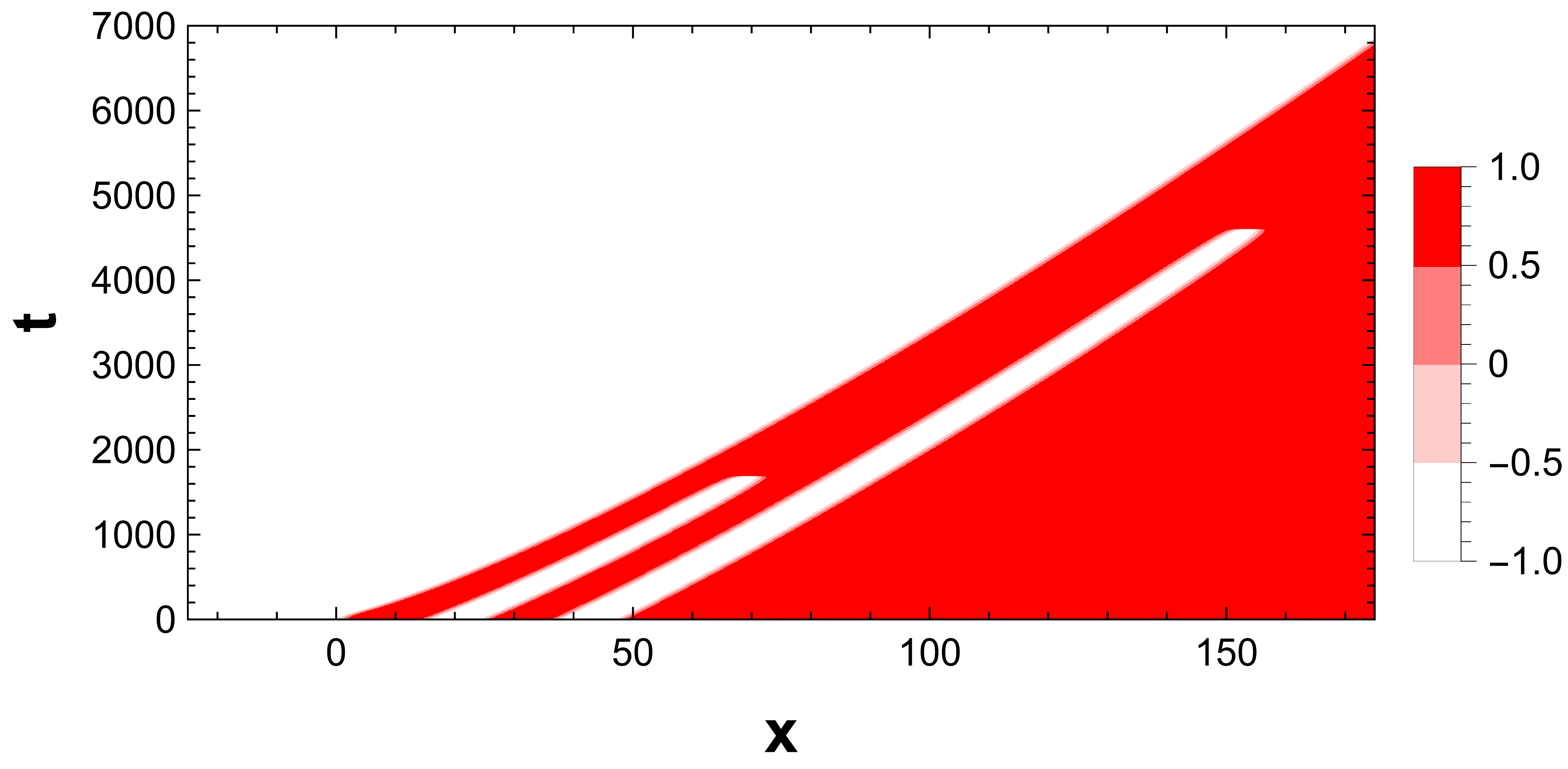}
\end{subfigure}
~
\begin{subfigure}[t]{0.32\textwidth}
\centering
\includegraphics[width = \textwidth]{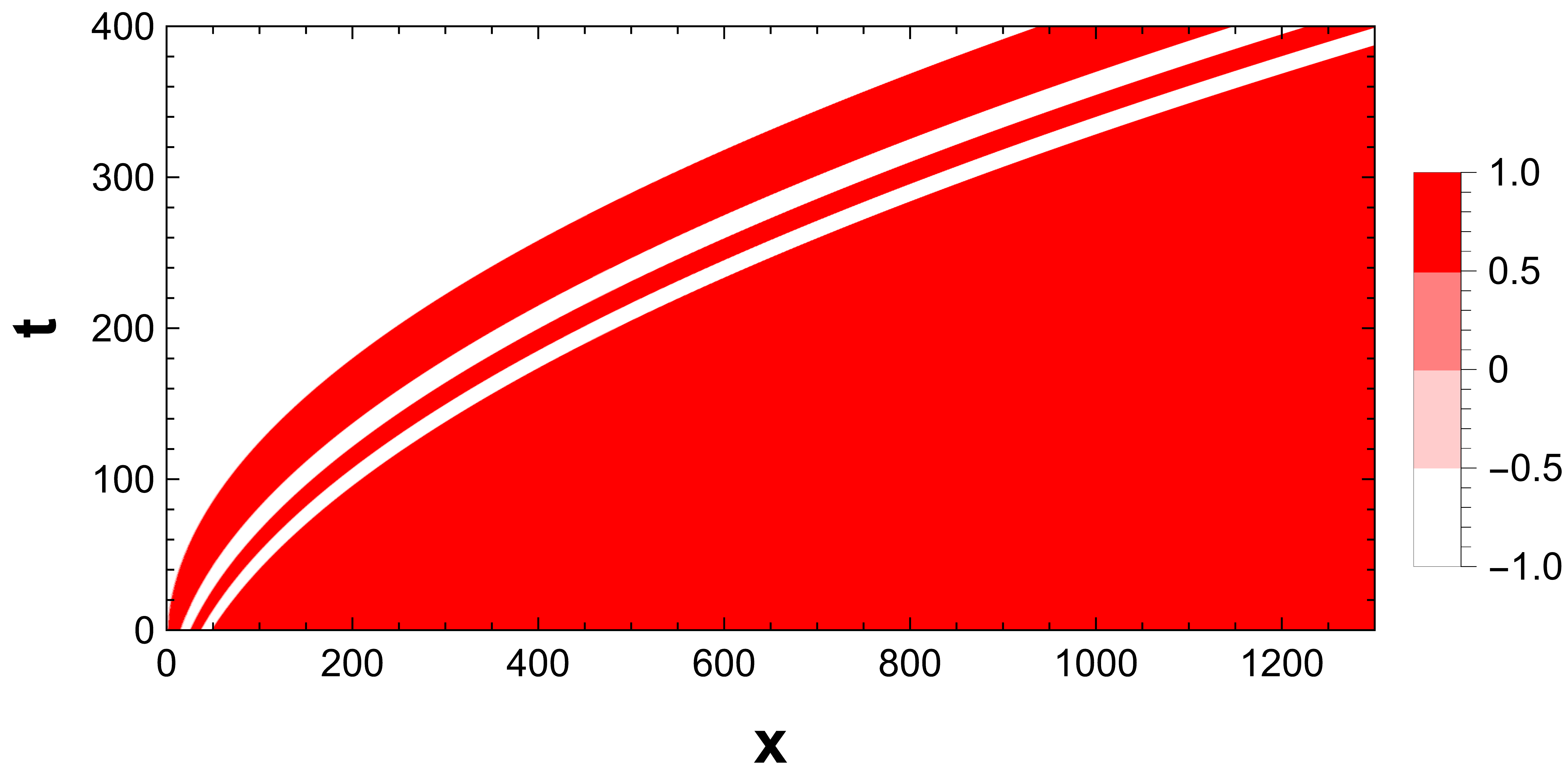}
\end{subfigure}
\caption{Simulations of \eqref{eq:mainEquation} with five-front initial condition $U(x,0) = \tanh(x - 1)  - \tanh(x - 15) + \tanh(x - 25) - \tanh(x - 37) + \tanh(x - 48)$ and  $\varepsilon = 0.1$ on the interval $(-500,1500)$ (plotted on relevant subintervals) with homogeneous Neumann boundary conditions (a) A spatially periodic inhomogeneity \eqref{eq:canonicalExample} with $f_1(x) = \sin(\frac15 x)$, $f_2(x) \equiv f_3(x) \equiv 0$. (b) An algebraically decaying spatially localized topographical inhomogeneity $H(x) = -H_{\rm alg}(x; p)$ with $0 < p = 0.25$ \eqref{eq:Ftopography}, \eqref{eq:defHuni-alg}. (c) The inhomogeneity $H(x) = -H_{\rm alg}(x; p)$ with $p = -0.50 \in (-1,0)$.}
\label{fig:5Fronts-Intro}
\end{figure}

Traditionally, the mathematical analysis of the formation and dynamics of patterns, {\it i.e.}, the dynamics of solutions of nonlinear systems of reaction-diffusion equations, has focused on systems of equations that do not explicitly vary in space (or time).
For the theory of the onset of near-homogeneous patterns, see for instance \cite{mielke2002ginzburg} and the references therein.
For the theory of localized, `far-from-equilibrium' patterns, see \cite{nishiura-book,sandstede2002stability, wei2013mathematical} and the references therein.
\\ \\
In recent years, there has been a growing interest in understanding pattern dynamics in biological, chemical, and physical problems which are spatially heterogeneous.
Patterns in biology are affected by the spatial variations of various rate constants, including for example in embryology where reaction or diffusion rates of various morphogen proteins may be modulated by preexisting spatial structure within the developing organism \cite{maini2003}.
Patterns in ecological models are influenced by the topography of the terrain, the soil inclination, the local soil quality, and the availability of nutrients; and these have all been shown to affect the appearance, magnitude, movement, and growth (or decline) of the patterns exhibited in models of ecosystems   \cite{bastiaansen2020pulse,bastiaansen2019dynamics,bastiaansen2018multistability,dunkerley2013vegetation,patterson2023spatial}, including also that non-trivial spatial patterns may increase the resilience of an ecosystem \cite{rietkerk2021evasion}.
Similarly, patterns in chemical reactions can  depend on the spatially heterogeneity of  reactors, spatially-dependent rate constants, or applied light sources \cite{VE2003}.
Some other recent references to spatial heterogeneities in reaction-diffusion models include \cite{JaramilloScheel2014,vangorder2021}.
Therefore, from a broad scientific point of view, a general approach needs to be developed to understand -- and preferably also to predict -- the impact of spatial inhomogeneities on the dynamics of localized structures in multi-component reaction-diffusion systems defined on large domains.
\\ \\
In this article, we study the impact of spatial heterogeneities in one of the simplest underlying pattern generating models: the scalar, bistable Allen-Cahn model in one space dimension.
In its standard form, this prototypical PDE is given by
\begin{equation}
\label{eq:standardAC}
	U_t = U_{xx} + U - U^3.
\end{equation}
We analyze the impact of small-amplitude, spatially heterogeneous terms on the dynamics exhibited by the Allen-Cahn equation, {\it i.e.}, we study the dynamics of patterns in
\begin{equation}
	U_t = U_{xx} + U - U^3 + \varepsilon F(U,U_x,x).
\label{eq:mainEquation}
\end{equation}
{\it A priori,} we impose the natural assumptions that $F(U,V,x)$ is continuous in all three arguments and that $\|x \mapsto F(U,V,x)\|_{C_b} \leq \mathcal{O}(1)$ for all $U, V \in \mathbb{R}$, though in the Discussion Section (Sec.~\ref{sec:Discussion}) we will also consider cases in which the boundedness assumption is violated.
\\ \\
The analysis in the article will focus largely on the case of general bounded heterogeneities $F$.
In addition, we provide concrete examples of the main analytical formulas -- and carry out numerical simulations -- for a variety of different explicit forms of the driving function $F$, including spatially-periodic and localized heterogeneities.
One particularly simple choice for $F(U,V,x)$ is
\begin{equation}
	F(U,V,x) = f_1(x) U + f_2(x) V + f_3(x).
\label{eq:canonicalExample}
\end{equation}
For example, in ecological models, $F(U,V,x)$ may represent the effect of a topography $H(x) \in C^2(\mathbb{R})$ with $f_2(x) = H'(x)$, $f_1(x) = f_2'(x)$, and $f_3(x) \equiv 0$, so that
\begin{equation}
F(U,V,x) =  H'(x) V + H''(x)U
\label{eq:Ftopography}
\end{equation}
(see \cite{bastiaansen2019dynamics} for a derivation and note that $H'(x) \equiv $ a constant in the classical Klausmeier model \cite{klausmeier1999regular}). This is the `topographical' case that we will also consider as one of the main examples throughout the article.
\\
\begin{figure}
\centering
	\begin{subfigure}[t]{0.325\textwidth}
		\centering
		\includegraphics[width=\linewidth]{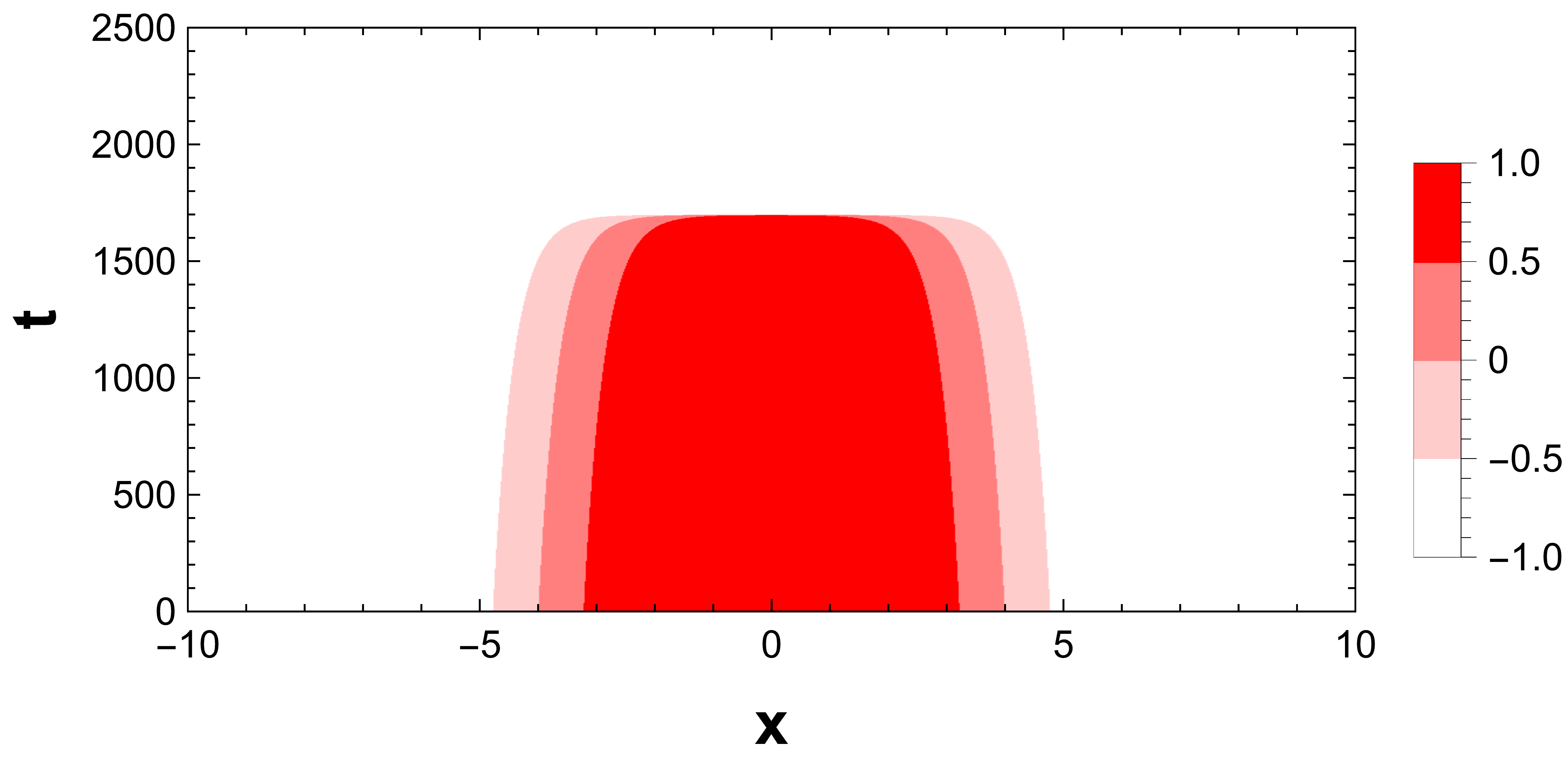}
		\caption{Homogeneous: $H(x) \equiv 0$.}
	\end{subfigure}
        \begin{subfigure}[t]{0.325\textwidth}
		\centering
		\includegraphics[width=\linewidth]{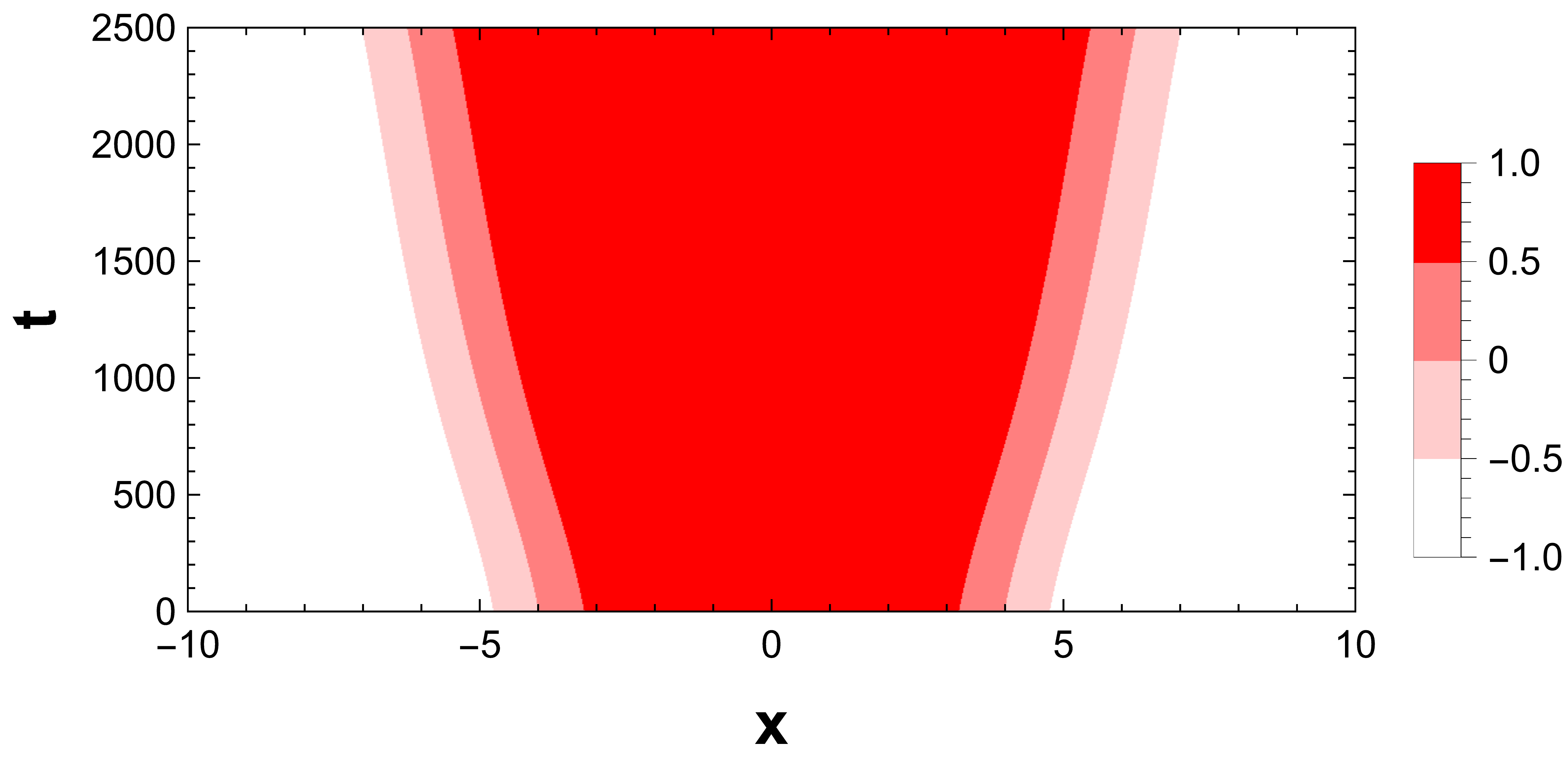}
        \caption{Localized: $H(x)=H_{\rm exp}(x; 1)$ \eqref{eq:defHuni-alg}.}
	\end{subfigure}
	\begin{subfigure}[t]{0.325\textwidth}
		\centering
		\includegraphics[width=\linewidth]{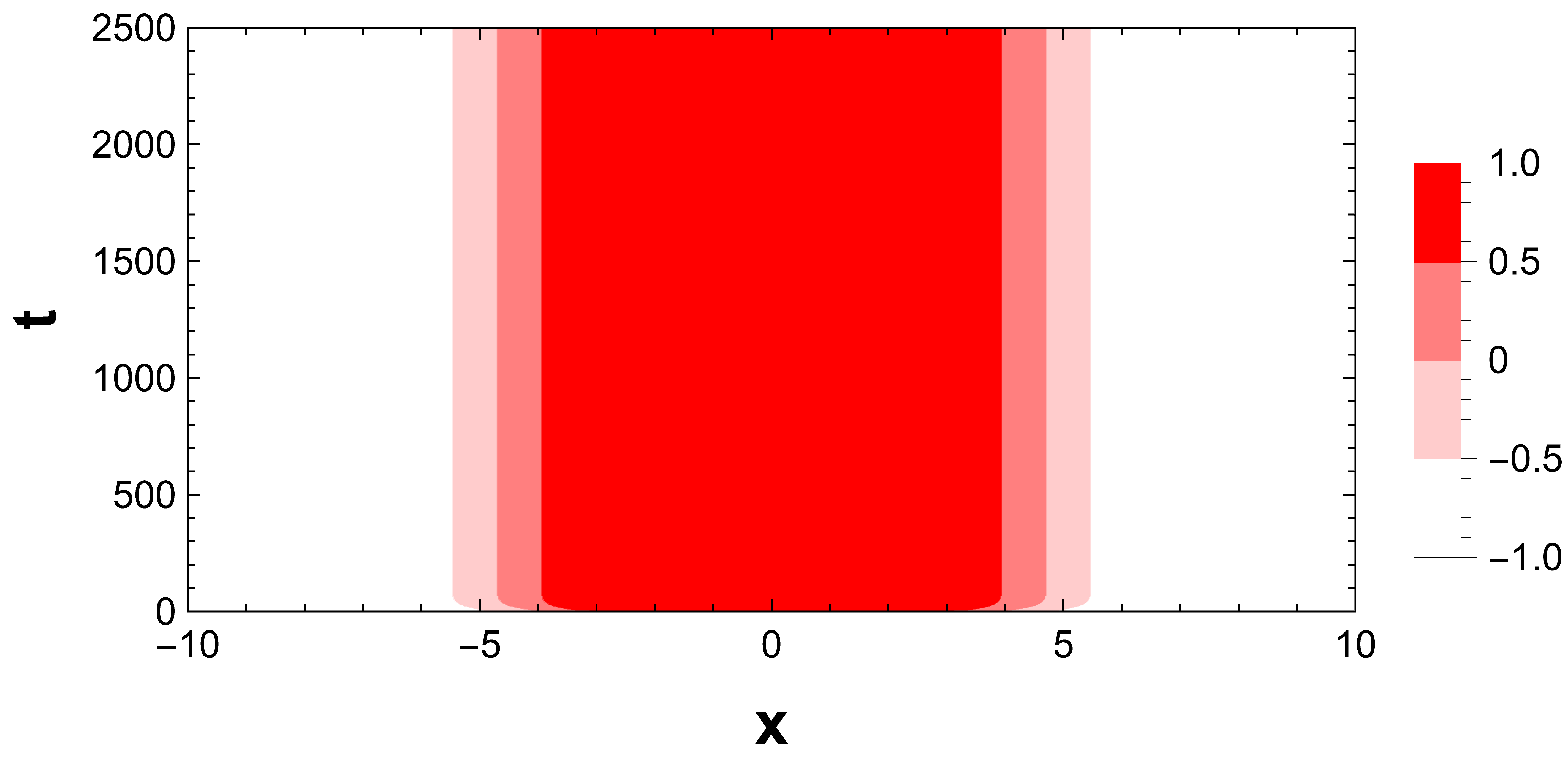}
	\caption{Periodic: $H(x) = \sin 2x$.}
	\end{subfigure}
	\begin{subfigure}[t]{0.325\textwidth}
		\centering
		\includegraphics[width=\linewidth]{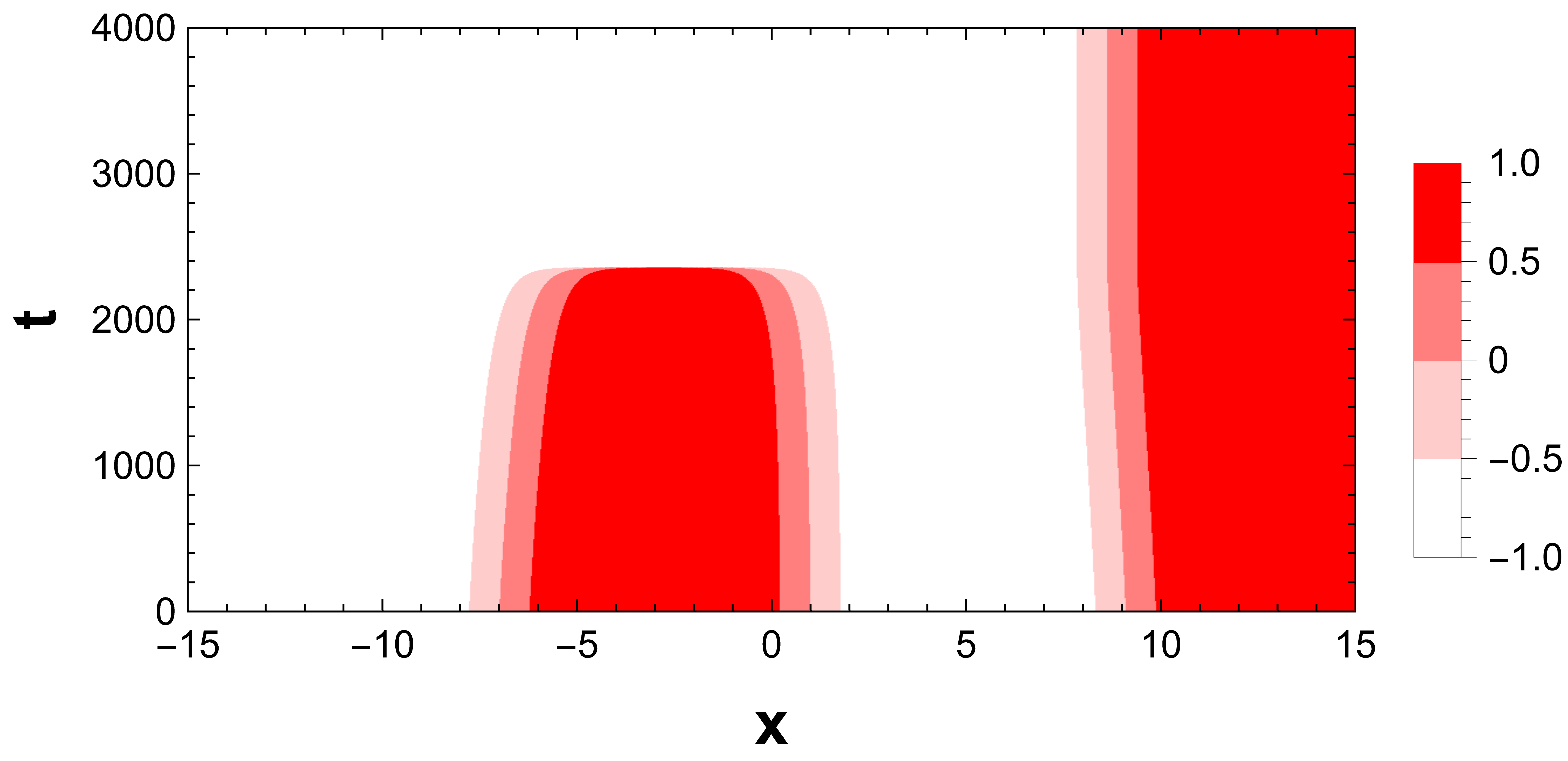}
	\caption{Homogeneous: $H(x) \equiv 0$.}
	\end{subfigure}
        \begin{subfigure}[t]{0.325\textwidth}
		\centering
		\includegraphics[width=\linewidth]{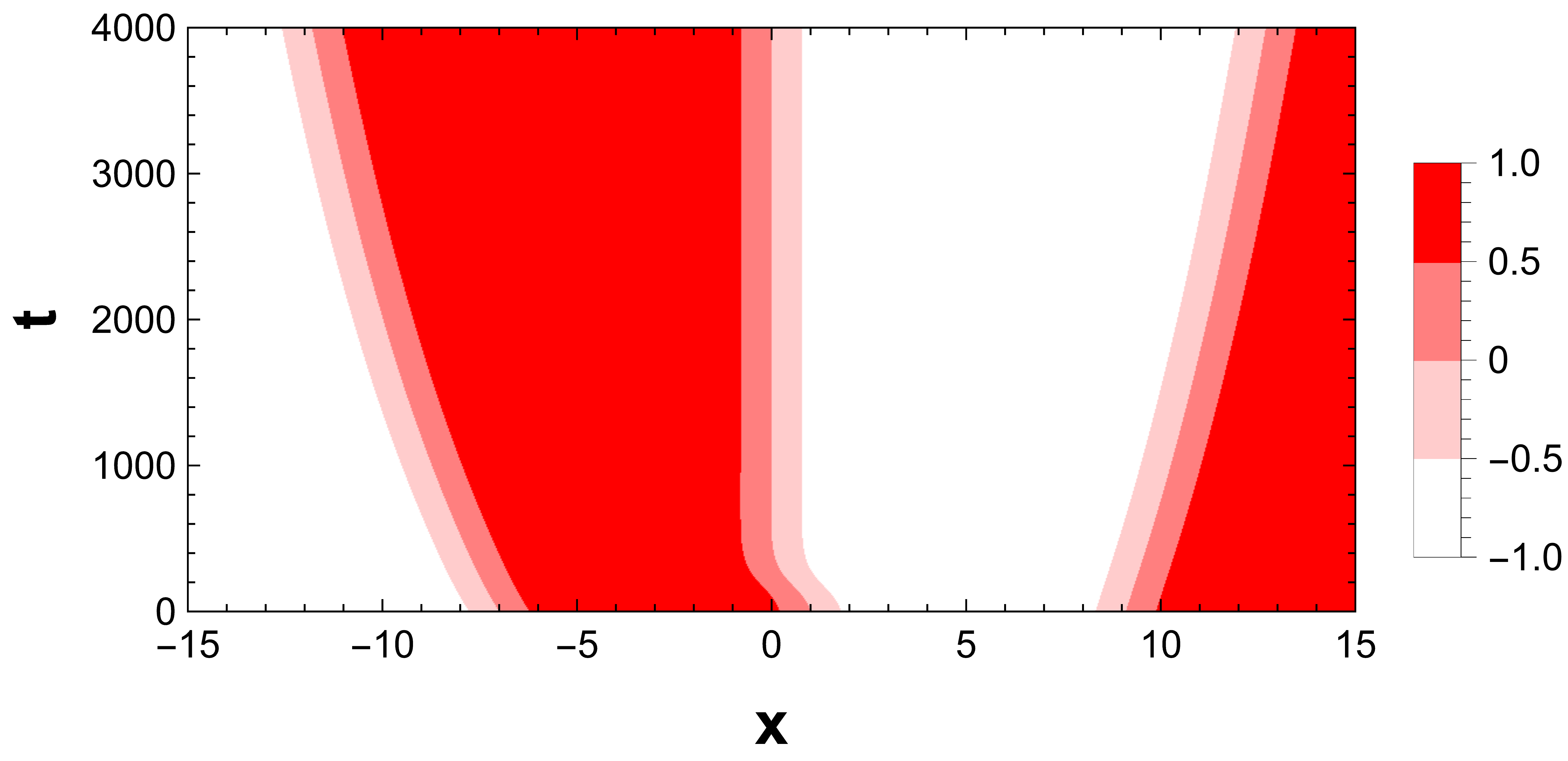}
	\caption{Localized: $H(x)=H_{\rm alg}(x; 2)$ \eqref{eq:defHuni-alg}.}
	\end{subfigure}
	\begin{subfigure}[t]{0.325\textwidth}
		\centering
		\includegraphics[width=\linewidth]{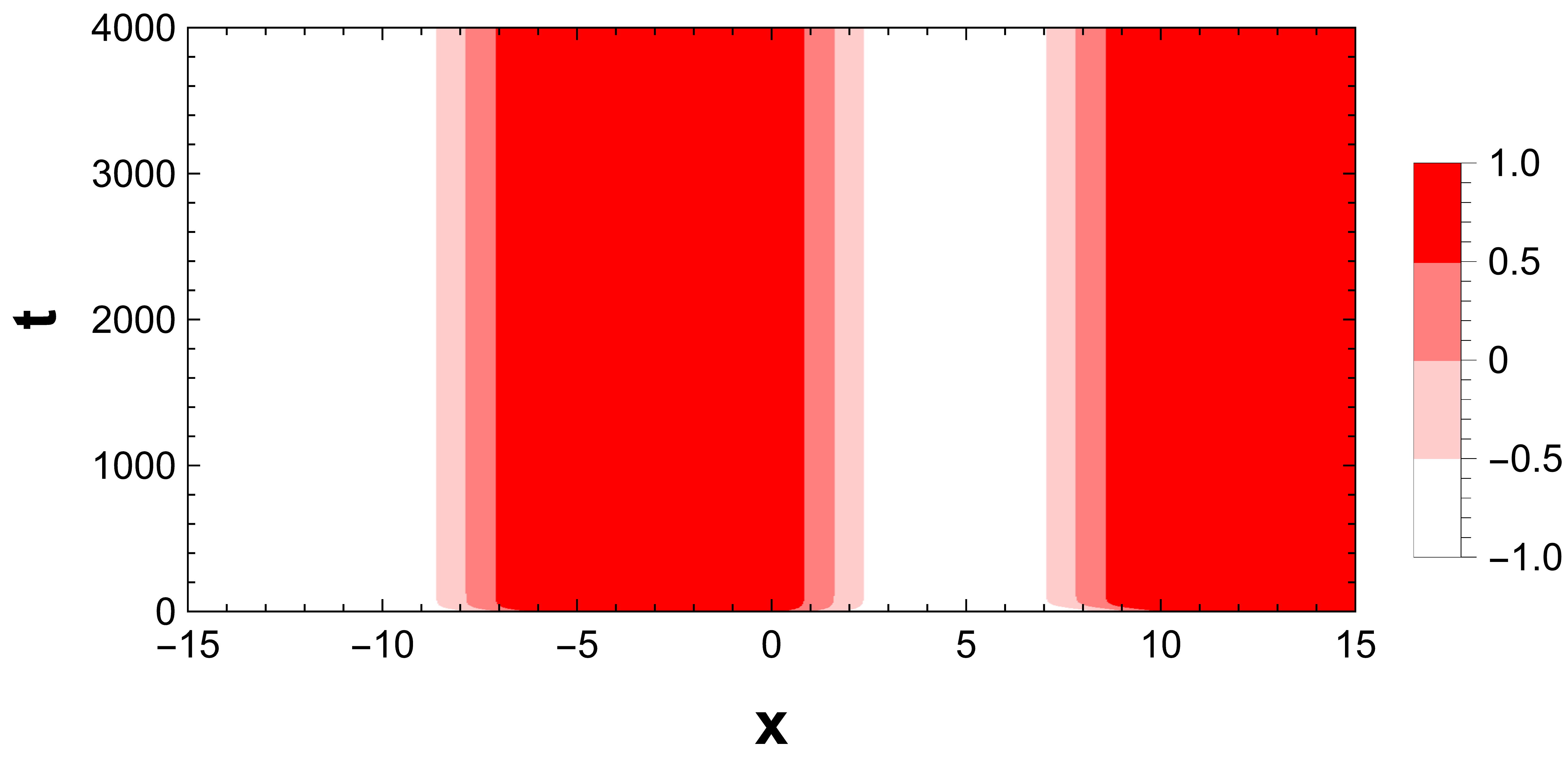}
        \caption{Periodic: $H(x) = \sin 2x$.}
	\end{subfigure}
\caption{Simulations of multi-front dynamics in \eqref{eq:mainEquation} (under homogeneous Neumann boundary conditions) with topographical inhomogeneities $H(x)$ given by \eqref{eq:Ftopography} (with $\varepsilon = 0.1$); (a), (b), (c): starting from the two-front initial condition $U(x,0) = \tanh(x + 4)  - \tanh(x - 4) - 1$ on the $x$-interval $(-10,10)$; (d), (e), (f): with three-front initial condition $U(x,0) = \tanh(x + 7)  - \tanh(x - 1) + \tanh(x - 9.1)$ on $(-25,25)$, plotted on $(-15,15)$.
The color indicates the magnitude of $u$.}
\label{fig:NumericsTwoThreeFronts-Intro}
\end{figure}
\\
To provide a first, brief quantification of the effects of small heterogeneous terms in \eqref{eq:mainEquation} on the front dynamics in the Allen-Cahn equation, {\it i.e.,} in order to fully motivate the main quantitative analysis of this article, we show in Fig. \ref{fig:NumericsTwoThreeFronts-Intro} the outcomes of some representative multi-front simulations.
In Figs. \ref{fig:NumericsTwoThreeFronts-Intro}(a-c), the evolution in time of \eqref{eq:mainEquation} from identical two-front initial data is shown for three kinds of inhomogeneities $F(U,U_x,x)$. Then, in Figs. \ref{fig:NumericsTwoThreeFronts-Intro}(d-f), the same is done for three-front initial data (again identical in all three cases). Frames (a) and (d) illustrate `classical' two-front and three-front dynamics in the homogeneous Allen-Cahn equation \eqref{eq:standardAC}.
Namely, two fronts collide, and the systems `coarsen' to either the trivial state $U \equiv -1$ (Fig. \ref{fig:NumericsTwoThreeFronts-Intro}(a)) or to a stationary one-front (Fig. \ref{fig:NumericsTwoThreeFronts-Intro}(d)),
as expected based on the classical literature applicable to equation \eqref{eq:standardAC} \cite{Carr1989,Chen2004,Fusco1989} (see also \cite{westdickenberg2021metastability} for a recent account).
The significant effects of spatially localized heterogeneities are illustrated in frames (b) and (e).
Instead of being attracted to each other (and eventually annihilating each other), the two fronts in Fig. \ref{fig:NumericsTwoThreeFronts-Intro}(b) are repelled by each other due to the localized heterogeneity.
Indeed, the final pattern consists of a dynamically evolving two-front state, with the fronts continuing to move outward, though with decreasing speed, and with the solution being  asymptotically close to $U=+1$ on the (growing) middle part between the fronts.
Thus, the attracting $U\equiv -1$ state of the homogeneous case (Fig. \ref{fig:NumericsTwoThreeFronts-Intro}(a)) is pushed further and further outward, to $-\infty$ and $+\infty$, and more and more of the domain
is covered by the other state.
In the three-front case, the localized spatial inhomogeneity also pushes the two outer fronts to $\mp \infty$, but here the middle front is attracted to a specific position at which it settles.
Hence, in both simulations, the final spatial patterns differ completely from those of the homogeneous cases,
and interesting new multi-front dynamics are created by the localized heterogeneities.
Next, Figs.~\ref{fig:NumericsTwoThreeFronts-Intro}(c) and (f) show that spatially periodic inhomogeneities can also have significant effects on two-front and three-front solutions.
Here, the fronts each get `pinned' at specific locations.
Hence, the spatially periodic heterogeneity also prevents the coarsening (seen in the homogeneous case) from happening, and it can generate new multi-front patterns with stationary, pinned fronts.
\\
\begin{figure}
\centering
        \begin{subfigure}[t]{0.4\textwidth}
		\centering
		\includegraphics[width=\linewidth]{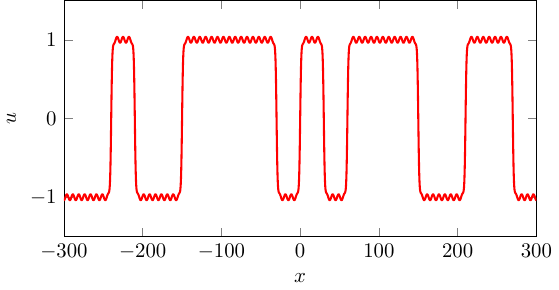}
	\end{subfigure}
	\begin{subfigure}[t]{0.4\textwidth}
		\centering
		\includegraphics[width=\linewidth]{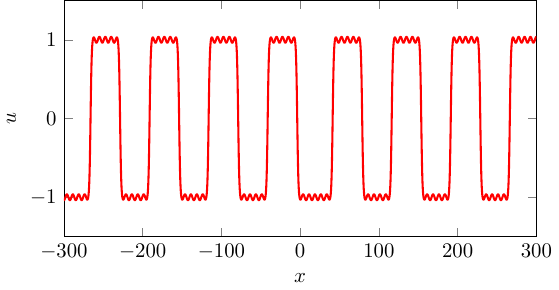}
	\end{subfigure}
	\begin{subfigure}[t]{0.4\textwidth}
		\centering
		\includegraphics[width=\linewidth]{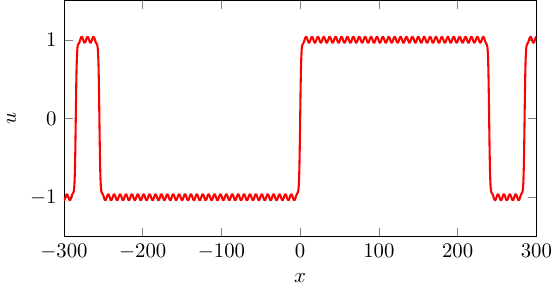}
	\end{subfigure}	
    \begin{subfigure}[t]{0.4\textwidth}
		\centering
		\includegraphics[width=\linewidth]{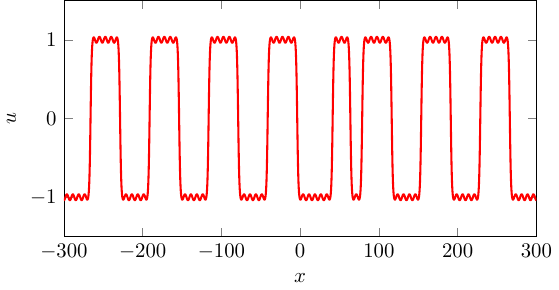}
	\end{subfigure}
	\caption{Examples of numerically attracting (and thus numerically stable) stationary multi-front patterns in the forced Allen-Cahn equation \eqref{eq:mainEquation} with spatially periodic inhomogeneous term $F(U,V,x)$ given by \eqref{eq:canonicalExample} with $f_1(x) = \cos(\frac{4}{15}\pi x)$, $f_2(x) \equiv f_3(x) \equiv 0$, and $\varepsilon = 0.1$. (Simulations done for $x \in (-300,300)$ with homogeneous Neumann boundary conditions; patterns are (numerically) stable based on simulation up to $T = 10^{12}$.)}
\label{fig:NumericsMultiFront-Intro}
\end{figure}
\\
The phenomenon of pinning (of fronts or pulses) is perhaps the most studied impact of spatial inhomogeneities on the dynamics of `far from equilibrium' localized structures -- see  \cite{bastiaansen2022fragmented, champneys2021bistability, ding2015transition, scheel2017depinning} and the references therein, where we add that typically only solitary fronts or pulses are considered in the mathematical literature.
Moreover, in the case of spatially periodic inhomogeneities, the dynamics of $N$-front patterns is typically driven by an interplay between coarsening and pinning.
Fronts that are `too close' attract and annihilate each other, while fronts that are `sufficiently far' apart get pinned, see the simulation shown in Fig. \ref{fig:5Fronts-Intro}(a) in which a $5$-front pattern eventually ends up in a pinned (and thus stationary) three-front state and see the pinned multi-front patterns of Fig. \ref{fig:NumericsMultiFront-Intro}. The situation is very different in the case of spatially localized inhomogeneities.
For these, at most one front can get pinned, see again Fig. \ref{fig:NumericsTwoThreeFronts-Intro}(e), and all other fronts necessarily start to travel.
However, the speed of fronts that move in the same direction depends on their distance from (the center of) the spatial inhomogeneity.
Fig.~\ref{fig:5Fronts-Intro}(b) shows the typical situation in which fronts catch up with those traveling directly ahead of them, which yields a stepwise process of successive annihilations, until there is at most one solitary traveling front. As we shall analyze, this is the expected behavior of multi-front patterns driven by localized inhomogeneities in the case in which the inhomogeneity is bounded as a function of $x$ ({\it i.e.}, if there is an $C=C(U,U_x) > 0$ such that $|F(U,U_x,x)| < C$ for all $x \in  \mathbb{R}$).
Finally, we observe that, although the analysis here requires that $F(U,U_x,x)$ is bounded, we can formally extrapolate beyond this and predict that it is possible that fronts may not catch up with their predecessors so that no mutual annihilation takes place, see Fig.~\ref{fig:5Fronts-Intro}(c).
\\
\\
The dynamics induced by various spatially heterogeneous terms, including those shown in Figs.~\ref{fig:NumericsTwoThreeFronts-Intro} and \ref{fig:5Fronts-Intro}, motivate the following fundamental questions about multi-front patterns: {\it What kind of novel (stable) multiple-front patterns may be initiated by the heterogeneous effects?}
{\it How do the different quantitative properties of the spatially heterogeneous terms affect the dynamics and stability of the interacting fronts in multi-front patterns?}
We address these analytically for the spatially heterogeneous Allen-Cahn model \eqref{eq:mainEquation}.
\\  \\
As a preparatory step for the analysis, we first study the persistence of two elementary (or basic) states of the homogeneous Allen-Cahn equation \eqref{eq:standardAC}: the trivial background states and the $\tanh$-type one-front solutions that connect different background states.
In particular, we obtain explicit results on the persistence of the background states and of the stationary one-fronts under heterogeneous effects.
These are obtained using a straightforward Melnikov analysis.
We also determine the linear (spectral) stability of the background states and one-front solutions of \eqref{eq:mainEquation}, see Theorem \ref{th:1frontExStab}.
\\  \\
Our analysis of the multi-front patterns proper then begins by establishing the existence and stability of stationary two-front solutions.
These are the simplest novel patterns that do not exist in the homogeneous case, see Theorem \ref{th:2frontExStab}.
Then, for $N$-front solutions with general $N$
and general spatially heterogeneous terms $F$, we derive the $N$-dimensional dynamical systems that govern the interactions of the dynamically-evolving fronts to leading order.
This system of $N$ ODEs is given by \eqref{eq:NFrontODE}.
(Here, we follow the approach also used in \cite{ei2002pulse, promislow2002renormalization}.)
\\ \\
The system \eqref{eq:NFrontODE} of ODEs has an insightful and explicitly tractable form for  problems in which the spatial heterogeneity has the `topographical' character of \eqref{eq:Ftopography}.
Therefore, after the general form \eqref{eq:NFrontODE} of the ODE system is derived for general $F$, we turn in Sec.~5 to present the main analysis of the multi-front dynamics for topographical $F$ (where the ODEs are given by \eqref{eq:dynNfronts}). For exponentially decaying topographies, we establish the existence of large families of stationary $N$-front patterns. These patterns correspond to equilibria of the ODEs and we prove that the maximal number is $(K+1)N - 1$, where $K$ is the number of zeroes of an associated Melnikov function. Also, the stability of these $N$-front patterns is determined by the stability type of the corresponding equilibria. See Theorem \ref{th:Nfrontslocexp}. Then, we also establish the existence of $N$-front solutions in the presence of algebraically decaying topographies. See Corollary \ref{cor:Nfrontslocalg}.
However, we show that, for $N > 1$, these $N$-front solutions are unstable.
This general result builds on a detailed analysis of the special cases $N=2,3,4$.
\\ \\
The analysis of these exponentially and algebraically decaying topographies also provides insights into the fate of patterns that start out from $N$-front configurations.
Namely, after an initial process of coarsening during which one front may get pinned at a special location (recall for example Fig. \ref{fig:NumericsTwoThreeFronts-Intro}(e)), `trains' of $M$-front patterns (with $M \leq N$) are formed that collectively travel to $+$ or $-\infty$.
We show that, as long as $H'(x)$ decays, fronts within such a train will eventually catch up with their predecessors, thereby driving a slow coarsening process -- recall also Fig.~\ref{fig:5Fronts-Intro}(b).
\\  \\
Next, we turn in Sec.~6 to investigate the system of $N$ ODEs \eqref{eq:NFrontODE} in the case of spatially periodic heterogeneities.
The dynamics of $N$-front patterns for these spatially-periodic problems are summarized in Theorems~\ref{th:Nfrontsper} and \ref{th:Nfrontsperattr}.
The former establishes the existence of stationary $N$-front patterns for patterns in which neighboring fronts are sufficiently far apart.
In fact, near any initial $N$-front configuration with sufficiently separated successive fronts,  there are $K^N$ (distinct) stationary $N$-fronts of which $(\frac12 K)^N$ are stable, where $K \geq 0$ is the number of zeros per period of the associated periodic Melnikov function.
This theorem also establishes the non-existence of stationary $N$-front patterns when the fronts are too close to each other.
The process by which stationary fronts are created (as the initial distance between fronts increases) is studied in detail in the specific setting of two-fronts.
It is shown how the underlying geometry of this process can be interpreted in terms of the well-established terminology of weakly driven nonlinear oscillators.
Then, Theorem \ref{th:Nfrontsperattr} presents results for dynamic $N$-fronts.
Again, these are all new patterns, which cannot exist in the homogeneous Allen-Cahn equation.
\\ \\
Throughout the article, we also illustrate our findings with numerical simulations of \eqref{eq:mainEquation}.
In the localized setting of \eqref{eq:Ftopography}, these simulations will be performed on specific topographies in which $H(x)$ is a `unimodal hill' (or `valley') with either exponentially or algebraically decaying tails,
\begin{equation}
\label{eq:defHuni-alg}
H_{\rm exp}(x; \mu) = \frac{1}{\cosh^2(\frac12 \sqrt{2} \mu x)}, \; \; \ \ \ \
H_{\rm alg}(x;p) = \frac{1}{(1+x^2)^{\frac12(p-1)}}.
\end{equation}
For a truly localized topography, we need to impose $\mu > 0$ for $H_{\rm exp}(x; \mu)$ and $p > 1$ for $H_{\rm alg}(x;p)$.
However, by (\ref{eq:Ftopography}) we note that only $H'(x)$ enters into equation (\ref{eq:mainEquation}), so that it is natural to consider $p > 0$, or even $p \geq 0$ since also $F(U,V,x)$ remains bounded for $p=0$ (which corresponds to the classical Klausmeier case \cite{klausmeier1999regular}).
In fact, we will find that our analysis may also be valid for $p \in (-1,0)$ for which $H''_{\rm alg}(x;p)$ remains bounded (but $H'_{\rm alg}(x;p)$ not), see Sec.~\ref{sec:Discussion} and Fig.~\ref{fig:5Fronts-Intro}(c). Since the `driving force' of the algebraically decaying tails of $H_{\rm alg}(x; p)$ is much stronger than that of the exponentially decaying tails of $H_{\rm exp}(x; \mu)$, almost all simulations will be performed with $H_{\rm alg}(x; p)$. However, $H_{\rm exp}(x; \mu)$ is more accessible to explicit analysis, therefore it will appear more often in the examples presented at various places in the article.
\\ \\
We conclude the article with an extensive discussion.
First, the impact of the spatial heterogeneities on the process of coarsening is discussed, followed by a brief exploration into the realm of topographies beyond those considered in the earlier sections (namely, beyond the periodic and localized (with monotonically decaying tails) cases).
Then, the main topic of this section is the behavior of $N$-front patterns driven by topographies $H(x)$ for which $H'(x)$ `decays' algebraically with rate $p < 0$ or even $p < 1$ ({\it i.e.}, $H'(x) \sim 1/x^p$ as $x \to \pm \infty$ -- see example $H_{\rm alg}(x;p)$ in \eqref{eq:defHuni-alg}).
For such $H(x)$, the methods of \cite{ei2002pulse, promislow2002renormalization} to settle the validity of the $N$-front interaction system cannot be applied directly, although there are indications that this should be possible for $p$ up to $-1$ ({\it i.e.}, for $p > -1$.
It follows from the interaction equations -- and is confirmed by numerical simulations (Fig. \ref{fig:5Fronts-Intro}(c), Fig. \ref{fig:5FrontDynVaryingP}) -- that for $-1 < p <0$, fronts within an `$M$-train' of jointly traveling fronts will travel slower than their predecessors.
This implies that the $M$-trains will persist and slowly `spread out'.
The situation becomes even more interesting as $p$ passes through $-1$.
For that value, the (formal) front interaction equations predict fronts with increasing speeds that blow up in finite time, which is also confirmed by simulations, see Fig.~\ref{fig:BeyondValidity}.
This final section is concluded by a brief discussion of the potential implications of our findings in the context of some scientific problems, especially ecosystem models which were a primary motivation.
\\
\begin{remark}
The stationary problem associated to \eqref{eq:mainEquation},
which is given by
\begin{equation}
\label{eq:stationaryODEsystem}
0 = u_{xx} + u - u^3 + \varepsilon F(u,u_x;x)
\; \; {\rm or} \; \;
\left\{
\begin{array}{r c l}
u_x & = & p, \\
p_x & = & u^3 - u - \varepsilon F(u,p;x),
\end{array}
\right.
\end{equation}
is of the form of a weakly driven (planar) nonlinear oscillator.
This form will play a central role in the upcoming analysis.
Since the classical theory of periodically driven oscillators predicts chaotic behavior (see \cite{Guckenheimer2002}), one may thus expect (stationary) `chaotic multi-front patterns' if $F(u,p;x)$ varies periodically in $x$ -- see Fig. \ref{fig:NumericsMultiFront-Intro} for stable patterns that indeed seem to be of this nature.
\end{remark}

\begin{remark}
\label{rem:slowlyvarying}
There is a more extensive literature on inhomogeneous Allen-Cahn-type equations in which the inhomogeneity varies slowly in space ({\it i.e.}, on long spatial scales) -- see \cite{ai2006layers, ei1990pattern, hale1999multiple} and the references therein.
In these references, \eqref{eq:mainEquation} is typically written as a singularly perturbed equation, $U_t = \varepsilon^2 U_{XX} + G(U,X)$ which is equivalent to $U_t = U_{xx} + G(U,\varepsilon x)$ with $x = X/\varepsilon$.
(In fact, in \cite{hale1999multiple} $G(U,\varepsilon x) = (1-U^2)(U - a(\varepsilon x))$ which corresponds to $\varepsilon F(U,U_x,x) = - a(\varepsilon x) (1-U^2)$ in \eqref{eq:mainEquation}.)
Here, we instead consider inhomogeneities that vary on the same spatial scale as $U(x,t)$.
To set up the asymptotic analysis, we however need to assume that the amplitude of the inhomogeneities is small (which is not necessary when the inhomogeneities vary on long spatial scales).
See also \cite{brena2014mathematical,tzou2018stability} and the references therein for the analysis of inhomogeneous systems of reaction-diffusion equations in which the inhomogeneity varies on a long spatial scale.
\end{remark}

\begin{remark}[numerical simulations] Most of the numerical PDE simulations presented in this article are made in Mathematica using the NDSolveValue and MethodOfLines routines. The exceptions are the results in Fig.~\ref{fig:NumericsMultiFront-Intro}, which are made in MATLAB using the pdepe routine. \label{remark:numerics}
\end{remark}

\subsection{Zero-front solutions: the persistence of the background states}
\label{s:uplusminus}
\begin{figure}
\centering
\begin{subfigure}[t]{0.32\textwidth}
\centering
\includegraphics[width=\textwidth]{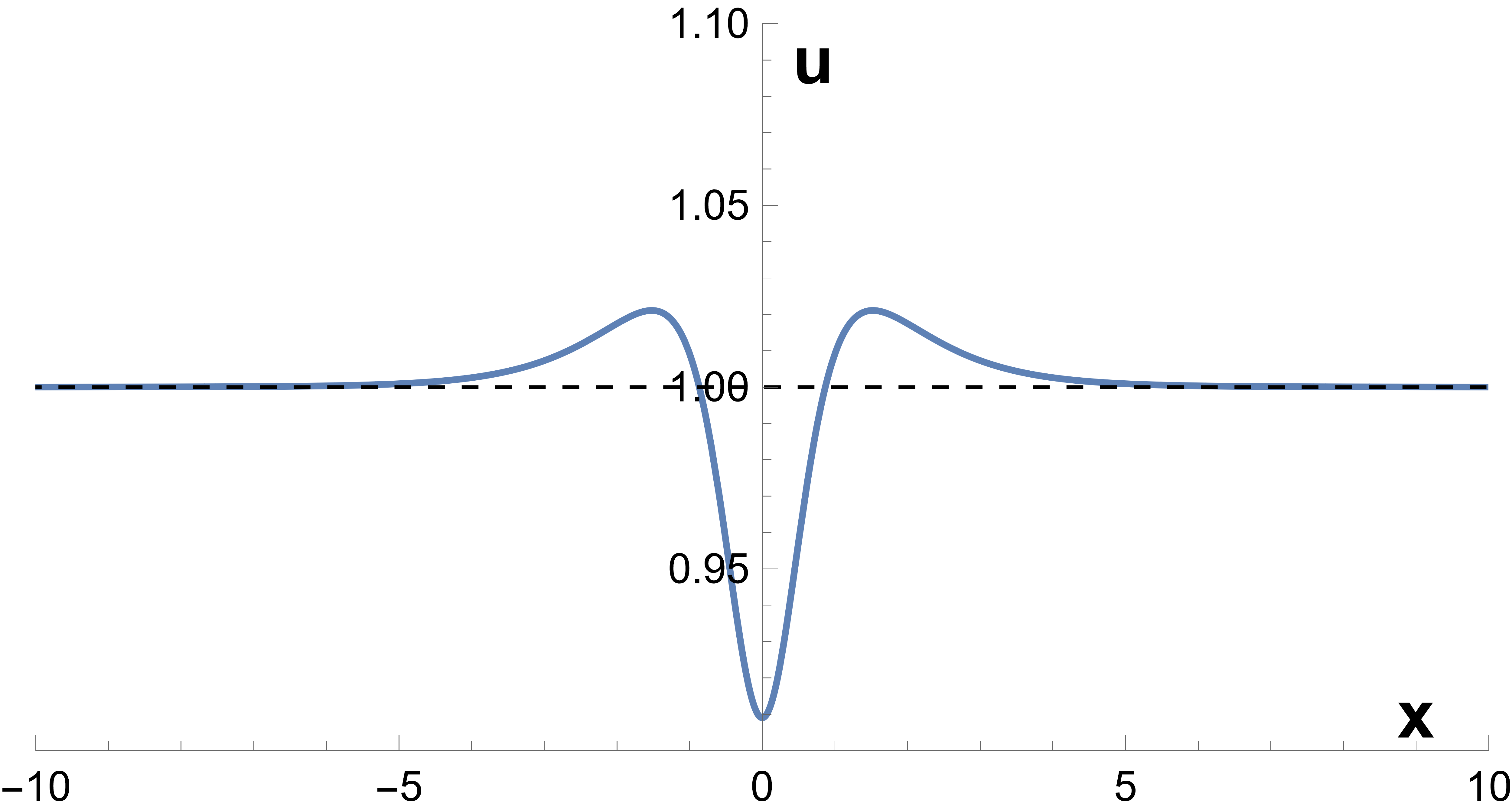}
\caption{}
\end{subfigure}
~
\begin{subfigure}[t]{0.32\textwidth}
\centering
\includegraphics[width = \textwidth]{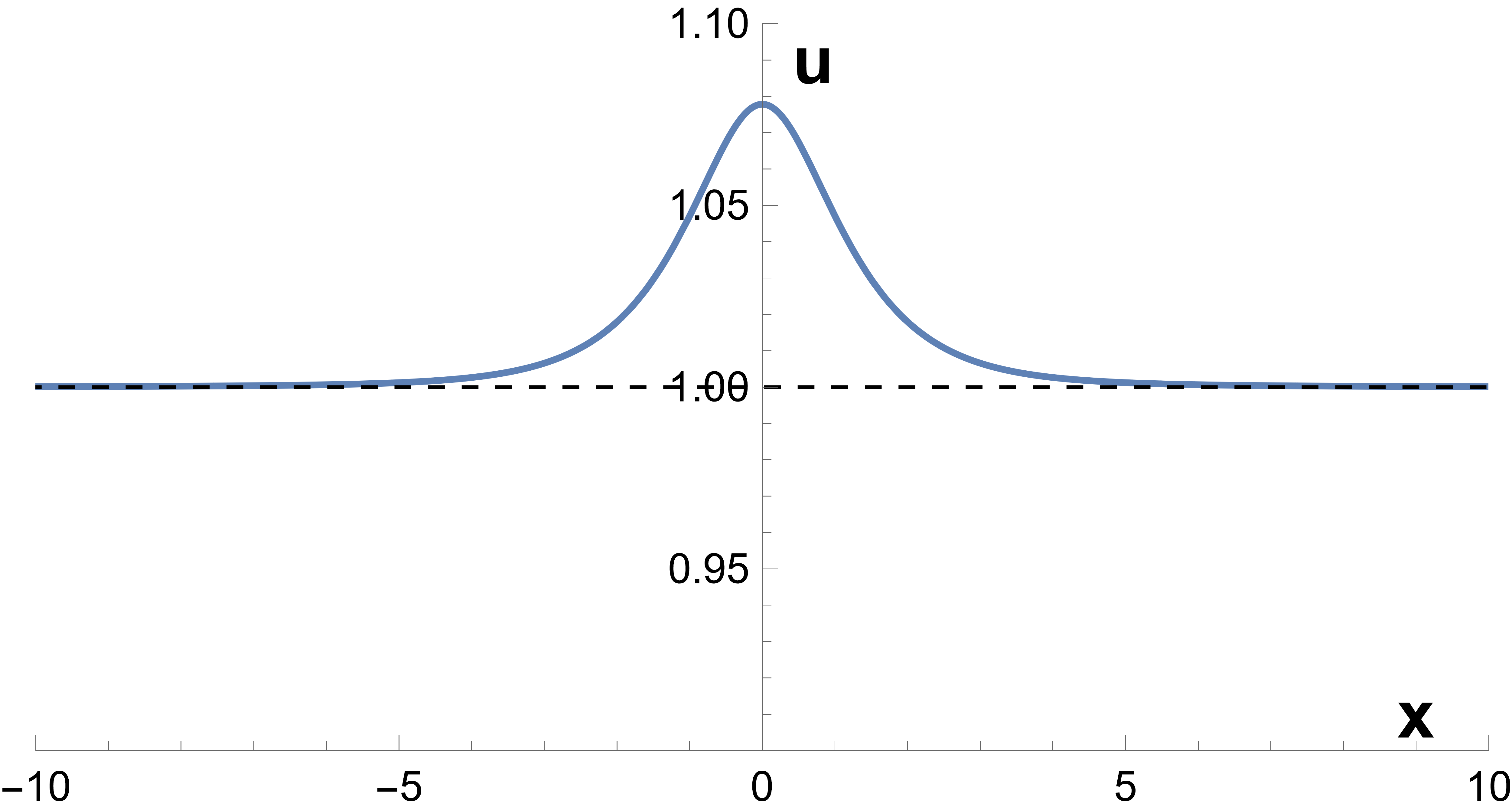}
\caption{}
\end{subfigure}
~
\begin{subfigure}[t]{0.32\textwidth}
\centering
\includegraphics[width = \textwidth]{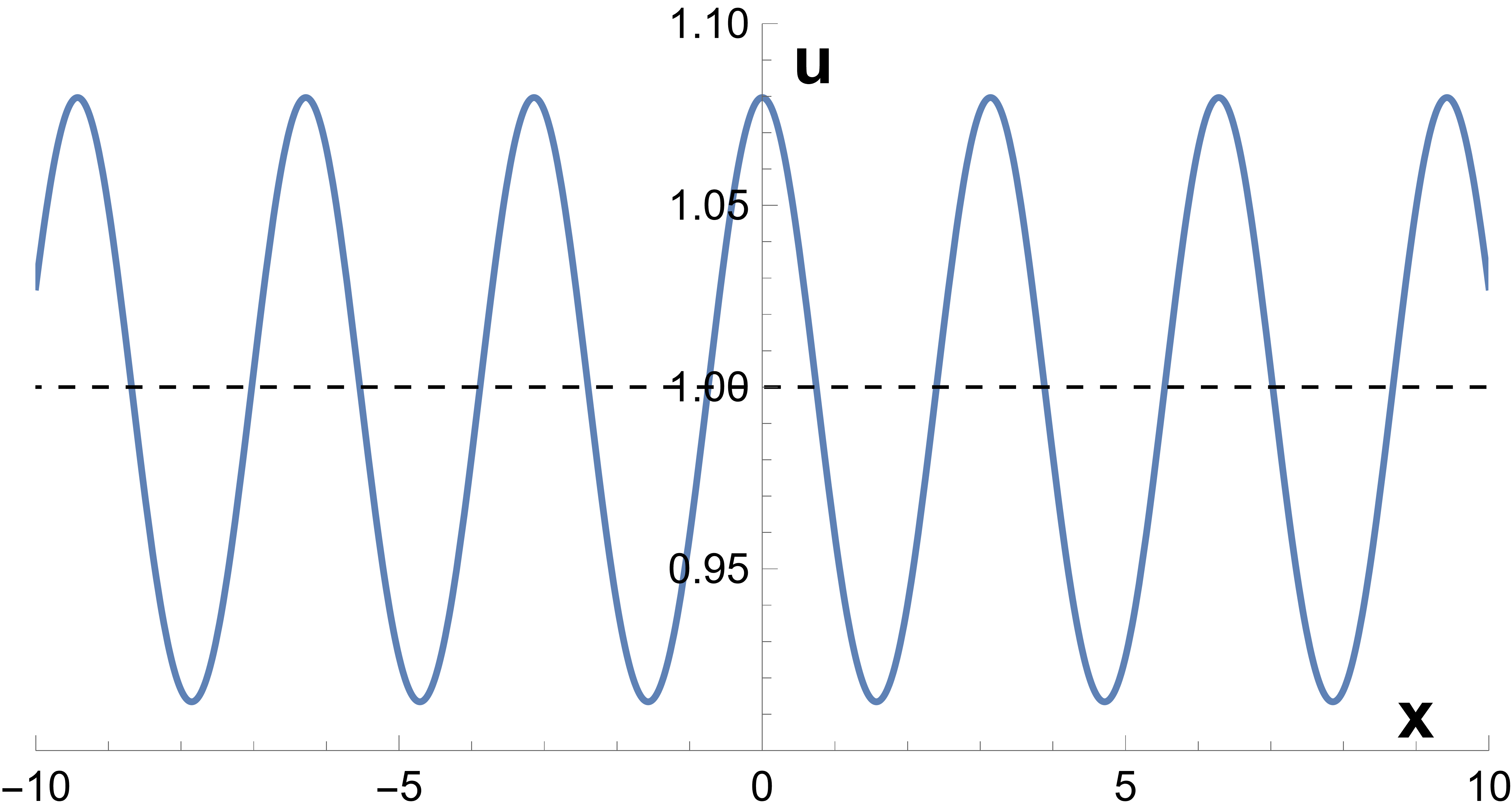}
\caption{}
\end{subfigure}
\caption{The zero-front background states $u_+^\varepsilon(x)$ for inhomogeneous equation \eqref{eq:mainEquation} with topographical inhomogeneities \eqref{eq:Ftopography}. (a) $H(x) = H_{\rm alg}(x;4.0)$, (b) $H(x) = H_{\rm alg}(x;0.0)$, (c) $H(x) = \sin 2x$ and $\varepsilon = 0.25$. Plots are obtained via numerical simulations; see also remark~\ref{remark:numerics}}
\label{fig:PertBackgroundStates}
\end{figure}

In this brief subsection, we consider the existence and stability of the background states of \eqref{eq:mainEquation} to which all multi-front patterns will converge in the limits $x \to \pm \infty$. In the homogeneous Allen-Cahn equation \eqref{eq:standardAC}, there are three such background states: $u_\pm (x) \equiv \pm 1$ and $u_0(x) \equiv 0$.
The following dispersion relations
define the location of the (essential) spectrum $\Sigma_\mathrm{ess}$ associated to their stability:
\begin{align*}
u_{\pm}:& \quad \lambda_{\pm} (k) = - k^2 + (1 - 3 u_{\pm}^2) = - k^2 - 2, \\
u_0:& \quad \lambda_0(k) = -k^2 + (1 - 3 u_0^2) = - k^2 + 1.
\end{align*}
Thus, the states $u_{\pm}$ are linearly stable against small-amplitude perturbations, whereas the state $u_0$ is unstable.
\\ \\
For \eqref{eq:mainEquation} with $\varepsilon > 0$, the classical solutions $u_-$, $u_0$, and $u_+$ are generally no longer stationary solutions. However, for sufficiently small $\varepsilon$, there are unique, bounded `zero-front' solutions, $u_{\pm}^\varepsilon$, which lie $\mathcal{O}(\varepsilon)$-close to the classical solutions $u_{\pm} = \pm 1$. These are given by
\begin{equation}
\label{eq:uMinusPlus}
u_\pm^\varepsilon(x) = \pm 1 + \frac{\varepsilon}{2\sqrt{2}} \left[ \int_x^\infty F(\pm 1,0,z) e^{\sqrt{2}(x-z)}\ dz + \int_{-\infty}^x F(\pm 1,0,z) e^{-\sqrt{2}(x-z)}\ dz \right] + \mathcal{O}(\varepsilon^2).
\end{equation}
For future reference, we note that, in the case of topographic inhomogeneities \eqref{eq:Ftopography},
\begin{equation}
\label{eq:uMinusPlusTopo}
u_\pm^\varepsilon(x) = \pm 1 \pm \frac{\varepsilon}{2\sqrt{2}} \left[ \int_x^\infty H''(z)e^{\sqrt{2}(x-z)}\ dz + \int_{-\infty}^x H''(z) e^{-\sqrt{2}(x-z)}\ dz \right] + \mathcal{O}(\varepsilon^2),
\end{equation}
{\it i.e.}, $u_\pm^\varepsilon(x)$ only depends on $H''(x)$ (to leading order) and is bounded (for $x \in \mathbb{R}$) under the condition that $H''(x)$ is (recall \eqref{eq:defHuni-alg}: $H'(x)$ may thus be unbounded here).
\\ \\
In the extended $(u,u_x,x)$-phase space of the stationary problem \eqref{eq:stationaryODEsystem}, the sets  $\{(u_{\pm}^\varepsilon(x), \frac{d}{dx} u_{\pm}^\varepsilon(x), x) | x \in \mathbb{R} \}$ representing the zero-front solutions are 1D saddle invariant manifolds, which we label $\mathcal{M}^\varepsilon_\pm$. These perturb for $0<\varepsilon\ll 1$ from the lines $\{(u_\pm, 0, x) | x \in \mathbb{R}\}$ of saddle points of the $\varepsilon=0$ system. 
Their persistence, for $\varepsilon$ sufficiently small and suitably smooth functions $F$, is guaranteed by classical theory for normally hyperbolic invariant manifolds, see for example \cite{Fenichel1971, Hirsch1970,HPS1977}.   Moreover, by the same persistence theory for normally hyperbolic manifolds, the normally hyperbolic invariant manifolds $\mathcal{M}^\varepsilon_\pm$ have 2D local stable and unstable manifolds, $W^{s,u}( \mathcal{M}^\varepsilon_\pm)$, for all $\varepsilon$ sufficiently small. These manifolds are the cornerstones of the theory for the front solutions presented in this article.
\\ \\
The stability of the zero-front solutions $u_\pm^\varepsilon(x)$
is determined by linearization,
\begin{equation}
	\left( \mathcal{L} - \lambda \right) \bar{u} = 0, \label{eq:eigenvalueProblem}
\end{equation}
where
\begin{equation}
	\mathcal{L} := \partial_x^2 + (1 - 3 (u_\pm^\varepsilon(x))^2) + \varepsilon \left[ F_u(u_\pm^\varepsilon(x),\partial_x u_\pm^\varepsilon(x), x) + F_{u_x}(u_\pm^\varepsilon(x),\partial_x u_\pm^\varepsilon(x), x) \partial_x \right].
\end{equation}
Hence, for $\varepsilon = 0$, the linear stability problem reduces to the classical case of the Allen-Cahn equation.
Then, for any $0<\varepsilon \ll 1$, the linear stability system can be rewritten (using \eqref{eq:uMinusPlus}) as
\begin{equation}
	\left( \mathcal{L}_0 + \mathcal{L}_\varepsilon  - \lambda \right) \bar{u} = 0,
\end{equation}
with
\begin{equation}
\label{eq:definitionL0Leps}
\mathcal{L}_0 = \partial_x^2 -2, \; \; \ \ \ \
\mathcal{L}_\varepsilon = - 3 \left[ (u_\pm^\varepsilon)^2 - 1 \right]+ \varepsilon \left[F_u(u^\varepsilon_\pm,\partial_x u^\varepsilon_\pm, x) + F_{u_x}(u^\varepsilon_\pm,\partial_x u^\varepsilon_\pm, x) \partial_x  \right] + \mathcal{O}(\varepsilon^2).
\end{equation}
Hence, $\mathcal{L}_\varepsilon = \mathcal{O}(\varepsilon)$.
By the theory of exponential dichotomies
\cite{coppel1978stability, palmer1988exponential}, it can be concluded that (for $\varepsilon$ sufficiently small), the spectrum of $u_\pm^\varepsilon$ lies $\varepsilon$-close to the spectrum of $u_\pm = \pm 1$. Therefore, $u_{\pm}^\varepsilon$ are stable as solutions of \eqref{eq:mainEquation}, see Fig.~\ref{fig:PertBackgroundStates}.

\section{Stationary one-front solutions in the inhomogeneous Allen-Cahn equation}
\label{sec:1fronts}
In this section, we consider the existence and stability of stationary one-front patterns in the inhomogeneous Allen-Cahn equation \eqref{eq:mainEquation}, starting from a brief preparation on (stable) one-fronts in the homogeneous equation \eqref{eq:standardAC}.
See Fig.~\ref{fig:1FrontSolutions}(b) for a homogeneous one-front pattern and Fig. \ref{fig:1FrontSolutions}(c), (d) for examples of stable stationary one-front patterns in the Allen-Cahn equation with spatially localized and spatially periodic inhomogeneities, respectively.

\subsection{One-front solutions in the homogeneous Allen-Cahn equation}
\label{sec:1fronts-class}
\begin{figure}
\centering
\begin{subfigure}[t]{0.32\textwidth}
\centering
\includegraphics[width = \textwidth]{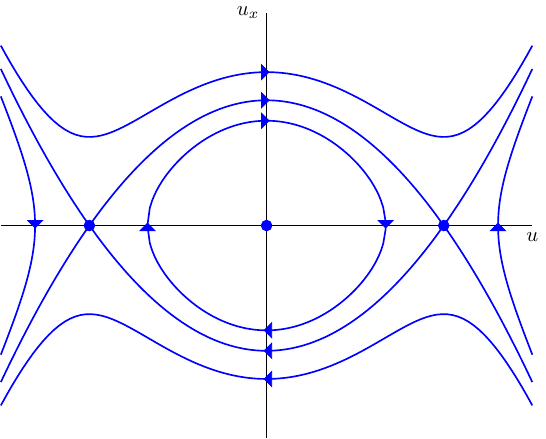}
\end{subfigure}
~
\begin{subfigure}[t]{0.205\textwidth}
\centering
\includegraphics[width = \textwidth]{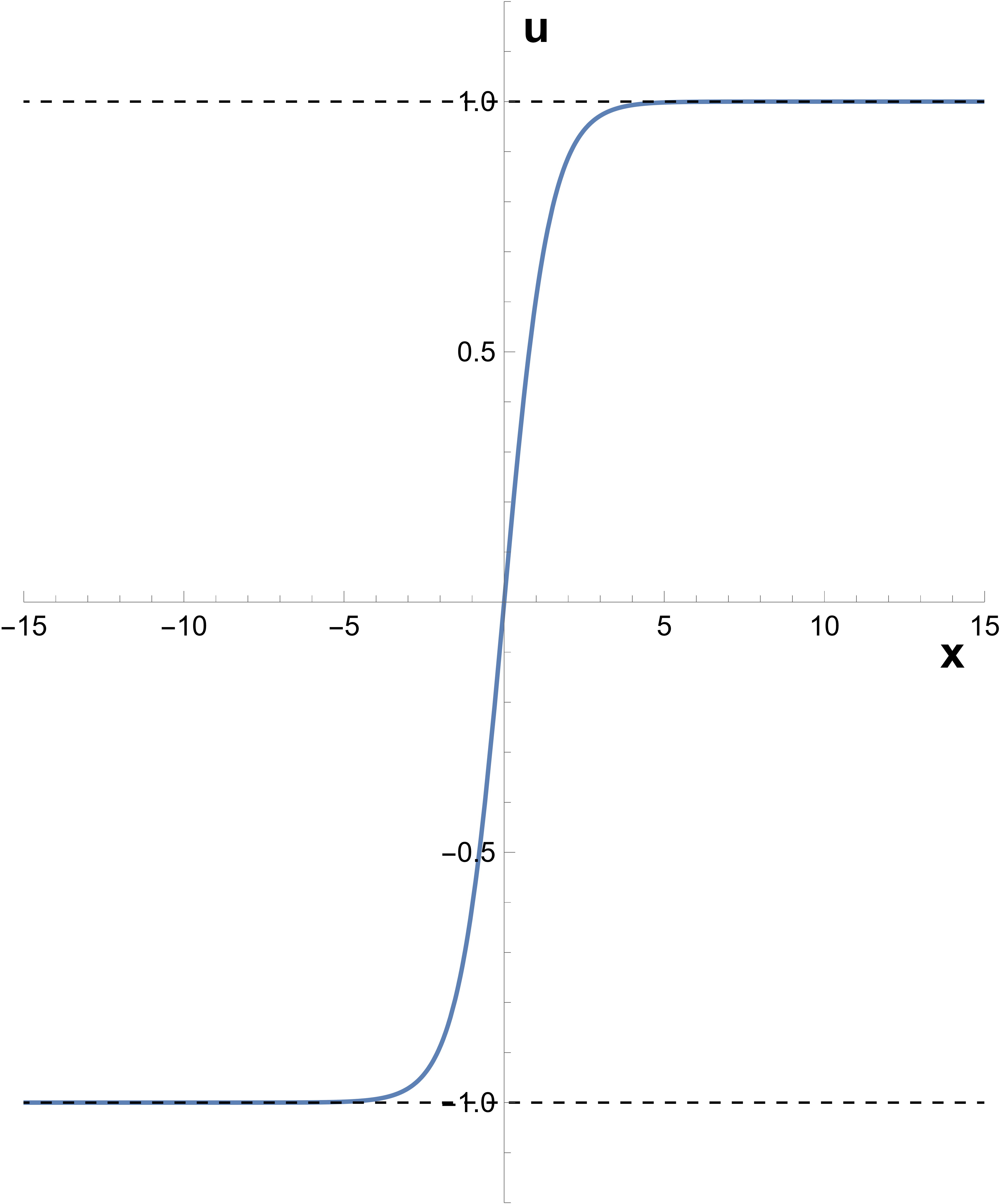}
\end{subfigure}
~
\begin{subfigure}[t]{0.205\textwidth}
\centering
\includegraphics[width = \textwidth]{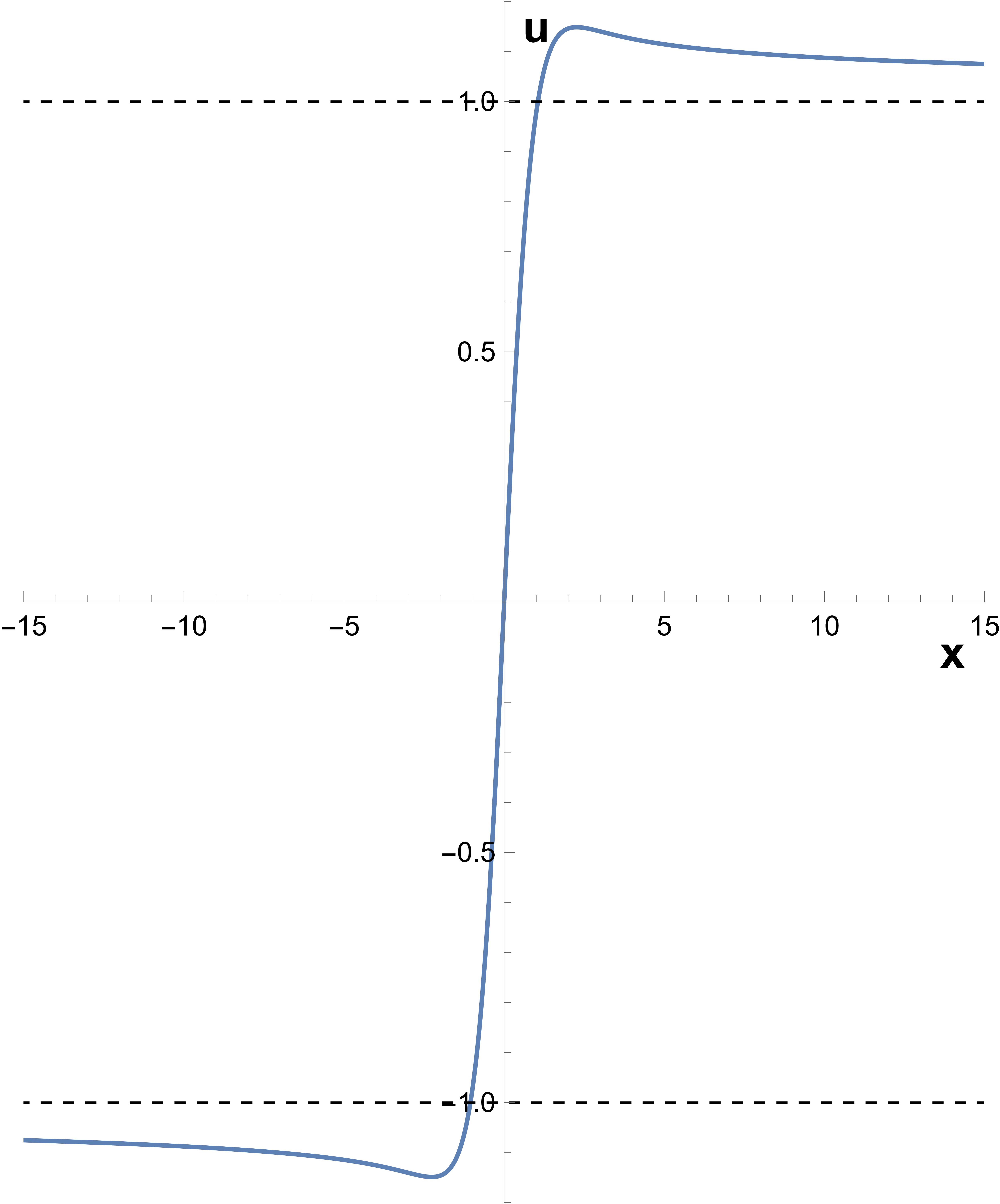}
\end{subfigure}
~
\begin{subfigure}[t]{0.205\textwidth}
\centering
\includegraphics[width = \textwidth]{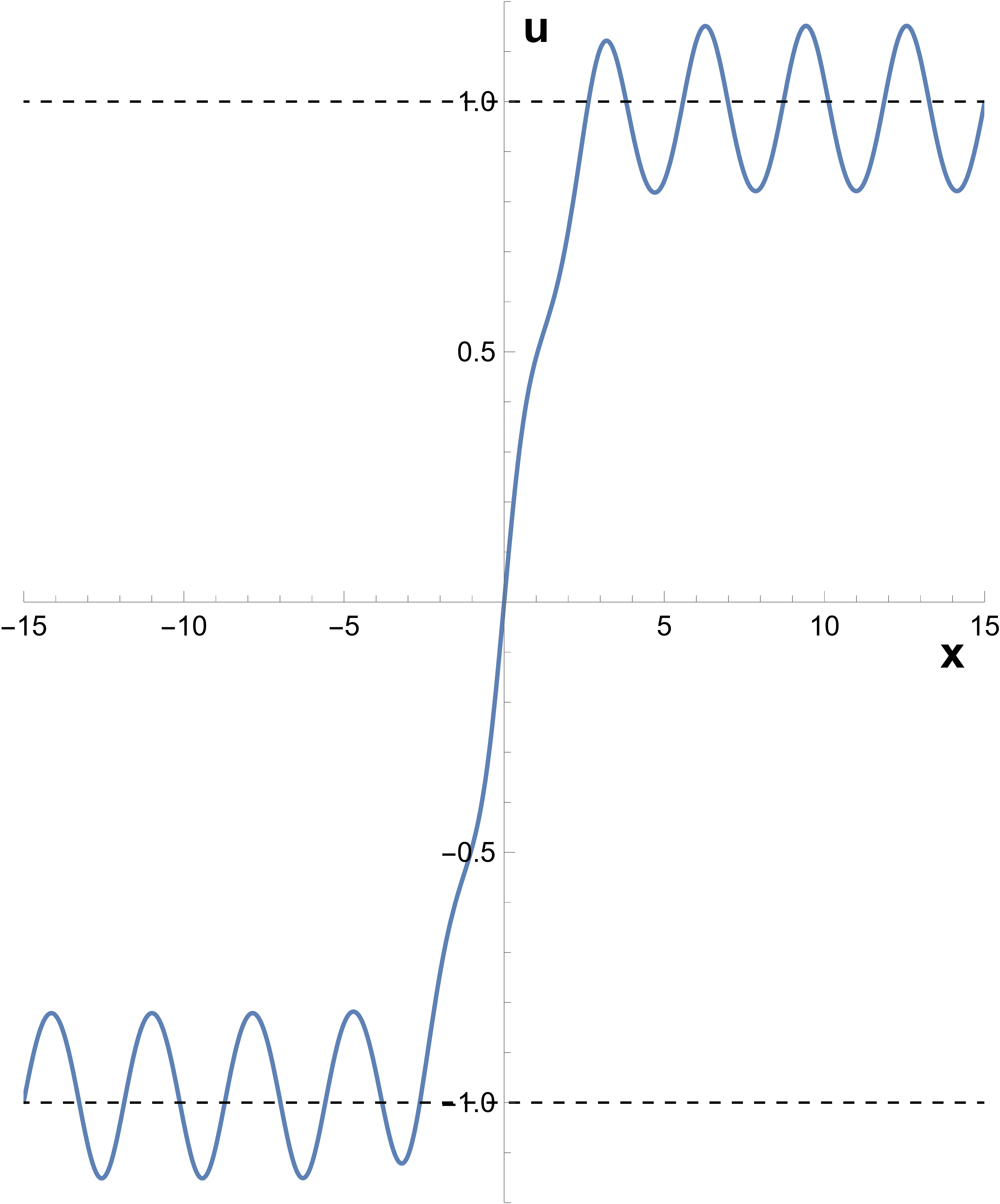}
\end{subfigure}
\caption{ (a) The phase portrait for the unperturbed ($\varepsilon =0$) planar system \eqref{eq:stationaryODEsystem} (determined by the level sets of $\mathcal{H}(u,u_x)$ \eqref{eq:defHamH}), (b) the unperturbed one-front pattern
$u_\textrm{up}(x;0)$ \eqref{eq:heteroclinicSolution} as stable solution of \eqref{eq:mainEquation} with $\varepsilon = 0$ and (c-d) (stable) stationary patterns in \eqref{eq:mainEquation} with $\varepsilon = 0.5$ with (c) a `localized'  and (d) a spatially periodic topography \eqref{eq:Ftopography} ((c) $H(x) = H_{\rm alg}(x;-0.6)$, (d) $H(x) = \sin 2x $). (Note that we have chosen $p < 0$ in (c) to better illustrate the impacts of the heterogeneity, whereas with $p < 0$ the topography is not really localised (see also Sec.~\ref{sec:Discussion}).}
\label{fig:1FrontSolutions}
\end{figure}

The stationary problem associated to the homogeneous Allen-Cahn equation \eqref{eq:standardAC}, {\it i.e.} \eqref{eq:stationaryODEsystem} with $\varepsilon = 0$, is an integrable system with Hamiltonian
\begin{equation}
\label{eq:defHamH}
	\mathcal{H}(u,u_x) = \frac12 u_x^2 + \frac12 u^2 - \frac14 u^4.
\end{equation}
The asymptotic states $u_\pm = \pm 1$ lie on the level set $\mathcal{H}(u,u_x) = \frac14$ so that there exist a heteroclinic connection $u_\textrm{up}(x)$ from $u_-$ to $u_+$ and a heteroclinic connection $u_\textrm{down}(x)$ from $u_+$ to $u_-$, see Fig. \ref{fig:1FrontSolutions}(a); $u_\textrm{up/down}(x)$ can be explicitly determined,
\begin{equation}
\label{eq:heteroclinicSolution}
u_\textrm{up}(x;\phi) = + \tanh\left(\frac12 \sqrt{2} [x - \phi] \right),  \; \;
u_\textrm{down}(x;\phi) = - \tanh\left(\frac12 \sqrt{2} [x - \phi] \right)
\end{equation}
(see Fig. \ref{fig:1FrontSolutions}(b)), where $\phi$ is an arbitrary phase shift. For notational clarity, we take $\phi = 0$ in the rest of this section. However, in the forced system, the phase shift is no longer arbitrary, and only specific phase shifts will be allowed. As a preparation for the upcoming analysis of the non-autonomously perturbed system \eqref{eq:stationaryODEsystem}, it is useful to view these heteroclinic orbits in the 3D $(u,u_x,x)$ extended phase space. For $\varepsilon = 0$, branches of the manifolds $W^u(\mathcal{M}^0_-)$ and $W^s(\mathcal{M}^0_+)$ coincide, as do branches of $W^u(\mathcal{M}^0_+)$ and $W^s(\mathcal{M}^0_-)$. These form 2D surfaces connecting $\mathcal{M}^0_-$ and $\mathcal{M}^0_+$, and they are foliated by the orbits \eqref{eq:heteroclinicSolution}, taken over all $x$ and $\phi$.
\\ \\
Since the one-fronts $u_\textrm{up/down}(x)$ are asymptotic to $u_\pm$, the essential spectrum associated to the stability of these one-front patterns is given by
\begin{equation}
\label{eq:SigmaEss}
	\Sigma_\mathrm{ess} = \{ \lambda \in \mathbb{C}: \lambda \leq -2\}.
\end{equation}
To find the eigenvalues in the point spectrum, $\Sigma_\mathrm{pt}$, we look for the $\lambda$ that solve the eigenvalue problem~\eqref{eq:eigenvalueProblem} with $\varepsilon =0$.
By the reversibility symmetry $x \rightarrow -x$ of \eqref{eq:mainEquation} (for $\varepsilon = 0$), both heteroclinic solutions have the same spectrum. Furthermore, in view of \eqref{eq:heteroclinicSolution}, the eigenvalue problems are given by
\begin{equation}
	\bar{u}_{xx} + \left(1 - 3 \tanh^2\left(\frac12 \sqrt{2} x\right) - \lambda\right)\bar{u} = 0.
\end{equation}
Also this equation can be solved explicitly (for any $\lambda \in \mathbb{C}$),
\begin{equation} \label{eq:solutionbasis}
\begin{array}{rcl}
\bar{u}_+(x;\lambda) & = & e^{-\sqrt{2+\lambda}x} \left(1 + \frac23 \lambda + \tanh^2\left(\frac12 \sqrt{2} x \right) + 2 \sqrt{1 + \frac12 \lambda} \tanh\left(\frac12 \sqrt{2} x\right) \right),
\\[1mm]
\bar{u}_-(x;\lambda) & = & e^{\sqrt{2+\lambda}x} \left(1 + \frac23 \lambda + \tanh^2\left(\frac12 \sqrt{2} x\right) - 2 \sqrt{1 + \frac12 \lambda} \tanh\left(\frac12 \sqrt{2} x\right) \right).
\end{array}
\end{equation}
For all $\lambda \in \mathbb{C} \backslash \Sigma_\mathrm{ess}$, one of these (Jost) solutions approaches $0$ as $x \rightarrow \infty$ and the other as $x \rightarrow -\infty$. The eigenvalues may be found as the roots of an Evans function (see e.g.~\cite{alexander1990topological, sandstede2002stability} and references therein), given by
\begin{equation}
	\mathcal{D}(\lambda) := \left| \begin{array}{cc} \bar{u}_+(0;\lambda) & \bar{u}_-(0;\lambda) \\ \bar{u}'_+(0;\lambda) & \bar{u}'_-(0;\lambda) \end{array} \right| = \frac{4}{9} \lambda (3 + 2 \lambda) \sqrt{2 + \lambda}.
\end{equation}
Therefore, the eigenvalues are $\lambda_0 = 0$, $\lambda_1 = - 3 / 2$ and $\lambda_2 = -2$, at the edge of $\Sigma_\mathrm{ess}$ \eqref{eq:SigmaEss}.
All eigenvalues have negative real parts, except the eigenvalue $\lambda_0 = 0$, which corresponds to the translational invariance of \eqref{eq:standardAC}. Hence, the one-front patterns $u_\textrm{up/down}(x)$ have a spectrum that is bounded away from the imaginary axis, except for a simple eigenvalue at the origin, which indicates that they are (nonlinearly) stable ({\it i.e.}, stable against all (small) perturbations up to phase shifts/translations of the front solution \cite{kapitula2013spectral}).

\subsection{The existence of stationary one-fronts} \label{sec:frontExistence}
For $0<\varepsilon \ll 1$, \eqref{eq:mainEquation} is no longer translation invariant. Thus, the phase shifts $\phi$ in \eqref{eq:heteroclinicSolution} can no longer be chosen arbitrarily like in the $\varepsilon = 0$ case. Instead, the heterogeneous terms force specific phases; or, perhaps more intuitively, specific locations for the heteroclinic connections. To determine which phases can be selected (if any), we need to determine which heteroclinic solutions persist for $\varepsilon > 0$.
\\ \\
Let $u_h^\varepsilon$ be a stationary solution for $\varepsilon \neq 0$ that lies close to a heteroclinic solution $u_h$ of the $\varepsilon = 0$ problem ({\it i.e.,} $u_h$ is either $u_\mathrm{up}$ or $u_\mathrm{down}$ as in~\eqref{eq:heteroclinicSolution}).
That is, let
\begin{equation}
	u_h^\varepsilon(x;\phi) = u_h(x;\phi) + \varepsilon u_1(x;\phi) + \mathcal{O}(\varepsilon^2).
\end{equation}
Substitution in \eqref{eq:mainEquation} yields
\begin{equation}
	\mathcal{L}_0 u_1 =-F(u_h(x;\phi),\partial_x u_h(x;\phi),x) + \mathcal{O}(\varepsilon),\label{eq:existenceProblem1Front}
\end{equation}
where (with slight abuse of notation)
\begin{equation}
\label{eq:defL0}
\mathcal{L}_0 = \left[ \partial_x^2 + 1 - 3 u_h^2(x;\phi) \right].
\end{equation}
The homogeneous equation $\mathcal{L}_0 \Psi = 0$ has two linearly independent solutions, $\Psi_\mathrm{b}$ and $\Psi_\mathrm{u}$. One of those, $\Psi_\mathrm{b}$, is bounded, while the other $\Psi_\mathrm{u}$ is unbounded, as neither the $x \rightarrow \infty$ nor the $x \rightarrow -\infty$ limits exist (see Appendix \ref{sec:orderEpsSystem} for explicit expressions $\Psi_\mathrm{b/u}(x; \phi)$). Hence, $\mathcal{L}_0$ has non-empty kernel, since $\Psi_\mathrm{b}(x; \phi) = \partial_x u_h(x;\phi)$ is contained in this kernel. The inhomogeneous equation~\eqref{eq:existenceProblem1Front} only has a bounded solution if a solvability condition is satisfied. Let
\begin{equation}
	\mathcal{R}(\phi) := \int_{-\infty}^\infty F(u_h(x;\phi),\partial_x u_h(x;\phi),x)\ \partial_x u_h(x;\phi)\ dx.
\label{eq:FredholmCondition1Front}
\end{equation}
The solvability condition is that $\mathcal{R}(\phi)$ has a simple zero.
It can be verified that solutions satisfying \eqref{eq:existenceProblem1Front} which stay bounded for $x \rightarrow -\infty$ have the following form (see Appendix~\ref{sec:orderEpsSystem}):
\begin{equation}
	u_1(x;\phi) = A(x;\phi) \Psi_b(x;\phi) + B(x;\phi) \Psi_u(x;\phi),
\label{eq:generalSolutionOrderEpsSystem}
\end{equation}
where
\begin{equation}
\label{eq:defABphi}
\begin{array}{rcl}
A(x;\phi) & = & \int_\phi^x F(u_h(z;\phi),u_h'(z;\phi),z)\ \Psi_u(z;\phi)\ dz,\\
B(x;\phi) & = & -\int_{-\infty}^x F(u_h(z;\phi),u_h'(z;\phi),z)\ \Psi_b(z;\phi)\ dz.
\end{array}
\end{equation}
Hence, $u_1$ stays bounded for $x \rightarrow \infty$ if and only if the solvability condition that $\mathcal{R}(\phi)$ has a simple zero is satisfied, where $\mathcal{R}$ is given by
\eqref{eq:FredholmCondition1Front} (see Appendix~\ref{sec:orderEpsSystem} for a comprehensive study of the inhomogeneous problem).
Thus, for suitable $\phi$, {\it i.e.,} for those that satisfy $\mathcal{R}(\phi) = 0$, a heteroclinic connection exists that connects the asymptotic state $u_+^\varepsilon$ to $u_-^\varepsilon$, or $u_-^\varepsilon$ to $u_+^\varepsilon$.
Examples are given in Fig.~\ref{fig:1FrontSolutions} and Sec.~2.4.
\\ \\
Finally, we observe that the solvability condition is possibly different for the persistence of front solutions connecting $u^\varepsilon_-$ to $u^\varepsilon_+$ and the persistence of front solutions connecting $u^\varepsilon_+$ to $u^\varepsilon_-$.
In particular, the former persist if
\begin{equation}
	\mathcal{R}_\mathrm{up}(\phi) := \int_{-\infty}^\infty F(u_\mathrm{up}(x;\phi),\partial_x u_\mathrm{up}(x;\phi),x)\ \partial_x u_\mathrm{up}(x;\phi)\ dx  \label{eq:FredholmConditionUp}
\end{equation}
has a simple zero,
and the latter if
\begin{equation}
	\mathcal{R}_\mathrm{down}(\phi) := \int_{-\infty}^\infty F(u_\mathrm{down}(x;\phi),\partial_x u_\mathrm{down}(x;\phi),x)\ \partial_x u_\mathrm{down}(x;\phi)\ dx  \label{eq:FredholmConditionDown}
\end{equation}
does.
These are the key conditions that determine the parameter values
for which stationary one-front solutions persist in \eqref{eq:mainEquation} with $0< \varepsilon \ll 1$.
\\ \\
Geometrically, these persistence conditions are equivalent to the condition that the Melnikov functions have simple zeros.
In particular, one measures the splitting distance between
$W^u(\mathcal{M}^\varepsilon_\mp)$ and
$W^s(\mathcal{M}^\varepsilon_\pm)$
on fixed $(u,u_x)$ cross-sections in the 3D extended phase space (using an outward pointing normal to the 2D surface foliated by the unforced heteroclinics).
For $\varepsilon$ sufficiently small, the perturbed local stable and unstable manifolds $W_\mathrm{loc}^s(\mathcal{M}^\varepsilon_\pm)$ and $W_\mathrm{loc}^u(\mathcal{M}^\varepsilon_\pm)$ lie $C^b$ $\mathcal{O}(\varepsilon)$ close to their unforced counterparts~\cite{Fenichel1979}.  By expanding solutions on these manifolds on appropriate semi-infinite intervals, the splitting distance can be measured, and the coefficient on the leading order ($\mathcal{O}(\varepsilon)$) term is, up to a normalization factor, the Melnikov function $\mathcal{R}(\phi)$.
For values of $\phi$ at which $\mathcal{R}$ is not equal to zero, there is a non-zero distance between the perturbed manifolds, as measured on that cross-section, and the unperturbed heteroclinics for those values of $\phi$ do not persist. In contrast, for values of $\phi$ that correspond to simple zeroes of $\mathcal{R}$, the perturbed stable and unstable manifolds intersect transversely for $\varepsilon$ sufficiently small (as measured along the normal on that cross-section), and each transverse intersection point corresponds to a perturbed heteroclinic connection in the full 3D extended phase space. For more detailed information about Melnikov theory, we refer the reader to \cite[Chapter 4]{Guckenheimer2002} in case of periodic perturbations and to \cite{Chow1992} for general $\mathcal{O}(\varepsilon)$ perturbations.
\\ \\
In Subsections~\ref{sec:top1fronts} and \ref{sec:per1fronts}, we illustrate the Melnikov function and the use of these conditions to find one-front solutions in some examples. However, before giving these examples,
we examine the linear stability of one-front
solutions for general $F$.

\subsection{Linear stability of stationary one-front solutions}
\label{sec:oneFrontStability}
\begin{figure}
\centering
\begin{subfigure}[c]{0.3\textwidth}
	\centering
	\includegraphics[width=\textwidth]{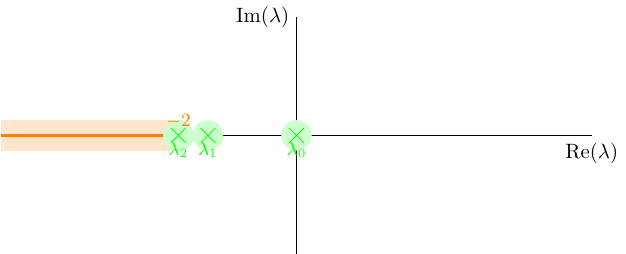}
	\caption{}
\end{subfigure}
~
\begin{subfigure}[c]{0.3\textwidth}
	\centering
	\includegraphics[width=\textwidth]{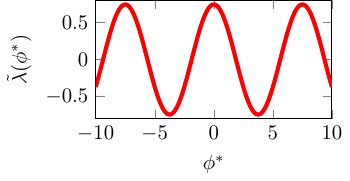}
	\caption{}
\end{subfigure}
~
\begin{subfigure}[c]{0.3\textwidth}
	\centering
	\includegraphics[width=\textwidth]{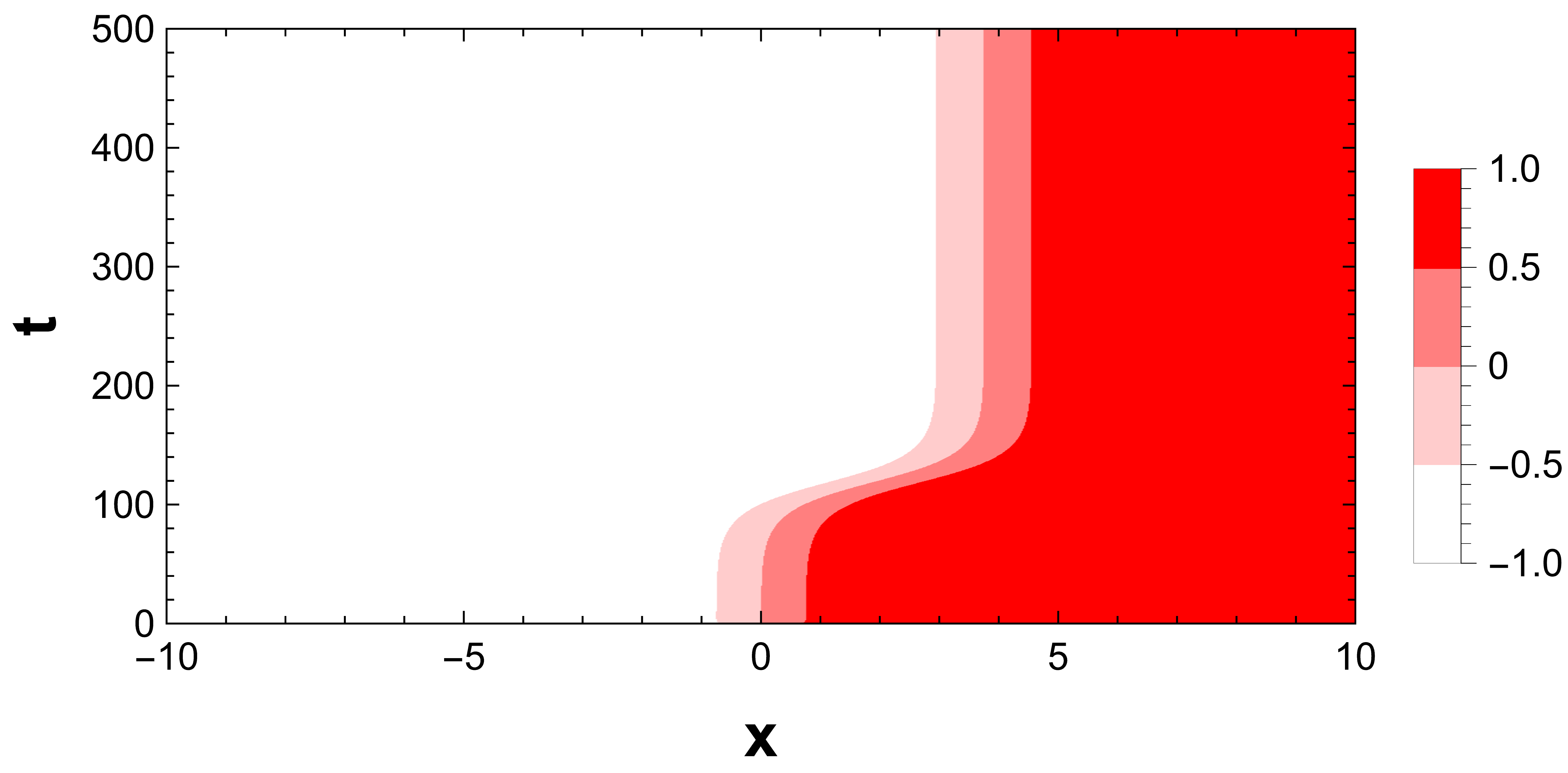}
	\caption{)}
\end{subfigure}
\caption{(a) Sketch of (bounds on) the spectrum of heteroclinic connections in the $\varepsilon \neq 0$ system. (b) A plot of $\tilde{\lambda}(\phi_{\ast})$ for $f_1(x) = \cos\left(\frac{4}{15}\pi x\right)$, $f_2(x) \equiv 0$, $f_3(x) \equiv 0$. Note $\tilde{\lambda}(0) \approx 0.75 > 0$, and $\tilde{\lambda}(15/4) \approx -0.75 < 0$. (c) Plot of $u(x,t)$ in the $(x,t)$ plane, with the indicated color code for the amplitude of $u$, obtained using direct simulation of the full PDE. The initial condition is very close to $u_h^0(x;0)$, i.e., the unstable front.
}
\label{fig:1FrontStability}
\end{figure}

For $0<\varepsilon\ll 1$, the spectral stability of the one-front patterns can be determined by the spectrum of the fronts in the $\varepsilon = 0$ system via the theory of exponential dichotomies. This theory shows that the spectrum lies $\varepsilon$-close to the spectrum of the fronts in the $\varepsilon = 0$ system. Therefore, for small enough $\varepsilon$, the essential spectrum is again to the left of the imaginary axis and bounded away from it. Concerning the eigenvalues, only the $\lambda_0 = 0$ eigenvalue can potentially move into the right-half of the complex plane and cause destabilization (see Fig.~\ref{fig:1FrontStability}).
\\ \\
To determine what happens to the $\lambda_0 = 0$ eigenvalue when $\varepsilon \neq 0$, we employ a formal regular expansion of the eigenvalue problem
\begin{equation}
    \left( \mathcal{L} - \lambda \right) \bar{u} = 0, \label{eq:eigenvalueProblemGeneral}
\end{equation}
where (again with a slight abuse of notation)
\begin{equation}
    \mathcal{L} := \partial_x^2 + (1 - 3 (u_*)^2) + \varepsilon \left[ F_u(u_*(x),\partial_x u_*(x), x) + F_{u_x}(u_*(x),\partial_x u_*(x), x) \partial_x \right].
\end{equation}
Here, $u_*$ is an exact one-front solution of the $\varepsilon \neq 0$ system. Since $\lambda_0 = 0$ is an eigenvalue with eigenfunction $\partial_x u_h(x;\phi)$ for the $\varepsilon = 0$ system, it is expected that the eigenfunction for the $\varepsilon \ll 1$ system is close to the derivative $\partial_x u_h^\varepsilon(x;\phi_{\ast})$ (where $\mathcal{R}(\phi_{\ast}) = 0$),
and that the associated eigenvalue is of size $\mathcal{O}(\varepsilon)$.
Hence, we set
\begin{equation}
	\begin{array}{rclcll}
		\bar{u}(x;\phi_{\ast}) & = & \partial_x u_h^\varepsilon(x;\phi_{\ast}) & + & \varepsilon \bar{u}_1(x;\phi_{\ast})  & +  \mathcal{O}(\varepsilon^2); \\
		\lambda & = & \varepsilon \tilde{\lambda} + \mathcal{O}(\varepsilon^2) & & & .
	\end{array}
\end{equation}
Substitution in \eqref{eq:eigenvalueProblemGeneral} yields the equation
\begin{equation}
	\varepsilon \tilde{\lambda} \partial_x u_h^\varepsilon(x;\phi_{\ast}) =  \mathcal{L}_0 \partial_x u_h^\varepsilon(x;\phi_{\ast}) + \varepsilon \mathcal{L}_0 \bar{u}_1 + \mathcal{L}_\varepsilon \partial_x u_h^\varepsilon(x;\phi_{\ast}) + \mathcal{O}(\varepsilon^2), \label{eq:eigenvalueProblem1Front}
\end{equation}
where $\mathcal{L}_0$ is defined in~\eqref{eq:defL0} and $\mathcal{L}_\varepsilon$ is (again with slight abuse of notation) given by
\begin{equation}
    \label{eq:defLeps}
    \mathcal{L}_\varepsilon = - 3 \left[ (u_*)^2 - 1 \right]+ \varepsilon \left[F_u(u_*,\partial_x u_*, x) + F_{u_x}(u_*,\partial_x u_*, x) \partial_x  \right] + \mathcal{O}(\varepsilon^2),
\end{equation}
and $u_* = u_h^\varepsilon(x;\phi_*)$. Since $u_h^\varepsilon(x;\phi_{\ast})$ is the exact solution to~\eqref{eq:mainEquation}, we can formally take the derivative of \eqref{eq:mainEquation} to find the (exact) relation
\begin{equation}
\mathcal{L}_0 \partial_x u_h^\varepsilon(x;\phi_{\ast}) + \mathcal{L}_\varepsilon \partial_x u_h^\varepsilon(x;\phi_{\ast}) = - \varepsilon F_x(u_h^\varepsilon(x;\phi_{\ast}), \partial_x u_h^\varepsilon(x;\phi_{\ast}),x).
\end{equation}
Upon using this relation in~\eqref{eq:eigenvalueProblem1Front}, we obtain
\begin{equation}
	\mathcal{L}_0 \bar{u}_1 = \tilde{\lambda} \partial_x u_h^\varepsilon(x;\phi_{\ast}) + F_x(u_h^\varepsilon(x;\phi_{\ast}), \partial_x u_h^\varepsilon(x;\phi_{\ast}),x) + \mathcal{O}(\varepsilon).
\end{equation}
Finally, we may approximate $u_h^\varepsilon = u_h + \mathcal{O}(\varepsilon)$ to obtain
\begin{equation}
	\mathcal{L}_0 \bar{u}_1 = \tilde{\lambda} \partial_x u_h(x;\phi_{\ast}) + F_x(u_h(x;\phi_{\ast}),\partial_x u_h(x;\phi_{\ast}),x) + \mathcal{O}(\varepsilon).
\end{equation}
As in the existence problem~\eqref{eq:existenceProblem1Front}, a bounded solution $\bar{u}_1$ only exists if the right-hand terms satisfy a (Fredholm) solvability condition, that now determines the eigenvalue $\tilde{\lambda}$. We find
\begin{equation}
\label{eq:eigenvalue1Front}
\begin{array}{rcl}
\tilde{\lambda}(\phi_{\ast})
& = &
-\frac{ \int_{-\infty}^\infty F_x(u_h(x;\phi_{\ast}),\partial_x u_h(x;\phi_{\ast}),x) \partial_x u_h(x;\phi_{\ast})\ dx }{\|\partial_x u_h(x;\phi_{\ast})\|_2^2}
\\
& = &
-\frac{3\sqrt{2}}{4} \int_{-\infty}^\infty F_x(u_h(x;\phi_{\ast}),\partial_x u_h(x;\phi_{\ast}),x) \partial_x u_h(x;\phi_{\ast})\ dx
=
- \frac{3\sqrt{2}}{4} \mathcal{R}'(\phi_{\ast}).
\end{array}
\end{equation}
(see \eqref{eq:FredholmCondition1Front}). Since the preceding analysis can be made rigorous in a straightforward way, we summarize our findings as follows.
\begin{theorem}
\label{th:1frontExStab}
Consider the weakly heterogeneous Allen-Cahn equation \eqref{eq:mainEquation} with bounded inhomogeneity $F(U,U_x,x)$.
Let the Melnikov function $\mathcal{R}(\phi)$ be given by \eqref{eq:FredholmCondition1Front} with $u_h = u_{\rm up}$ and assume that $\mathcal{R}(\phi)$ has a non-degenerate zero at $\phi = \phi_\ast$ (i.e., $\mathcal{R}(\phi_\ast) = 0$ and $\mathcal{R}'(\phi_\ast) \neq 0$).
Then, for $0 < \varepsilon \ll 1$ sufficiently small, \eqref{eq:mainEquation} has a stationary one-front solution $U_{\rm 1-f}(x)$ that is heteroclinic to $u^\varepsilon_\pm(x) = \pm 1 + \mathcal{O}(\varepsilon)$ \eqref{eq:uMinusPlus} (i.e., $\lim_{x \to \pm \infty} U_{\rm 1-f}(x) = u^\varepsilon_\pm(x)$) and that is $\mathcal{O}(\varepsilon)$ close to $u_{\rm up}(x;\phi_\ast)$ \eqref{eq:heteroclinicSolution} (i.e., ${\rm sup}_{x \in \mathbb{R}} |U_{\rm 1-f}(x) - u_{\rm up}(x;\phi_\ast)| =  \mathcal{O}(\varepsilon)$). The one-front pattern $U_{\rm 1-f}(x)$ is stable if $\mathcal{R}'(\phi_\ast) > 0$.
\end{theorem}
A similar theorem holds for the persistence of the stationary one-front solution $u_{\rm down}$.

\subsection{Stationary one-front solutions in topographically driven and in spatially periodic systems}
\label{s:pertop1fronts}
In this section, we apply the general theory to two main explicit classes of forcing functions $F(U,V,x)$ identified in the introduction: the  `topographically driven' case in which $F(U,V,x)$ is determined by a topography $H(x)$ \eqref{eq:Ftopography}, which is assumed here to be spatially localized, and the spatially periodically driven case in which $F(U,V,x)$ is given by \eqref{eq:canonicalExample}, with $f_j(x)$ and $X > 0$ such that $f_j(x+X) \equiv f_j(x)$.

\subsubsection{One-front solutions arising from weak topographical effects}
\label{sec:top1fronts}
\begin{figure}
\centering
\begin{subfigure}[t]{0.3\textwidth}
\centering
\includegraphics[width=\textwidth]{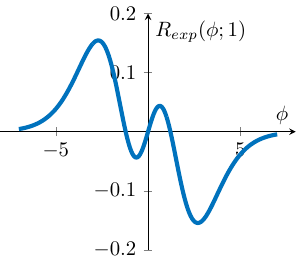}
\end{subfigure}
~
\begin{subfigure}[t]{0.3\textwidth}
\centering
\includegraphics[width = \textwidth]{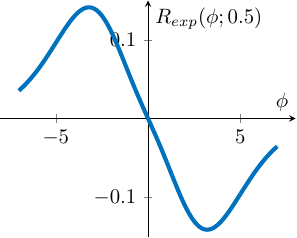}
\end{subfigure}
~
\begin{subfigure}[t]{0.3\textwidth}
\centering
\includegraphics[width = \textwidth]{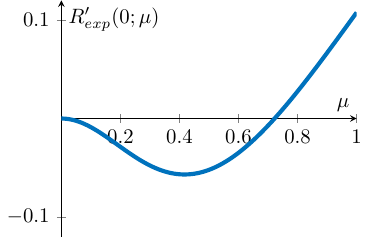}
\end{subfigure}

\caption{(a) The Melnikov function $\mathcal{R}_{\rm exp}(\phi;\mu)$ for a topography with a unimodal hill given by \eqref{eq:defHuni-alg} with $\mu = 1$, see (\ref{eq:Rsolhill1}). There are three positions for which one-front solutions exist: the front at $\phi =0$ is stable (since $\mathcal{R}_{\rm exp}'(0;1) > 0$), and the other two (with $\phi\ne 0$) are unstable. (b) The curve $\mathcal{R}_{\rm exp}(\phi;\frac{1}{2})$. By decreasing $\mu$ from $1$ to $\frac12$, we see that the two unstable zeroes have disappeared in a pitchfork bifurcation at the value $\mu_{\rm PF} = 0.722133...$, and the remaining one-front solution, which corresponds to the phase $\phi=0$, has become unstable. (c) The graph of $\mathcal{R}'_H(0;\mu)$ as a function of $\mu \in [0,1]$ reveals the value of
$\mu_{\rm PF} $.}
\label{fig:SolHillPitchfork}
\end{figure}

With \eqref{eq:Ftopography}, the spatial inhomogeneity $F(U,V,x)$ may be interpreted as topographically-driven.
The Melnikov function $\mathcal{R}(\phi)$ \eqref{eq:FredholmCondition1Front} is calculated explicitly by translating the variable of integration to $y=x-\phi$,
\begin{equation}
\begin{array}{rcl}
\mathcal{R}(\phi) & = & \int_{-\infty}^\infty \left( H'(y+\phi)\partial_y u_h(y;0) + H''(y+\phi) u_h(y;0)\right) \partial_y u_h(y;0) dy
\\
& = &
\int_{-\infty}^\infty \left( \partial_y\left[H'(y+\phi)u_h(y;0) \partial_y u_h(y;0) \right] - H'(y+\phi) u_h(y;0) \partial_{yy} u_h(y;0) \right) dy
\\
& = &
\int_{-\infty}^\infty H'(y+\phi) W_h(y) dy
\end{array}
\label{eq:Ex2R}
\end{equation}
with
\begin{equation}
W_h(y) :=  u_h^2(y;0) (1-u_h^2(y;0)) \in \left[0,\frac14\right].
\label{eq:defWh}
\end{equation}
See also Fig.~\ref{fig:1FrontSolutions}.
Since  $u_\textrm{up}(x;\phi) = -u_\textrm{down}(x;\phi)$ by \eqref{eq:heteroclinicSolution}, it
follows that
\begin{equation}
\begin{array}{rclcl}
\mathcal{R}_\mathrm{up}(\phi) & = & \int_{-\infty}^\infty H'(y+\phi) u_\mathrm{up}^2(y;0) (1-u_\mathrm{up}^2(y;0)) dy & &
\\
& = & \int_{-\infty}^\infty H'(y+\phi) u_\mathrm{down}^2(y;0) (1-u_\mathrm{down}^2(y;0)) dy & = & \mathcal{R}_\mathrm{down}(\phi).
\end{array}
\label{eq:Rtopoupdown}
\end{equation}
This is a very useful explicit form of the Melnikov conditions for the topographical case.
\\ \\
Here, we highlight three direct consequences of these explicit formulas for $\mathcal{R}_\mathrm{up}$ and $\mathcal{R}_\mathrm{down}$.
First, it follows directly from \eqref{eq:Ex2R} (and \eqref{eq:Rtopoupdown}) that if the topography $H(x)$ is spatially periodic, {\it i.e.}, if $H(x+X) \equiv H(x)$ for some (minimal) $X>0$, then also the Melnikov function $\mathcal{R}(\phi)$ is spatially periodic (with period $X$):
\begin{equation}
\label{eq:Rtopoper}
\mathcal{R}(\phi+X) =
\int_{-\infty}^\infty H'(y+\phi + X) W_h(y) dy =
\int_{-\infty}^\infty H'(y+\phi) W_h(y) dy =
\mathcal{R}(\phi).
\end{equation}
Second, for topographies $H(x)$ that are symmetric under $x \to -x$, {\it i.e.}, for which $H(x) \equiv H(-x)$, one has that $H'(x) \equiv - H'(-x)$ and, hence, that
\begin{equation}
\mathcal{R}(0) = \left(\int_{0}^\infty + \int_{-\infty}^0 \right) H'(y) W_h(y) dy = \int_{0}^\infty H'(y) W_h(y) dy - \int_{0}^\infty H'(z) W_h(z) dz = 0.
\end{equation}
Therefore, in the case of heterogeneities in which the topography $H(x)$ is symmetric around $x=0$, there must be a stationary one-front of PDE \eqref{eq:mainEquation} at $\phi_{\ast} =0$.
Also, these one-front patterns are linearly stable for $\mathcal{R}'(0)>0$, where
\begin{equation}
\mathcal{R}'(0) = 2 \int_{0}^\infty H''(y) W_h(y) dy.
\end{equation}
\\ \\
Third, and of most interest, the equation $\mathcal{R}(\phi) = 0$ may have more solutions.
As an explicit example, we consider the `unimodal hill' with exponentially decaying tails $H_{\rm exp}(x; \mu)$ \eqref{eq:defHuni-alg}, which gives for $\mu = 1$, with $\psi = \sqrt{2} \phi$ and $\mathcal{S}_{\rm exp}(\psi;\mu) = \mathcal{R}_{\rm exp}(\phi;\mu)$,
\begin{equation}
\label{eq:Rsolhill1}
\mathcal{S}_{\rm exp}(\psi;1) = - \frac{16 e^{\psi} [
(3\psi-13)e^{4\psi} +
2(27\psi-47)e^{3\psi} +
126\psi e^{2\psi} +
2(27\psi+47)e^{\psi} +
(3\psi+13)
]}{3(e^{\psi}-1)^6}.
\end{equation}
As shown in Fig.~\ref{fig:SolHillPitchfork}, $\mathcal{R}_{\rm exp}(\phi;1)$ has three zeroes. As the hill becomes less steep, {\it i.e.}, as $\mu$ decreases, two of these merge into a pitchfork bifurcation (at $\mu_{\rm PF} = 0.722133 ...$). The one-front at $\phi=0$ is only stable for solitary hills that are sufficiently steep ({\it i.e.}, for $\mu > \mu_{\rm PF}$). Note on the other hand that if we change the hill into a valley, {\it i.e.}, if we take $V_{\rm exp}(x;\mu) = -H_{\rm exp}(x;\mu)$ so that $\mathcal{R}_V(\phi;\mu) = - \mathcal{R}_{\rm exp}(\phi;\mu)$, then the one-front at $\phi=0$ is stable only if the slope of the ground is not too steep ({\it i.e.}, for $0 < \mu < \mu_{\rm PF}$). In this case, two stable states that are not centered at the deepest part of the value appear as the steepness of the valley increases through $\mu_{\rm PF}$.

\subsubsection{One-front solutions in systems with spatially-periodic heterogeneities}
\label{sec:per1fronts}

For general functions $F$ (and/or $f_j$) that are $X$-periodic in $x$, {\it i.e.}, $F(U,V,x+X) = F(U,V,x)$ for all $U,V,x \in \mathbb{R}$, the phase selection condition $\mathcal{R}(\phi) = 0$ is $X$-periodic, as well.
Hence, generically, there will be an even number of $\phi \in [0,X]$ satisfying the condition $\mathcal{R}(\phi) = 0$.
Moreover, the function $\tilde{\lambda}(\phi_{\ast})$ is $X$-periodic in $\phi_{\ast}$. Therefore, generically, half of the phase shifts $\phi_{\ast}$ that satisfy $\mathcal{R}(\phi_{\ast}) = 0$ will have $\tilde{\lambda}(\phi_{\ast}) < 0$ (making these stable), and the other half have $\tilde{\lambda}(\phi_{\ast}) > 0$ (making these unstable).
In the degenerate case that $\mathcal{R}(\phi) = 0$ has an odd number of solutions $\phi_{\ast}$, there will be a $\phi_{\ast}$ for which $\tilde{\lambda}(\phi_{\ast}) = 0$ (i.e. in that case $\phi_{\ast}$ is an extreme point of $\mathcal{R}$), signaling a (saddle-node) bifurcation.
\\ \\
For explicit calculations, we focus on the case of \eqref{eq:canonicalExample} with
\begin{equation}
\label{eq:fgh-example1}
f_1(x) = \alpha_1 \cos(k x), \ \ \ \
f_2(x) = \alpha_2 \sin(k x), \ \ \ \
f_3(x) = \alpha_3 \sin(k x),
\end{equation}
with $\alpha_j, k \in \mathbb{R}$ ($j=1,2,3$) and period $X = 2 \pi/k$.
This specific form enables one to observe some of the typical behavior in the most accessible and concrete way.
In fact, the Melnikov functions $\mathcal{R}_\mathrm{up/down}$ in \eqref{eq:FredholmConditionUp} and \eqref{eq:FredholmConditionDown} can be computed in closed form, as follows:
\begin{equation}
\label{Rupdownperex}
\begin{array}{rcl}
\mathcal{R}_\mathrm{up}(\phi)
& = & \int_{-\infty}^\infty \left[ \alpha_1 \cos(k x) u_\mathrm{up}(x;\phi) + \alpha_2 \sin(k x) u_\mathrm{up}'(x;\phi) + \alpha_3 \sin(k x) \right] u_\mathrm{up}'(x;\phi)\ dx
\\
& = & 16 \left( A + B \right) \sin(k \phi);
\\
\mathcal{R}_\mathrm{down}(\phi)
& = & \int_{-\infty}^\infty \left[ \alpha_1 \cos(k x) u_\mathrm{up}(x;\phi) + \alpha_2 \sin(k x) u_\mathrm{up}'(x;\phi) - \alpha_3 \sin(k x) \right] u_\mathrm{up}'(x;\phi)\ dx
\\
& = & 16 \left( A - B \right) \sin(k \phi),
\end{array}
\end{equation}
where $A$ and $B$ are defined by
\begin{equation}
\label{perexAB}
16 A  := \frac{k \pi}{3 \sinh( k \pi / \sqrt{2} )} \left[ - 3 \alpha_1 k + \alpha_2 (2 + k^2) \right], \; \; \ \ \
16 B := \frac{k \pi}{3 \sinh( k \pi / \sqrt{2} )} \left[ 3 \sqrt{2} \alpha_3 \right]. \end{equation}
Therefore, persistent one-front solutions exist at $k \phi=\ell \pi$ (to leading order) for all integers $\ell$.
Fig.~\ref{fig:1FrontStability} illustrates a particular case of this example, with $\alpha_1=1$, $\alpha_2=\alpha_3=0$, and $k=4\pi/15$.
(Also, for completeness, we note that $B=0$ in the topographic case, so that we recover $\mathcal{R}_\mathrm{up}(\phi) \equiv \mathcal{R}_\mathrm{up}(\phi)$, recall \eqref{eq:Rtopoupdown}).

\section{Stationary two-front solutions in the inhomogeneous Allen-Cahn equation}
\label{sec:2fronts}
As we have seen in the introduction,
the classical homogeneous Allen-Cahn PDE \eqref{eq:mainEquation} with $\varepsilon=0$
(and $x \in \mathbb{R}$) has four primary stable solutions: the homogeneous background states $u_\pm$ and the one-front solutions (heteroclinics) $u_\mathrm{up}$ and $u_\mathrm{down}$.
These are illustrated in the phase portrait for $\varepsilon = 0$, recall Fig.~\ref{fig:1FrontSolutions}.
All other solutions are attracted to one of these four solutions.
In particular, stationary multi-front solutions, with two or more fronts, are not present in \eqref{eq:mainEquation} with $\varepsilon = 0$.
\\ \\
However, with $0<\varepsilon \ll 1$,
the symmetry is broken, and system \eqref{eq:mainEquation} may possess stationary two-front solutions.
Geometrically, stationary two-front solutions will be found in the (transversal) intersections of the stable and unstable manifold of the same asymptotic state, {\it i.e.} intersections of $W^u(\mathcal{M}_-^\varepsilon)$ and $W^s(\mathcal{M}_-^\varepsilon)$, or of $W^u(\mathcal{M}_+^\varepsilon)$ and $W^s(\mathcal{M}_+^\varepsilon)$.
The fronts in a stationary two-front solution will lie sufficiently far away from each other ({\it i.e.} the difference in their phases will be sufficiently large),
since the invariant manifolds pass close by the other asymptotic state,
that is, solutions on them pass asymptotically close to a saddle orbit of the system.
\\ \\
In this section, we construct stationary two-front solutions using regular expansions of solutions on invariant stable and unstable manifolds, determine their linear stability, and apply the general theory to the same two cases (topographical and spatially periodic heterogeneities) studied in the previous section.
For clarity of notation, we focus on two-front solutions limiting on the asymptotic state $u_-^\varepsilon$; the case of two-front solutions limiting on $u_+^\varepsilon$ is similar.

\begin{remark}
Beyond the stationary two-fronts whose existence is established here, the spatially heterogeneous Allen-Cahn PDE \eqref{eq:mainEquation} also has dynamically-evolving two-front solutions.
We refer to the upcoming Secs.~\ref{sec:ODE2frontloc} and \ref{sec:ODE2frontper} -- and various simulations shown throughout this article -- for more on the dynamics of two-front patterns and on the role played by the stationary patterns.
\end{remark}

\begin{figure}
\centering
\begin{subfigure}[t]{0.45\textwidth}
	\includegraphics[width=\textwidth]{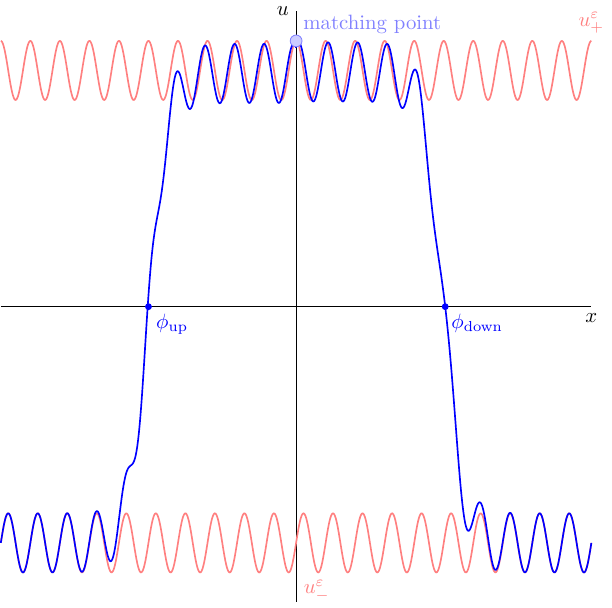}
\end{subfigure}
\caption{Sketch of the construction of a two-front solution.}
\label{fig:twoFrontSketch}
\end{figure}

\subsection{The existence of stationary two-front solutions}
\label{sec:2frontExistence}
As the name suggests, two-front solutions consist, to leading order, of two fronts: one which goes (approximately) from the state $u_-^\varepsilon$ to $u_+^\varepsilon$
and the other which goes (approximately) from $u_+^\varepsilon$ to $u_-^\varepsilon$.
We denote the phase/location of the first front by $\phi_\mathrm{up}$ and that of the second one by $\phi_\mathrm{down}$. We start with the following Ansatz for two-front solutions,
which is based on an expansion around the former jump on the domain $x < 0$ and an expansion around the latter jump on the domain $x > 0$:
\begin{equation}
	u_\mathrm{2f}(x;\phi_\mathrm{up},\phi_\mathrm{down}) =
	\begin{cases}
		u_\mathrm{up}(x;\phi_\mathrm{up}) + \varepsilon u_1^-(x;\phi_\mathrm{up}) + \mathcal{O}(\varepsilon^2) & \mbox{for $x < 0$};\\
		u_\mathrm{down}(x;\phi_\mathrm{down}) + \varepsilon u_1^+(x;\phi_\mathrm{down}) + \mathcal{O}(\varepsilon^2) & \mbox{for $x > 0$},
	\end{cases}
	\label{eq:twoFrontExpansion}
\end{equation}
where $u_\mathrm{up/down}$ are as defined in~\eqref{eq:heteroclinicSolution}.
See also Fig.~\ref{fig:twoFrontSketch} for a sketch of a stationary two-front solution and Sec.~3.3 for examples.
A key element of the existence proof is to impose continuity of the two components of $u_{2f}$ at the point $x=0$, as well
as continuity of their first derivatives. One notes that, by themselves, the terms which are {\it prima facie} $\mathcal{O}(1)$ cannot be joined continuously at $x = 0$ for any $\mathcal{O}(1)$ choice of $\phi_\mathrm{up}$ and $\phi_\mathrm{down}$, since otherwise two-front solutions would exist for the unperturbed Allen-Cahn equation.
Therefore, in order to insure continuity of $u_{2f}$ and its first derivative at $x=0$,
one needs to examine simultaneously the terms which are {\it prima facie} $\mathcal{O}(1)$ and those which are {\it prima facie} $\mathcal{O}(\varepsilon)$.
\\ \\
Continuity of the components of \eqref{eq:twoFrontExpansion} at $x=0$
and continuity of their first derivatives at $x = 0$ hold when
\begin{equation}
\label{eq:continuity-of-u+uprime}
\begin{array}{rcl}
u_\mathrm{up}(0,\phi_\mathrm{up}) + \varepsilon u_1^-(0,\phi_\mathrm{up})
+\mathcal{O}(\varepsilon^2)
& = & u_\mathrm{down}(0,\phi_\mathrm{down}) + \varepsilon u_1^+(0,\phi_\mathrm{down})
+\mathcal{O}(\varepsilon^2)
\\[1mm]
u_\mathrm{up}'(0,\phi_\mathrm{up}) + \varepsilon \partial_x u_1^-(0,\phi_\mathrm{up})
+\mathcal{O}(\varepsilon^2)
& = & u_\mathrm{down}'(0,\phi_\mathrm{down}) + \varepsilon \partial_x u_1^+(0,\phi_\mathrm{down})
+\mathcal{O}(\varepsilon^2).
\end{array}
\end{equation}
We set
\begin{equation}
\label{eq:phiupdown}
\phi_\mathrm{up}  = -\frac{|\log \varepsilon|}{2\sqrt{2}}  - \ell_\mathrm{up}, \; \; \ \ \ \ \
\phi_\mathrm{down}  =  \frac{|\log \varepsilon|}{2\sqrt{2}} + \ell_\mathrm{down},
\end{equation}
where $\ell_\mathrm{up}$ and $\ell_\mathrm{down}$ are $\mathcal{O}(1)$ constants to be determined (noting that only some choices will lead to two-front solutions), and we examine each of the terms in \eqref{eq:continuity-of-u+uprime}. By~\eqref{eq:phiupdown}, the terms which are {\it prima facie} $\mathcal{O}(1)$ in~\eqref{eq:continuity-of-u+uprime} are
\begin{equation}
\label{eq:asymptoticResultsLO12}
\begin{array}{rcl}
u_\mathrm{up,down}(0;\phi_\mathrm{up,down}) &=& 1 - 2 e^{\pm \sqrt{2} \phi_\mathrm{up,down}} + \mathcal{O}(\varepsilon) = 1 - 2  \sqrt{\varepsilon} e^{- \sqrt{2} \ell_\mathrm{up,down}} + \mathcal{O}(\varepsilon)
\\[1mm]
u'_\mathrm{up,down}(0;\phi_\mathrm{up,down}) &=& \pm 2\sqrt{2}  \sqrt{\varepsilon} e^{-\sqrt{2} \ell_\mathrm{up,down}} + \mathcal{O}(\varepsilon).
\end{array}
\end{equation}
Next, the terms which are {\it prima facie} $\mathcal{O}(\varepsilon)$ in~\eqref{eq:continuity-of-u+uprime} satisfy equation~\eqref{eq:existenceProblem1Front} with $u_h$ being given
by the corresponding one-front $u_\mathrm{up}$ or $u_\mathrm{down}$.
A comprehensive study of this inhomogeneous problem is given in Appendix~\ref{sec:orderEpsSystem}. In fact, by stipulating boundedness for $u_1^-$ in the limit $x \rightarrow -\infty$ and $u_1^+$ in the limit $x \rightarrow \infty$, we can determine the limit behavior for $u_1^\mp$ at $x = 0$. Using~\eqref{eq:u1limitbehavior} and \eqref{eq:phiupdown}, we obtain
\begin{equation}
\begin{array}{rclcl}
u_1^-(0;\phi_\mathrm{up}) & = & - \frac{\mathcal{R}_\mathrm{up}(\phi_\mathrm{up})}{8}  e^{-\sqrt{2} \phi_\mathrm{up}} + \mathcal{O}(1) & = & - \frac{1}{\sqrt{\varepsilon}} \frac{\mathcal{R}_\mathrm{up}(\phi_\mathrm{up})}{8} e^{\sqrt{2} \ell_\mathrm{up}} + \mathcal{O}(1);
\\[1mm]
\partial_x u_1^- (0;\phi_\mathrm{up}) & = & - \frac{\sqrt{2}}{\sqrt{\varepsilon}} \frac{\mathcal{R}_\mathrm{up}(\phi_\mathrm{up})}{8} e^{\sqrt{2}\ell_\mathrm{up}} + \mathcal{O}(1). & &
 \end{array}
 \end{equation}
Similarly, using~\eqref{eq:u1limitbehavior} and \eqref{eq:phiupdown}, one finds
\begin{equation}
\begin{array}{rclcl}
u_1^+(0;\phi_\mathrm{down}) & = &
\frac{\mathcal{R}_\mathrm{down}(\phi_\mathrm{down})}{8} e^{\sqrt{2} \phi_\mathrm{down}} + \mathcal{O}(1) & = & \frac{1}{\sqrt{\varepsilon}}
\frac{\mathcal{R}_\mathrm{down}(\phi_\mathrm{down})}{8} e^{\sqrt{2} \ell_\mathrm{down}} + \mathcal{O}(1),
\\[1mm]
\partial_x u_1^+ (0;\phi_\mathrm{down}) & = & - \frac{\sqrt{2}}{\sqrt{\varepsilon}} \frac{\mathcal{R}_\mathrm{down}(\phi_\mathrm{down})}{8} e^{\sqrt{2}\ell_\mathrm{down}} + \mathcal{O}(1), & &
\end{array}
\end{equation}
where $\mathcal{R}_\mathrm{up,down}$ are as defined in~\eqref{eq:FredholmConditionUp} and \eqref{eq:FredholmConditionDown}.
Hence, the continuity condition on $u$ implies
\begin{equation*}
1 - 2 \sqrt{\varepsilon} e^{-\sqrt{2}\ell_\mathrm{up}}
-\sqrt{\varepsilon} \frac{\mathcal{R}_\mathrm{up} (\phi_\mathrm{up})}{8} e^{\sqrt{2}\ell_\mathrm{up}}
+ \mathcal{O}(\varepsilon)
= 1 - 2 \sqrt{\varepsilon} e^{-\sqrt{2}\ell_\mathrm{down}}
+\sqrt{\varepsilon} \frac{\mathcal{R}_\mathrm{down} (\phi_\mathrm{down})}{8} e^{\sqrt{2}\ell_\mathrm{down}}
+ \mathcal{O}(\varepsilon).
\end{equation*}
Putting this together with a similar condition for the continuity of the derivatives at $x = 0$ (per \eqref{eq:continuity-of-u+uprime}), one finds the following conditions at
$\mathcal{O}(\sqrt{\varepsilon})$:
\begin{equation}
\begin{array}{rcl}
2 e^{-\sqrt{2} \ell_\mathrm{up}} + \frac{\mathcal{R}_\mathrm{up}(\phi_\mathrm{up})}{8} e^{\sqrt{2} \ell_\mathrm{up}} & = & 2 e^{-\sqrt{2} \ell_\mathrm{down}} - \frac{\mathcal{R}_\mathrm{down}(\phi_\mathrm{down})}{8} e^{\sqrt{2} \ell_\mathrm{down}}
\\[1mm]
2 e^{-\sqrt{2} \ell_\mathrm{up}} - \frac{\mathcal{R}_\mathrm{up}(\phi_\mathrm{up})}{8} e^{\sqrt{2} \ell_\mathrm{up}} & = & -2 e^{-\sqrt{2} \ell_\mathrm{down}} - \frac{\mathcal{R}_\mathrm{down}(\phi_\mathrm{down})}{8} e^{\sqrt{2} \ell_\mathrm{down}}
\end{array}
\end{equation}
Finally, addition and subtraction of these two conditions and recalling~\eqref{eq:phiupdown} allows us to rearrange them into the following condensed and insightful way:
\begin{equation}
\mathcal{R}_\mathrm{up}(\phi_\mathrm{up}) = -\mathcal{R}_\mathrm{down}(\phi_\mathrm{down})
= 16 e^{-\sqrt{2}\left[ \ell_\mathrm{up} + \ell_\mathrm{down} \right]}
= \frac{16}{\varepsilon} e^{-\sqrt{2}\left[ \phi_\mathrm{down} - \phi_\mathrm{up}\right]}.
\label{eq:twoFrontCondition}
\end{equation}
Thus, we have shown that the PDE \eqref{eq:mainEquation} possesses stationary two-front solutions
for all parameter values that satisfy the conditions \eqref{eq:twoFrontCondition} on the Melnikov functions, and the fronts are located near $\phi_\mathrm{up,down}$.

\begin{remark}
	In the limit $\ell_\mathrm{up} + \ell_\mathrm{down} \rightarrow \infty$ ({\it i.e.,} in the limit in which the two fronts are infinitely far apart from each other), condition~\eqref{eq:twoFrontCondition} reduces to $\mathcal{R}_\mathrm{up}(\phi_\mathrm{up}) = 0$ and $\mathcal{R}_\mathrm{down}(\phi_\mathrm{down}) = 0$. This is the same as the existence condition for one-front solutions, recall \eqref{eq:FredholmCondition1Front}. Thus, in this limit, interactions between fronts are only through their exponentially small tails, and their locations/phases are determined to leading order by the standard phase selection criterion.
\end{remark}

\subsection{Linear stability of stationary two-front solutions} \label{sec:2frontStability}
To determine the linear stability of stationary two-front solutions, we observe that eigenvalues in the point spectrum lie close to those of the one-front solutions. Hence, for sufficiently small $\varepsilon$, only the $\mathcal{O}(\varepsilon)$ eigenvalues (which are near $\lambda_0=0$) can possibly lie in the right-half plane, {\it i.e.,} can possibly make the two-front solution unstable. In general, two such eigenvalues are expected here for two-front solutions.
We determine the locations of these small eigenvalues using regular perturbation theory and a combination of techniques from the existence problem in Sec.~\ref{sec:2frontExistence} and the stability analysis of one-front solutions in Sec.~\ref{sec:oneFrontStability}. Specifically, we use regular perturbation theory with one expansion on $x < 0$ and another on $x > 0$, and we impose continuity of the eigenfunctions and their derivatives at $x = 0$.
Furthermore, we observe that, for these small eigenvalues, the solutions to the eigenvalue problem~\eqref{eq:eigenvalueProblem} in both regions should, to leading order, be multiples of the derivative of the (exact) two-front solutions to~\eqref{eq:mainEquation}. However, the multiples can be different in the two regions.
\\
\\
Specifically,
we set
\begin{align}
	\bar{u}(x;\phi_\mathrm{up},\phi_\mathrm{down})
& = \begin{cases}
		\gamma_-\ u_\mathrm{2f}'(x;\phi_\mathrm{up},\phi_\mathrm{down}) + \varepsilon \bar{u}_1^-(x;\phi_\mathrm{up},\phi_\mathrm{down}) + \mathcal{O}(\varepsilon^2) & \mbox{ for $x < 0$};\\
		\gamma_+\ u_\mathrm{2f}'(x;\phi_\mathrm{up},\phi_\mathrm{down}) + \varepsilon \bar{u}_1^+(x;\phi_\mathrm{up},\phi_\mathrm{down}) + \mathcal{O}(\varepsilon^2) & \mbox{ for $x > 0$};
	\end{cases} \nonumber \\
	\lambda & = \varepsilon \bar{\lambda} + h.o.t. \label{eq:defbarlamb2fronts}
\end{align}
where $\gamma_\pm$ are constants. Because of the linearity of the eigenvalue problem~\eqref{eq:eigenvalueProblem}, we are free to choose one of these constants, while the other needs to be determined via continuity. However, we should beware of the possibility that one of the constants is zero, and therefore we wait to specify the value of one of these constants.
\\ \\
Using a similar approach as in Sec.~\ref{sec:oneFrontStability}, we combine equation~\eqref{eq:mainEquation} and the $x$-derivative of~\eqref{eq:mainEquation} to find to leading order
\begin{equation}
	\mathcal{L}_0 \bar{u}_1^\mp = \gamma_\mp \left[ \pm \bar{\lambda} u_\mathrm{up}' + F_x(\pm u_\mathrm{up}, \pm u_\mathrm{up}', x) \right].
\end{equation}
Stipulating boundedness for $\bar{u}_1^-$ (respectively $\bar{u}_1^+$) in the limit $x \rightarrow -\infty$ (respectively in the limit $x \rightarrow +\infty$), we find $\bar{u}^\mp_1$ using formulae from Appendix~\ref{sec:orderEpsSystem} (after replacing the inhomogeneous terms). As in the existence problem, only the leading order behavior at the point $x = 0$ is important.
Here, also by Appendix~\ref{sec:orderEpsSystem}, straightforward computations yield the leading order $\mathcal{O}(\sqrt{\varepsilon})$ terms at the point $x=0$.
In particular, imposing continuity of the two components and their derivatives at $x = 0$, we find the following conditions:
\[
\begin{array}{llcl}
& \sqrt{2}\ \gamma_-  \left[2 e^{-\sqrt{2} \ell_\mathrm{up}}-\frac{\mathcal{R}_\mathrm{up}(\phi_\mathrm{up})}{8}e^{\sqrt{2}\ell_\mathrm{up}}\right] &+& \frac{\gamma_-}{8} \left[ \bar{\lambda}\ \|u_\mathrm{up}'\|_2^2 + \mathcal{R}'_\mathrm{up}(\phi_\mathrm{up})\right] e^{\sqrt{2} \ell_\mathrm{up}}
\\[1mm]
= & \sqrt{2}\ \gamma_+  \left[-2 e^{-\sqrt{2}\ell_\mathrm{down}} - \frac{\mathcal{R}_\mathrm{down}(\phi_\mathrm{down})}{8} e^{\sqrt{2}\ell_\mathrm{down}}\right] &-& \frac{\gamma_+}{8} \left[ \bar{\lambda}\ \|u_\mathrm{up}'\|_2^2 + \mathcal{R}_\mathrm{down}'(\phi_\mathrm{down})\right] e^{\sqrt{2} \ell_\mathrm{down}},
\end{array}
\]
and
\[
\begin{array}{llcl}
& -\sqrt{2}\ \gamma_- \left[ 2 e^{-\sqrt{2}\ell_\mathrm{up}} + \frac{\mathcal{R}_\mathrm{up}(\phi_\mathrm{up})}{8}e^{\sqrt{2} \ell_\mathrm{up}}\right] &+& \frac{\gamma_-}{8} \left[ \bar{\lambda}\ \|u_\mathrm{up}'\|_2^2 + \mathcal{R}'_\mathrm{up}(\phi_\mathrm{up})\right] e^{\sqrt{2} \ell_\mathrm{up}}
\\[1mm]
= & -\sqrt{2}\ \gamma_+ \left[ 2 e^{-\sqrt{2}\ell_\mathrm{down}} - \frac{\mathcal{R}_\mathrm{down}(\phi_\mathrm{down})}{8} e^{\sqrt{2}\ell_\mathrm{down}}\right] &+& \frac{\gamma_+}{8} \left[ \bar{\lambda}\ \|u_\mathrm{up}'\|_2^2 + \mathcal{R}_\mathrm{down}'(\phi_\mathrm{down})\right] e^{\sqrt{2} \ell_\mathrm{down}}.
\end{array}
\]
Hence, using \eqref{eq:twoFrontCondition} which are the conditions for existence of two-front solutions, and adding and subtracting the above two equations, we find that these conditions reduce to
\begin{equation}
	\gamma_- \left[ \bar{\lambda}\ \|u_\mathrm{up}'\|_2^2 + \mathcal{R}'_\mathrm{up}(\phi_\mathrm{up})\right] = - \gamma_+ \left[ \bar{\lambda}\ \|u_\mathrm{up}'\|_2^2 + \mathcal{R}'_\mathrm{down}(\phi_\mathrm{down})\right] = 16 \sqrt{2} (\gamma_- - \gamma_+) e^{-\sqrt{2}( \ell_\mathrm{up} + \ell_\mathrm{down} ) }. \label{eq:twoFrontStabilityCondition}
\end{equation}
These are the conditions on the constants $\gamma_-$ and $\gamma_+$ that must be satisfied for nontrivial eigenfunctions to exist.
\\
\\
Analysis of these conditions shows that there are no solutions of~\eqref{eq:twoFrontStabilityCondition} with $\gamma_- = 0$ (except for in the limit $\ell_\mathrm{up}+\ell_\mathrm{down} \rightarrow \infty$; see also Remark~\ref{remark:twoFrontStabilityConditionLimit}).
Hence, without loss of generality, we may scale the eigenfunctions by setting $\gamma_- = 1$. Then, the system of linear equations has two solutions given by
\begin{align}
	(\gamma_+)_{1,2} & = \frac{-\sqrt{2}}{64}\ e^{\sqrt{2}  \left( \ell_\mathrm{up} + \ell_\mathrm{down} \right) } \left[ \mathcal{R}_\mathrm{up}'(\phi_\mathrm{up})-\mathcal{R}_\mathrm{down}'(\phi_\mathrm{down}) \pm \ \sqrt{ D\ } \right]; \nonumber \\
	\bar{\lambda}_{1,2} & = \frac{1}{\|u_\mathrm{up}'\|_2^2} \left[ 16 \sqrt{2}\ e^{-\sqrt{2}  \left( \ell_\mathrm{up} + \ell_\mathrm{down} \right)} - \frac{\mathcal{R}_\mathrm{up}'(\phi_\mathrm{up})+\mathcal{R}_\mathrm{down}'(\phi_\mathrm{down})}{2} \mp \frac{\sqrt{ D\ }}{2} \right]; \label{eq:twoFrontEigenvalues}
\intertext{where}
D & := \left( \mathcal{R}_\mathrm{up}'(\phi_\mathrm{up})-\mathcal{R}_\mathrm{down}'(\phi_\mathrm{down})  \right)^2 + 2048\ e^{-2 \sqrt{2}  \left( \ell_\mathrm{up} + \ell_\mathrm{down} \right) }.\label{eq:eigenvalues2Front}
\end{align}
This completes the derivation of the eigenfunctions and eigenvalues.
\\ \\
We summarize the preceding analysis in a theorem (that can again be proven by standard methods) and illustrate our findings in two examples in the upcoming  subsection.
\begin{theorem}
\label{th:2frontExStab}
Consider the weakly heterogeneous Allen-Cahn equation \eqref{eq:mainEquation} with bounded inhomogeneity $F(U,U_x,x)$. Let the Melnikov functions $\mathcal{R}_{\rm up}(\phi)$ and $\mathcal{R}_{\rm down}(\phi)$ be given by \eqref{eq:FredholmConditionUp} and \eqref{eq:FredholmConditionDown} assume that $(\phi_{\rm up}, \phi_{\rm down})$, with $\phi_{\rm up} < \phi_{\rm down}$, is a non-degenerate solution of \eqref{eq:twoFrontCondition} (with $\ell_{\rm up, down}$ as in \eqref{eq:phiupdown}). Then, for $0 < \varepsilon \ll 1$ sufficiently small, \eqref{eq:mainEquation} has a stationary two-front solution $U_{\rm 2-f}(x)$ that is biasymptotic to $u^\varepsilon_-(x) = - 1 + \mathcal{O}(\varepsilon)$ \eqref{eq:uMinusPlus} and that is $\mathcal{O}(\sqrt{\varepsilon})$ close to $u_{\rm up}(x;\phi_{\rm up})  + u_{\rm down}(x;\phi_{\rm down}) - 1$ \eqref{eq:heteroclinicSolution} (i.e., ${\rm sup}_{x \in \mathbb{R}} |U_{\rm 2-f}(x) - (u_{\rm up}(x;\phi_{\rm up})  + u_{\rm down}(x;\phi_{\rm down}) - 1)| =  \mathcal{O}(\sqrt{\varepsilon})$). The two-front pattern $U_{\rm 2-f}(x)$ is stable if $\lambda_{1,2} = \varepsilon \bar{\lambda}_{1,2} < 0$ \eqref{eq:defbarlamb2fronts}, \eqref{eq:twoFrontEigenvalues}.
\end{theorem}
\begin{remark}
\label{rem:symmR'}
In the special case in which
$\mathcal{R}_\mathrm{up}'(\phi_\mathrm{up}) = \mathcal{R}_\mathrm{down}'(\phi_\mathrm{down})$,
one of the eigenfunctions is symmetric, and the other anti-symmetric. Specifically, in this case, we have
\begin{equation}
\label{eq:al12la12symmR'}
\begin{array}{lcl}
(\gamma_+)_1 = + 1, & &
\bar{\lambda}_1 = - \frac{\mathcal{R}_\mathrm{up}'(\phi_\mathrm{up})}{\|u_\mathrm{up}'\|_2^2},
\\[1mm]
(\gamma_+)_2 = - 1, & &
\bar{\lambda}_2 = - \frac{ \mathcal{R}_\mathrm{up}'(\phi_\mathrm{up}) - 32 \sqrt{2} e^{-\sqrt{2}\left(\ell_\mathrm{up}+\ell_\mathrm{down}\right)}}{\|u_\mathrm{up}'\|_2^2}.
\end{array}
\end{equation}

\end{remark}

\begin{remark}\label{remark:twoFrontStabilityConditionLimit}
In the limit $\ell_\mathrm{up} + \ell_\mathrm{down} \rightarrow \infty$ ({\it i.e.,} in the limit in which the two fronts are infinitely far apart from each other), conditions~\eqref{eq:twoFrontStabilityCondition} reduce to
$\gamma_-\left[ \bar{\lambda} \|u_\mathrm{up}'\|_2^2 + \mathcal{R}_\mathrm{up}'(\phi_\mathrm{up})\right] = \gamma_+ \left[ \bar{\lambda} \|u_\mathrm{up}'\|_2^2 + \mathcal{R}_\mathrm{down}'(\phi_\mathrm{down})\right] = 0.$
Thus, in this limit, there are two eigenfunctions: one corresponds to the one-front eigenfunction of the left front (and to leading order nothing happens at the other front, i.e. $\gamma_+ = 0$), and the other to the one-front eigenvalue of the right front (and $\gamma_- = 0$). Moreover, the eigenvalues also reduce to those found in the stability analysis of the one-front solutions~\eqref{eq:eigenvalue1Front}.
Furthermore, \eqref{eq:twoFrontEigenvalues} reduces to these values in the limit $\ell_\mathrm{up} + \ell_\mathrm{down} \rightarrow \infty$, as expected.
\end{remark}

\subsection{Stationary two-fronts in topographically driven systems and spatially periodic systems}
\label{s:pertop2fronts}

\subsubsection{Two-front solutions solutions arising from weak topographical effects}
\label{sec:top2fronts}
In Subsec.~\ref{sec:top1fronts}, we derived that $\mathcal{R}_\mathrm{up}(\phi) = \mathcal{R}_\mathrm{down}(\phi)$ in case \eqref{eq:Ftopography} in which driving term $F(U,V,x)$ models the impact of a topography $H(x)$.
As a consequence, the existence condition (\ref{eq:twoFrontCondition}) reduces to
\begin{equation}
\int_{-\infty}^\infty H'(y+\phi_\mathrm{up}) W_h(y) dy =
- \int_{-\infty}^\infty H'(y+\phi_\mathrm{down}) W_h(y) dy =  \frac{16}{\varepsilon} e^{-\sqrt{2}\left[\phi_\mathrm{down} - \phi_\mathrm{up} \right]},
\label{eq:Ex2TwoFrontCondition}
\end{equation}
((\ref{eq:Ex2R}), (\ref{eq:defWh}), (\ref{eq:phiupdown})).
If we again also assume that $H(x) = H(-x)$ and additionally restrict our search for two-front patterns to symmetric patterns, {\it i.e.}, assume that $\phi = \phi_\mathrm{down} = -\phi_\mathrm{up} > 0$, then we have
\[
\mathcal{R}(\phi) = \int_{-\infty}^\infty H'(y+\phi) W_h(y) dy = - \int_{-\infty}^\infty H'(-y-\phi) W_h(y) dy = - \int_{-\infty}^\infty H'(z-\phi) W_h(z) dz.
\]
As a result, (\ref{eq:Ex2TwoFrontCondition}) simplifies to
\begin{equation}
\mathcal{R}(\phi) = - \frac{16}{\varepsilon} e^{-2\sqrt{2}\phi}.
\label{eq:Ex2TwoFrontCondSymm}
\end{equation}
In the case of an $X$-periodic topography $H(x)$, with also $\mathcal{R}(\phi+X) = \mathcal{R}(\phi)$ \eqref{eq:Rtopoper}, it immediately follows that there is a critical $\bar{\phi}_\mathrm{min} = \mathcal{O}\left(\frac{|\log \varepsilon|}{2\sqrt{2}}\right)$ such that (\ref{eq:Ex2TwoFrontCondSymm}) has countably many zeroes $\phi_{\ast,j} \geq \bar{\phi}_\mathrm{min}$ (see also Theorem \ref{th:Nfrontsper} for a generalization to $N$-front patterns). For the unimodal hill $H_{\rm exp}(x;\mu)$ \eqref{eq:defHuni-alg}, we obtain for $\psi = \sqrt2 \phi \gg 1$ by the leading order behavior of (\ref{eq:Rsolhill1}) -- see also Fig. \ref{fig:SolHillPitchfork}(a) -- that for $\mu = 1$ there is a symmetric two-front solution at $-\phi_\mathrm{up}=\phi_\mathrm{down}=\phi^\ast(1) \gg 1$ (the $1$ referring to the $\mu$-value) with $\phi^\ast(1)$ the (unique) solution of
\[
\varepsilon \sqrt{2} \phi e^{\sqrt{2} \phi} = 1,
\]
so that $\phi_{\ast} \in (\frac12 \nu  \sqrt{2} |\log \varepsilon|, \frac12 \sqrt{2} |\log \varepsilon|)$ for any $\nu \in (0,1)$.
In fact, it will follow from the more general asymptotic analysis of Sec.~\ref{sec:ODE2frontloc} that $\phi^\ast(\mu) > 0$ exists for any $\mu > 0$ (see Fig.~\ref{fig:SolHillPitchfork}(b) for $\mu = \frac12$). There is a symmetric two-front pattern for all unimodal hill topographies $H_{\rm exp}(x;\mu)$. Moreover, it follows by the same arguments that there cannot be symmetric two-front patterns in the case of the unimodal valleys $V_{\rm exp}(x; 1) = - H_{\rm exp}(x,1)$.
\\ \\
To determine the stability of the two-front patterns, we again first consider the more general setting of \eqref{eq:Ftopography}. Again assume that $H(x)$ is symmetric and that the two-front is symmetric, i.e. $\phi = \phi_\mathrm{down} = - \phi_\mathrm{up} > 0$, we find
\[
\mathcal{R}'_\mathrm{up}(\phi_\mathrm{up}) = \int_{-\infty}^\infty H''(y-\phi) W_h(y) dy
= \int_{-\infty}^\infty H''(-z - \phi) W_h(-z) dz = \int_{-\infty}^\infty H''(z + \phi) W_h(z) dz =
\mathcal{R}'_\mathrm{down}(\phi_\mathrm{down})
\]
Thus, this is the situation of the special case discussed in Remark \ref{rem:symmR'}.
It follows by (\ref{eq:twoFrontCondition}), (\ref{eq:al12la12symmR'}), and (\ref{eq:Ex2TwoFrontCondSymm}) that the stability of symmetric two-front patterns at $\phi_\mathrm{down} = - \phi_\mathrm{up} = \phi^\ast$ is determined by
\begin{equation}
\bar{\lambda}_1 = - \frac{\mathcal{R}'(\phi^\ast)}{\|u_h'\|_2^2},
\; \; \ \ \ \ \
\bar{\lambda}_2 = - \frac{\mathcal{R}'(\phi^\ast) + 2\sqrt{2} \mathcal{R}(\phi^\ast)}{\|u_h'\|_2^2}.    \label{eq:la12top2front}
\end{equation}
In the case that $H(x)$ is localized with monotonically decaying tails -- as with the explicit hill and valley examples discussed above -- it will be shown in Sec.~\ref{sec:ODElocalized} that both $\mathcal{R}(\phi)$ and $\mathcal{R}'(\phi) \to 0$ as $\phi \to \infty$.
In fact, it follows in this case that if $\mathcal{R}(\phi) > 0$, respectively $\mathcal{R}(\phi) < 0$, for $\phi \gg 1 $, then $\mathcal{R}(\phi)$ decreases to 0, i.e. $\mathcal{R}'(\phi) < 0$, resp. $\mathcal{R}(\phi)$ increases with $\mathcal{R}'(\phi) > 0$ (see Fig. \ref{fig:SolHillPitchfork}).
Therefore, since $\mathcal{R}(\phi^\ast) < 0$ by (\ref{eq:Ex2TwoFrontCondSymm}), it follows that $\bar{\lambda}_1 < 0$ and that the symmetric two-front pattern is stable if also $\bar{\lambda}_2<0$, that is by (\ref{eq:la12top2front}), if also
\begin{equation}
\label{symm2frontstab}
\mathcal{R}'(\phi^\ast) + 2\sqrt{2} \mathcal{R}(\phi^\ast) > 0.
\end{equation}
Let us analyze this condition.
First, we recall that $\mathcal{R}(\phi^\ast) < 0 < \mathcal{R}'(\phi^\ast)$.
Hence, this condition may hold for certain localized, symmetric topographies $H(x)$. However, for the specific unimodal hill $H_{\rm exp}(x;\mu)$ \eqref{eq:defHuni-alg}, we find explicitly from (\ref{eq:Rsolhill1}) for $\mu = 1$ that
\[
\mathcal{R}(\phi^\ast) = - 16 \sqrt{2} \phi^\ast e^{-\sqrt{2} \phi^\ast} + \; {\rm h.o.t}, \; \;
\mathcal{R}'(\phi^\ast) = 32 \phi^\ast e^{-\sqrt{2} \phi^\ast} + \; {\rm h.o.t}.
\]
This implies that $\bar{\lambda}_2 > 0$ and that the two-front pattern is unstable.
Moreover, the $\mu=1$ is not exceptional.
The analysis in Sec.~\ref{sec:ODE2frontloc} will show that the symmetric two-front pattern associated to the $H_{\rm exp}(x;\mu)$ topography is unstable for all $0 < \mu < 1$, as well (see Theorem \ref{th:Nfrontslocexp}).

\subsubsection{Two-front solutions in systems with spatially-periodic heterogeneities}
\label{sec:per2fronts}
In this subsection, we consider the explicit case introduced in Sec.~\ref{sec:per1fronts}:
\eqref{eq:canonicalExample} with $f_j(x)$ ($j=1,2,3$) as in \eqref{eq:fgh-example1}.
Recalling the Melnikov functions
$\mathcal{R}_\mathrm{up}(\phi)$ and $\mathcal{R}_\mathrm{down}(\phi)$ \eqref{Rupdownperex} with $A$ and $B$ given by \eqref{perexAB}, we observe that \eqref{eq:twoFrontCondition} provides an explicit condition on the existence of two-front patterns in this case,
\begin{equation}
\label{eq:exper2front-exist}
\left( A + B \right) \sin(k \phi_\mathrm{up}) =
-\left( A - B \right) \sin(k \phi_\mathrm{down}) =
\frac{1}{\varepsilon} e^{-\sqrt{2}\left[ \phi_\mathrm{down} - \phi_\mathrm{up}\right]} + \mathcal{O}(\varepsilon).
\end{equation}
Here, the correction term represents the higher order terms so far neglected in this section (see the discussions in Secs.~\ref{sec:ODElocalized} and \ref{sec:Discussion}).
Assuming the two fronts are logarithmically far apart, {\it i.e.}, assuming that there is an $\rho > 0$ such that $0<\sqrt{2}(\phi_\mathrm{down} - \phi_\mathrm{up})/|\log \varepsilon| \to \rho$ as $\varepsilon \to 0$ (see Theorem \ref{th:Nfrontsper}), it follows that
\[
\left(A + B \right) \sin(k \phi_\mathrm{up}) =
-\left(A - B \right) \sin(k \phi_\mathrm{down}) =
\varepsilon^{\rho-1} + \mathcal{O}(\varepsilon),
\]
which implies that there cannot be any stationary two fronts for $0<\rho <1$. On the other hand, we also see that if $\rho > 1$ (and $A\pm B \ne 0$), there is a stationary two-front $\mathcal{O}(\varepsilon^{\rho-1}, \varepsilon)$ close to any pair $(n_\mathrm{up}(\varepsilon) \pi/k, n_\mathrm{down}(\varepsilon)\pi/k)$ with $(n_\mathrm{up}(\varepsilon),n_\mathrm{down}(\varepsilon)) \in \mathbb{Z}^2$ such that $0 < n_\mathrm{up}(\varepsilon) - n_\mathrm{down}(\varepsilon) \to k |\log \varepsilon|/\pi \sqrt{2}$ as $\varepsilon \to 0$. We refer to Theorem \ref{th:Nfrontsper} for a similar statement on the more general case of $N$-front patterns in \eqref{eq:mainEquation} with a periodic inhomogeneity (restricted to the topographic subcase \eqref{eq:Ftopography}). The leading order expressions for the two eigenvalues \eqref{eq:twoFrontEigenvalues} that determine the stability of these two-front patterns ({\it i.e.}, with $\rho > 1$) also simplify considerably to
\[
\lambda_{1,2} = -\frac{8 k \varepsilon}{\|u_\mathrm{up}'\|_2^2} \left[
(A + B) \cos(n_\mathrm{up} \pi) +
(A - B) \cos(n_\mathrm{down} \pi) \pm \sqrt{\left((A + B) \cos(n_\mathrm{up} \pi) - (A - B) \cos(n_\mathrm{down} \pi)\right)^2}\right]
\]
(\eqref{eq:defbarlamb2fronts}, \eqref{eq:eigenvalues2Front}), so that
\begin{equation}
\label{eq:exper2front-lambda12}
\lambda_{1} = -\frac{16 k (A+B) \varepsilon}{\|u_\mathrm{up}'\|_2^2} \cos(n_\mathrm{up} \pi), \; \; \ \ \ \
\lambda_{2} = -\frac{16 k (A-B) \varepsilon}{\|u_\mathrm{up}'\|_2^2}  \cos(n_\mathrm{down} \pi)
\end{equation}
(to leading order in $\varepsilon$). We conclude that a two-front pattern with $\rho > 1$ is stable if both $n_\mathrm{up}$ and $n_\mathrm{down}$ are even, {\it i.e.}, that one out of every quartet of these two-front patterns is stable (see again Theorem \ref{th:Nfrontsper}).
\\ \\
In Sec.~\ref{sec:ODE2frontper} -- and again in Sec.~\ref{sec:geom-int} -- we consider the many bifurcations associated to increasing $\rho$ from below $\rho <1$, for which there are no stationary two-fronts, to $\rho >1$, for which there are countably infinite.
In Sec.~\ref{sec:ODE2frontper}, this will be done in the setting of the (reduced) two-dimensional ODE that governs the dynamics of the positions $(\phi_1(t),\phi_2(t))$ of an evolving two-front pattern solution of \eqref{eq:mainEquation} (that will be derived in Sec.~\ref{sec:interactiondynamics}).
It will also follow from this analysis that any two-front pattern $(\phi_1(t),\phi_2(t))$ will be attracted to a nearby stationary two-front pattern, {\it i.e.}, that the limiting behavior of all possible $(\phi_1(t),\phi_2(t))$-front patterns is determined by the family of stable stationary two-front patterns (see Theorem \ref{th:Nfrontsperattr} for the $N$-front case).

\section{The (slow) dynamics of interacting fronts}
\label{sec:interactiondynamics}
In the homogeneous ($\varepsilon = 0$) case, the dynamics of multi-front solutions are similar to that for the two-front solutions: fronts attract each other and consequently annihilate in pairs.
The only possible stable end states are the trivial states $u_-$ and $u_+$ or the (stationary) one-fronts $u_{\rm up}$ and $u_{\rm down}$.
By contrast, for $\varepsilon \neq 0$, there are many stable and unstable stationary multi-front patterns possible. These stationary $N$-front solutions are given to leading order by a superposition of $N$ fronts, where now each subsequent front should be of different type. Both the construction and the stability of these stationary multi-fronts can be established along the lines set out in Sec.~\ref{sec:2fronts} for two-front patterns. We refrain from following this approach here (however, see Sec.~\ref{sec:ODENfrontloc}). The richness of the types of stable multi-front patterns exhibited by \eqref{eq:mainEquation} is indicated for example by the attracting (and thus stable) stationary multi-front patterns presented in Figs.~\ref{fig:NumericsMultiFront-Intro}, as obtained by numerical simulations.
\\
\\
In the present section, we take a different approach. We have seen that (stationary) front solutions have small, $\mathcal{O}(\varepsilon)$ eigenvalues in their associated point spectrum. These small eigenvalues originate from the translational invariance of $1$-front solutions to \eqref{eq:mainEquation} for $\varepsilon = 0$ (which leads to an eigenvalue at $\lambda = 0$ with $u_h'(x)$ -- the derivative of the unperturbed one-front pattern -- as corresponding eigenfunction).
When multiple fronts and/or heterogeneous terms (i.e. $\varepsilon > 0$) are present, the small eigenvalues in the point spectrum have eigenfunctions that are to leading order determined by these derivatives, see Sec.~\ref{sec:2frontStability}.
The dynamics induced by this so-called small spectrum can be captured by projecting the full state space onto a finite-dimensional space spanned by the associated eigenspaces.
Doing so for $N$-front solutions leads to an $N$-th order system of ODEs that describes the (slow) evolution of multi-front solutions to \eqref{eq:mainEquation}.
\\
\\
Indeed, the main results of this section is the derivation of the ODE system \eqref{eq:NFrontODE}.
The ODE system \eqref{eq:NFrontODE} governs the slow migratory movements of the multi-front solutions by tracking changes in the locations $\phi_i(t)$, $i=1,\ldots,N$, (or phase shifts) of the fronts to leading order in some small quantity -- typically associated to the inverse of the distance between two successive fronts. The validity of this procedure, thus both of the low-dimensional approximation and of the ODE governing the low-dimensional dynamics of the {\it weakly} interacting fronts, has a rigorous basis and can be established by following the methods developed in the literature, see for instance \cite{ei2002pulse,promislow2002renormalization}.
\\ \\
In this article, we refrain from going into the details of the validity analysis (see also Remark \ref{rem:Fbounded}).
Instead, we focus on the derivation of the governing ODEs describing the interactions of fronts to leading order and on their implications for the dynamics of multi-front patterns.
More specifically, we derive the ODEs governing the interactions in $N$-front patterns in upcoming Sec.~\ref{sec:ODEderivation}.
Then, in the next sections of this article, we study various aspects of $N$-front interactions in the two explicit cases we already considered in the previous sections.
Namely, in Sec.~\ref{sec:ODElocalized}, we study the impacts of localized topographies, and in Sec.~\ref{sec:ODEperiodic} we study the impacts of spatially periodic driving terms.
Moreover, we will return to the above-mentioned stationary multi-front patterns, since they reappear as critical points of the interaction ODEs.
\begin{remark}
\label{rem:Fbounded}
The assumption that the inhomogeneous term $F(U,V,x)$ in \eqref{eq:mainEquation} is bounded for all $x$ is essential to directly apply the methods of \cite{ei2002pulse,promislow2002renormalization} to establish the finite dimensional character of front interactions and the governing ODEs.
In Sec.~\ref{sec:Discussion}, we will discuss observations indicating that the derived  interaction equations may also be valid beyond this setting, see also Fig.~\ref{fig:5Fronts-Intro}(c) in which $F(U,V,x)$ is given by \eqref{eq:Ftopography} with unbounded $H'(x)$. Thus, although we refrain from going into the details of the standard case (in which $F(U,V,x)$ is bounded), we motivate at several places in this work  that setting up a detailed validity analysis that goes beyond that of the standard case would be an interesting and relevant mathematical enterprise.
\end{remark}

\subsection{Derivation of the governing ODEs}
\label{sec:ODEderivation}

\begin{figure}
\centering
\includegraphics[width=\textwidth]{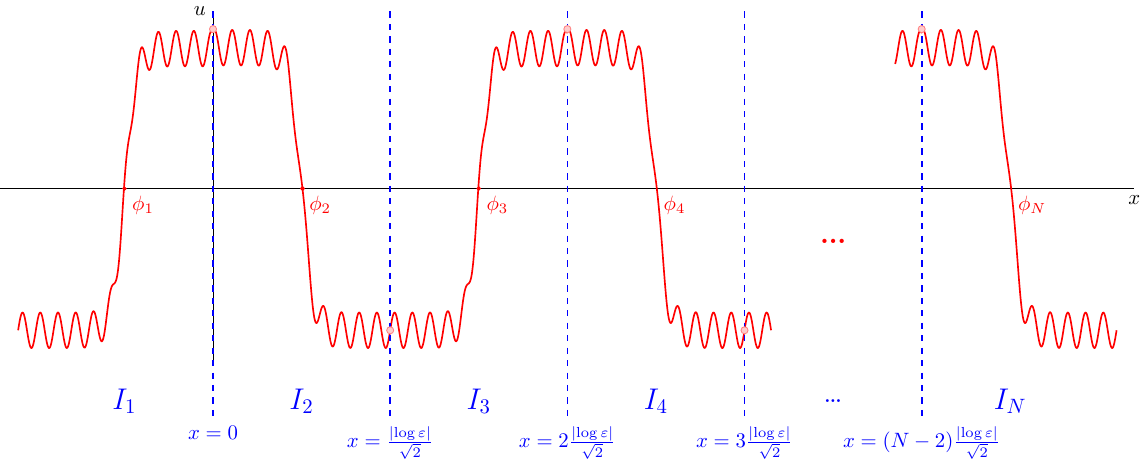}

\caption{Sketch of am $N$-front pattern anchored by \eqref{eq:anchor} with $p_j = j - \frac{3}{2}$, a standard choice. Here, $\phi_1, \ldots, \phi_N$ denote the locations -- or phase shifts -- of the fronts. The domain ($\mathbb{R}$) is divided into the regions $I_1, \ldots, I_N$ as given in \eqref{eq:defIjs} in and solutions are matched at the boundaries of these regions.}
\label{fig:NFrontSketch}
\end{figure}

We consider an $N$-front pattern with front locations (or phase shifts) $\phi_1(t) < \phi_2(t)  < \ldots < \phi_N(t)$, see Fig.~\ref{fig:NFrontSketch}.
We anchor the asymptotic analysis by introducing $p_j$ ($j = 1, ...,N$), with $p_j < p_{j+1}$, and $\ell_j(t)$ similar as in \eqref{eq:phiupdown}, using
\begin{equation}
	\phi_j(t) = p_j \frac{|\log \varepsilon|}{\sqrt{2}} + \ell_j(t).
\label{eq:anchor}
\end{equation}
The exact values of the $p_j$'s will not show up explicitly in the outcome of the analysis. However, the $p_j$'s determine the relative positions of the fronts and the driving terms, which will be important (among others) in the case of localized topographies (Sec.~\ref{sec:ODElocalized}).
By \eqref{eq:anchor}, the distance between two successive fronts
$\Delta \phi_j$ is `anchored' by
\begin{equation}
 \label{eq:Deltaphij}
\Delta \phi_j := \phi_{j+1} - \phi_j =  (p_{j+1}-p_j) \frac{|\log \varepsilon|}{\sqrt{2}} + \ell_{j+1} - \ell_j, \; \; (j = 1, \ldots, N-1).
\end{equation}
This is our chosen analytic expression for the requirement that the fronts are `sufficiently far apart'.
As in \eqref{eq:phiupdown}, we assume in the upcoming analysis that $p_{j+1}-p_j = 1$ and $\ell_j(t) = \mathcal{O}(1)$, so that $\Delta \phi_j =  \frac{|\log \varepsilon|}{\sqrt{2}}$ to leading order in $\varepsilon$ (see Fig. \ref{fig:NFrontSketch}). We note however that we will go beyond that in our analysis of the ODE dynamics driven by localized topographies in  Sec.~\ref{sec:ODElocalized}.
\\ \\
By \eqref{eq:Deltaphij} the evolution of $\phi_j(t)$ will be slow and -- anticipating the magnitude of $\frac{d\phi_j(t)}{dt}$ for $\ell_j(t) = \mathcal{O}(1)$ -- we introduce $q_j(t)$ by
\begin{equation}
	\frac{d\phi_j(t)}{dt} = \varepsilon q_j.
\end{equation}
To study the migratory movement of $N$-front patterns, {\it i.e.,} to determine $q_j$ to leading order in $\varepsilon$, we subdivide the domain $\mathbb{R}$ into $N$ domains that each contains one front. (See Fig.~\ref{fig:NFrontSketch}.)
Specifically, we define the regions
\begin{equation}
\label{eq:defIjs}
\begin{array}{rccr}
I_1 & := & \left( -\infty, \frac{1}{2}(p_1 + p_2)\frac{|\log \varepsilon|}{\sqrt{2}} \right); & \\
I_j & := & \left(\frac{1}{2}(p_{j-1} + p_{j}) \frac{|\log \varepsilon|}{\sqrt{2}}, \frac{1}{2}(p_{j} + p_{j+1}) \frac{|\log \varepsilon|}{\sqrt{2}} \right), & (j = 2, \ldots, N-1);\\
I_N & := & \left(\frac{1}{2}(p_{N-1} + p_{N}) \frac{|\log \varepsilon|}{\sqrt{2}}, \infty \right). &
\end{array}
\end{equation}
Following the procedure used in section~\ref{sec:2frontExistence}, we now solve~\eqref{eq:mainEquation} separately in each region and match these expressions at the boundaries of the region. Doing so imposes conditions on $q_j$, which leads to a $N$-dimensional ODE describing the slow migratory movements of the front locations.

As before, we employ regular expansions to find approximate solutions in each of the regions. In regions with a front connecting (to leading order) $u_-^\varepsilon$ to $u_+^\varepsilon$ we set $u_j = u_\mathrm{up} + u_{j,1}$ and in regions with a front connecting (to leading order) $u_+^\varepsilon$ to $u_-^\varepsilon$ we set $u_j = u_\mathrm{down} + u_{j,1}$. For notational clarity, we assume that the odd regions are of the former type and the even regions of the latter type and we comment on other situation at the end of this section. So, we expand $u$ as follows
\begin{equation}
	u(x;\phi_1(t),\ldots,\phi_N(t)) =
\begin{cases}
	u_\mathrm{up}(x;\phi_j(t)) + \varepsilon u_{1,j}(x;\phi_j(t)) + h.o.t. & x \in I_j, j \mbox{ odd}; \\
	u_\mathrm{down}(x;\phi_j(t)) + \varepsilon u_{1,j}(x;\phi_j(t)) + h.o.t.& x \in I_j, j \mbox{ even}.
\end{cases}
\end{equation}
Substitution in~\eqref{eq:mainEquation} then leads to the following equation at $\mathcal{O}(\varepsilon)$:
\begin{equation}
	\mathcal{L}_0 u_{1,j} =
\begin{cases}
	- q_j u_\mathrm{up}' - F(u_\mathrm{up}, u_\mathrm{up}'; x) & j \mbox{ odd}; \\
	- q_j u_\mathrm{down}' - F(u_\mathrm{down}, u_\mathrm{down}'; x) & j \mbox{ even},
\end{cases}
\end{equation}
where $\mathcal{L}_0$, as defined in~\eqref{eq:defL0}, is now given by
\begin{equation*}
	\mathcal{L}_0 = \partial_x^2 + 1 - 3 \tanh^2\left(\frac{\sqrt{2}}{2}[x-\phi_j(t)]\right).
\end{equation*}
As this is essentially the same equation as before -- only with slightly different inhomogeneous terms -- the formulas from Appendix~\ref{sec:orderEpsSystem} can be used by replacing the inhomogenous terms.
Using the limit behavior of $u_\mathrm{up,down}$, computed in \eqref{eq:asymptoticResultsLO12} and of $u_{1,j}$, computed in \eqref{eq:u1limitbehavior}, we can now match solutions at the boundaries $x = \frac{1}{2}(p_{j+1}+p_j) \frac{|\log \varepsilon|}{\sqrt{2}}$ at order $\mathcal{O}\left(\sqrt{\varepsilon}\right)$ -- similar to the procedure in Sec.~\ref{sec:2frontExistence}.
\\ \\
Performing this procedure leads, for $j$ odd, to the following set of conditions:
\[
\begin{array}{rcl}
- 2 e^{\sqrt{2} \ell_j} + \frac{1}{8} (B_+)_j e^{-\sqrt{2} \ell_j} & = & - 2 e^{-\sqrt{2} \ell_{j+1}} - \frac{1}{8} (B_-)_{j+1} e^{\sqrt{2} \ell_{j+1}};
\\
2 e^{\sqrt{2} \ell_j} + \frac{1}{8} (B_+)_j e^{-\sqrt{2} \ell_j} & = & - 2 e^{-\sqrt{2} \ell_{j+1}} + \frac{1}{8} (B_-)_{j+1} e^{\sqrt{2} \ell_{j+1}}.
\end{array}
\]
Here, $(B_\pm)_j$ are defined consistently with~\eqref{eq:BpmLimits}; that is,
\begin{equation}
	(B_\pm)_j =
\begin{cases}
	B_{0,j} - \int_{\phi_j}^{\pm \infty} \left[ q_j u_\mathrm{up}(z;\phi_j) + F(u_\mathrm{up}(z;\phi_j),u_\mathrm{up}'(z;\phi_j);z)\right] u_\mathrm{up}'(z;\phi_j) dz & j \mbox{ odd}; \\
	B_{0,j} - \int_{\phi_j}^{\pm \infty} \left[-q_j u_\mathrm{up}(z;\phi_j) + F(-u_\mathrm{up}(z;\phi_j),-u_\mathrm{up}'(z;\phi_j);z)\right] u_\mathrm{up}'(z;\phi_j) dz & j \mbox{ even},
\end{cases}
\label{eq:BplusminODE}
\end{equation}
where $B_{0,j}$ are (unknown) constants.
By adding and subtracting these conditions, we arrive at the equivalent conditions -- that hold for $j$ odd:
\[
(B_+)_j = - 16 e^{-\sqrt{2} \left( \ell_{j+1} - \ell_{j} \right)}, \; \; \ \ \ \ \
(B_-)_{j+1} =16 e^{-\sqrt{2} \left( \ell_{j+1} - \ell_{j}\right)}.
\]
Similarly, for $j$ even, the same procedure results -- for $j$ even -- in the conditions
\[
(B_+)_j = 16 e^{-\sqrt{2} \left( \ell_{j+1} - \ell_{j} \right)}, \; \; \ \ \ \ \
(B_-)_{j+1}=-16 e^{-\sqrt{2} \left( \ell_{j+1} - \ell_{j}\right)}.
\]
Combining these sets of conditions, we thus have the following $2N$ conditions:
\begin{equation}
(B_+)_j =
\begin{cases}
	- 16 e^{-\sqrt{2} \left( \ell_{j+1} - \ell_{j} \right)} & j \mbox{ odd}\\
	+ 16 e^{-\sqrt{2} \left( \ell_{j+1} - \ell_{j} \right)} & j \mbox{ even},
\end{cases}
\; \; \ \ \ \ \
	(B_-)_j =
\begin{cases}
	- 16 e^{-\sqrt{2} \left( \ell_{j} - \ell_{j-1} \right)} & j \mbox{ odd}\\
	+ 16 e^{-\sqrt{2} \left( \ell_{j} - \ell_{j-1} \right)} & j \mbox{ even}.
\end{cases}
\label{eq:BConditions}
\end{equation}
Note that we have $2N$ conditions here, since we also have $2N$ unknowns ($q_1,\ldots,q_N$) and $(B_{0,1},\ldots,B_{0,N}$). Since we are only interested in the values for $q_j$, we can subtract $(B_-)_j$ and $(B_+)_j$ to eliminate the other constants.
That is, by~\eqref{eq:BplusminODE}, we have
\begin{equation}
	(B_-)_j - (B_+)_j =
\begin{cases}
\int_{-\infty}^\infty \left[ q_j u_\mathrm{up}(z;\phi_j) + F(u_\mathrm{up}(z;\phi_j),u_\mathrm{up}'(z;\phi_j);z)\right] u_\mathrm{up}'(z;\phi_j) dz & j \mbox{ odd} \\
\int_{-\infty}^\infty \left[-q_j u_\mathrm{up}(z;\phi_j) + F(-u_\mathrm{up}(z;\phi_j),-u_\mathrm{up}'(z;\phi_j);z)\right] u_\mathrm{up}'(z;\phi_j) dz & j \mbox{ even}.
\end{cases}
\end{equation}
Also, by~\eqref{eq:BConditions}, we have
\begin{equation}
	(B_-)_j - (B_+)_j = (-1)^j 16 \left[ e^{-\sqrt{2}\left(\ell_j-\ell_{j-1}\right)} - e^{-\sqrt{2}\left(\ell_{j+1}-\ell_j\right)}\right].
\end{equation}
Hence, we can combine these two expressions to find $N$ conditions that determine $q_1,\ldots,q_N$. Moreover, the integrals can be worked out and related to the Melnikov integrals $\mathcal{R}_\mathrm{up/down}$ from~\eqref{eq:FredholmConditionUp}-\eqref{eq:FredholmConditionDown}. This results in
\begin{equation}
	q_j =
\begin{cases}
	\frac{1}{\|u_\mathrm{up}'\|_2^2} \left[ - \mathcal{R}_\mathrm{up}(\phi_j) + 16 \left( e^{-\sqrt{2}(\ell_{j+1}-\ell_j)} - e^{-\sqrt{2}(\ell_j-\ell_{j-1})} \right) \right], & j \mbox{ odd} \\
	\frac{1}{\|u_\mathrm{up}'\|_2^2} \left[ - \mathcal{R}_\mathrm{down}(\phi_j) + 16 \left( e^{-\sqrt{2}(\ell_{j+1}-\ell_j)} - e^{-\sqrt{2}(\ell_j-\ell_{j-1})} \right) \right], & j \mbox{ even}.
\end{cases}
\end{equation}
Recalling $\frac{d\phi(t)}{dt} = \varepsilon q_j$ and noting $\varepsilon e^{-\sqrt{2} (\ell_{j+1}-\ell_j)} = e^{-\sqrt{2}(\phi_{j+1}-\phi_j)} = e^{-\sqrt{2} \Delta \phi_j}$, we find that the slow dynamics of $N$-front patterns is determined to leading order in $\varepsilon$ by the ODE,
\begin{equation}
	\frac{d\phi_j}{dt} = \frac{1}{\|u_\mathrm{up}'\|_2^2} \left[ - \varepsilon \mathcal{R}_j(\phi_j) + 16 \left( e^{-\sqrt{2}\Delta \phi_j} - e^{-\sqrt{2}\Delta \phi_{j-1}}\right) \right],
\label{eq:NFrontODE}
\end{equation}
where $\Delta \phi_j := \phi_{j+1}-\phi_j$ ($j = 1, \ldots, N-1$), $\Delta \phi_0 \rightarrow \infty$, $\Delta \phi_N \rightarrow \infty$ , and
\begin{equation}
\label{eq:RjRupdown}
	\mathcal{R}_j =
    \begin{cases} \mathcal{R}_\mathrm{up} & \mbox{region $I_j$ contains a front going up} \\ \mathcal{R}_\mathrm{down} &\mbox{region $I_j$ contains a front going down}.
    \end{cases}
\end{equation}
We now recall that the above derivation has been performed on $N$-front configurations in which the odd regions have fronts going up (from $u_-^\varepsilon$ to $u_+^\varepsilon$). The above procedure can also be used on the other kind of $N$-front configurations -- with odd regions that have fronts going down (from $u_-^\varepsilon$ to $u_+^\varepsilon$). The end result is identical.
\\ \\
Hence, the interactions of N-fronts anchored by \eqref{eq:anchor} are in general governed to leading order by the $N$-dimensional system of ODEs \eqref{eq:NFrontODE}.
Also, system \eqref{eq:NFrontODE} confirms that the choice of $p_j$ (in \eqref{eq:anchor}) does not show up explicitly, however, it does implicitly by \eqref{eq:BplusminODE}.
Furthermore, the (neglected) correction terms to \eqref{eq:NFrontODE} are $\mathcal{O}(\varepsilon^2)$, as can be established by a straightforward procedure, under the assumption that the inhomogeneity $F(U,V,x)$ in \eqref{eq:mainEquation} is bounded.
This implies that \eqref{eq:NFrontODE} loses its validity for multi-fronts $(\phi_1, \phi_2,...,\phi_N)$ that are spread out so much that the magnitude of the right hand side of \eqref{eq:NFrontODE} has decreased to $\mathcal{O}(\varepsilon^2)$ -- see also the extended discussions in Secs.~\ref{sec:ODElocalized} and \ref{sec:Discussion}.
\\ \\
Finally, since the system \eqref{eq:NFrontODE} of ODEs captures the dynamics of localized $N$-front patterns in the PDE \eqref{eq:mainEquation} governed by the small spectrum, the small eigenvalues associated to specific $N$-front patterns can also be found via \eqref{eq:NFrontODE}.
More specifically, stationary $N$-front patterns in \eqref{eq:mainEquation} directly correspond to fixed points of \eqref{eq:NFrontODE}.
By construction, the associated small eigenvalues of the $N$-front patterns correspond to the eigenvalues of these critical points (see Remark \ref{rem:1fronts} for the case of one-front patterns). If the critical point is stable in \eqref{eq:NFrontODE}, then the stability of the corresponding pattern as solution of \eqref{eq:mainEquation} immediately follows because all other elements of the spectrum necessarily lie to the left of the imaginary axis and are bounded away from it (see \ref{sec:oneFrontStability}).
\\ \\
We refrain from summarizing the results of this section in the form of a general theorem that captures the validity of $N$-front interaction equations \eqref{eq:NFrontODE} -- like we did with Theorems~\ref{th:1frontExStab} and \ref{th:2frontExStab} for the main results of Secs.~\ref{sec:1fronts} and \ref{sec:2fronts}. Indeed, it can be shown by the methods of \cite{ei2002pulse, promislow2002renormalization} that there is an attracting, approximate $N$-dimensional manifold $\mathcal{M}_N^0$ for the (semi-)flow generated by \eqref{eq:mainEquation} and that the flow on $\mathcal{M}_N^0$ is governed by \eqref{eq:NFrontODE} to leading order.
However, the manifold $\mathcal{M}_N^0$ will have boundaries, and the nature of (parts of) these boundaries will vary.
Indeed, the approximation scheme breaks down if the distance between (some) neighboring fronts will either become too small or too large. Therefore, a precisely formulated validity result will ask for many technical details -- see \cite{van2010front} for an example -- and not add to the transparency of the present text, see also Sec.~\ref{sec:Discussion}.

\begin{remark}
\label{rem:homAC}
In the classical (unforced) case of Allen-Cahn equation~\eqref{eq:mainEquation}, $\mathcal{R}_\mathrm{up/down}(\phi) \equiv 0$. Hence, the slow migration of fronts in a $N$-front configuration is given by the ODE
\begin{equation}
	\frac{d\phi_j}{dt} = \frac{16}{\|u_\mathrm{up}'\|_2^2} \left( e^{-\sqrt{2} \Delta \phi_j} - e^{-\sqrt{2} \Delta \phi_{j-1}}\right)
\end{equation}
(see \cite{Carr1989,Chen2004,Fusco1989}).
It is clear that this ODE cannot have any fixed points for $N > 1$ (recalling that $\Delta \phi_{0} = \infty$, we for instance note that $\frac{d\phi_1}{dt}$ cannot vanish). In fact,
\begin{equation}
\label{ddtphiNphi1}
\frac{d}{dt}(\phi_N - \phi_1) = - \frac{16}{\|u_\mathrm{up}'\|_2^2} \left[e^{-\sqrt{2}(\phi_2-\phi_1)} + e^{-\sqrt{2}(\phi_N-\phi_{N-1})} \right] < 0,
\end{equation}
which implies that the extent of an $N$-front pattern -- {\it i.e.}, the distance $\phi_N - \phi_1$ between the two outermost fronts -- always decreases.
By contrast, in the case of the inhomogeneous Allen-Cahn equation \eqref{eq:mainEquation}, the forcing terms may strengthen or counteract the effects that fronts have on each other. Thus, system \eqref{eq:NFrontODE} can have critical points corresponding to stationary multi-front patterns in \eqref{eq:mainEquation} -- as we have already seen in the explicit one-front and two-front examples in the preceding sections. This observation forms the foundation for some of the main results in the upcoming sections (see Theorems \ref{th:Nfrontslocexp} and \ref{th:Nfrontsper}). Moreover, the inhomogeneity may `spread out' an $N$-front pattern, {\it i.e.},  may force $\frac{d}{dt}(\phi_N - \phi_1)$ to be positive (as in Fig. \ref{fig:5Fronts-Intro} -- see also Sec.~\ref{sec:Discussion}).
\end{remark}

\begin{remark}
\label{rem:1fronts}
The dynamics of a solitary front is described by
\begin{equation}
\label{eq:1frontdyn}
\frac{d\phi_1}{dt} = -\frac{\varepsilon}{\|u_\mathrm{up}'\|_2^2} \mathcal{R}_\mathrm{up}(\phi_1),
\end{equation}
where we note that we considered the case of an `up-front' connecting $u_-^\varepsilon(x) = - 1 + \mathcal{O}(\varepsilon)$ as $x \to - \infty$ to $u_{1}^\varepsilon(x) = 1 + \mathcal{O}(\varepsilon)$ as $x \to +\infty$ \eqref{eq:uMinusPlus}, so that
$\mathcal{R}_1(\phi_1) = \mathcal{R}_\mathrm{up}(\phi_1)$ \eqref{eq:RjRupdown}.
We observe that this agrees fully with the (stationary) one-fronts analysis in  Sec.~\ref{sec:1fronts}.
Namely, the existence of a stationary one-front is determined by the zeroes $\phi_1^\ast$ of $\mathcal{R}_\mathrm{up}(\phi_1)$, and their stability by $\lambda(\phi_1^\ast) = \varepsilon \tilde{\lambda}(\phi_1^\ast) = - \varepsilon \mathcal{R}'_\mathrm{up}(\phi_1^\ast)/\|u_\mathrm{up}'\|_2^2$ (see  \eqref{eq:eigenvalue1Front} which uses that $\|u_\mathrm{up}'\|_2^2 = \frac{2}{3}\sqrt{2}$).
By construction,  equation \eqref{eq:1frontdyn} provides more information: it also shows that a one-front pattern will slowly move towards a zero $\phi_1^\ast$ of $\mathcal{R}_\mathrm{up}(\phi_1)$ that has $\mathcal{R}'_\mathrm{up}(\phi^\ast) > 0$ (if such a zero exists), {\it i.e.}, it will move toward and `settle down' (get pinned) at a nearby stable stationary one-front state (if such a state exists). (See also Theorem \ref{th:Nfrontsperattr} for a more precise statement in the setting of interacting $N$-front patterns driven by a spatially periodic topography \eqref{eq:Ftopography}.) Moreover, there can also be stable traveling one-front patterns that do not approach a stable stationary state, as shown in Fig.~\ref{fig:5Fronts-Intro}(b) (for $t > 4300$) and below in Fig.~\ref{fig:3FrontLocPer}(a) (for $t > 3400$). The speed at which these fronts travel will in general not be constant; in fact, in the case of localized topographies the fronts typically slow down, see however Fig.~\ref{fig:5FrontDynVaryingP} and the discussion in Sec.~\ref{sec:Discussion}.
\end{remark}

\begin{remark}
\label{rem:gradient}
System \eqref{eq:NFrontODE} can be written as a gradient flow,
\[
\frac{d \phi_j}{d \tau} = \frac{1}{\|u_\mathrm{up}'\|_2^2} \frac{\partial}{\partial \phi_j} \left[ -\varepsilon \sum_{i=1}^{N}\int_0^{\phi_i} \mathcal{R}_j(\varphi) \, d \varphi + 8 \sqrt{2} \sum_{i=1}^{N-1} e^{-\sqrt{2}(\phi_{i+1} - \phi_i)}\right], \; \; (j=1,2,...,N).
\]
This is common for the interaction equations governing the dynamics of localized structures (see \cite{van2010front}).
Thus, we may conclude for instance that the system of ODEs \eqref{eq:NFrontODE} cannot exhibit periodic dynamics.
\end{remark}

\section{Front dynamics driven by localized topographies}
\label{sec:ODElocalized}
In this section, we consider the impact of a localized topography on the dynamics of \eqref{eq:mainEquation}, {\it i.e.} we consider \eqref{eq:mainEquation} with $F(U,U_x,x)$ prescribed by a topography $H(x)$ \eqref{eq:Ftopography} that is localized. As was noted in Sec.~\ref{sec:top1fronts}, the fact that we assume that $F(U,U_x,x)$ is prescribed by a topography $H(x)$ implies that $\mathcal{R}_\mathrm{up}(\phi) \equiv \mathcal{R}_\mathrm{down}(\phi) \equiv \mathcal{R}(\phi)$ (see  \eqref{eq:Rtopoupdown}) with $\mathcal{R}(\phi)$ given by \eqref{eq:Ex2R}. This is in fact the most important reason to restrict ourselves to this class of inhomogeneous perturbations $\varepsilon F(U,U_x,x)$ in \eqref{eq:mainEquation}: continuing to work with in principle distinct Melnikov functions $\mathcal{R}_\mathrm{up}(\phi)$ and $\mathcal{R}_\mathrm{down}(\phi)$ would make the `bookkeeping' in the upcoming analysis more complicated. In other words, this choice is -- again -- made to simplify the upcoming analysis  as much as possible. The nature of our findings indicate that very similar results will hold for the general case.
\\ \\
The simplification of the topographic setting can be made explicit by (re-)introducing
\begin{equation}
\label{eq:defpsiStau}
\psi_j(t) = \sqrt{2} \phi_j(t), \; \; \mathcal{S}(\psi_j) = \mathcal{R}(\phi_j), \; \; \tau = \sqrt{2}t/\|u_\mathrm{up}'\|_2^2.
\end{equation}
(see \eqref{eq:Rsolhill1}).
By this rescaling, the system of $N$ ODEs \eqref{eq:NFrontODE} can be written to leading order in $\varepsilon$ as
\begin{equation}
\label{eq:dynNfronts}
\left\{
\begin{array}{r c r c l}
\frac{d \psi_1}{d \tau} & = &
- \varepsilon \mathcal{S}(\psi_1) & + & 16 e^{-(\psi_2 - \psi_1)} \\
\frac{d \psi_2}{d \tau} & = & - \varepsilon \mathcal{S}(\psi_2) & + & 16 \left(e^{-(\psi_3 - \psi_2)} - e^{-(\psi_2 - \psi_1)}\right) \\
& \vdots & \\
\frac{d \psi_j}{d \tau} & = & - \varepsilon \mathcal{S}(\psi_j) & + & 16 \left(e^{-(\psi_{j+1} - \psi_j)} - e^{-(\psi_j - \psi_{j-1})}\right) \\
& \vdots & \\
\frac{d \psi_N}{d \tau} & = & - \varepsilon \mathcal{S}(\psi_N) & - & 16 e^{-(\psi_N - \bar\psi_{N-1})}.
\end{array}
\right.
\end{equation}
We recall that this system has been derived under the condition that $\psi_{j} < \psi_{j+1}$ with $\psi_{j+1} - \psi_j = \mathcal{O}(|\log \varepsilon|)$ ($j=1,..,N-1$) and that the (neglected) leading order correction terms are of $\mathcal{O}(\varepsilon^2)$.
\\ \\
The main analytical result of this section is Theorem \ref{th:Nfrontslocexp} on the existence and stability of stationary $N$-front patterns in the inhomogeneous Allen-Cahn equation \eqref{eq:mainEquation} driven by a general localized topography $H(x)$ -- see Sec.~\ref{sec:ODENfrontloc}.
The analysis in essence reduces to an investigation of the existence and stability of critical points of system \eqref{eq:dynNfronts} (with $\mathcal{S}(\psi)$ given by \eqref{eq:Ex2R} and \eqref{eq:defpsiStau}). The nature of the decay of topography $H(x)$ plays a crucial role in both the statement and the proof of Theorem \ref{th:Nfrontslocexp}, not only directly but also through its interaction with the leading order corrections to \eqref{eq:dynNfronts}. To carefully build our understanding, and to make the analysis more transparent, we will first study the special cases $N=2$ (Subsec.~\ref{sec:ODE2frontloc}) and $N=3,4$ (Subsec.~\ref{sec:ODE34frontloc}). In these sections, we will simultaneously consider the interacting front dynamics of two-, three- and four-front patterns, {\it i.e.}, the dynamics exhibited by system \eqref{eq:dynNfronts} (for $N=2,3,4$) beyond just its critical points. This will also yield general insights and conjectures about the dynamics of $N$-front patterns -- as will also be discussed in  Sec.~\ref{sec:ODENfrontloc}.
\\ \\
The reason why the nature of the decay of topography $H(x)$ plays such a crucial role in the upcoming analysis is that the limiting behavior of $H(x)$ -- or more precisely of $H'(x)$ -- as $x \to \pm \infty$ directly affects the limiting behavior of $\mathcal{S}(\psi)$ as $\psi \to \pm \infty$.
In fact, we first provide two technical lemmas in which the relation between the decay behavior of $H'(x)$ and that of $\mathcal{S}(\psi)$ is precisely described.
These are valid for general topographies $H(x)$ (that attain $\mathcal{O}(1)$ values for $x = \mathcal{O}(1)$), {\it i.e.} $H(x)$ need not be symmetric here, nor given by (\ref{eq:defHuni-alg}).
However, we do assume that the derivative of the topography $H'(x)$ has uniformly decaying tails and that the rates of decay are the same for $x \to \pm \infty$ (see \eqref{eq:expdecayHx}, \eqref{eq:algdecayHx} for more explicit statements).
Both assumptions are made for simplicity,
and we focus on two specific types of decay: exponential decay (Lemma \ref{lem:expdecayH}) and algebraic decay (Lemma \ref{lem:algdecayH}).

\begin{lemma}[Exponential Decay]
\label{lem:expdecayH}
Assume that there is a constant $\mu > 0$ and expressions $h_\pm = h_\pm(\mu) \neq 0$ such that
\begin{equation}
\lim_{x \to \pm \infty} e^{\pm \mu \sqrt{2} x} H'(x) = h_\pm(\mu).
\label{eq:expdecayHx}
\end{equation}
If $0 < \mu < 1$, then
\begin{equation}
\label{eq:expdecayR0mu1}
\lim_{\psi \to \pm \infty} e^{\pm \mu \psi} \mathcal{S}(\psi) = h_\pm(\mu) w_\pm(\mu),
\end{equation}
with
\begin{equation}
w_\pm(\mu) := \int_{-\infty}^{\infty} e^{\mp \mu \sqrt{2} y}W_h(y) dy  \; \; (> 0);
\label{eq:defwpm}
\end{equation}
recall \eqref{eq:defWh}.
Finally, if $\mu > 1$, then
\begin{equation}
\label{eq:expdecayRmug1}
\lim_{\psi \to \pm \infty} e^{\pm \psi} \mathcal{S}(\psi) = 4 \hat{h}_\pm(\mu) = 4 \int_{-\infty}^\infty e^{\sqrt{2} z} H'(z) dz,
\end{equation}
which implies that $\mathcal{S}(\psi)$ cannot decay faster than $e^{\mp \psi}$ as $\psi \to \pm \infty$, independent of the decay rate of $H'(x)$.
\end{lemma}
As the parameter $\mu$ tends to $\mu = 1$, the decay behavior as indicated by the lemma breaks down in the following way: the function $w_\pm(\mu) \rightarrow \infty$ as $\mu \uparrow 1$, since by \eqref{eq:heteroclinicSolution} and \eqref{eq:defWh}
\begin{equation}
\label{eq:decayWh}
W_h(y) = 4 \frac{(1-e^{-\sqrt{2}y})^2}{(1+e^{-\sqrt{2}y})^4} e^{-\sqrt{2}y} =
4 \frac{(1-e^{\sqrt{2}y})^2}{(1+e^{\sqrt{2}y})^4} e^{\sqrt{2}y}.
\end{equation}
Moreover, $|\hat{h}_\pm(\mu)| \uparrow \infty$ as $\mu \downarrow 1$, since the integral in~\eqref{eq:expdecayRmug1} diverges as $\mu \downarrow 1$ by the assumed decay behavior of $H'(x)$. Further, for the explicit example of the unimodal hill $H_{\rm exp}(x;\mu)$~\eqref{eq:defHuni-alg}, the parameter $\mu$ in the definition of $H_{\rm exp}$ plays the same role as defined in \eqref{eq:expdecayHx} -- for this example $h_\pm(\mu) = \mp 4 \sqrt{2} \mu$. Moreover, as $\mu = 1$, the explicit expression \eqref{eq:Rsolhill1} associated to $H_{\rm exp}(x;1)$ -- with $\mathcal{R}_{\rm exp}(\phi;\mu) = \mathcal{S}_{\rm exp}(\psi;\mu)$ \eqref{eq:defpsiStau} -- indicates that
\[
\lim_{\psi \to \infty}\frac{e^{\psi} \mathcal{S}_{\rm exp}(\psi;1)}{\psi} = - 16\sqrt{2} = h_+(1) w_+(1) = (- 4 \sqrt{2}) \times 4,
\]
{\it i.e.}, for $\mu = 1$, $\mathcal{S}_{\rm exp}(\phi)$ decays as $-16 \sqrt2 \psi \, e^{-\psi}$ for $\psi \to \infty$.
\\
\begin{figure}
\centering
\begin{subfigure}[t]{0.32\textwidth}
\centering
\includegraphics[width=\textwidth]{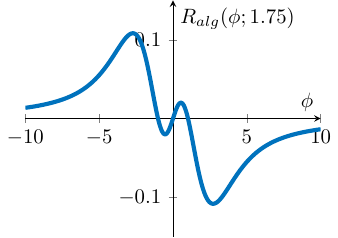}
\end{subfigure}
~
\begin{subfigure}[t]{0.32\textwidth}
\centering
\includegraphics[width = \textwidth]{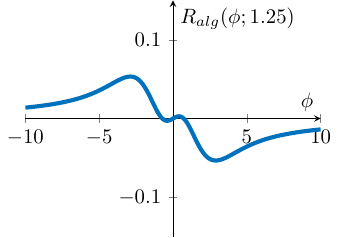}
\end{subfigure}
~
\begin{subfigure}[t]{0.32\textwidth}
\centering
\includegraphics[width = \textwidth]{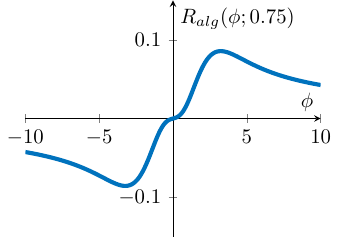}
\end{subfigure}
\caption{The Melnikov function $\mathcal{R}(\phi;p)$ for the topography $H_{\rm alg}(x; p)$ of a symmetrical hill with algebraically decaying tails \eqref{eq:defHuni-alg} (see \eqref{eq:Ex2R}). By construction, $p$ determines the leading order algebraic decay as in \eqref{eq:algdecayHx} (a) $p = 1.75$, (b) $p = 1.25$, (c) $p = 0.75$ (see Fig. \ref{fig:SolHillPitchfork}).}
\label{fig:AlgHill}
\end{figure}
\\
The case of algebraic decay is simpler.

\begin{lemma}[Algebraic Decay]
\label{lem:algdecayH}
Assume that there is a constant $p \in \mathbb{R}$ and $\tilde{h}_\pm = \tilde{h}_\pm(p) \neq 0$ such that
\begin{equation}
\lim_{x \to \pm \infty} |x|^p H'(x) = \tilde{h}_\pm(p),
\label{eq:algdecayHx}
\end{equation}
then
\begin{equation}
\label{eq:algdecayR}
\lim_{\psi \to \pm \infty} |\psi|^p \mathcal{S}(\psi) = \frac13 2^{\frac12 (p+3)}\, \tilde{h}_\pm(p).
\end{equation}
\end{lemma}
Note that we have thus allowed for topographies $H(x)$ that may also grow algebraically -- see especially Sec.~\ref{sec:Discussion}.
Further, for the unimodal topography with algebraically decaying/growing tails $H_{\rm alg}(x;p)$ \eqref{eq:defHuni-alg}, we have $\tilde{h}_\pm(p) = \pm (1-p)$ -- recall also \eqref{eq:defHuni-alg}.
\\ \\
{\bf Proofs of Lemmas \ref{lem:expdecayH} and \ref{lem:algdecayH}.} The derivations of \eqref{eq:expdecayR0mu1}, \eqref{eq:expdecayRmug1}, and \eqref{eq:algdecayR} follow from a number of straightforward observations, of which we only present the main lines. Moreover, we only consider the limits $x \to +\infty$, noting that the limits $x \rightarrow -\infty$ follow in a similar vein.
\\ \\
Firstly, for $0 < \mu < 1$, we obtain, using  \eqref{eq:Ex2R}, that
\[
\begin{array}{rcl}
\lim_{\psi \to \infty} e^{\psi} \mathcal{S}(\psi)
& = &
\lim_{\psi \to \infty} \int_{-\infty}^\infty e^{\psi} H'(y+\frac12 \sqrt2 \psi) W_h(y) dy
\\
& = &
\int_{-\infty}^\infty \left(\lim_{\psi \to \infty} \left[e^{\mu \sqrt{2} (y + \frac12 \sqrt2 \psi)} H'(y+\frac12 \sqrt2\psi)\right] e^{-\mu \sqrt{2} y}W_h(y) \right)dy
\\
& = &
h_+(\mu) \int_{-\infty}^{\infty} e^{- \mu \sqrt{2} y} \, W_h(y) dy.
\end{array}
\]
Then, for $ \mu > 1$, we obtain by \eqref{eq:decayWh},
\[
\begin{array}{rcl}
\lim_{\psi \to \infty} e^{\psi} \mathcal{S}(\phi) & = &
\int_{-\infty}^\infty \left(e^{\sqrt{2} z} H'(z) \lim_{\psi \to \infty} \left[e^{- \sqrt{2}(z-\frac12 \sqrt2 \psi)}W_h(z-\frac12 \sqrt2 \psi) \right] \right)dz
\\
& = & 4 \int_{-\infty}^\infty e^{\sqrt{2} z} H'(z) dz.
\end{array}
\]
Secondly, for the case of algebraic decay, we have
\[
\begin{array}{rcl}
\lim_{\psi \to \infty} \psi^p \mathcal{S}(\psi) & = &
\int_{-\infty}^\infty \left(\lim_{\phi \to \infty} \left[(y+\frac12 \sqrt2\psi)^p H'(y+\frac12 \sqrt2 \psi) \frac{\psi^p}{(y + \frac12 \sqrt2\psi)^p} W_h(y) \right] \right) \, dy
\\
& = &
\tilde{h}_+(p) \int_{-\infty}^\infty \left(\lim_{\psi \to \infty} \left[\frac{\sqrt2^p}{(1 + \frac{\sqrt{2}y}{\psi})^p} W_h(y) \right] \right) \, dy
\\
& = & 2^{\frac12 p}\tilde{h}_+(p) \int_{-\infty}^\infty W_h(y) \, dy = \frac13 2^{\frac12 (p+3)} \, \tilde{h}_\pm(p).
\end{array}
\]
\hfill $\Box$

\subsection{Interacting two-fronts}
\label{sec:ODE2frontloc}
The leading order interaction ODE for two-front patterns in \eqref{eq:mainEquation} with driving term \eqref{eq:Ftopography} (that approach $u_-^\varepsilon(x) = -1 + \mathcal{O}(\varepsilon)$ \eqref{eq:uMinusPlus} as $x \to \pm \infty$) is given by
\begin{equation}
\left\{
\begin{array}{r c l}
\frac{d\psi_1}{d\tau} & = & - \varepsilon \mathcal{S}(\psi_1) + 16 e^{-(\psi_2 - \psi_1)} \\
\frac{d\psi_2}{d\tau} & = & - \varepsilon \mathcal{S}(\psi_2) - 16 e^{-(\psi_2 - \psi_1)},
\end{array}
\right.
\label{eq:intODE2loc}
\end{equation}
with the additional condition that $\psi_1(t) < \psi_2(t)$.
A priori, one expects that the dynamics of two-dimensional system \eqref{eq:intODE2loc} are dominated by its critical points $(\bar{\phi}_1, \bar{\phi}_2)$.
These satisfy
\begin{equation}
\mathcal{S}(\bar{\psi}_1) = - \mathcal{S}(\bar{\psi}_2) = \frac{16}{\varepsilon} e^{-(\bar{\psi}_2 - \bar{\psi}_1)},
\label{eq:crit2loc}
\end{equation}
which coincides with \eqref{eq:twoFrontCondition}. The $\mathcal{O}(1/\varepsilon)$ factor on the right hand side implies that $\bar{\psi}_2 - \bar{\psi}_1$ must be asymptotically large. By the assumptions on $H(x)$, $\mathcal{S}(\psi)$ does not explicitly depend on $\varepsilon$ and thus varies on the $\mathcal{O}(1)$ scale. This implies that either both $|\bar{\psi}_1|, |\bar{\psi}_2| \gg 1$ or only one of these. In other words, there are three kinds of (potential) critical points in \eqref{eq:intODE2loc}:
\\ \\
$\bullet$ critical points of the first kind have $|\bar{\psi}_j| \gg 1$ and, by \eqref{eq:expdecayR0mu1}, $|\mathcal{S}(\bar{\psi}_j)| \ll 1$ (j= 1,2).
\\
$\bullet$ critical points of the second kind have $\bar{\psi}_1 = \mathcal{O}(1)$ and $\bar{\psi}_2 \gg 1$.
\\
$\bullet$ critical points of the third kind have $\bar{\psi}_1 \ll - 1$ and $\bar{\psi}_2 = \mathcal{O}(1)$.
\\ \\
In Sec.~\ref{sec:ODE34frontloc}, we will introduce a fourth kind, and in Sec.~\ref{sec:ODENfrontloc} we will establish that all possible stationary $N$-front patterns can be seen as extensions of these four fundamental kinds (see Theorem \ref{th:Nfrontslocexp}).
\\ \\
We start by searching for critical points of the first kind. For this, we introduce $\nu_\pm > 0$ and $\bar{\ell}_\pm$ and set
\begin{equation}
\bar{\psi}_1 = -\nu_- |\log \varepsilon| + \bar{\ell}_- \ll -1, \; \;
\bar{\psi}_2 = \nu_+ |\log \varepsilon| + \bar{\ell}_+ \gg 1
\label{eq:Ansat2critpts1}
\end{equation}
so that (for topographies $H(x)$ that decay exponentially), we have to leading order (by \eqref{eq:expdecayR0mu1}) that
\begin{equation}
\label{eq:exp21Rpsi12}
\begin{array}{rlll}
0 < \mu < 1: &
\mathcal{S}(\bar{\psi}_{1}) = \varepsilon^{\mu \nu_-} h_- w_- e^{\mu\bar{\ell}_-}, &
\mathcal{S}(\bar{\psi}_{2}) = \varepsilon^{\mu \nu_+} h_+ w_+ e^{-\mu\bar{\ell}_+}, &
e^{-(\bar{\psi}_2 - \bar{\psi}_1)} = \varepsilon^{\nu_- + \nu_+} e^{-(\bar{\ell}_+ - \bar{\ell}_-)};
\\
\mu > 1: &
\mathcal{S}(\bar{\psi}_{1}) = 4 \varepsilon^{\nu_-} \hat{h}_- e^{\bar{\ell}_-}, &
\mathcal{S}(\bar{\psi}_{2}) = 4 \varepsilon^{\nu_+} \hat{h}_+ e^{-\bar{\ell}_+}, &
e^{-(\bar{\psi}_2 - \bar{\psi}_1)} = \varepsilon^{\nu_- + \nu_+} e^{-(\bar{\ell}_+ - \bar{\ell}_-)},
\end{array}
\end{equation}
where we have suppressed explicit notation of dependence on $\mu$ of the functions $h_\pm$, $w_\pm$ and $\hat{h}_\pm$.
Hence, by \eqref{eq:crit2loc} this implies that for the critical points we need to take $\nu_- = \nu_+ := \nu$ and for $0 < \mu < 1$, $1 + \mu \nu = 2 \nu$, {\it i.e.}, $\nu = 1/(2-\mu)$ while $\nu = 1$ for $\mu > 1$ (see below for a brief discussion of the case $\mu =1$). It also follows that it is necessary for the existence of these  kinds of critical points that $h_+, \hat{h}_+ < 0 < h_-, \hat{h}_-$ and that there exists a unique critical point if this conditions holds. Note that the fact that $h_+, \hat{h}_+ < 0 < h_-, \hat{h}_-$ is necessary also implies that there cannot be stationary points $(\bar{\phi}_1, \bar{\phi}_2)$ with $(\bar{\phi}_1, \bar{\phi}_2)$ on the same tail of $\mathcal{S}(\psi)$, {\it i.e.}, with $\bar{\phi}_1 \ll \bar{\phi}_2 \ll -1$ or $1 \ll \bar{\phi}_1 \ll \bar{\phi}_2$. Note also that the symmetric critical point associated to the symmetrical topographies $H(x)$ considered in Sec.~\ref{sec:top2fronts} coincides with the present critical point of the first kind.
\\ \\
At this point it also becomes clear that it may be natural to adapt the `anchor' points  \eqref{eq:anchor} of the analysis by which system \eqref{eq:intODE2loc} was derived to the (interacting) patterns under consideration.
For instance, for $0 < \mu < 1$, the choices $p_2 = \nu_+ = 1/(2-\mu) = - p_1$ are self-evident (where we note that the factor $1/\sqrt{2}$ is incorporated in the definition of $\psi$).
\\ \\
More importantly, this is the first time we need to realize the potential impact of the neglected higher order effects in the derivation of the interaction ODEs \eqref{eq:NFrontODE}. These higher order terms can be represented by the correction terms $\varepsilon^2 R^{j,2}_{2}(\psi_1, \psi_2;\varepsilon)$, with $j=1,2$ corresponding to the correction to the $\frac{d \psi_j}{d \tau}$ equation.
Moreover, it follows by standard techniques that if the inhomogeneous term $F(U,V,x)$ in \eqref{eq:mainEquation} indeed is bounded (for $x \in \mathbb{R}$) and if $\psi_{1,2}$ remain bounded (in the appropriate norms), then there is a $C > 0$ such that $|R^{j,2}_2(\psi_1, \psi_2;\varepsilon)| < C$ ($j=1,2$). In other words, the higher order terms are of $\mathcal{O}(\varepsilon^2)$, which is the expected natural behavior of a regular perturbation expansion. However, since we did (and will) not evaluate (the leading order approximation of) $R^{j,2}_2(\psi_1, \psi_2;\varepsilon)$, it is necessary to only consider critical points $(\bar{\psi}_1,\bar{\psi}_2)$ with $\varepsilon |\mathcal{S}(\bar{\psi}_{1,2})|, e^{-(\bar{\psi}_2 - \bar{\psi}_1)} \gg \varepsilon^2$ (since the terms $R^{j,2}_2(\psi_1, \psi_2;\varepsilon)$ may have a leading order effect if these conditions do not hold). Hence, we need to impose for $0 < \mu < 1$ that $1 + \mu \nu = 2 \nu = \frac{2}{2-\mu} < 2$, {\it i.e.}, that $\mu < 1$, which is consistent. However, for $\mu > 1$ we find $1 + \nu = 2 \nu = 2 < 2$ and conclude that the $R^{j,2}_2(\psi_1, \psi_2;\varepsilon)$ terms cannot be neglected.
By working out the details of the case $\mu = 1$ in \eqref{eq:expdecayHx}, we find that exactly the same holds for $\mu =1$. Since we do not intend here to go deeper into the analysis of the impact of these higher order terms -- see also Sec.~\ref{sec:Discussion} -- we therefore impose the assumption that $0 < \mu < 1$ in our analysis of the critical points of the first kind.
\\ \\
The second kind of critical points have $\bar{\psi}_1 = \mathcal{O}(1)$ and $\bar{\psi}_2 \gg 1$, which implies that $|\mathcal{S}(\bar{\psi}_2)| \ll 1$. We note by \eqref{eq:crit2loc} that $\bar{\psi}_1$ must be asymptotically close to a zero of $\mathcal{S}(\psi)$. Thus, we assume that there is a $\psi_\ast$ such that $\mathcal{S}(\psi_\ast) = 0$ and $\mathcal{S}'(\psi_\ast) \neq 0$ (i.e. we assume that the zero $\psi_\ast$ is non-degenerate); as example, see Fig. \ref{fig:SolHillPitchfork} for $\mathcal{S}(\psi)$ associated to the unimodal hill \eqref{eq:defHuni-alg}, which has three such zeroes for $\mu \in (\mu_{\rm PF}, 1)$. Therefore, we change \eqref{eq:Ansat2critpts1} into
\begin{equation}
\bar{\psi}_1 = \psi_\ast + \varepsilon^{\nu_1} \tilde{\psi}_1, \; \; \ \ \ \
\bar{\psi}_2 = \nu_2 |\log \varepsilon| + \bar{\ell}_2.
\label{eq:Ansat2critpts2}
\end{equation}
We deduce for $0 < \mu < 1$ (using similar reasoning as in the previous case) the critical points are
\begin{equation}
\label{sol:nu12}
\nu_1(\mu) = \mu/(1-\mu), \; \; \nu_2(\mu) = 1/(1-\mu)
\end{equation}
and that $\tilde{\psi}_1$ and $\bar{\ell}_2$ can be solved uniquely (in terms of $\mathcal{S}'(\psi_\ast) \neq 0$) if $h_+ < 0$ (note that for critical points of this second kind it is thus \emph{not} necessary that $h_- > 0$ as is necessary for critical points of the first kind). For $\mu > 1$, we observe from \eqref{eq:expdecayR0mu1} and the second identity of \eqref{eq:crit2loc} that,
\[
\mathcal{O}(\varepsilon^{1 + \nu_2}) = |\varepsilon \mathcal{S}(\bar{\psi}_2)| = 16e^{-(\bar{\psi}_2 - \bar{\psi}_1)} = \mathcal{O}(\varepsilon^{\nu_2}),
\]
which implies that there cannot be critical points of the second kind in case $\mu > 1$.
\\ \\
We again impose that the impact of the correction terms $\varepsilon^2 R^{j,2}_2(\psi_1, \psi_2;\varepsilon)$ ($j=1,2$) is indeed of higher order. Thus we find that this leading order outcome only is valid if $1+\nu_1 = 1 + \mu \nu_2 < 2$, {\it i.e.}, if $0 < \mu < \frac{1}{2}$. Note that this does not necessarily imply that there cannot be stationary two-front patterns of the second kind for $\mu \geq \frac{1}{2}$. In fact, it is not unlikely that these patterns do exists for (some) $\mu \in  [\frac{1}{2}, 1)$, since it is natural to expect that $|R^{j,2}_2(\psi_1, \psi_2;\varepsilon)| \ll 1$ for $\bar{\psi_2} \gg 1$. However, to investigate this, one would need to develop a higher analysis of \eqref{eq:intODE2loc}, which we refrain from doing here.
\\ \\
Determining the third kind of critical points goes completely along the same lines. We find that $0 < \mu < \frac12$.
However, in this case, (only) $h_- > 0$ appears as necessary existence condition.
\\ \\
We conclude that in a case like that of Fig. \ref{fig:SolHillPitchfork}(a) for the unimodal hill $H_{\rm exp}(x; \mu)$ \eqref{eq:defHuni-alg} with $\mu \in (\mu_{\rm PF}, 1)$) -- that also has $h_+ < 0 < h_-$ -- system \eqref{eq:intODE2loc} has in total $\mathcal{N}(2) = 7$ critical points: one of the first kind, three of the second and three of the third kind. Moreover, we also conclude that there can be no critical points if neither $h_+ < 0$ nor $h_- > 0$ and that there are as many stationary two-fronts as there are stationary one-fronts, {\it i.e.}, zeroes of $\mathcal{S}(\psi)$, if only $h_+ < 0$ -- in which case they are all of the second kind -- or only $h_- > 0$ (third kind).
\\ \\
However, all critical points are unstable as we shall show in the rest of this subsection. For this, we define
\begin{equation}
E_j := 16 e^{-(\bar{\psi}_{j+1} - \bar{\psi}_{j})}.
\label{eq:defbarEj}
\end{equation}
Hence, in the first case,
\[
E_1 = 16 \varepsilon^{\frac{2}{2-\mu}} e^{-(\bar{\ell}_{+} - \bar{\ell}_{-})} =
\varepsilon^{\frac{2}{2-\mu}} |h_+| w_+ e^{-\mu \bar{\ell}_+}
\]
(\eqref{eq:crit2loc}, \eqref{eq:Ansat2critpts1}, \eqref{eq:exp21Rpsi12}). By the exponential nature of $\mathcal{S}(\psi)$ (for $|\psi| \gg 1$) \eqref{eq:expdecayR0mu1}, we conclude that the spectral stability of the critical points $(\bar{\psi}_1,\bar{\psi}_2)$ is (to leading order) determined by the eigenvalues of the Jacobian matrix
\begin{equation}
\label{eq:stab2loc1}
\left(
\begin{array}{cc}
- \varepsilon \mathcal{S}'(\bar{\psi}_1) + E_1  & - E_1
\\
- E_1 & - \varepsilon \mathcal{S}'(\bar{\psi}_2) + E_1
\end{array}
\right)
=
\varepsilon^{\frac{2}{2-\mu}} |h_+| w_+ e^{-\mu \bar{\ell}_+}
\left(
\begin{array}{cc}
1 - \mu & - 1
\\
- 1 &  1 - \mu
\end{array}
\right).
\end{equation}
Thus, $(\bar{\psi}_1,\bar{\psi}_2)$ is a saddle point with (leading order) eigenvalues,
\begin{equation}
\label{eq:lamba21loc}
\lambda_1(\mu) =  -\mu |h_+| w_+ e^{-\mu \bar{\ell}_+}  = -\mu  E_1 < 0,
\; \;
\lambda_2(\mu) = (2-\mu)|h_+| w_+ e^{-\mu \bar{\ell}_+} = (2 - \mu) E_1 > 0
\end{equation}
(where we recall that we needed to assume in the existence analysis that $0 < \mu < 1$). A very similar approach yields that the stability of a critical point $(\bar{\psi}_1,\bar{\psi}_2)$ of the second kind is determined by the eigenvalues of the Jacobian matrix
\begin{equation}
\label{eq:stab2loc2}
E_1
\left(
\begin{array}{cc}
1 - \frac{1}{\tilde{\psi}_1} \varepsilon^{-\frac{\mu}{1-\mu}} & - 1
\\
- 1 & 1 - \mu
\end{array}
\right),
\end{equation}
with
\[
E_1 = 16 \varepsilon^{\frac{1}{1-\mu}} e^{\psi_\ast -\bar{\ell}_2}.
\]
The eigenvalues are given to leading order by
\begin{equation}
\label{eq:eigenv2front2nd}
\lambda_1 = - 16 \varepsilon \frac{e^{\psi_\ast -\bar{\ell}_2}}{\tilde{\psi}_1} = - \varepsilon \mathcal{S}'(\psi_\ast), \; \;
\lambda_2 = 16 (1-\mu) \varepsilon^{\frac{1}{1 - \mu}} e^{\psi_\ast -\bar{\ell}_2} = (1-\mu) E_1
\end{equation}
\eqref{eq:crit2loc}, \eqref{eq:Ansat2critpts2}.
Hence, $\lambda_2 >0$ (again), so that $(\bar{\psi}_1,\bar{\psi}_2)$ is unstable. Moreover, the sign of $\lambda_1$ is determined by (the sign of) $\mathcal{S}'(\psi_\ast)$. Hence, by \eqref{eq:eigenvalue1Front}, $(\bar{\psi}_1,\bar{\psi}_2)$ again is a saddle of the flow induced by \eqref{eq:intODE2loc} if the one-front pattern located at $\psi = \psi_\ast$ is stable; otherwise it is an unstable node.
Essentially the same analysis yields the same instability results for the critical points of the third kind.
\\ \\
We refer to Sec.~\ref{sec:ODENfrontloc}, Theorem \ref{th:Nfrontslocexp}, and Appendix \ref{ap:ProofTh} (applied to the case of $N=2$-front patterns) for a precise statement of this existence and stability result (and its proof).
\\ \\
We conclude that there are no stable stationary two-fronts in this case. As we have already seen in Fig. \ref{fig:NumericsTwoThreeFronts-Intro}(b), this does not imply that the two-front dynamics of \eqref{eq:mainEquation} with a localized topography is similar to that of the standard homogeneous Allen-Cahn equation \eqref{eq:standardAC}.
Fronts in \eqref{eq:mainEquation} with localized topography may move (exponentially) slowly to $\pm \infty$, while the two fronts of a two-front pattern in \eqref{eq:standardAC} attract each other (see \eqref{eq:intODE2loc} with $\mathcal{S}(\psi) \equiv 0$) and subsequently annihilate each other, see Fig. \ref{fig:NumericsTwoThreeFronts-Intro}(a).
\\ \\
In the richest case with $h_- > 0$, $h_+ < 0$, and $\mathcal{S}'(\psi_\ast) > 0$ (so that the  stationary one-front pattern associated to $\psi_\ast$ is stable), there are various scenarios -- depending on the initial conditions in \eqref{eq:intODE2loc}, or, equivalently, the initial positions of the two fronts:
\\ \\
{\it (i)} The $\psi_1$-front travels towards $-\infty$, the $\psi_2$-front travels towards $+\infty$ (see Figs. \ref{fig:NumericsTwoThreeFronts-Intro}(b), \ref{fig:4FrontDynamics}(b)).
\\
{\it (ii)} One of the fronts moves to a stable $\psi_\ast$ one-front position, the other to $+ \infty$ or to $- \infty$ (see Fig. \ref{fig:Coarsening}(d)).
\\
{\it (iii)} The two fronts attract and annihilate each other (see Fig. \ref{fig:4FrontDynamics}(d)).
\\
{\it (iv)} The two fronts travel as a(n evolving) pair -- or `train' -- to either $+\infty$ or $-\infty$.
\\ \\
(Note that the above given references to numerical simulations mostly concern simulations of \eqref{eq:mainEquation} with localized topographies (\ref{eq:Ftopography}) with algebraically decaying tails and that the prescribed behavior often only occurs after several pairs of fronts of the initial $N$-front patterns (with $N > 2$) have merged.) By considering the evolution of a solitary front $\psi(\tau)$ that starts out with $\psi(0) = \psi_0 \gg 1$, {\it i.e.}, at the positive tail of $\mathcal{S}(\psi)$, it is natural to conjecture that the fourth scenario cannot occur. Assuming that $h_+ < 0$, we see that the evolution of $\psi(\tau)$ is (to leading order) given by
\begin{equation}
\label{eq:1frontfareqexp}
\frac{d\psi}{d\tau} =  - \varepsilon \mathcal{S}(\psi) = \varepsilon |h_+| w_+ e^{-\mu \psi}
\end{equation}
so that,
\begin{equation}
\label{eq:1frontfarsolexp}
\psi(\tau) = \frac{1}{\mu} \log \left(\varepsilon |h_+| w_+ \tau + e^{\mu \psi_0} \right).
\end{equation}
As a consequence, the leading order distance between two non-interacting fronts, $\tilde{\psi}(\tau) < \psi(\tau)$ that start out sufficiently far from each other, evolves as
\begin{equation}
\label{eq:deltafrontsexp}
\psi(\tau) - \tilde{\psi}(\tau) =
\frac{1}{\mu} \log \left(\frac{\varepsilon |h_+| w_+ \tau + e^{\mu \psi_0}}{\varepsilon |h_+| w_+ \tau  + e^{\mu \tilde{\psi}_0}} \right) = \frac{e^{\mu \psi_0} - e^{\mu \tilde{\psi}_0}}{\varepsilon \mu |h_+| w_+ \tau} + \; \; {\rm h.o.t.~for} \; \; \tau \gg \frac{e^{\mu \psi_0}}{\varepsilon}.
\end{equation}
Thus, on a (very) long timescale, the $\tilde{\psi}$-front is expected to catch up with its forerunner $\psi$. Moreover, when the distance between $\psi$ and $\tilde{\psi}$ becomes sufficiently `small', the additional attractive front-interaction effect of the homogeneous Allen-Cahn equation will kick in, so that the distance between $\psi$ and $\tilde{\psi}$ decreases even faster.
\\ \\
However, this reasoning implicitly assumes that the impact of the  higher order terms $\varepsilon^2 R^{j,2}_2(\psi_1, \psi_2;\varepsilon)$ will always remain higher order, and thus negligible. Moreover, the simulations shown in Fig. \ref{fig:5Fronts-Intro}(c) and in (upcoming) Fig. \ref{fig:5FrontDynVaryingP}(c) indicate that -- on long temporal and spatial scales -- persistent traveling `trains' consisting of slowly dispersing two-, three- and five-front patterns may exist. However, these simulations are done for \eqref{eq:mainEquation} with topographical inhomogeneities \eqref{eq:Ftopography} $H(x) = H_{\rm alg}(x; p)$ with $p \in (0,1)$, so that the $H'(x)$ term in \eqref{eq:Ftopography} is not bounded (recall \eqref{eq:defHuni-alg}). Thus, it is not clear whether or not there are circumstances for which the `trains' mentioned under {\it (iv)} can exist. We refer to Sec.~\ref{sec:Discussion} for a discussion.
\\ \\
Finally, we note that if one of the conditions $h_- > 0$, $h_+ < 0$, $\mathcal{S}'(\psi_\ast) > 0$ does not hold, some of the above scenarios are no longer possible, or have to be adapted. For instance, if $h_+ > 0$ (but still $h_- > 0$, $\mathcal{S}'(\psi_\ast) > 0$), fronts cannot travel towards $+ \infty$, so that scenario {\it (i)} will not occur and scenarios {\it (ii)} and {\it (iv)} only take place with a front(s) moving to $-\infty$. Moreover, if neither of these conditions hold, the localized topography will not have impact on the final outcome of the front interactions, but it may speed up the coarsening process -- see Sec.~\ref{sec:Discussion} and Fig.~\ref{fig:Coarsening}.

\subsection{Three- and four-front dynamics}
\label{sec:ODE34frontloc}
As a bridge towards the general case of interacting $N$-fronts of Sec.~\ref{sec:ODENfrontloc}, we now consider several aspects of the interactions of three-fronts and four-fronts (where we again for simplicity only consider exponentially decaying topographies $H(x)$). Our main goal is to indicate by way of the $N=3,4$-example systems how all stationary $N$-front patterns can be seen as extensions of four kinds of basic patterns, the three introduced in the previous section and a fourth one (that can be seen as a combination of the second and third). Like in the previous section, it can be shown hat the higher order terms of \eqref{eq:dynNfronts} can (a priori) not be neglected for $\mu > 1$. Therefore, we will only consider the case $0 < \mu < 1$ in this section.
\\  \\
By \eqref{eq:dynNfronts}, the (leading order) interaction ODE for three-front patterns (that approach $u_\pm^\varepsilon(x) = \pm 1 + \mathcal{O}(\varepsilon)$ \eqref{eq:uMinusPlus} as $x \to \pm \infty$) in \eqref{eq:mainEquation} with driving term \eqref{eq:Ftopography} is given by
\begin{equation}
\left\{
\begin{array}{r c c c l}
\frac{d\psi_1}{d\tau} & = & - \varepsilon \mathcal{S}(\psi_1) & + & 16 e^{-(\psi_2 - \psi_1)} \\
\frac{d\psi_2}{d\tau} & = & - \varepsilon \mathcal{S}(\psi_2) & + & 16 \left(e^{-(\psi_3 - \psi_2)} - e^{-(\psi_2 - \psi_1)}\right) \\
\frac{d\psi_3}{d\tau} & = & - \varepsilon \mathcal{S}(\psi_3) & - & 16 e^{-(\psi_3 - \psi_2)}.
\end{array}
\right.
\label{eq:intODE3loc}
\end{equation}
Similarly, the system of ODEs governing the front locations in four-front patterns (homoclinic to $u_-^\varepsilon(x) = -1 + \mathcal{O}(\varepsilon)$) is given by
\begin{equation}
\left\{
\begin{array}{r c c c l}
\frac{d\psi_1}{d\tau} & = & - \varepsilon \mathcal{S}(\psi_1) & + & 16 e^{-(\psi_2 - \psi_1)} \\
\frac{d\psi_2}{d\tau} & = & - \varepsilon \mathcal{S}(\psi_2) & + & 16 \left(e^{-(\psi_3 - \psi_2)} - e^{-(\psi_2 - \psi_1)}\right) \\
\frac{d\psi_3}{d\tau} & = & - \varepsilon \mathcal{S}(\psi_3) & + & 16 \left(e^{-(\psi_4 - \psi_3)} - e^{-(\psi_3 - \psi_2)}\right) \\
\frac{d\psi_4}{d\tau} & = & - \varepsilon \mathcal{S}(\psi_4) & - & 16 e^{-(\psi_4 - \psi_3)}.
\end{array}
\right.
\label{eq:intODE4loc}
\end{equation}
An immediate consequence of the form of these ODEs is that critical points $(\bar{\psi}_1,\bar{\psi}_2, \bar{\psi}_3)$ of \eqref{eq:intODE3loc}, respectively $(\bar{\psi}_1,\bar{\psi}_2, \bar{\psi}_3,\bar{\psi}_4)$ of \eqref{eq:intODE4loc}, must have $\bar{\psi}_1 \ll \bar{\psi}_2 \ll \bar{\psi}_3$, resp. $\bar{\psi}_1 \ll \bar{\psi}_2 \ll \bar{\psi}_3 \ll \bar{\psi}_4$, and that it is thus natural to adapt the `anchors' $p_j$ of the fronts \eqref{eq:anchor} in a suitable fashion, according to the character of the critical point under consideration (see \eqref{eq:anchors4ex} below).
\\ \\
Since the number of fronts is also even in the case of four-front patterns, the 4-fronts case is somewhat more similar to the two-fronts case, and we therefore start by studying the four-dimensional system \eqref{eq:intODE4loc}. To determine the critical point $(\bar{\psi}_1,\bar{\psi}_2, \bar{\psi}_3,\bar{\psi}_4)$ that is a natural extension of the critical point of the first kind of \eqref{eq:intODE2loc}, we first look for critical points with with two fronts located at the left and two at the right tail of the heterogeneity, {\it i.e}. with $\bar{\psi_1} \ll \bar{\psi_2} \ll -1 < 1 \ll \bar{\psi_3} \ll \bar{\psi_4}$. For this, we follow \eqref{eq:Ansat2critpts1}, introduce $\nu^\pm_{1,2}$ -- with $0 < \nu^\pm_{2} < \nu^\pm_{1}$ -- and $\bar{\ell}^\pm_{1,2}$, and set
\begin{equation}
\bar{\psi}_1 = -\nu_1^- |\log \varepsilon| + \bar{\ell}_1^- \ll
\bar{\psi}_2 = -\nu_2^- |\log \varepsilon| + \bar{\ell}_2^- \ll -1 < 1 \ll
\bar{\psi}_3 = \nu_2^+ |\log \varepsilon| + \bar{\ell}_2^+ \ll
\bar{\psi}_4 = \nu_1^+ |\log \varepsilon| + \bar{\ell}_1^+,
\label{eq:Ansat4critpts1}
\end{equation}
where we note that indeed $\bar{\psi}_2 - \bar{\psi}_1 \gg 1$, $\bar{\psi}_3 - \bar{\psi}_2 \gg 1$, $\bar{\psi}_4 - \bar{\psi}_3 \gg 1$. It thus follows that (to leading order),
\begin{equation}
\begin{array}{lcr}
|\mathcal{S}(\bar{\psi}_1)| = \varepsilon^{\mu \nu_1^-} |h_-| w_- e^{\mu \bar{\ell}_1^-} & \ll &
\varepsilon^{\mu \nu_2^-} |h_-| w_- e^{\mu \bar{\ell}_2^-} =  |\mathcal{S}(\bar{\psi}_2)| \\
|\mathcal{S}(\bar{\psi}_4)| = \varepsilon^{\mu \nu_1^+} |h_+| w_+ e^{-\mu \bar{\ell}_1^+} & \ll &
\varepsilon^{\mu \nu_2^+} |h_+| w_+ e^{-\mu \bar{\ell}_2^+} =  |\mathcal{S}(\bar{\psi}_3)|
\end{array}
\label{eq:exp41R}
\end{equation}
(see \eqref{eq:exp21Rpsi12}). Since, by \eqref{eq:intODE4loc}, $\left|16 e^{-(\bar{\psi_2}-\bar{\psi_1})} \right| = |\varepsilon S(\bar{\psi}_1)| \ll |\varepsilon S(\bar{\psi}_2) |$ and $\left| 16 e^{-(\bar{\psi_4}-\bar{\psi_3})} \right| = |\varepsilon S(\bar{\psi}_4)| \ll |\varepsilon S(\bar{\psi}_3)|$, and $\mathcal{S}(\bar{\psi}_1) + \mathcal{S}(\bar{\psi}_2) + \mathcal{S}(\bar{\psi}_3) + \mathcal{S}(\bar{\psi}_4) = 0$, we deduce that to leading order
\begin{equation}
\mathcal{S}(\bar{\psi}_2) = - \mathcal{S}(\bar{\psi}_3) =
\frac{16}{\varepsilon} e^{-(\bar{\psi}_3 - \bar{\psi}_2)},
\label{eq:crit423loc}
\end{equation}
which is identical to \eqref{eq:crit2loc}. Hence, we conclude that $\nu_2^- = \nu_2^+ = \nu_2 = 1/(2-\mu)$ and that $\bar{\ell}_2^-$ and $\bar{\ell}_2^+$ can be determined uniquely if (and only if) $h_+ < 0 < h_-$. It should be noted however that the two-front (sub)pattern thus formed by $\bar{\psi}_2$ and $\bar{\psi}_3$ is a reflected version of that in Sec.~\ref{sec:ODE2frontloc}: at $\bar{\psi}_2$, the pattern `jumps' from $+1$ to $-1$ and back to $+1$ at $\bar{\psi}_2$, while all two-front patterns in Sec.~\ref{sec:ODE2frontloc} jump from $-1$ to $+1$ back to $-1$ (by construction). This reflection is a completely symmetric operation and has no influence on existence nor stability. Under the assumption that indeed $h_+ < 0 < h_-$, we next substitute the outcome of \eqref{eq:crit423loc} in the $\frac{d\psi_1}{d\tau} =0$ and $\frac{d\psi_4}{d\tau}=0$ equations, to obtain
\begin{equation}
\label{eq:psi14eqs}
\begin{array}{r c c c c}
0 & = & -\varepsilon \mathcal{S}(\bar{\psi}_1) + 16 e^{-\bar{\psi_2}} e^{\bar{\psi_1}}
& = & - \varepsilon^{1+\mu \nu_1^-} h_- w_- e^{\mu \bar{\ell}_1^-} + 16 \varepsilon^{\nu_1^- - \frac{1}{2-\mu}} e^{-(\bar{\ell}_2^- - \bar{\ell}_1^-)}
\\
0 & = &  -\varepsilon \mathcal{S}(\bar{\psi}_4) - 16 e^{\bar{\psi_3}} e^{-\bar{\psi_4}}
& = & - \varepsilon^{1+\mu \nu_1^+} h_+ w_+ e^{-\mu \bar{\ell}_1^+} - 16 \varepsilon^{\nu_1^+ - \frac{1}{2-\mu}} e^{(\bar{\ell}_2^+ - \bar{\ell}_1^+)},
\end{array}
\end{equation}
find
\begin{equation}
\nu_1^- = \nu_1^+ := \nu_1 = \frac{3-\mu}{(1-\mu)(2-\mu)},
\label{eq:loc4nu1}
\end{equation}
and conclude that $\bar{\ell}_{1,2}^-$ can again be determined uniquely (since $h_+ < 0 < h_-$). Once again, we note that it is natural -- and even necessary -- to adapt our choices for the anchors $p_j$ \eqref{eq:anchor} to the patterns under consideration:
\begin{equation}
\label{eq:anchors4ex}
p_1 = - \nu_ 1 < p_2 = -\nu_2  < p_3 = \nu_2 < p_4 = \nu_1,
\end{equation}
which also implies that the lengths of two bounded interval $I_2$ and $I_3$ becomes $\frac{1}{2}(\nu_1 + \nu_2) = 1/(1-\mu)$ \eqref{eq:defIjs}, thus longer than the standard length $1$ (see \eqref{eq:phiupdown}).
\\ \\
Especially since the $\psi_1$- and $\psi_4$-fronts are now deeper into the tails of $\mathcal{S}(\psi)$ compared to those of the stationary two-front (of the first kind), we again need to check whether we can still neglect the higher order contributions to \eqref{eq:intODE4loc}. Following the approach of Sec.~\ref{sec:ODE2frontloc}, we denote the correction terms to the $\frac{d \psi_j}{d \tau}$ equation by $\varepsilon^2 R^{j,4}_{2}(\psi_1, \psi_2, \psi_3, \psi_4;\varepsilon)$ ($j=1,2,3,4$). We need to ensure that the stationary four-fronts can be constructed in circumstances under which the $\mathcal{O}(\varepsilon^2)$ correction terms $\varepsilon^2 R^{j,4}_{2}(\psi_1, \psi_2, \psi_3, \psi_4;\varepsilon)$ -- i.e. that satisfy $\varepsilon^2 |R^{j,4}_{2}(\psi_1, \psi_2, \psi_1, \psi_2;\varepsilon)| < \varepsilon^2 C$ for some $C > 0$ -- are indeed of higher order. Thus, we need to check whether $\nu_1 - \nu_2 =  2/((1-\mu)(2-\mu)) < 2$: if this is the case we know that there is no possible leading order contribution by correction terms we did not take into account. It follows that the above leading order analysis is only valid for $\mu \in (0, \frac12(3-\sqrt{5})) \approx (0,0.3819..)$. In other words, like for the critical points of the second and third kind in \eqref{eq:intODE2loc}, the original $\mu$-interval $(0,1)$ must be shrunk into a smaller interval in order to establish the existence of this kind of stationary four-front patterns. As for the critical points of the second kind discussed in the previous section, it is natural to expect that this domain may be extended, but a full higher order analysis would need to be performed to (possibly) come to such a conclusion (which, again, we refrain from doing in this study).
\\ \\
Like the existence, also the stability of the stationary four-front pattern is determined by the two-front that underlies it. By the relative magnitudes \eqref{eq:exp41R}, we can write down the leading order part of $4\times4$ matrix associated to the stability of $(\bar{\psi}_1,\bar{\psi}_2, \bar{\psi}_3,\bar{\psi}_4)$,
\begin{equation}
\label{eq:leadingmatrix4}
\left(
\begin{array}{cccc}
0 & 0 & 0 & 0
\\
0 & - \varepsilon \mathcal{S}'(\bar{\psi}_2) + E_2  & - E_2 & 0
\\
0 & - E_2 & - \varepsilon \mathcal{S}'(\bar{\psi}_3) + E_2 & 0
\\
0 & 0 & 0 & 0
\end{array}
\right)
\end{equation}
(see \eqref{eq:defbarEj}, \eqref{eq:stab2loc1}). Hence, the two (relatively) large eigenvalues associated to the stability of $(\bar{\psi}_1,\bar{\psi}_2, \bar{\psi}_3,\bar{\psi}_4)$ are given by
\begin{equation}
\label{eq:lamba41loc}
\lambda_1 =  -\mu \varepsilon^{\frac{2}{2-\mu}} |h_+| w_+ e^{-\mu \bar{\ell}_2^+} < 0, \; \; \ \ \ \ \
\lambda_2 =  (2-\mu)\varepsilon^{\frac{2}{2-\mu}} |h_+| w_+ e^{-\mu \bar{\ell}_2^+} > 0
\end{equation}
(see \eqref{eq:lamba21loc}).
The associated stationary four-front pattern in (\ref{eq:mainEquation}) is unstable. (We refrain from determining the other two smaller eigenvalues associated to the $\psi_1$-front and the $\psi_4$-front.)
\\ \\
Thus, we conclude that the kind of stationary four-front patterns can be regarded as straightforward extensions of the stationary two-front pattern of the first kind considered in Sec.~\ref{sec:ODE34frontloc}: both the existence and the stability of these four-front patterns is determined by its two-front core.
\\ \\
Instead of considering critical points of \eqref{eq:intODE4loc} that correspond to the critical points of \eqref{eq:intODE2loc} of the second and third kind, we now first consider a new class of critical points $(\bar{\psi}_1,\bar{\psi}_2, \bar{\psi}_3)$ of \eqref{eq:intODE3loc} -- the critical points of the fourth kind --  that can be seen as combinations of critical points of the second and third kind,
\\ \\
$\bullet$ critical points of the fourth kind have $\bar{\psi}_1 \ll -1 < 1 \ll \bar{\psi}_3$ and $\bar{\psi}_2 = \mathcal{O}(1)$.
\\ \\
Thus, we combine \eqref{eq:Ansat2critpts1} and \eqref{eq:Ansat2critpts2} and set
\begin{equation}
\bar{\psi}_1 = -\nu_- |\log \varepsilon| + \bar{\ell}_-, \; \; \ \ \ \
\bar{\psi}_2 = \psi_\ast + \varepsilon^{\nu_0} \tilde{\psi}_0, \; \;  \ \ \ \
\bar{\psi}_3 = \nu_+ |\log \varepsilon| + \bar{\ell}_+.
\label{eq:Ansat3critpts1}
\end{equation}
Since also here $\mathcal{S}(\bar{\psi}_1) + \mathcal{S}(\bar{\psi}_2) + \mathcal{S}(\bar{\psi}_3) = 0$, we conclude by balancing the asymptotic magnitudes (in $\varepsilon$)
that $\nu_- = \nu_+ = \nu$ and $\nu_0 = \mu \nu$. It then follows from (for instance) the $\frac{d\psi_3}{d\tau}=0$ equation that
\[
\nu_-= \nu_+ = \nu = \frac{1}{1-\mu}, \; \; \ \ \ \
\nu_0 = \frac{\mu}{1-\mu},
\]
where we note that this leading order analysis is only valid for $0 < \mu < \frac12$ (as for the critical points of \eqref{eq:intODE3loc} of second and third kind). If (and only if) again $h_+ < 0 < h_-$ and assuming that $\mathcal{S}'(\psi_\ast) \neq 0$, then the resulting equations for $\bar{\ell}_{\pm}$ and $\tilde{\psi}_0$ can again be solved uniquely.
\\ \\
By defining
\[
M_\pm = \mp h_\pm w_\pm e^{\mp \mu \bar{\ell}_\pm} > 0,
\]
(by \eqref{eq:defwpm} and since $h_+ < 0 < h_-$), we can (after some algebra) write the $3\times3$ matrix associated to the stability of this critical point $(\bar{\psi}_1,\bar{\psi}_2, \bar{\psi}_3)$ as,
\[
\varepsilon^\frac{1}{1-\mu}
\left(
\begin{array}{ccc}
-\mu M_- + M_- & - M_- & 0 \\
-M_- & -\varepsilon^{-\frac{\mu}{1-\mu}} \mathcal{S}'(\psi_\ast) + M_- + M_+ & -M_+ \\
0 & -M_+ & -\mu M_+ + M_+
\end{array}
\right)
\]
This matrix has one (relatively) large eigenvalue and two smaller ones,
\[
\lambda_0 = - \mathcal{S}'(\psi_\ast) \varepsilon, \; \; \ \ \ \ \
\lambda_\pm = (1-\mu) M_\pm \varepsilon^\frac{1}{1-\mu} > 0
\]
(all to leading order and under the assumptions $h_+ < 0 < h_-$ and $\mathcal{S}'(\psi_\ast) \neq 0$). Thus, also this stationary multi-front pattern is unstable (with two or three unstable eigenvalues according to the stability of the stationary one-front associated to $\psi_\ast$).
\\ \\
We may now proceed and determine all possible stationary three-fronts, using the stationary two-fronts of the first, second and third kind introduced in Sec.~\ref{sec:ODE2frontloc} as foundation. We start out with the extended fronts of the first kind -- as we already did above for the four-front systems -- and consider critical points $(\bar{\psi}_1,\bar{\psi}_2, \bar{\psi}_3)$ of \eqref{eq:intODE3loc} that have all three fronts in the tails of $\mathcal{S}(\psi)$, {\it i.e.}, that either have $\bar{\psi}_1 \ll - 1 < 1 \ll \bar{\psi}_2 \ll \bar{\psi}_3$ or $\bar{\psi}_1 \ll \bar{\psi}_2 \ll - 1 < 1 \ll \bar{\psi}_3$. As for the above studied four-fronts, it follows from $\mathcal{S}(\bar{\psi}_1) + \mathcal{S}(\bar{\psi}_2) + \mathcal{S}(\bar{\psi}_3) = 0$ that (to leading order) $\mathcal{S}(\bar{\psi}_1) = -\mathcal{S}(\bar{\psi}_2)$ with $|\mathcal{S}(\bar{\psi}_3)| \ll  |\mathcal{S}(\bar{\psi}_{1,2})|$ in the former case and $\mathcal{S}(\bar{\psi}_2) = -\mathcal{S}(\bar{\psi}_3)$ with $|\mathcal{S}(\bar{\psi}_1)| \ll  |\mathcal{S}(\bar{\psi}_{2,3})|$ in the latter. Thus, in the first case, the positions of the $\bar{\psi}_1$- and $\bar{\psi}_2$-fronts coincide to leading order with those of the two-front of the first kind, with $\bar{\psi}_3$ determined by
\[
\varepsilon \mathcal{S}(\bar{\psi}_3) = 16 e^{\bar{\psi}_2} e^{-\bar{\psi}_3},
\]
similar to \eqref{eq:psi14eqs}. Under the condition $h_+ < 0 < h_-$ we thus find a uniquely determined three-front pattern, that can be seen as a three-front extension of the two-front pattern of the first kind. The second case goes along very similar lines, the three-front pattern is built around a two-front pattern of the first kind -- although this again is a reflected version of the two-front patterns considered in Sec.~\ref{sec:ODE2frontloc} (i.e. it connects $+1$ to $-1$ back to $+1$).
The fact that the `added' third front goes deeper into the tail of $\mathcal{S}(\psi)$ implies that the $\mathcal{O}(\varepsilon^2)$ correction term are only of higher order for $\mu$ in a certain subinterval if $(0,1)$ -- see Theorem \ref{th:Nfrontslocexp} and especially its proof in Appendix \ref{ap:ProofTh}.
\\ \\
As before, also the stability of these three-fronts is completely determined by that of the underlying two-front: the $3\times3$-matrices associated to the stability of the present three-front patterns have two (relatively) large eigenvalues that coincide (to leading order) with those of the underlying two-front pattern of the first kind (see \eqref{eq:leadingmatrix4}), \eqref{eq:lamba41loc}): both stationary three-front patterns are unstable.
\\ \\
The final kind of critical points $(\bar{\psi}_1,\bar{\psi}_2, \bar{\psi}_3)$ of \eqref{eq:intODE3loc} are points that have a two-front pattern of the second or third kind as underlying pattern (or a reflected version of these), {\it i.e.}, that either have $\bar{\psi}_1 \approx \psi_\ast \ll \bar{\psi}_2$ as in \eqref{eq:Ansat2critpts2} and $\bar{\psi}_2 \ll \bar{\psi}_3$ or $\bar{\psi}_1 \ll \bar{\psi}_2 \ll \bar{\psi}_3 \approx \psi_\ast$. Again we refer to Sec.~\ref{sec:ODENfrontloc}, Theorem \ref{th:Nfrontslocexp} and Appendix \ref{ap:ProofTh} for more systematic constructions, statements and proofs. It also again follows that these patterns must be unstable since the two large eigenvalues of the $3\times3$ stability matrices are determined to leading order by the underlying two-front patterns (of the second or third kind).
\\ \\
Thus, apart from the new three-front patterns of the fourth kind, all other three-front patterns can be seen as (relatively straightforward) extensions of the two-front patterns. We wrap up the analysis of stationary three-front patterns by concluding that if $\mathcal{S}(\psi)$ has one, respectively three, zeroes and $h_+ < 0 < h_-$ -- as for the unimodal hill \eqref{eq:defHuni-alg} with $\mu \in (0,\mu_{\rm PF})$, resp. $\mu \in (\mu_{\rm PF}, 1)$  (Fig. \ref{fig:SolHillPitchfork}) -- system \eqref{eq:intODE3loc} has $\mathcal{N}(3) = 5$, resp. $\mathcal{N}(3) = 11$, critical points. More specifically, there are two critical points that are extended two-fronts of the first kind and do not have a front near a zero of $\mathcal{S}(\psi)$ and three critical points that have a front near given a zero $\psi_\ast$ of $\mathcal{S}(\psi)$ -- the three-fronts that are extended two-fronts of the second and third kind and the new three-front of the fourth kind (so that $5 = 2 + 3 \times 1$ and $11 = 2 + 3 \times 3$ -- see \eqref{eq:NN}). All critical points/all possible stationary three-front patterns are unstable.
\\ \\
The stationary four-front patterns can all be constructed along the same lines, {\it i.e.}, as `(relatively straightforward) extensions' of the two-front patterns of the first, second and third kind and of the three-front pattern of the fourth kind. Without going into the details -- see Sec.~\ref{sec:ODENfrontloc} and Appendix \ref{ap:ProofTh} for the general set-up -- we note that for $h_+ < 0 < h_-$ there are three critical points of four-dimensional system \eqref{eq:intODE4loc} that do not have a front near a zero of $\mathcal{S}(\psi)$ and thus can be seen as extensions of the two-front patterns of the first kind (above we worked out the case $\bar{\psi}_1 \ll \bar{\psi}_2 \ll -1 < 1 \ll \bar{\psi}_3 \ll \bar{\psi}_4$, the cases $\bar{\psi}_1 \ll -1 < 1 \ll \bar{\psi}_2 \ll \bar{\psi}_3 \ll \bar{\psi}_4$ and $\bar{\psi}_1 \ll \bar{\psi}_2 \ll \bar{\psi}_3 \ll -1 < 1 \ll \bar{\psi}_4$ follow along similar lines). Any two-front of the second kind has a unique extension as stationary four-front (with $\bar{\psi}_1 = \mathcal{O}(1) \ll \bar{\psi}_2 \ll \bar{\psi}_3 \ll \bar{\psi}_4$), and the same holds (as usual) for the two-fronts of the third kind. Finally, each three-front of the fourth kind can be extended to a four-front in two ways (either by placing the fourth `far away' front to the left or to the right af the original three-front pattern). Thus, if $\mathcal{S}(\psi)$ has 1, respectively 3, zeroes then \eqref{eq:intODE4loc} has $\mathcal{N}(4) = 7$, resp. $\mathcal{N}(4) = 15$, critical points (see \eqref{eq:NN}). The stability of the associated stationary four-front patterns in \eqref{eq:mainEquation} is again fully determined by the underlying `basic' two- or three-front patterns. Since these are all unstable, it follows again that all critical points of system \eqref{eq:intODE4loc} -- and thus all stationary four-front patterns in \eqref{eq:mainEquation} -- are unstable.
\\ \\
As for the two-front case, if only one of the conditions $h_+ < 0$ or $h_- > 0$ hold, the number of critical points $\mathcal{N}(4)$ is equal to the number of zeroes of $\mathcal{S}(\psi)$ -- both for the three-front patterns as well as for the four-front patterns (and the three- and four-front patterns all are extensions of the two-front patterns of the second kind -- for $h_+ < 0$ -- or the third kind -- for $h_- > 0$). There are no critical points if neither of these conditions hold.
\\
\begin{figure}
\centering
\begin{subfigure}[t]{0.23\textwidth}
\centering
\includegraphics[width=\textwidth]{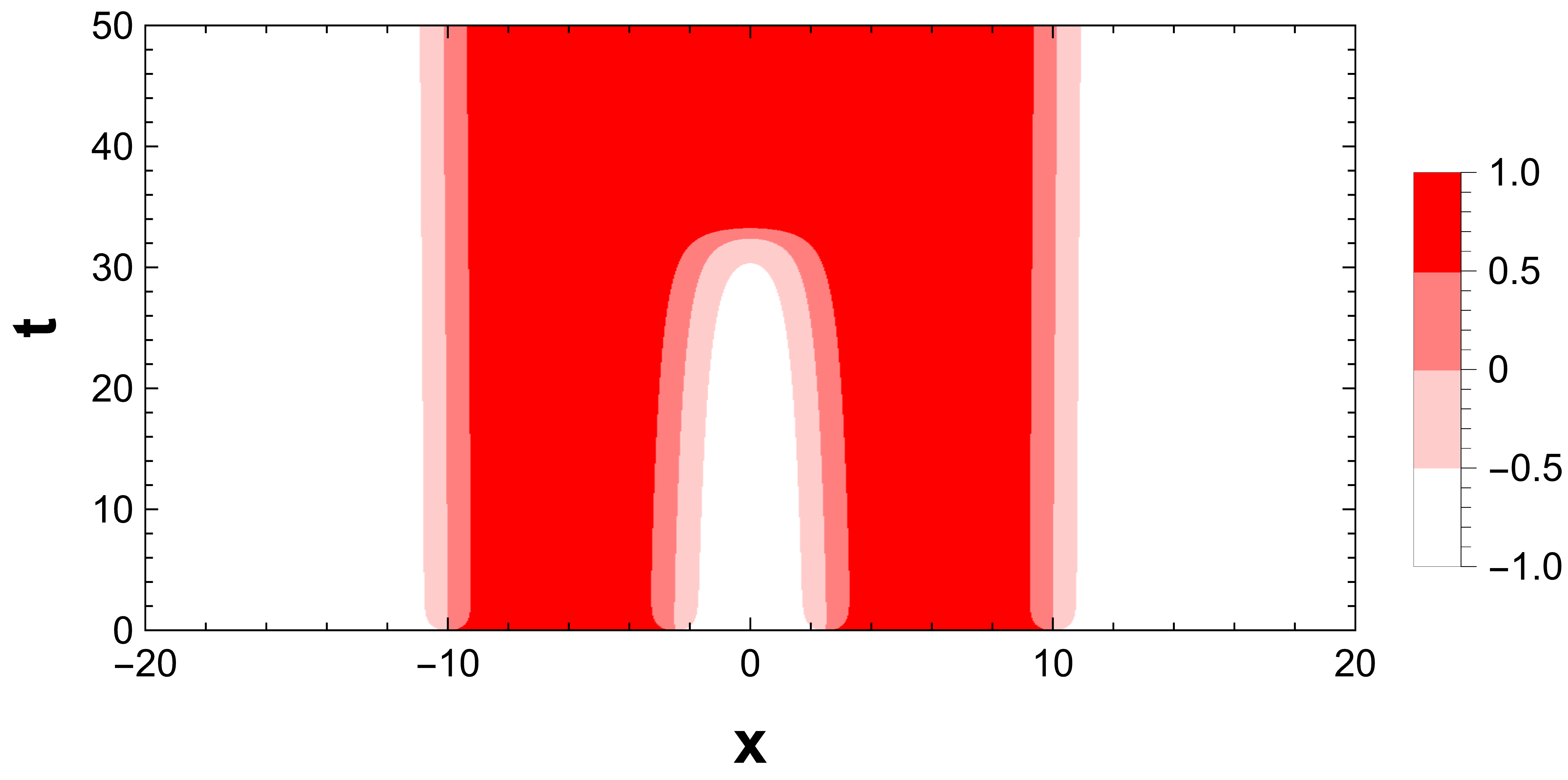}
\end{subfigure}
~
\begin{subfigure}[t]{0.23\textwidth}
\centering
\includegraphics[width = \textwidth]{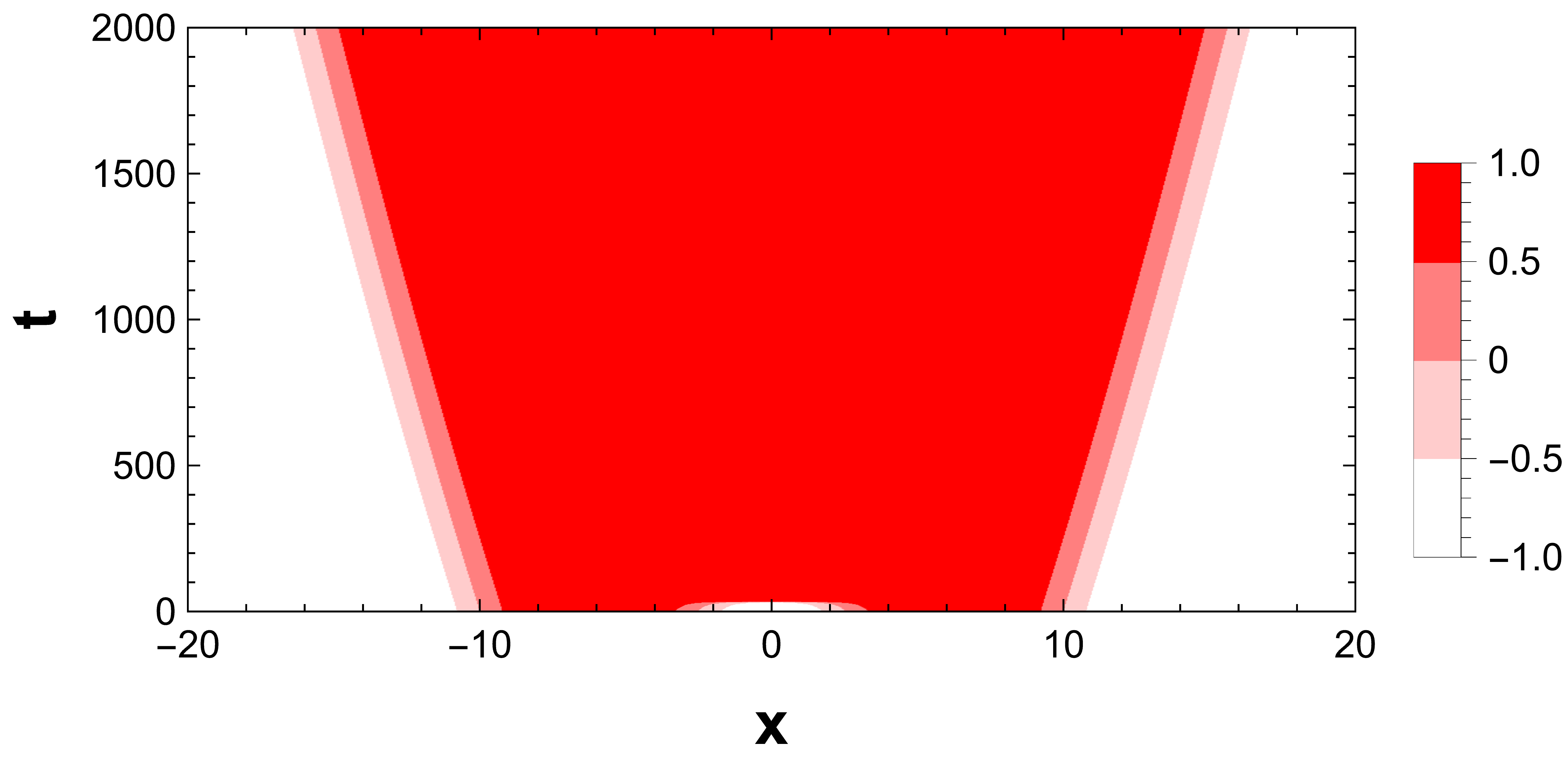}
\end{subfigure}
~
\begin{subfigure}[t]{0.23\textwidth}
\centering
\includegraphics[width = \textwidth]{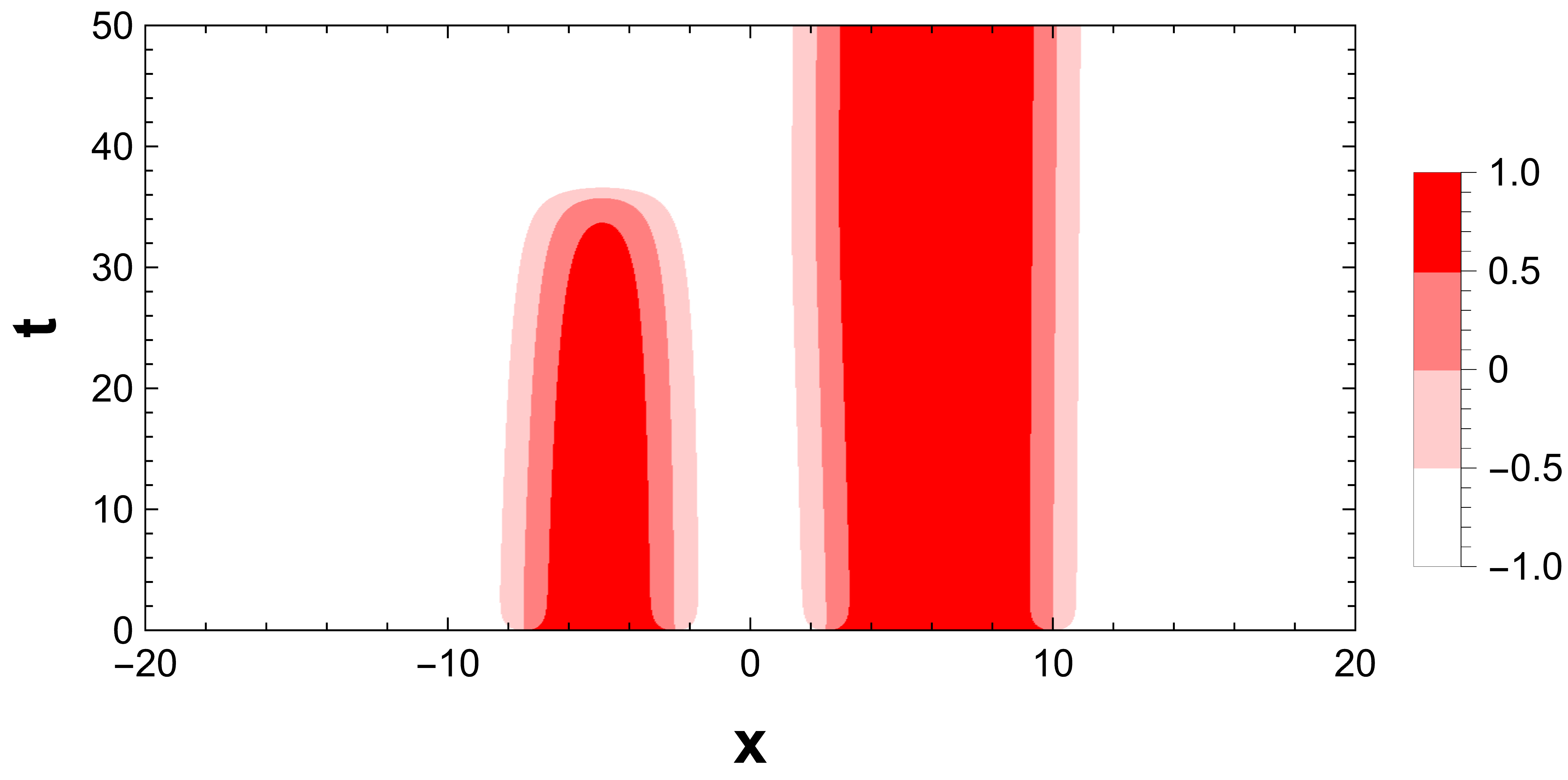}
\end{subfigure}
~
\begin{subfigure}[t]{0.23\textwidth}
\centering
\includegraphics[width = \textwidth]{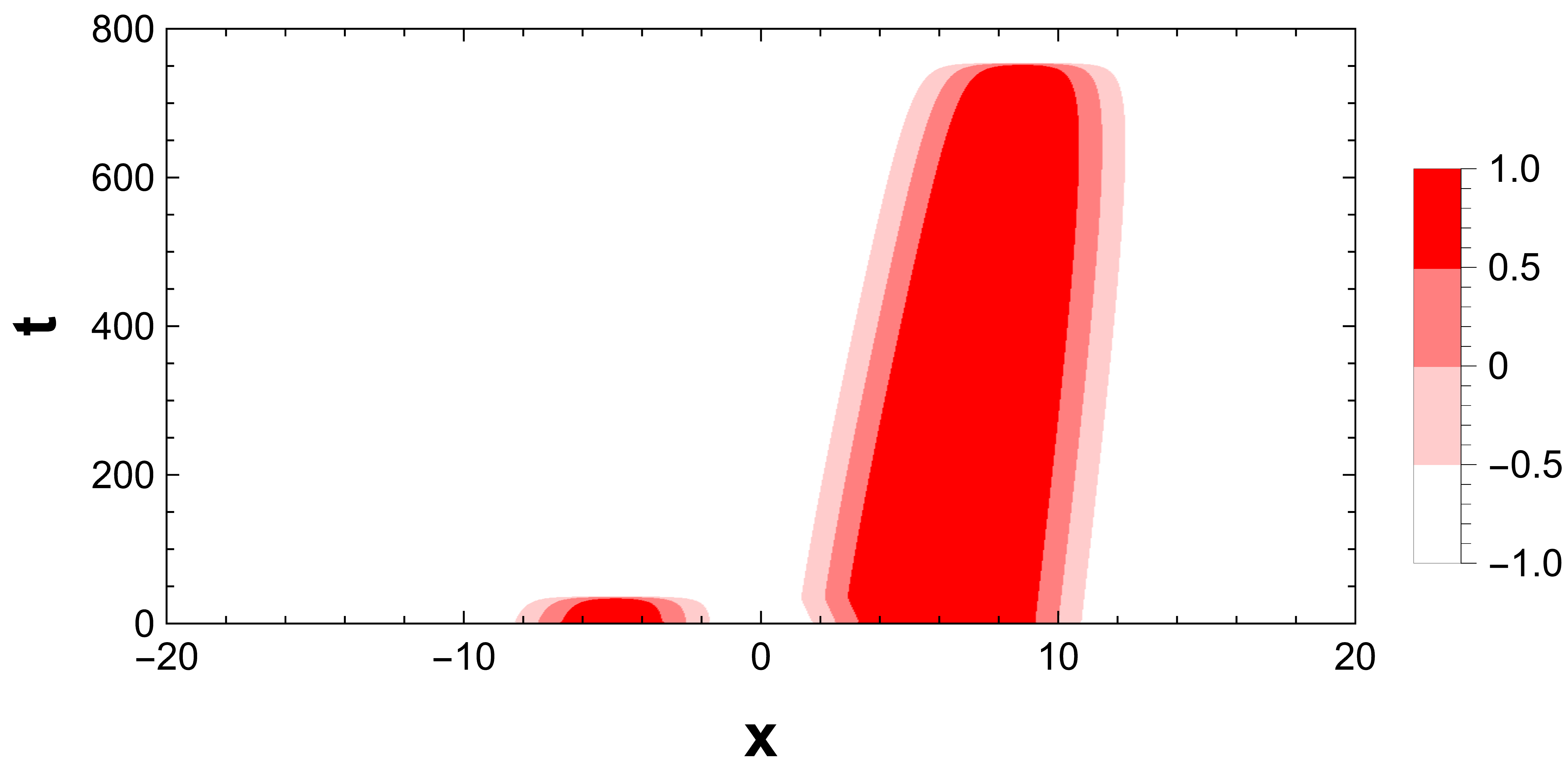}
\end{subfigure}
\caption{Simulations of \eqref{eq:mainEquation} for $x \in (-20,20)$ (and homogeneous Neumann boundary conditions) with topographical inhomogeneity \eqref{eq:Ftopography} given by $-H_{\rm alg}(x; 0.5)$ \eqref{eq:defHuni-alg} and $\varepsilon = 0.02$. (a)/(b) A four-front initial condition: $U(x,0) = -1 + \tanh(5 (x + 10))  - \tanh(5 (x + 2.5)) + \tanh(5 (x - 2.5)) - \tanh(5 (x - 10))$: the two middle fronts $\psi_2(\tau) < 0 < \psi_3(\tau)$ attract each other and merge, while the outer fronts travel away, {\it i.e.}, the initial four-front pattern $(\psi_1(\tau),\psi_2(\tau),\psi_3(\tau),\psi_4(\tau))$ evolves into a two-front pattern $(\psi_1(\tau),\psi_4(\tau))$ with $\psi_1(\tau) \to -\infty$ and $\psi_4(\tau) \to +\infty$ as $\tau \to \infty$ -- $t \in (0,50)$ in (a), $t \in (0,2000)$ in (b). (c)/(d) A four-front initial condition: $U(x,0) = -1 + \tanh(5 (x + 7.5))  - \tanh(5 (x + 2.5)) + \tanh(5 (x - 2.5)) - \tanh(5 (x - 10))$ that is equal to that of (a)/(b) except that $\psi_1(0)$ is slightly shifted to the right (from $-10$ to $-7.5$): now the two left hand side fronts $\psi_1(\tau)$ and $\psi_2(\tau)$ merge and both $\psi_{3,4}(\tau)$ of the resulting two-front pattern $(\psi_3(\tau),\psi_4(\tau))$ initially travel to the right until $\psi_{3}(\tau)$ catches up with $\psi_{4}(\tau)$ and also this pair annihilates each other -- $t \in (0,50)$ in (c), $t \in (0,800)$ in (d).}
\label{fig:4FrontDynamics}
\end{figure}
\\
For both the interacting three-fronts and the interacting four-fronts, the different possible scenarios can be worked out as was done for the case of two-front solutions in Sec.~\ref{sec:ODE2frontloc}.
The scenarios depend on the initial conditions, and also on the parity of the number of fronts considered.
We come back to this discussion in more detail in the general setting of the upcoming section.
We refer to Fig. \ref{fig:4FrontDynamics} for two four-front simulations of \eqref{eq:mainEquation} with only slightly different initial conditions that show very different behavior (where we used a localized topography $H_{\rm alg}(x; p)$ with algebraically decaying tails, instead of a $H_{\rm exp}(x; \mu)$ with exponentially decaying tails -- recall \eqref{eq:defHuni-alg}).

\subsection{N-front dynamics}
\label{sec:ODENfrontloc}

Our first general result on the existence and stability of stationary $N$-front dynamics -- that is a direct generalization of the insights of the previous sections -- can be stated as follows.

\begin{theorem}
\label{th:Nfrontslocexp}
Consider the weakly heterogeneous Allen-Cahn equation \eqref{eq:mainEquation} with inhomogeneity $F(U,U_x,x)$ given by the topography $H(x)$ \eqref{eq:Ftopography}. Assume that $H(x)$ is localized and that there is a $\mu > 0$ so that for $x \to \pm \infty$ $H'(x)$ decays exponentially like $h_\pm(\mu) e^{\mp \mu \sqrt{2} x}$ \eqref{eq:expdecayHx}. Assume furthermore that the Melnikov function $\mathcal{S}(\psi)$ associated to $H(x)$ \eqref{eq:Ex2R}, \eqref{eq:defpsiStau} has $K \geq 0$ non-degenerate zeroes $\psi_{\ast,k}$ ({\it i.e.}, $\mathcal{S}'(\psi_{\ast,k}) \neq 0$) and that $0 < \varepsilon \ll 1$ is sufficiently small. Let $\mathcal{N}(N)$ be the number of stationary $N$-front patterns of \eqref{eq:mainEquation} that approach $u_-^\varepsilon(x) = - 1 + \mathcal{O}(\varepsilon)$ as $x \to - \infty$ and $u_{(-1)^{N+1}}^\varepsilon(x) = (-1)^{N+1} + \mathcal{O}(\varepsilon)$ as $x \to +\infty$ \eqref{eq:uMinusPlus}.\\
Then $\mathcal{N}(1) = K$, and for $N \geq 2$ and $K$ odd,
\\
$\bullet$ if $h_+(\mu) > 0$ and $h_-(\mu) < 0$, then $\mathcal{N}(N) = 0$;
\\
$\bullet$ if $h_+(\mu) < 0 < h_-(\mu)$, then there is a $\mu_\ast(N) \in (0,1)$ such that for $\mu \in (0, \mu_\ast(N))$, $\mathcal{N}(N) = (K+1)N - 1$ (with $\mu_{\ast}(N) = 1 - 2^{-\frac{1}{N-1}} \in (0,1)$);
\\
$\bullet$ if either $h_+(\mu) < 0$ or $h_-(\mu) > 0$, then $\mathcal{N}(N) = K$ for $\mu \in (0, \mu_\ast(N))$.
\\
For $N \geq 2$ and $K$ even,
\\
$\bullet$  $\mathcal{N}(N) = K$ for $\mu \in (0, \mu_\ast(N))$.
\\
The one-front patterns are associated to the zeroes $\psi_{\ast,k}$ and are stable if $\mathcal{S}'(\psi_{\ast,k}) > 0$. All $N$-front patterns are unstable for $N \geq 2$.
\end{theorem}
An essentially identical result holds for stationary $N$-front patterns that connect $u_+^\varepsilon(x) = 1 + \mathcal{O}(\varepsilon)$ (as $x \to - \infty$) to $u_{(-1)^{N}}^\varepsilon(x) = (-1)^{N} + \mathcal{O}(\varepsilon)$ (as $x \to +\infty$). Note also that the condition that $\mathcal{S}(\psi)$ only has non-degenerate zeroes implies that the richest case $h_+ < 0 < h_-$ can only occur for $K$ odd: the signs of $h_+$ and $h_-$ must be equal if $\mathcal{S}(\psi)$ has an even number of zeroes.
\\ \\
The total number of stationary $N$-front patterns $\mathcal{N}(N)$ can be deduced from the observation in the previous sections that all multi-front patterns can be seen as direct extensions of only 4 kinds of basic patterns -- the first, second and third kinds two-front patterns of Sec.~\ref{sec:ODE2frontloc}, and the fourth three-front pattern of Sec.~\ref{sec:ODE34frontloc}. Starting out in the richest case, {\it i.e.}, assuming that $h_+ < 0 < h_-$ we note -- either by direct observations or by inductive arguments -- that there are $N-1$ distinct extensions of the two-front pattern of the first kind to a (stationary) $N$-front pattern (where it may be necessary to flip the original two-front pattern from one that approaches $-1$ as $x \to \pm \infty$ to one that approaches $+1$ as $x \to \pm \infty$). These patterns are different from all other stationary $N$-front patterns in the sense that they do not have a front asymptotically close to a zero of $\mathcal{S}(\psi)$. As was the case for the four-front patterns: each of the two-front patterns of the second and third kind can be extended uniquely to an $N$-front pattern (these patterns are studied in most analytical detail in Appendix \ref{ap:ProofTh}, since they give rise to the critical value $\mu_\ast(N)$). Similar to the two-front patterns of the first kind, all three-front patterns of the fourth kind can be extended in $N-2$ ways to $N$-front patterns. Thus, we conclude that for $h_+ < 0 < h_-$ (and $N \geq 2$), there can be in total
\begin{equation}
\label{eq:NN}
\mathcal{N}(N) = (N-1) + K[1 + 1 + (N-2)] = (K+1)N -1
\end{equation}
distinct $N$-front patterns. As before, $\mathcal{N}(N)$ reduces to only $K$ if either only $h_+ < 0$ or only $h_- > 0$, and to 0 if $h_+ > 0$ and $h_- < 0$.
\\ \\
The core of the proof of Theorem \ref{th:Nfrontslocexp} lies in establishing that each of these stationary $N$-front patterns indeed will exist if $\mu$ is sufficiently close to 0. As before, this is for a large part a study of the potential impact of the $\mathcal{O}(\varepsilon^2)$ terms $\varepsilon^2 R^{j,N}_{2}(\psi_1,...,\psi_N;\varepsilon)$ ($j=1,...,N$) that describe the higher order corrections to each of the $\frac{d \psi_j}{d \tau}$ equations. As in the previous sections, we only consider the leading order terms of (rewritten) $N$-front interaction ODE \eqref{eq:dynNfronts} in the proof of Theorem \ref{th:Nfrontslocexp}: the occurrence of the bounds $\mu_\ast(N)$ on $\mu$ is caused by the fact that the magnitudes of the leading order terms in the asymptotic analysis decrease towards the upperbounds $\varepsilon^2 C$ that can be established for each of these correction terms. Only by going deeper into the analysis of the correction terms we can draw conclusions on the existence of (all) stationary $N$-fronts for $\mu \geq \mu_\ast(N)$ -- which we have chosen not to do. (As follows from the proof of Theorem \ref{th:Nfrontslocexp} in Appendix \ref{ap:ProofTh}, the upperbound $\mu_\ast(N)$ on $\mu$ typically varies with the character of the $N$-front pattern, $\mu_\ast(N)$ is determined explicitly as minimum over each of these character-dependent upperbounds.) However, since $\mu_\ast(N) > 0$, we know by Theorem \ref{th:Nfrontslocexp} that stationary $N$-front patterns exist for exponentially decaying topographies $H(x)$ -- as long as this exponential decay is sufficiently weak. This observation directly motivates the next corollary, that considers $N$-front patterns for topographies with algebraically decaying tails.

\begin{figure}
\centering
\begin{subfigure}[t]{0.23\textwidth}
\centering
\includegraphics[width=\textwidth]{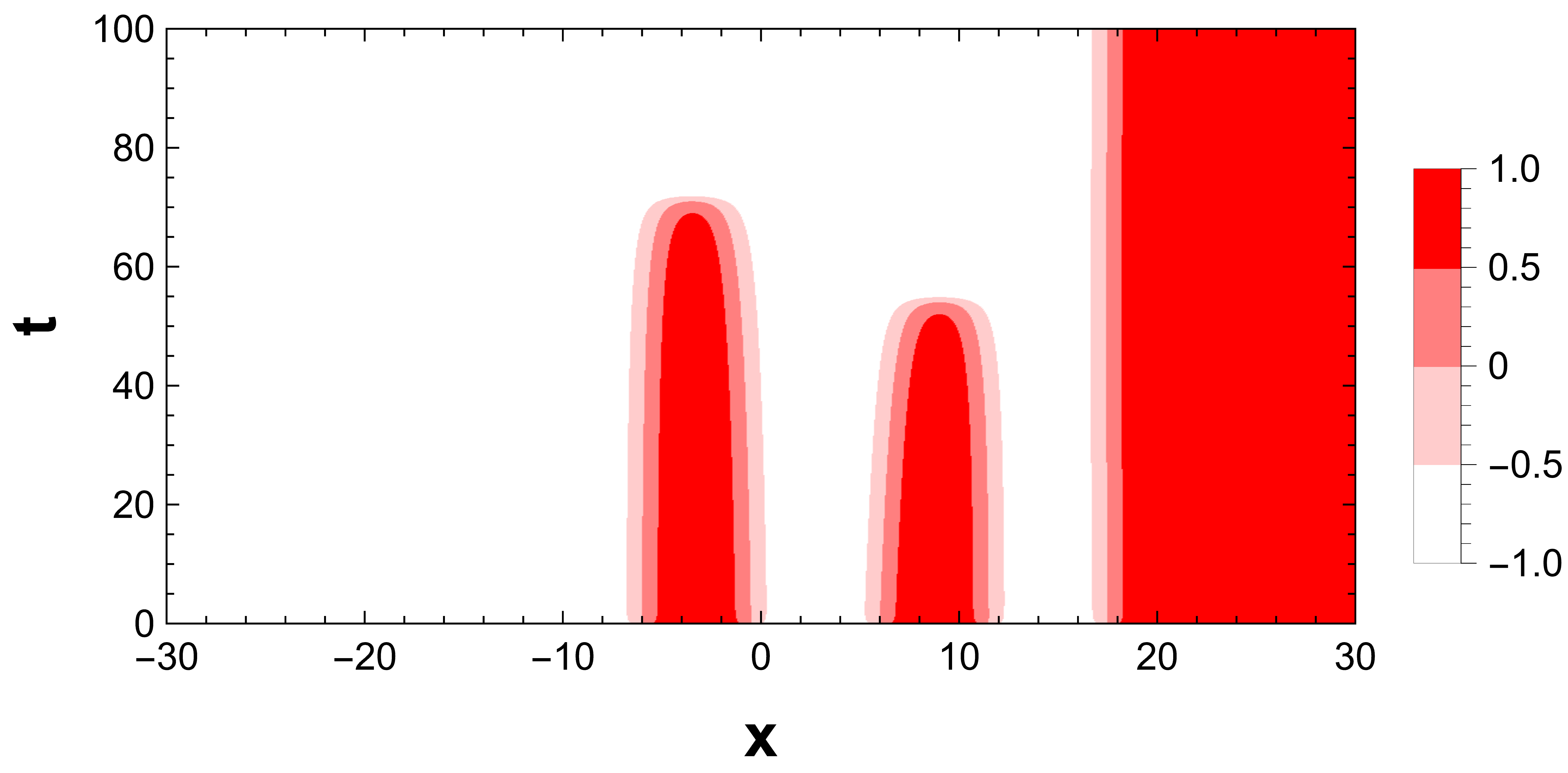}
\end{subfigure}
~
\begin{subfigure}[t]{0.23\textwidth}
\centering
\includegraphics[width = \textwidth]{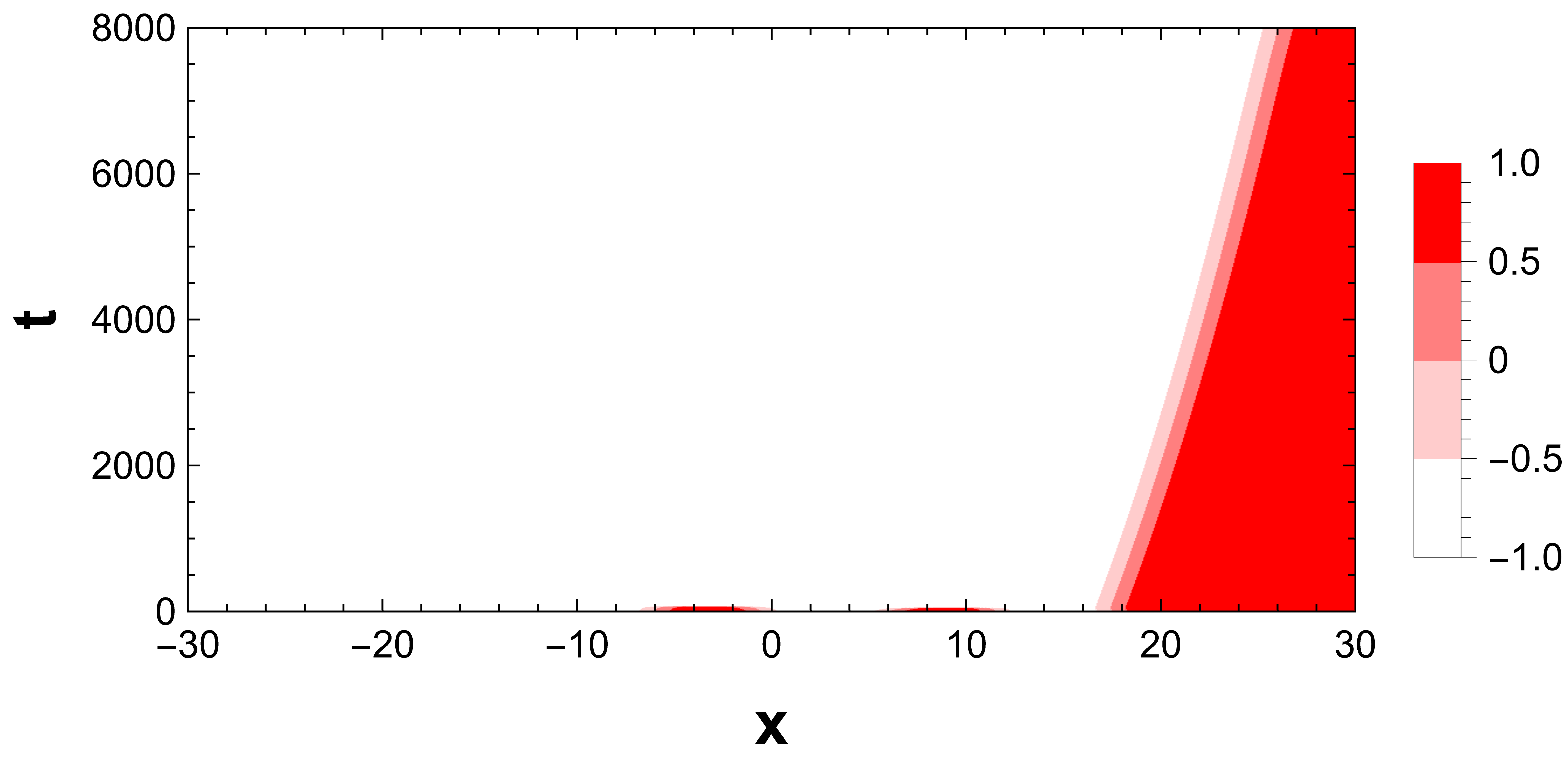}
\end{subfigure}
~
\begin{subfigure}[t]{0.23\textwidth}
\centering
\includegraphics[width = \textwidth]{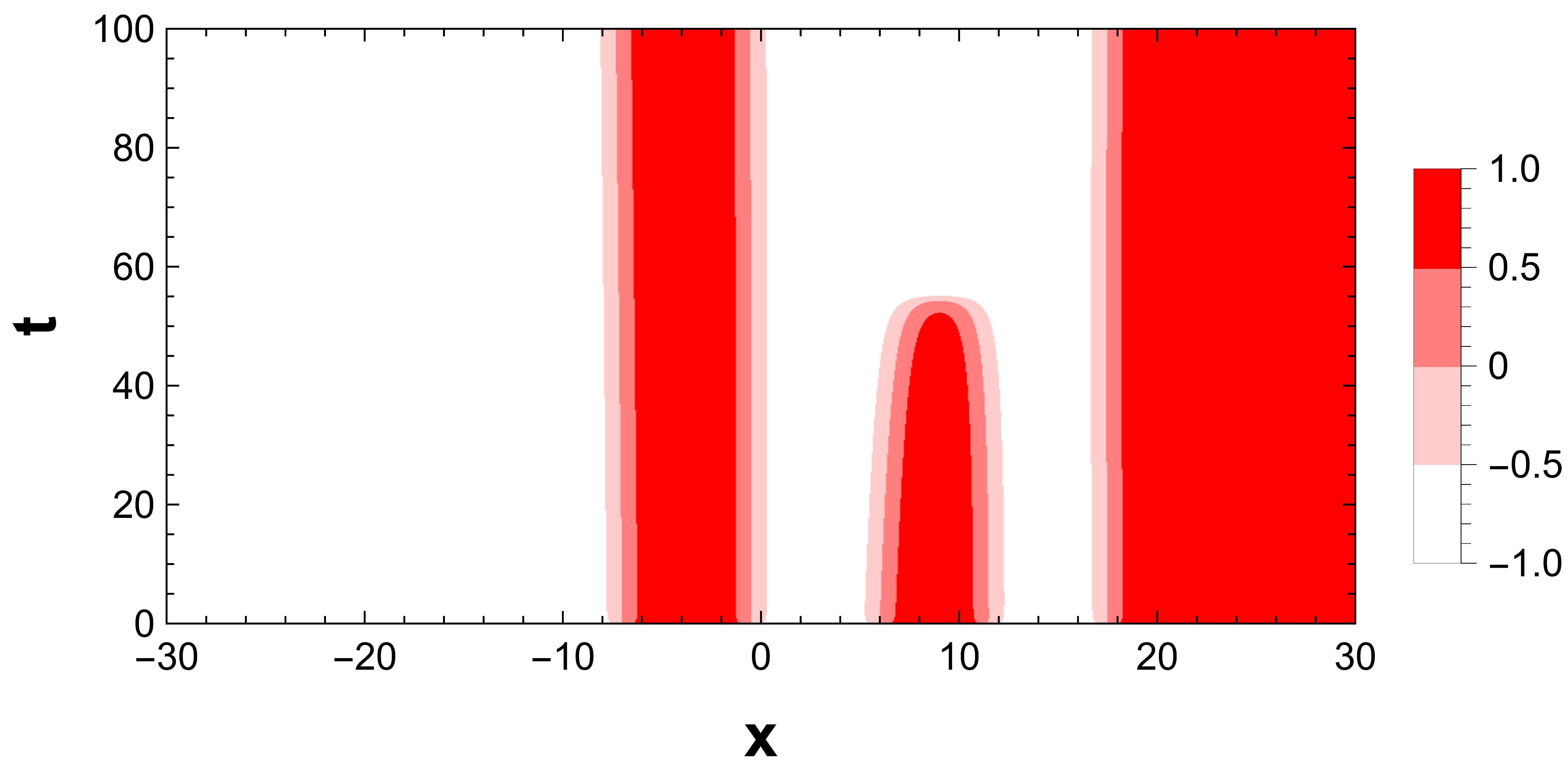}
\end{subfigure}
~
\begin{subfigure}[t]{0.23\textwidth}
\centering
\includegraphics[width = \textwidth]{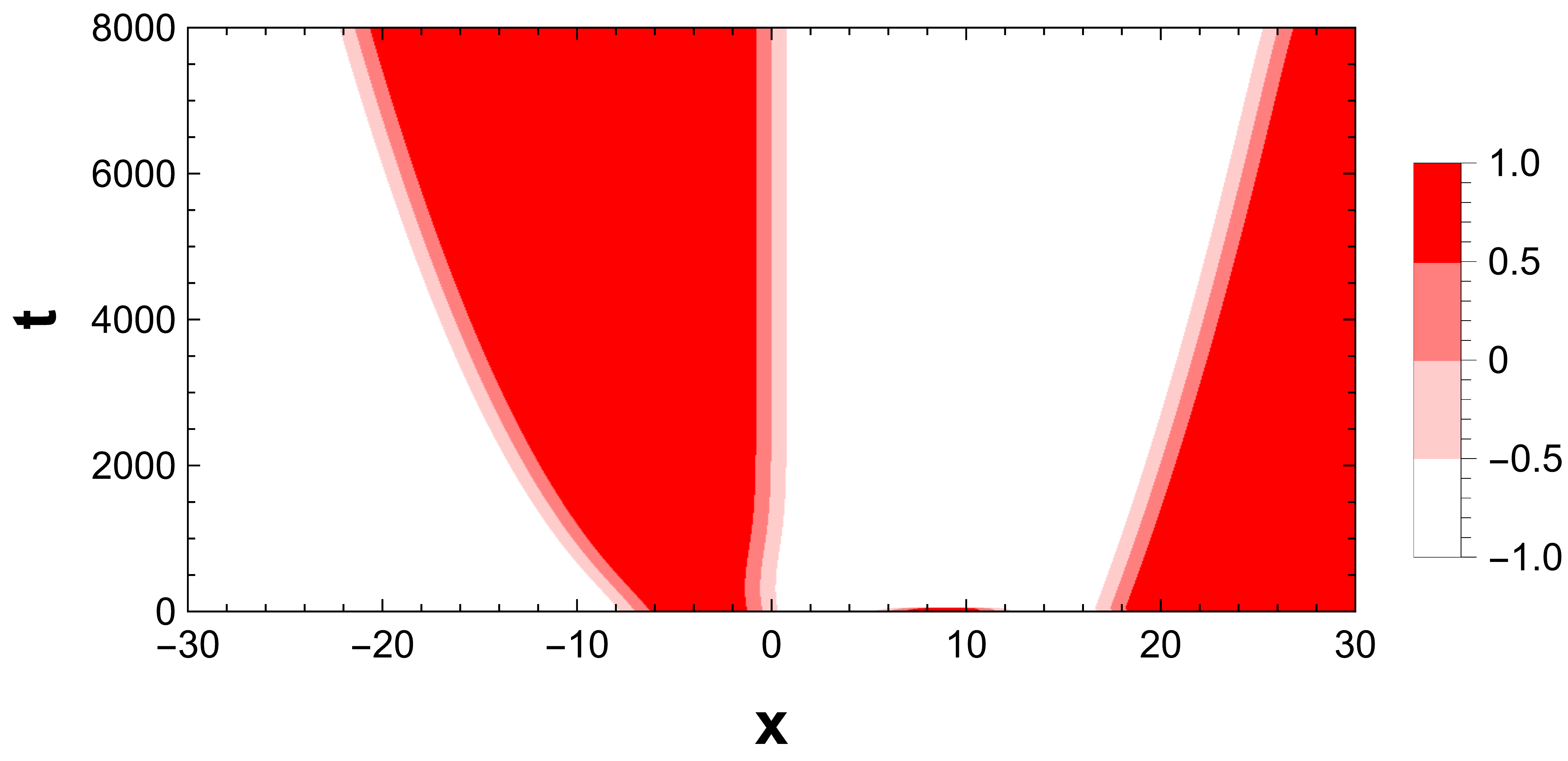}
\end{subfigure}
\caption{Simulations of \eqref{eq:mainEquation} for $x \in (-50,50)$ (and homogeneous Neumann boundary conditions) with topographical inhomogeneity \eqref{eq:Ftopography} given by $H_{\rm alg}(x; 1.25)$ \eqref{eq:defHuni-alg} and $\varepsilon = 0.2$. (a)/(b) A five-front initial condition: $U(x,0) = \tanh(x + 6)  - \tanh(x + 0.5) + \tanh(x - 6) - \tanh(x - 11.5) + \tanh(x - 17.5)$: first the $\psi_3$- and $\psi_4$-fronts merge, followed by the merging of $\psi_1$ and $\psi_2$; eventually, only the $\psi_5$-front persists, it travels with decaying speed as a solitary front to $+\infty$ -- $t \in (0,100)$ in (a), $t \in (0,8000)$ in (b). (c)/(d) Initially, the only difference with (a)/(b) is the position of $\psi_1$ -- that is placed more to the left compared to (a)/(b) -- $\psi_1(0) = -7$ instead of $-6$ in (a)/(b); after the $\psi_3$- and $\psi_4$-fronts have merged, $\psi_2$ settles at the stable stationary one-front location at $\psi = 0$, while $\psi_1$ travels (slowly) to $-\infty$ and $\psi_5$ to $+\infty$, both with decaying speeds -- $t \in (0,100)$ in (c)/$t \in (0,8000)$ in (d). Thus, changing the initial condition of $\phi_1$ changed the end configuration from one single front traveling to $+\infty$, to a three-front pattern with the middle one fixed at $\psi=0$ and the other two traveling to $\pm \infty$.}
\label{fig:5FrontDynamics}
\end{figure}

\begin{corollary}
\label{cor:Nfrontslocalg}
Consider the situation of Theorem \ref{th:Nfrontslocexp} but now for topographies $H(x)$ for which $H'(x)$ decays algebraically as $x \to \pm \infty$. Then, all statements of Theorem \ref{th:Nfrontslocexp} hold true -- with $\tilde{h}_\pm(p)$ \eqref{eq:algdecayHx} in the roles of $h_\pm(\mu)$ \eqref{eq:expdecayHx} -- without any further conditions on the rates of decay.
\end{corollary}
The proofs of Theorem \ref{th:Nfrontslocexp} and Corollary \ref{cor:Nfrontslocalg} are given in Appendix \ref{ap:ProofTh}. Here, we note that although the statement of Corollary \ref{cor:Nfrontslocalg} is less involved than that of Theorem \ref{th:Nfrontslocexp}, the asymptotic nature of the main properties of the $N$-front patterns -- {\it i.e.}, the positions of the fronts -- is technically more involved in the case of algebraic decay. Therefore, the proof of Corollary \ref{cor:Nfrontslocalg} builds directly on that of Theorem \ref{th:Nfrontslocexp}.
\\ \\
Even for the homogeneous Allen-Cahn model \eqref{eq:standardAC}, long term simulations of multi-front patterns on large domains are challenging, since the interaction between fronts is only exponentially weak so that the fronts typically get pinned (even in high precision simulations), {\it i.e.} the fronts stop traveling while the analysis shows that they must keep on moving (albeit exponentially slow). In the case of the inhomogeneous Allen-Cahn equation \eqref{eq:mainEquation} with localized topography $H(x)$ (see \eqref{eq:Ftopography}) this effect persists if it is assumed that $H(x)$ decays exponentially. However, the topographical driving forces are much stronger when $H(x)$ only decays algebraically slow (in fact, our theory applies if only $H'(x)$ decays, which implies that $H(x)$ may even grow as $x \to \pm \infty$, see \eqref{eq:algdecayHx} and Fig. \ref{fig:AlgHill}). Thus, additional to the `clean' result of Corollary \ref{cor:Nfrontslocalg}, there is a second reason to focus on algebraically decaying localized topographies: it is more tractable -- easier -- to  numerically explore the nature of multi-front dynamics if it is assumed that $H(x)$ decays algebraically.
\\ \\
In Fig. \ref{fig:5FrontDynamics} we present simulations of five-front dynamics in \eqref{eq:mainEquation} that run along lines similar to those of scenarios {\it (i)} - {\it (iii)} in Sec.~\ref{sec:ODE2frontloc} for general two-front dynamics.
In other words, like the two- and three-front dynamics of Fig. \ref{fig:NumericsTwoThreeFronts-Intro}(b),(e) and the four-front dynamics of Fig. \ref{fig:4FrontDynamics}, the five-front dynamics shown in Fig. \ref{fig:5FrontDynamics} corroborate the leading order natures of interaction ODEs \eqref{eq:dynNfronts} (for $N=2,3,4$ and $5$). This indicates that it will be possible to use \eqref{eq:dynNfronts} to give an approximate description of the expected long term dynamics of a given initial $N$-front configuration, as we did in Sec.~\ref{sec:ODE2frontloc} for $N=2$.
At least, this can be done as long as the distance between successive fronts does not become too small, or too large.
\\ \\
In the former case, the front interactions can no longer be considered weak (the assumption on which the derivation (and validation) of system \eqref{eq:dynNfronts} is based). In fact, it is clear from \eqref{eq:dynNfronts} that before distances between successive fronts become this small, the front evolution is already dominated by that of the homogeneous Allen-Cahn equation. This case is described and thoroughly analyzed in the classical works \cite{Carr1989,Chen2004,Fusco1989} (noting that the pairwise front-annihilation process that is especially relevant in the present setting was studied in \cite{Chen2004}). Thus, we may claim that the dynamics of multi-front patterns in inhomogeneous PDE \eqref{eq:mainEquation} for which the front interactions can no longer be considered weak can be analytically `controlled' by existing methods. However, this is a priori not clear for the other limit situation in which the distances between successive fronts become so large that it can no longer be assumed that the higher order correction terms $\varepsilon^2 R^{j,N}_{2}(\psi_1, \psi_2, ..., \psi_N;\varepsilon)$ to \eqref{eq:dynNfronts} -- that have not been determined -- can be neglected, {\it i.e.} that the asymptotic magnitudes of the $\varepsilon^2 R^{j,N}_{2}(\psi_1, \psi_2, ..., \psi_N;\varepsilon)$ terms indeed are small compared to the magnitudes of the terms $\varepsilon \mathcal{S}(\psi_j)$ ($j=1,...,N$). According to \eqref{eq:deltafrontsexp} -- that is based on the assumption that $\varepsilon \mathcal{S}(\psi_j)$ remains the leading order term in each $\frac{d \phi_j}{d \tau}$ equation -- the distances between successive fronts (that are sufficiently far apart) must always decrease. The simulations shown in Fig. \ref{fig:5Fronts-Intro}(c) seem to violate this prediction, which may suggest that $\varepsilon \mathcal{S}(\psi_j)$ indeed cannot determine the leading order dynamics of a front (that is sufficiently far removed from its neighbors). However, there is no solid foundation for such a conclusion, for instance because the simulations of Fig. \ref{fig:5Fronts-Intro}(c) are done for a topography $H(x)$ with algebraic behavior in its tails, while \eqref{eq:deltafrontsexp} has been derived for topographies that decay exponentially.
We refer to Sec.~\ref{sec:Discussion} for a further deliberation.
\section{The spatially periodically forced Allen-Cahn Equation}
\label{sec:ODEperiodic}
In this section, we consider the case of multi-front dynamics \eqref{eq:mainEquation} with a spatially periodic inhomogeneity, {\it i.e.}, we assume that there is an $X>0$ such that $F(U,V,x+X) \equiv F(U,V,x)$. First, we consider the dynamics of two-front patterns, again by focusing on the special case of $F(U,V,x)$
given by \eqref{eq:canonicalExample} with $f_{1,2,3}(x)$ as in \eqref{eq:fgh-example1}.
Next, we further develop the geometrical point of view; although this approach is general, we especially apply it to provide a deeper understanding of the appearance of stationary two-fronts in periodically driven systems. As in the localized setting of Sec.~\ref{sec:ODElocalized}, we conclude this section by deriving general insights in the dynamics of $N$-front patterns in periodically driven systems (focusing on those that differ essentially from those in Sec.~\ref{sec:ODElocalized}).

\subsection{Interacting two-fronts in the example system}
\label{sec:ODE2frontper}
\begin{figure}
	\centering
	\begin{subfigure}[t]{0.13\textwidth}
		\centering
		\includegraphics[width=\textwidth]{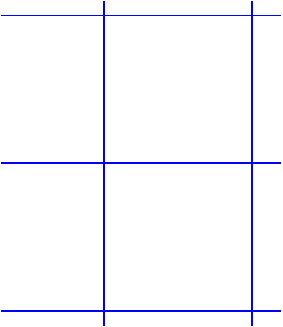}
		\caption{}
	\end{subfigure}\,
	\begin{subfigure}[t]{0.13\textwidth}
		\centering
		\includegraphics[width=\textwidth]{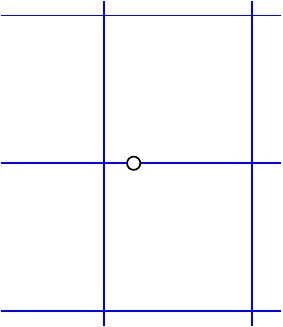}
		\caption{}
	\end{subfigure}\,
	\begin{subfigure}[t]{0.13\textwidth}
		\centering
		\includegraphics[width=\textwidth]{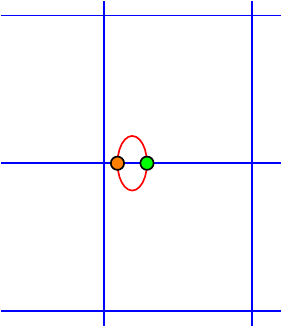}
		\caption{}
	\end{subfigure}\,
	\begin{subfigure}[t]{0.13\textwidth}
		\centering
		\includegraphics[width=\textwidth]{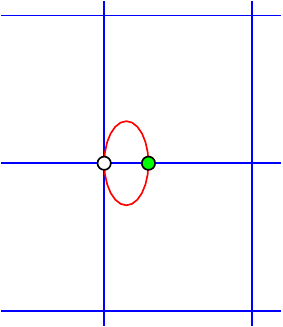}
		\caption{}
	\end{subfigure}\,
	\begin{subfigure}[t]{0.13\textwidth}
		\centering
		\includegraphics[width=\textwidth]{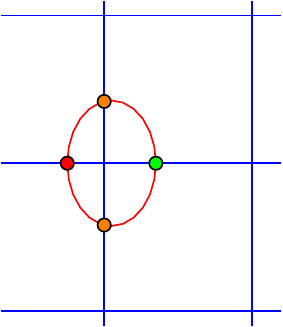}
		\caption{}
	\end{subfigure}\,
	\begin{subfigure}[t]{0.13\textwidth}
		\centering
		\includegraphics[width=\textwidth]{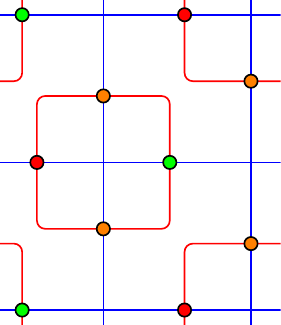}
		\caption{}
	\end{subfigure}\,
	\begin{subfigure}[t]{0.13\textwidth}
		\centering
		\includegraphics[width=\textwidth]{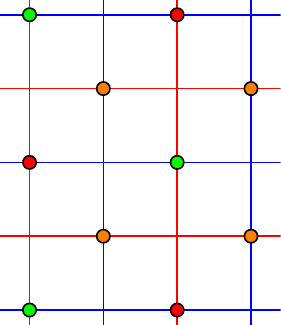}
		\caption{}
	\end{subfigure}
	\begin{subfigure}[t]{0.13\textwidth}
		\centering
		\includegraphics[width=\textwidth]{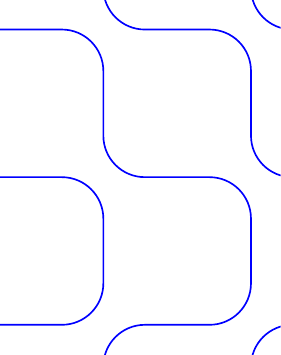}
		\caption{}
	\end{subfigure}\,
	\begin{subfigure}[t]{0.13\textwidth}
		\centering
		\includegraphics[width=\textwidth]{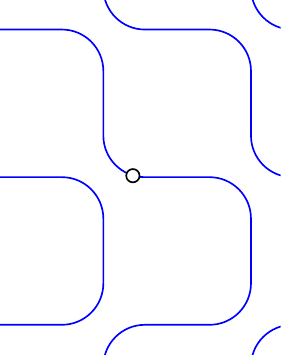}
		\caption{}
	\end{subfigure}\,
	\begin{subfigure}[t]{0.13\textwidth}
		\centering
		\includegraphics[width=\textwidth]{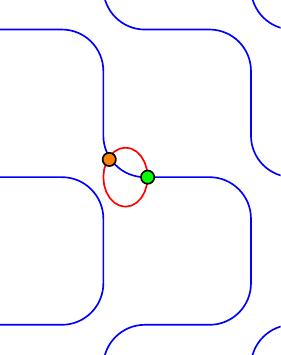}
		\caption{}
	\end{subfigure}\,
	\begin{subfigure}[t]{0.13\textwidth}
		\centering
		\includegraphics[width=\textwidth]{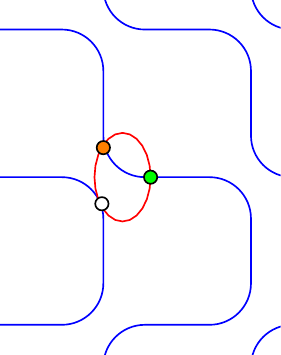}
		\caption{}
	\end{subfigure}\,
	\begin{subfigure}[t]{0.13\textwidth}
		\centering
		\includegraphics[width=\textwidth]{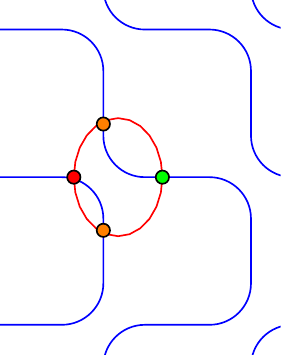}
		\caption{}
	\end{subfigure}\,
	\begin{subfigure}[t]{0.13\textwidth}
		\centering
		\includegraphics[width=\textwidth]{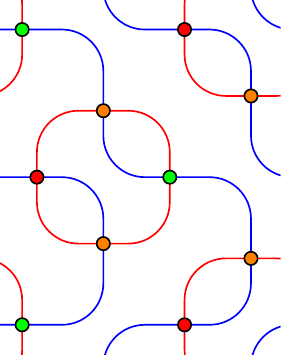}
		\caption{}
	\end{subfigure}\,
	\begin{subfigure}[t]{0.13\textwidth}
		\centering
		\includegraphics[width=\textwidth]{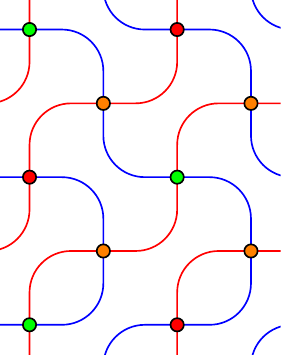}
		\caption{}
	\end{subfigure}
	\caption{Sketches of (possible) parts of the phase portrait of \eqref{eq:ODE2frontExample}/\eqref{eq:experfront-ODEsd} in the symmetric case (with the topographic case $B = 0$ as example) -- (a)-(g) -- and with broken symmetry ($B \neq 0$) -- (h)-(n). Blue lines denote the nullclines associated to the $S/s$-equation, red lines the nullclines associated to $D/d$. Circles indicate fixed points; green ones are stable nodes, orange ones saddles, red ones unstable nodes, white ones have one zero eigenvalue. (a)/(h) No fixed points ({\it i.e.} $D$/$d$ is sufficiently low). (b)/(i) A saddle-node bifurcation. (c)/(j) A saddle and stable node created by the saddle node bifurcation. (d)/(k) A pitchfork bifurcation in the symmetric case/a second saddle node bifurcation in the asymmetric case. (e)/(l) Both cases have two saddles, a stable node, and an unstable node (f)/(m) As $D/d$ increases several similar structures approach each other. (g)/(n) For $D/d \rightarrow \infty$, the structures merge and the nullclines form a regular grid.}
\label{fig:ODEbifurcations}
\end{figure}

\begin{figure}
	\centering
	\begin{subfigure}[t]{0.49\textwidth}
		\centering
		\includegraphics[width=\textwidth]{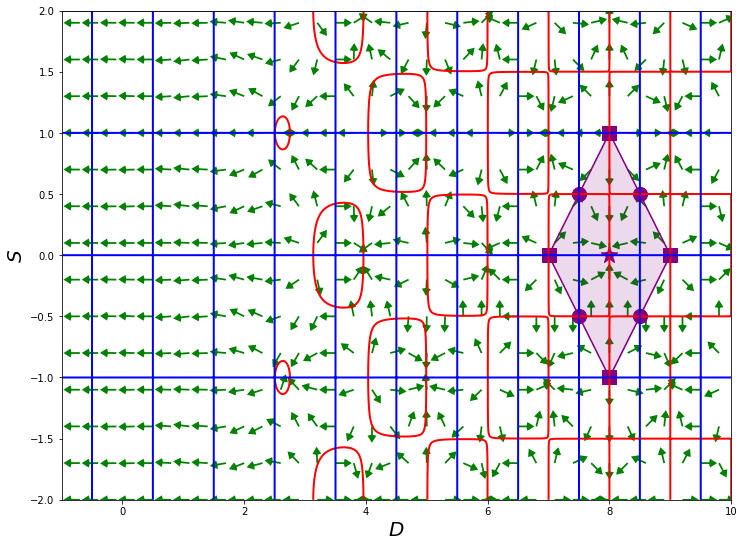}
		\caption{}
	\end{subfigure}\,
	\begin{subfigure}[t]{0.49\textwidth}
		\centering
		\includegraphics[width=\textwidth]{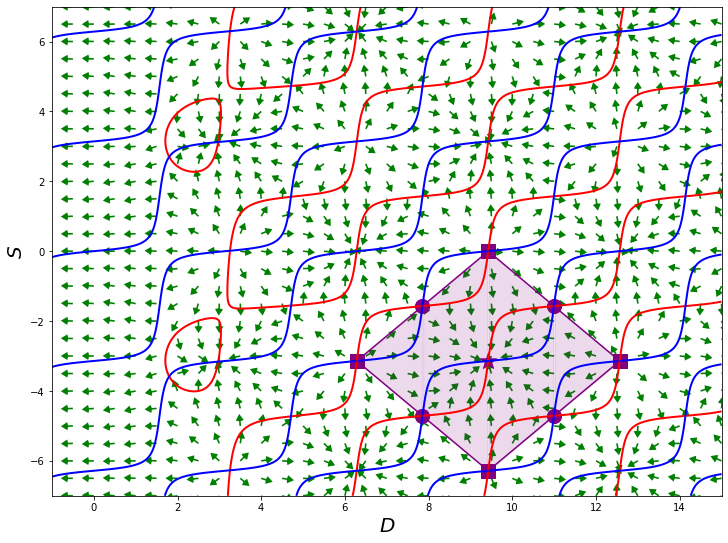}
		\caption{}
	\end{subfigure}\,
\caption{(a) Phase portrait of \eqref{eq:ODE2frontExample} in a symmetric case: $\varepsilon A = 0.029$, $\varepsilon B = 0$, and $k = 2 \pi$. (b) An asymmetric case: $\varepsilon A = 0.09$, $\varepsilon B = 0.01$ and $k = 2$. The red, respectively blue, lines denote the nullclines associated to the $D$, resp. $S$, equation. Green arrows denote the direction of the flow. The purple shaded area is the basin of attraction for the fixed point within that region (indicated with a purple star); the purple lines forming its boundary are the stable manifolds of the surrounding fixed points of saddle type (indicated with purple circles), which limit to the unstable nodes (indicated with purple squares). (Numerical computations are made using Python's pyplot's plotting routines quiver and contour and unstable manifolds are computed using backward time-integration with scipy's odeint.)}
\label{fig:PhasePortrait2FrontSymmASymm}
\end{figure}

\begin{figure}
\centering
	\begin{subfigure}[t]{0.24\textwidth}
		\centering
		\includegraphics[width=\linewidth]{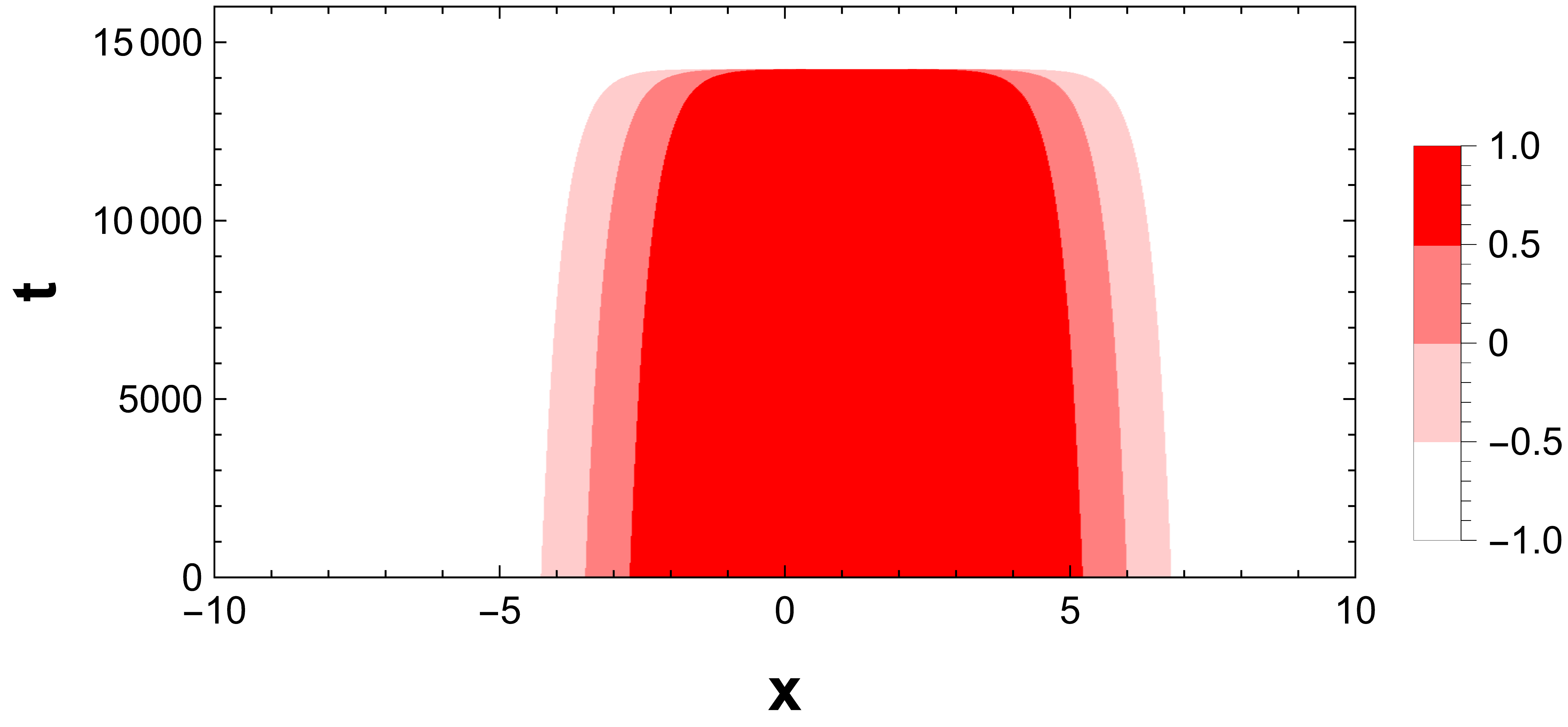}
	\end{subfigure}
	\begin{subfigure}[t]{0.24\textwidth}
		\centering
		\includegraphics[width=\linewidth]{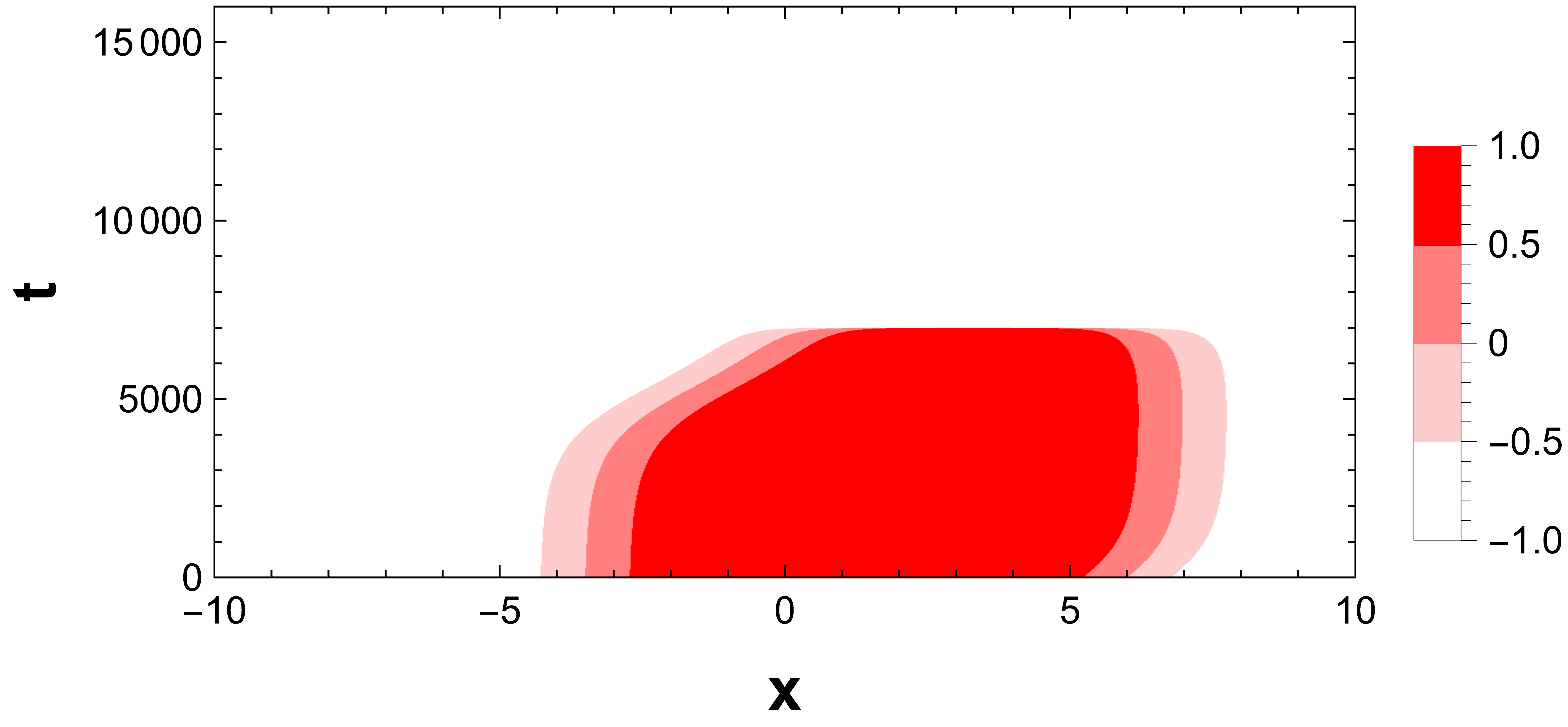}
	\end{subfigure}
	\begin{subfigure}[t]{0.24\textwidth}
		\centering
		\includegraphics[width=\linewidth]{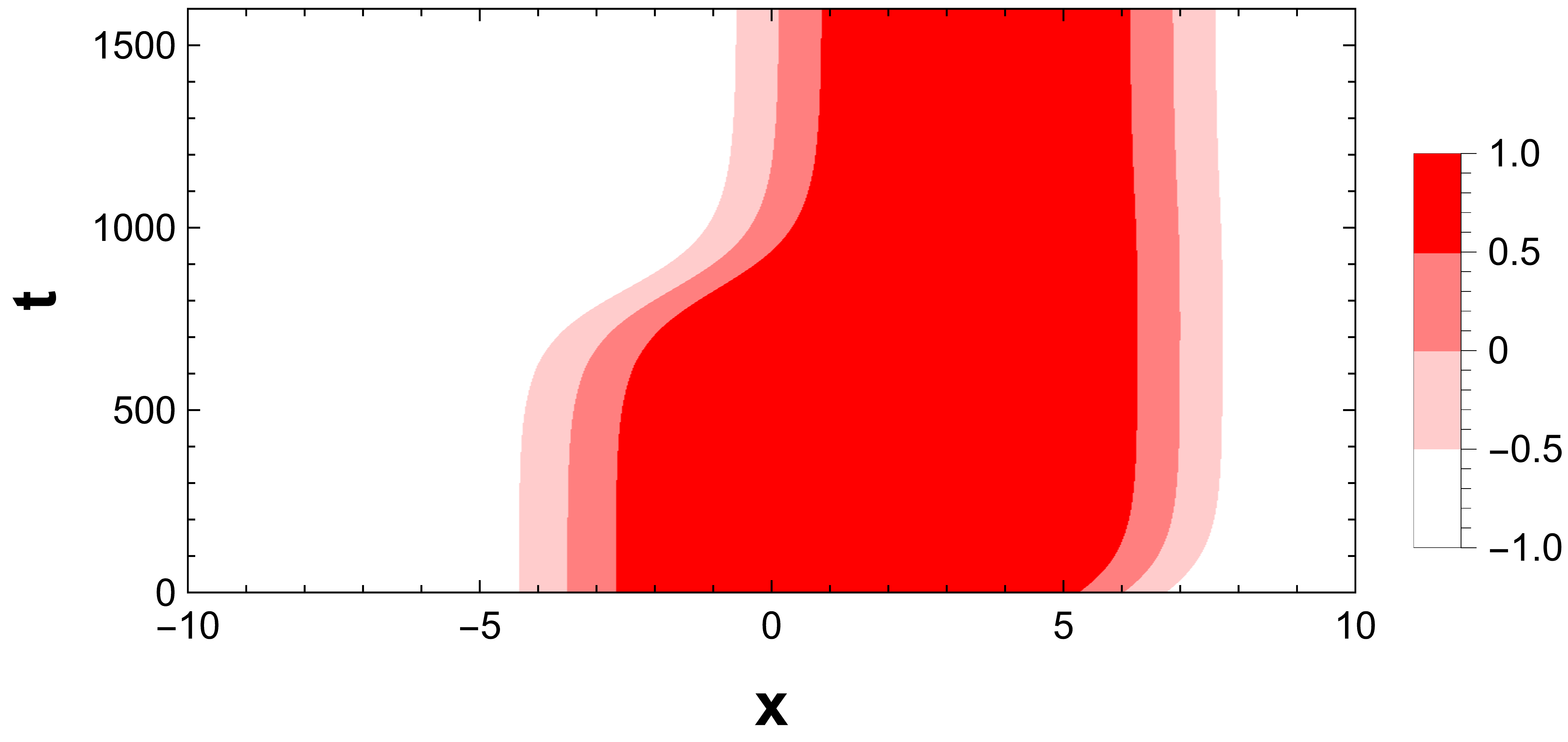}
	\end{subfigure}
	\begin{subfigure}[t]{0.24\textwidth}
		\centering
		\includegraphics[width=\linewidth]{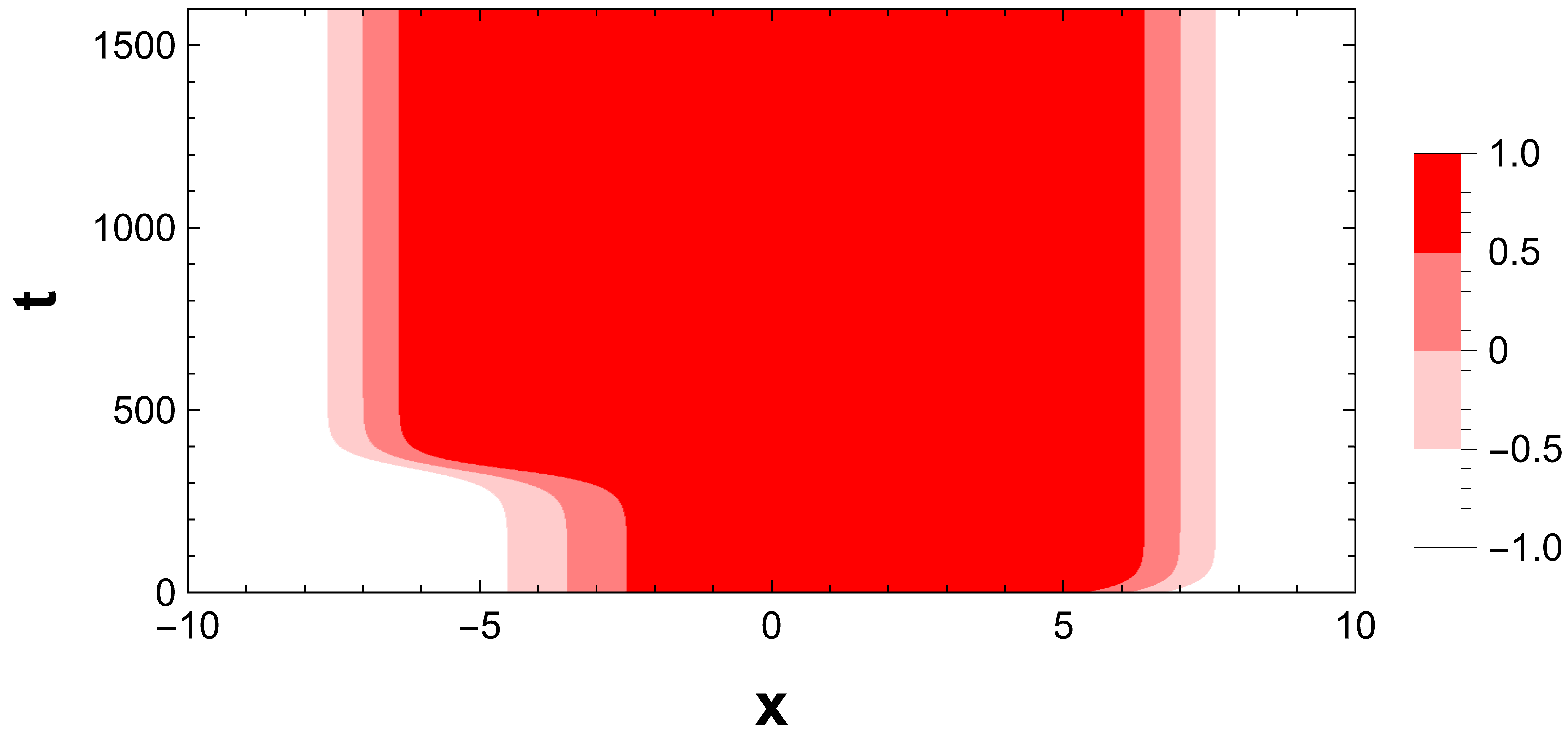}
	\end{subfigure}
	\caption{Simulations of \eqref{eq:mainEquation} on the $x$-interval $(-40,40)$ (with homogeneous Neumann boundary conditions) plotted on (-10,10) with spatially periodic topographic inhomogeneity given by \eqref{eq:Ftopography} with $H(x) = \sin kx$ with $k=0.8977$ starting from the two-front initial condition $U(x,0) =   -1 + \tanh(x + 3.5) -  \tanh(x - 6)$. (a) $\varepsilon = 0$, (b) $\varepsilon = 0.01$, (c) $\varepsilon = 0.1$ and (d) $\varepsilon = 0.4$.
 }
\label{fig:Numerics2Front}
\end{figure}

In this subsection, we consider the strongly simplified case of front-interactions (again) in the setting of explicit example \eqref{eq:canonicalExample} with (\ref{eq:fgh-example1}), and choose $N = 2$. We use this explicit form to be able to showcase the typical behavior in the most accessible and concrete way; the results are expected to be qualitatively similar for different periodically forced $F(U,V,x)$ (or $f_{1,2,3}(x)$) -- if the coefficients have two roots per period. (Also, if they have more, then the dynamics are slightly more complex, but follow in a similar vein.)
Thus, the Melnikov functions $\mathcal{R}_\mathrm{up/down}(\phi)$ are given by \eqref{Rupdownperex} with \eqref{perexAB}, where we recall that the choice $\alpha_1 = k \alpha_2$, $\alpha_3=0$ in (\ref{eq:fgh-example1}) represents the topographic subcase \eqref{eq:Ftopography} of \eqref{eq:canonicalExample}.
\\ \\
Substitution in \eqref{eq:NFrontODE}  gives the following ODE that describes the displacement of fronts over time. By rescaling time ($t' = 16t/\|u_\mathrm{up}'\|_2^2$, so that $t'= 8 \sqrt{2} \tau$ \eqref{eq:defpsiStau}), we find by \eqref{eq:NFrontODE} that the dynamics of the two-front pattern $(\phi_1(t'), \phi_1(t'))$ in \eqref{eq:mainEquation} -- with the first front going up -- is governed to leading order by
\begin{equation}
\left\{ \begin{array}{r c l l}
	\frac{d\phi_1}{dt'} & = & -\ \varepsilon (A+B) \sin(k\phi_1) & +\ e^{-\sqrt{2} (\phi_2-\phi_1)}; \\
	\frac{d\phi_2}{dt'} & = & -\ \varepsilon (A-B) \sin(k\phi_2) & -\ e^{-\sqrt{2} (\phi_2-\phi_1)},
\end{array}\right.
\label{eq:exper2front-ODE}
\end{equation}
with $A$ and $B$ determined by \eqref{perexAB}. Note that the critical points of \eqref{eq:exper2front-ODE} indeed coincide with the configurations determined by \eqref{eq:exper2front-exist}. This system simplifies (a bit) by introducing the new coordinates $(D(t'),S(t')) = (\phi_2(t') - \phi_1(t'),\phi_1(t') + \phi_2(t'))$,
\begin{equation}
\left\{ \begin{array}{r c l l l}
	\frac{dD}{dt'} & = & -\ 2 \varepsilon A \sin\left( \frac{1}{2} kD \right)\cos\left( \frac{1}{2} kS \right) & +\ 2 \varepsilon B \cos\left( \frac{1}{2} kD \right) \sin\left( \frac{1}{2} kS \right) & -\ 2 e^{-\sqrt{2} D}; \\
	\frac{dS}{dt'} & = & -\ 2 \varepsilon A \cos\left( \frac{1}{2} kD\right) \sin\left( \frac{1}{2} kS\right) & +\ 2 \varepsilon B \sin\left( \frac{1}{2} kD\right)\cos\left( \frac{1}{2} kS\right). &
\end{array}\right.
\label{eq:ODE2frontExample}
\end{equation}
(where we have used some standard trigonometric identities).
Observe that for (asymptotically large) $D$ such that $0<\sqrt{2}D/|\log \varepsilon| \to \rho$ as $\varepsilon \to 0$ with $\rho > 1$, the fixed points $(\bar{D},\bar{S})$ of \eqref{eq:ODE2frontExample} are determined to leading order by
\begin{equation}
\label{eq:exper2front-critptsfar}
\sin\left( \frac{1}{2} k \bar{D}\right) = \sin\left( \frac{1}{2} k \bar{S} \right) = 0 \; \; {\rm or} \; \; \cos\left( \frac{1}{2} k\bar{D} \right) = \cos\left( \frac{1}{2} k\bar{S} \right) = 0.
\end{equation}
This fully agrees with $(\phi_1,\phi_2) = (\phi_\mathrm{up},\phi_\mathrm{down}) = (n_\mathrm{up}(\varepsilon) \pi/k, n_\mathrm{up}(\varepsilon) \pi/k)$, as found in Sec.~\ref{sec:per2fronts}. To investigate the (dis)appearance of critical points as $\rho$ passes through $1$, we first consider the special case $B=0$, so that $S(t')$ is symmetric under $S \rightarrow -S$ in \eqref{eq:ODE2frontExample} -- which is the case in the topographic setting \eqref{eq:Ftopography}. Note that the upcoming analysis is valid for all $f_{1,2,3}(x)$ with $B=0$, {\it i.e.}, with $\alpha_3 = 0$ \eqref{perexAB}, not specifically for a topographic choice of $f_{1,2,3}(x)$. In fact, the underlying simplification comes from the symmetry $F(-U,-V;x) = -F(U,V;x)$ in \eqref{eq:mainEquation}, which is trivially satisfied by \eqref{eq:Ftopography}. A similar analysis can be performed for $F(-U,-V;x) = F(U,V;x)$ that implies $A=0$, so that \eqref{eq:ODE2frontExample} is symmetric under $S \rightarrow -S + \pi$. (In fact, there is a third special case, $A=B$, that we do not discuss.)
\\ \\
In the symmetric/topographic case with $B=0$, the $\frac{dS}{dt'} = 0$ nullclines of \eqref{eq:ODE2frontExample} are given by the straight lines for which $D = (2m_{S,D}+1)\pi/k$ or $S = 2m_{S,S}\pi/k$ ($m_{S,D/S} \in \mathbb{Z}$).
The $\frac{dD}{dt'} = 0$ nullclines satisfy
\[
\varepsilon A \sin\left( \frac{1}{2} kD \right) \cos\left( \frac{1}{2} kS \right) = - e^{-\sqrt{2} D}.
\]
Since the left hand side can only take on values between $+\varepsilon A$ and $-\varepsilon A$, necessarily only points with $D > D_c := (|\log \varepsilon| - \log |A|)/\sqrt{2}$ can lie on this nullcline, which confirms the critical $\rho =1$ growth rate of $D_c=D_c(\varepsilon)$. For $D(\varepsilon)$ with $\rho > 1$, these nullclines tend to the straight lines $D = 2 m_{D,D} \pi/k$, $S = (2m_{D,S}+1)\pi/k$ ($m_{D,D/S} \in \mathbb{Z}$) -- (again) reconfirming Sec.~\ref{sec:per2fronts} and \eqref{eq:exper2front-critptsfar}.
Decreasing $D$ so that $\rho$ becomes 1, while keeping $D> D_c$, causes these lines to break up into closed curves. These closed `circles' contract as $D$ decreases further and subsequently disappear, see Fig. \ref{fig:ODEbifurcations}(a)-(g) for a sketch.
\\
\\
To understand the precise mechanisms driving the appearance of the critical points of \eqref{eq:ODE2frontExample} and/or the associated bifurcations, we focus on the case $\rho = 1$. Therefore, we introduce $(d(\tau'),s(\tau'))$, $\tau'$, and constants $N=N_{\varepsilon,k} \in \mathbb{Z}$ and $R=R_{\varepsilon,k} \in [0,2 \pi)$ by
\begin{equation}
\label{eq:defsdtauNR}
d(\tau') = D(t') - \frac{|\log \varepsilon|}{\sqrt{2}}, \; \;  \ \ \ \
s(\tau') = S(t'), \; \; \ \ \ \tau' = 2 \varepsilon t',
\; \; \ \ \ \ \ 2 N \pi + R = \frac{k |\log \varepsilon|}{2 \sqrt{2}},
\end{equation}
where we note that $(d,s)$ correspond to $(\ell_\mathrm{down} + \ell_\mathrm{up}, \ell_\mathrm{down} - \ell_\mathrm{up})$ as introduced in
\eqref{eq:phiupdown} (and used throughout Sec.~\ref{sec:2fronts}). In the symmetric topographic ($B=0$) case, \eqref{eq:ODE2frontExample} becomes
\begin{equation}
\left\{ \begin{array}{r c l l }
	d' & = & - A \sin\left(R + \frac{1}{2} kd \right)\cos\left( \frac{1}{2} ks \right) & -\  e^{- d \sqrt{2}} \\[1mm]
	s' & = & - A \cos\left(R + \frac{1}{2} kd\right) \sin\left( \frac{1}{2} ks\right). &
\end{array}\right.
\label{eq:experfront-ODEsd}
\end{equation}
This ODE is periodic under $s \to s + 4 \pi/k$.
We may consider it on a cylinder $(d,s) \in \mathbb{R} \times S^1$, with $s \in [0,4\pi/k)$. We can distinguish between four families of fixed points, $(\bar{d}^{\mathrm{a}}_{\pm},\bar{s}^{\mathrm{a}}_{\pm})$ and $(\bar{d}^{\mathrm{b}}_{\pm},\bar{s}^{\mathrm{b}}_{\pm})$, distinguishable by setting either the cosine or the sine in the $s$-equation equal to zero:
\begin{equation}
\label{eq:exper2front-4famscritpts}
\begin{array}{rlcl}
{\rm (a)} & \bar{s}^{\mathrm{a}}_{+} = 0, \bar{s}^{\mathrm{a}}_{-} = \frac{2\pi}{k}, & & A \sin\left(R + \frac{1}{2} k\bar{d}^{\mathrm{a}}_{\pm} \right) = \mp e^{-\bar{d}^{\mathrm{a}}_{\pm}\sqrt{2}},
\\
{\rm (b)} & \bar{d}^{\mathrm{b}}_{\pm} = \frac{(4 m_\mathrm{b} \pm 1)\pi - 2R}{k}, m_\mathrm{b} \in \mathbb{Z} & & A \cos\left(\frac{1}{2} k \bar{s}^{\mathrm{b}}_{\pm} \right) = \mp e^{-\frac{(4 m_\mathrm{b} \pm 1)\pi - 2R}{k}\sqrt{2}}.
\end{array}
\end{equation}
(to leading order in $\varepsilon$), where we immediately note that each family $(\bar{d}^{\mathrm{a/b}}_{\pm},\bar{s}^{\mathrm{a/b}}_{\pm})$ contains countably many elements. The associated eigenvalues are given by
\begin{equation}
\label{eq:exper2front-evscritpts}
\begin{array}{rl}
{\rm (a)} & \lambda^{\mathrm{a}}_{\pm, 1} = \mp \frac12 k A \cos\left(R + \frac{1}{2} k\bar{d}^{\mathrm{a}}_{\pm}\right) +  \sqrt{2} e^{-\sqrt{2} \bar{d}^{\mathrm{a}}_{\pm}}, \; \; \lambda^{\mathrm{a}}_{\pm, 2} = \mp \frac12 k A \cos\left(R + \frac{1}{2} k\bar{d}^{\mathrm{a}}_{\pm} \right),
\\
{\rm (b)} &
\lambda^{\mathrm{b}}_{\pm, 1/2} = \frac{1}{2} \left(\sqrt{2} e^{-\sqrt{2} \bar{d}^{\mathrm{b}}_{\pm}} \pm \sqrt{ 2 e^{-2 \sqrt{2} \bar{d}^{\mathrm{b}}_{\pm}} + k^2 A^2 \sin^2\left(\frac{1}{2} k\bar{s}^{\mathrm{b}}_{\pm}\right)} \right)
\end{array}
\end{equation}
(which for $\bar{d}^{\mathrm{a/b}}_{\pm} \gg 1$ agrees with \eqref{eq:exper2front-lambda12} (with $B=0$) -- after undoing all transformations and scalings). Hence, it immediately follows that all (b)-type critical points $(\bar{d}^{\mathrm{b}}_{\pm},\bar{s}^{\mathrm{b}}_{\pm})$ are saddles (for all parameters values considered).
\\
\\
To study the nature of the (a)-type critical points $(\bar{d}^{\mathrm{a}}_{\pm},\bar{s}^{\mathrm{a}}_{\pm})$ ,we consider $A$ as bifurcation parameter. It follows from \eqref{eq:exper2front-4famscritpts}(a) that as $A$ increases new critical points are created at the bifurcational values $A=A^{\mathrm{SN}}_{\pm}$ at which the curve  $A \sin\left(R + \frac{1}{2} k\bar{d} \right)$ becomes tangent to either one of the curves $\mp e^{-\bar{d}\sqrt{2}}$ (as function of $\bar{d}$). Note that these new critical points have $\bar{d}$-coordinates smaller than the $\bar{d}$-coordinates of all (countably many) critical points that already exist before $A$ passes through $A=A^{\mathrm{SN}}_{\pm}$, {\it i.e.}, the value of the minimal $\bar{d}$-coordinate over all existing the critical points decreases as $A$ increases.
For $\bar{d}^{\mathrm{a}}_{+}$, respectively $\bar{d}^{\mathrm{a}}_{-}$, the saddle node bifurcation happens with a positive, resp. negative, derivative, which implies that $\cos\left(R + \frac{1}{2} k\bar{d}^{\mathrm{a}}_{+} \right) > 0$, resp. $\cos\left(R + \frac{1}{2} k\bar{d}^{\mathrm{a}}_{-} \right) < 0$. Hence, by \eqref{eq:exper2front-evscritpts}(a), both $\lambda^{\mathrm{a}}_{\pm, 1} < 0$ at the tangency (and by construction $\lambda^{\mathrm{a}}_{\pm, 2} = 0$). It follows that as $A$ increases through each of the values $A=A^{\mathrm{SN}}_{\pm}$, two new critical points are created: $(\bar{d}^{\mathrm{a},1}_{\pm},\bar{s}^{\mathrm{a},1}_{\pm})$ always appears as a saddle and $(\bar{d}^{\mathrm{a},2}_{\pm},\bar{s}^{\mathrm{a},2}_{\pm})$ as a stable node (with $\bar{d}^{\mathrm{a},1}_{\pm} < \bar{d}^{\mathrm{a},2}_{\pm}$).
\\
\\
As $A$ increases further, the stable node $(\bar{d}^{\mathrm{a},2}_{\pm},\bar{s}^{\mathrm{a},2}_{\pm})$ does not undergo any further bifurcations (and $\bar{d}^{\mathrm{a},2}_{\pm}$ increases towards the next zero of $\sin\left(R + \frac{1}{2} k\bar{d}\right)$ as $A \to \infty$ \eqref{eq:exper2front-4famscritpts}). However, $(\bar{d}^{\mathrm{a},1}_{\pm},\bar{s}^{\mathrm{a},1}_{\pm})$ does bifurcate: as $A$ increases beyond $A^{\mathrm{SN}}_{\pm}$, $\bar{d}^{\mathrm{a},1}_{\pm}$ decreases and reaches at  $A=A^{\mathrm{PF}}_{\pm} (> A^{\mathrm{SN}}_{\pm})$ the point where the curve $\mp e^{-\bar{d}\sqrt{2}}$ intersects $A \sin\left(R + \frac{1}{2} k\bar{d} \right)$ in its nearest minimum ($+$) or maximum ($-$).
As a consequence, $\lambda^{\mathrm{a}}_{\pm, 2}$ becomes zero at $A^{\mathrm{PF}}_{\pm}$ (while $\lambda^{\mathrm{a}}_{\pm, 1} = \sqrt{2}\,{\rm exp}(-\sqrt{2} \bar{d}^{\mathrm{a}}_{\pm}) > 0$) \eqref{eq:exper2front-evscritpts}. Moreover, at this point, $(\bar{d}^{\mathrm{a},1}_{\pm},\bar{s}^{\mathrm{a},1}_{\pm})$ coincides with two appearing (symmetrical) saddle points of (b)-type, $(\bar{d}^{\mathrm{b},1/2}_{\pm},\bar{s}^{\mathrm{b},1/2}_{\pm})$
\eqref{eq:exper2front-4famscritpts} that split off from $(\bar{d}^{\mathrm{a},1}_{\pm},\bar{s}^{\mathrm{a},1}_{\pm})$ as $A$ increases further beyond $A^{\mathrm{PF}}_{\pm}$.
The saddle  $(\bar{d}^{\mathrm{a},1}_{\pm},\bar{s}^{\mathrm{a},1}_{\pm})$ changes into an unstable node by this pitchfork bifurcation, its nature no longer changes as $A$ increases further (in fact, $\bar{d}^{\mathrm{a},1}_{\pm}$ decreases towards the preceding zero of $\sin\left(R + \frac{1}{2} k\bar{d}\right)$ as $A \to \infty$ \eqref{eq:exper2front-4famscritpts}). We refer to Fig. 
\ref{fig:ODEbifurcations}(a)-(g) for graphical visualizations of this process and to Fig.\ref{fig:PhasePortrait2FrontSymmASymm}(a) for a phase portrait of the flow governed by \eqref{eq:ODE2frontExample} (or equivalently \eqref{eq:experfront-ODEsd}) for $A$ just below an $A^{\mathrm{PF}}_{\pm}$.
\\
\\
Now, if system \eqref{eq:ODE2frontExample} does not possess a symmetry, {\it i.e.}, if $B \neq 0$ (or $A\neq 0$), the above scenario breaks down. However, the situation is that of a standard broken symmetry scenario: the saddle node bifurcations persist and the pitchfork bifurcations break open into saddle-node bifurcations with a persisting critical point. A precise distinction between (a)- and (b)-type critical points can no longer be made, but the perturbed, asymmetric scenario can be described based on these notions. For increasing $A$, first a saddle and a stable node corresponding to  $(\bar{d}^{\mathrm{a},1}_{\pm},\bar{s}^{\mathrm{a},1}_{\pm})$ $(\bar{d}^{\mathrm{a},2}_{\pm},\bar{s}^{\mathrm{a},2}_{\pm})$ appear (with $\bar{d}^{\mathrm{a},1}_{\pm} < \bar{d}^{\mathrm{a},2}_{\pm}$). Neither of these undergo any additional bifurcations. However, as $A$ increases further, two unstable critical points appear in a second saddle node bifurcation, a saddle and an unstable node. This saddle node bifurcation takes over the role of $A^{\mathrm{PF}}_{\pm}$.
The unstable node corresponds to $(\bar{d}^{\mathrm{a},1}_{\pm},\bar{s}^{\mathrm{a},1}_{\pm})$ in the symmetrical case: to maintain the relation with the symmetric case, the persisting saddle that appeared at the first saddle node bifurcation has to be relabeled as a (b)-type saddle point as $A$ passes through the second saddle node bifurcation -- see Fig.~\ref{fig:ODEbifurcations}(h)-(n) for a schematic sketch and Fig.~\ref{fig:PhasePortrait2FrontSymmASymm}(b) for a phase portrait.
\\
\\
Our insights so far on the existence, nature and bifurcations of the critical points are a priori only the first steps towards understanding the dynamics of (two) interacting fronts governed by equations \eqref{eq:exper2front-ODE}/\eqref{eq:ODE2frontExample}/\eqref{eq:experfront-ODEsd}. However, it is not hard to see from both the symmetrical and asymmetrical phase portraits of Fig. \ref{fig:PhasePortrait2FrontSymmASymm} that for fronts that are sufficiently far apart, {\it i.e.}, for $D$ (or $d$) sufficiently large, each stable ((a)-type) node $(\bar{d}^{\mathrm{a}}_{\pm}, \bar{a}^{\mathrm{b}}_{\pm})$ is enclosed by a network formed by the unstable manifolds of 4 ((b)-type) saddle points $(\bar{d}^{\mathrm{b}}_{\pm}, \bar{s}^{\mathrm{b}}_{\pm})$ -- with 4 ((a)-type) unstable nodes $(\bar{d}^{\mathrm{a}}_{\pm}, \bar{s}^{\mathrm{a}}_{\pm})$ at its  vertices ({\it i.e.}, the boundary of the basin of attraction of each stable ((a)-type) node $(\bar{d}^{\mathrm{a}}_{\pm}, \bar{a}^{\mathrm{b}}_{\pm})$ is formed by the stable manifolds of the surrounding saddle-type fixed points; see Fig~\ref{fig:PhasePortrait2FrontSymmASymm}) Hence, for any initial condition that is not chosen exactly on this network, the two-front pattern will converge to the nearby (enclosed) attractor. In other words, two-front patterns that are sufficiently far apart will always converge to a nearby stationary two-front pattern. Here, `sufficiently far apart' is as usual measured with respect to the standard distance scale $|\log \varepsilon|$ -- see the simulations of Fig. \ref{fig:Numerics2Front} in which the initial position of a two-front pattern is kept fixed, while the magnitude of $\varepsilon$ increases from $0$, through $0.01$, to $0.1$ and $0.4$, so that the distance between the fronts increases relative to the decreasing magnitude $|\log \varepsilon|$. For $\varepsilon = 0.001$ (Fig. \ref{fig:Numerics2Front}(b)), the system behaves as the homogeneous Allen-Cahn case (with $\varepsilon = 0$, Fig. \ref{fig:Numerics2Front}(a)): the fronts are attracted to each other, merge and disappear (nevertheless, the influence of the spatial inhomogeneity can clearly not be completely neglected for $\varepsilon = 0.01$). As $\varepsilon$ increases further, the two-front pattern indeed settles in a nearby stationary pattern (Figs. \ref{fig:Numerics2Front}(c),(d)). However, the value of $\varepsilon$ also influences the position of the initial condition with respect to the above network of unstable manifolds and thus has impact on the stationary stable two-front state the PDE dynamics is attracted to -- as can be seen by comparing Fig. \ref{fig:Numerics2Front}(c) with $\varepsilon=0.1$ to Fig. \ref{fig:Numerics2Front}(d) with $\varepsilon=0.4$.
\\
\\
In Sec.~\ref{sec:ODENfrontper}, a general result like this -- {\it i.e.}, on the existence of nearby stationary attractors for $N$-front patterns -- will be proven (for the case of a topographic spatially periodic driving term \eqref{eq:Ftopography}, under the condition that the fronts are sufficiently separated -- see Theorem \ref{th:Nfrontsperattr}). The proof is completely based on the observation that if (in the setting of special case \eqref{eq:exper2front-ODE})  $\sqrt{2} (\phi_2-\phi_1)/|\log \varepsilon| > 1$ (as $\varepsilon \to 0$), the front equations decouple (to leading order) and both $\phi_{1/2}$ evolve towards a zero $\phi_{*,1/2}$ of $\sin(k \phi)$ with $\cos(k \phi_{*,1/2}) < 0$ (so that $(\phi_{*,1},\phi_{*,2})$ represents the stable attracting stationary two-front pattern to leading order). On the other hand, for fronts that are relatively close to each other, {\it i.e.}, with $\sqrt{2} (\phi_2-\phi_1)/|\log \varepsilon| < 1$, the terms associated to the driving heterogeneity become of higher order relevance and the two-front interaction equation is to leading order given by that for the classical homogeneous Allen-Cahn equation case \cite{Carr1989,Chen2004, Fusco1989}, which explains why the $\varepsilon = 0.001$ dynamics of Fig. \ref{fig:Numerics2Front}(b) are similar to that of Fig. \ref{fig:Numerics2Front}(a) (in which  $\varepsilon = 0$).

\subsection{The existence problem as a weakly driven nonlinear oscillator: the geometric approach}
\label{sec:geom-int}
\begin{figure}[ht!]
	\centering
		\begin{subfigure}[t]{0.325 \textwidth}
			\centering
			\includegraphics[width=\textwidth]{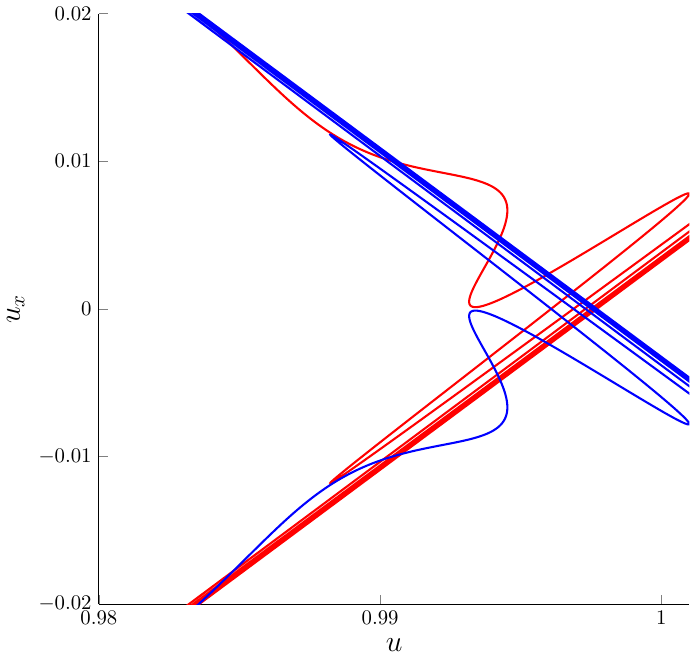}
			\caption{$\alpha_1 = -0.95$}
		\end{subfigure}
		\begin{subfigure}[t]{0.325 \textwidth}
			\centering
			\includegraphics[width=\textwidth]{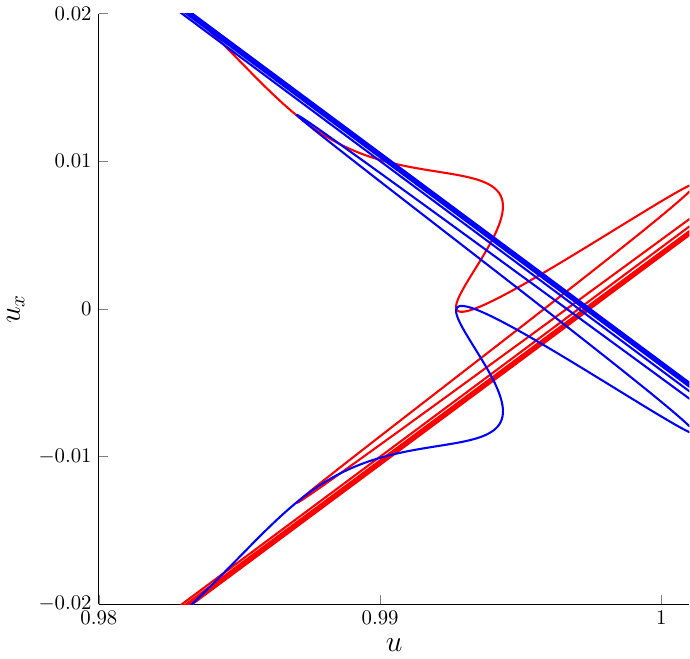}
			\caption{$\alpha_1 = -1.05$}
		\end{subfigure}
		\begin{subfigure}[t]{0.325 \textwidth}
			\centering
			\includegraphics[width=\textwidth]{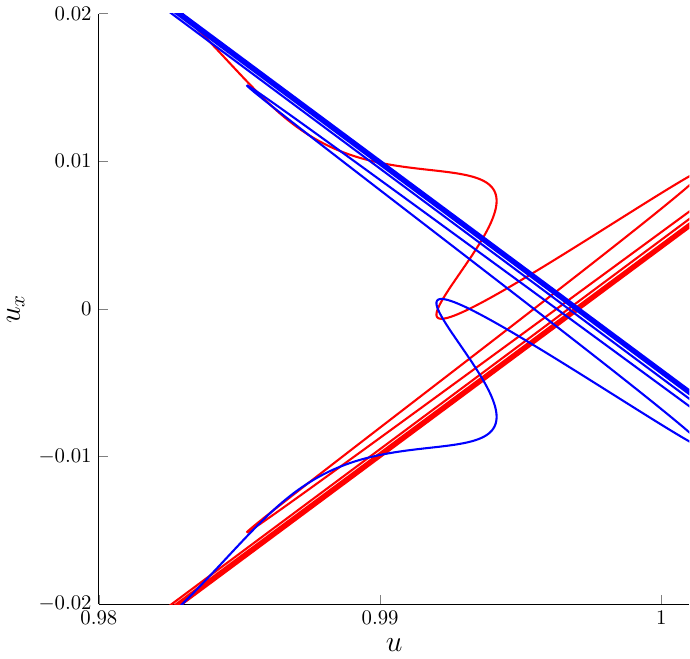}
			\caption{$\alpha_1 = -1.20$}
		\end{subfigure}
	\caption{Approximations of $W_0^{s/u}(\mathcal{M}_-^\varepsilon)$ with $\varepsilon = 0.1$ and $F(U,V,x) = \alpha_1 \cos(2 \pi x ) U$ for various values for $\alpha_1$. The red lines denote the unstable manifold $W_0^u(\mathcal{M}_-^\varepsilon)$, and the blue lines the stable manifolds $W_0^s(\mathcal{M}_-^\varepsilon)$.
    We focus on the small region in the middle where a small segment of the red manifold dips right to left down toward a small (symmetrically-diposed) segment of the blue manifold. (a) these local segments do not intersect, (b) these local segments have two intersections that are created in a saddle-node bifurcation, (c) these local segments have four intersections after a pitchfork bifurcation.}
	\label{fig:manifolds2FrontsSym2}
\end{figure}

\begin{figure}
	\centering
		\begin{subfigure}[t]{0.325 \textwidth}
			\centering
			\includegraphics[width=\textwidth]{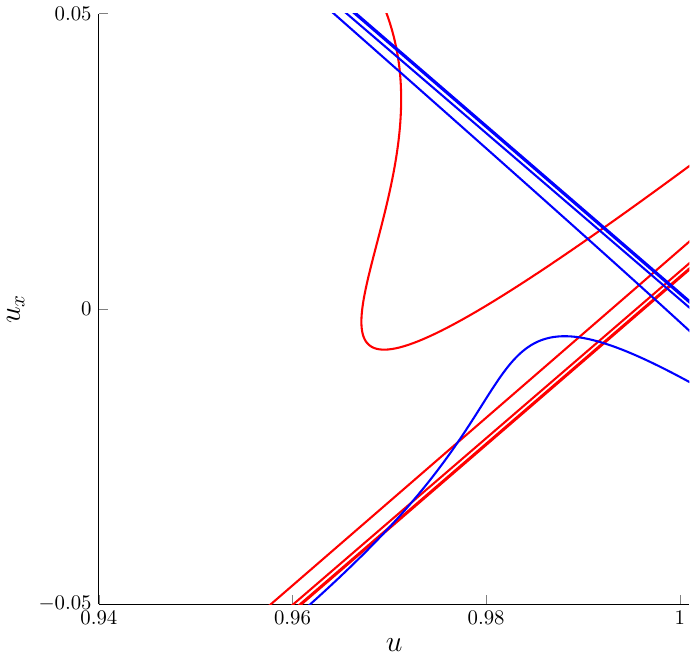}
			\caption{$\alpha_1 = -0.1$}
		\end{subfigure}
		\begin{subfigure}[t]{0.325 \textwidth}
			\centering
			\includegraphics[width=\textwidth]{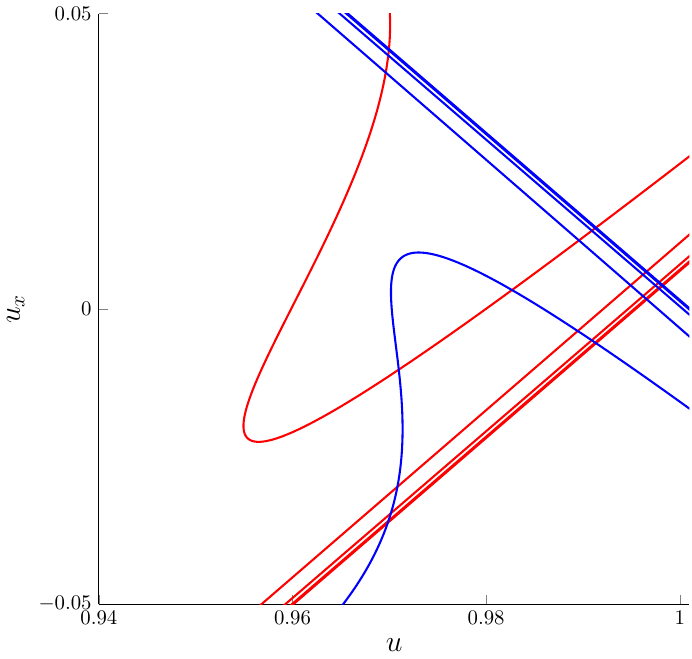}
			\caption{$\alpha_1 = -0.2$}
		\end{subfigure}
		\begin{subfigure}[t]{0.325 \textwidth}
			\centering
			\includegraphics[width=\textwidth]{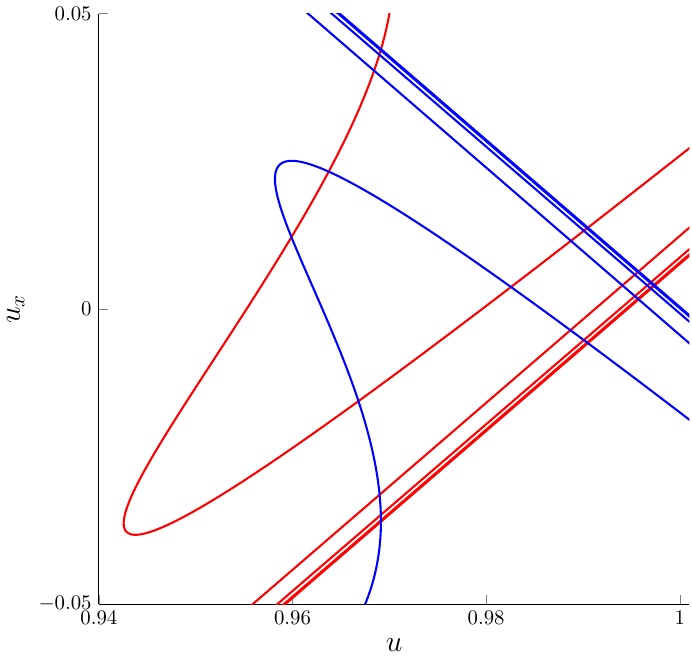}
			\caption{$\alpha_1 = -0.3$}
		\end{subfigure}
	\caption{Approximations of $W_0^{s/u}(\mathcal{M}_-^\varepsilon)$ with $\varepsilon = 0.1$ and $F(U,V,x) = \alpha_1 \cos(\pi x ) U + 0.15 \sin(\pi x)$ for various values for $\alpha_1$. The red lines denote the unstable manifold $W_0^u(\mathcal{M}_-^\varepsilon)$ and the blue lines the stable manifolds $W_0^s(\mathcal{M}_-^\varepsilon)$. (a) No intersection between manifolds, (b) two fixed points appeared by a first saddle-node bifurcation, (c) four fixed points after a second saddle-node bifurcation.}
	\label{fig:manifolds2FrontsNonSym}
\end{figure}

In this section, we again restrict ourselves to the issue of the existence of stationary multi-front patterns, {\it i.e.}, we do not consider the dynamics of multi-front patterns in \eqref{eq:mainEquation} as we did at the end of the previous section (and will again in  Sec.~\ref{sec:ODENfrontper} below). Here, we take a geometric point of view by going deeper into the correspondence between the existence question for localized stationary multi-front patterns \eqref{eq:stationaryODEsystem} and the (classical) theory of weakly driven (planar) nonlinear oscillators (see for example \cite{Guckenheimer2002}). This point of view provides an alternative way to study the appearance of multi-front patterns -- {\it i.e.}, homoclinic and heteroclinic orbits of \eqref{eq:stationaryODEsystem} -- as functions of changing parameters. In fact, especially in the case that $F(u,v,x)$ is periodic in $x$, it provides a natural way to embed this analysis in that of the classical subject of `lobe dynamics' of perturbed nonlinear oscillators.
\\ \\
Recapitulating the observations of Sec.~\ref{sec:1fronts} (especially Sec.~\ref{sec:1fronts-class}), we note that for $\varepsilon = 0$, \eqref{eq:stationaryODEsystem} possesses three fixed points, the two saddles $\mathcal{M}_\pm^0 := (u_\pm,0) = (\pm1,0)$, with stable/unstable manifolds $W^{s/u}(\mathcal{M}_\pm^0)$, and the center $\mathcal{M}_0^- := (u_0,0) = (0,0)$ (Fig. \ref{fig:1FrontSolutions}). There are two (unperturbed) heteroclinic orbits $W^u(\mathcal{M}_-^0) \cap W^s(\mathcal{M}_+^0)$ and $W^u(\mathcal{M}_+^0) \cap W^s(\mathcal{M}_-^0)$ between the saddles $\mathcal{M}_\pm^0$ that are explicitly given by the one-front solutions $u_\mathrm{up}(x;\phi)$ and $u_\mathrm{down}(x;\phi)$ \eqref{eq:heteroclinicSolution}.
\\ \\
For $0<\varepsilon \ll 1$, \eqref{eq:stationaryODEsystem} is non-autonomous, so that it is natural to consider the system in the extended $(u,p,x)$-phase space. Here, the saddles $\mathcal{M}_\pm^0$ persist (for $0 < \varepsilon$ sufficiently small) as unique, bounded solutions $u_{\pm}^\varepsilon(x)$ (given explicitly in \eqref{eq:uMinusPlus}) which both have $2D$ stable/unstable manifolds $W^{s/u}(\mathcal{M}_\pm^\varepsilon) \subset \mathbb{R}^3$. These manifolds can be explicitly approximated by constructing the solutions $(u^\varepsilon_{\pm,s/u},\frac{d}{dx} u^\varepsilon_{\pm,u/s},x) \in W^{s/u}(\mathcal{M}_\pm^\varepsilon)$. Focusing on $W^{u}(\mathcal{M}_-^\varepsilon)$ we set
\begin{equation}
\label{eq:regexpuux}
	u^\varepsilon_{-,u}(x;\phi) = u_\mathrm{up}(x;\phi) + \varepsilon u_{u,1}(x;\phi) + \varepsilon^2 u_{u,2}(x;\phi) + \; {\rm h.o.t.}
\end{equation}
and find, once again (see ~\eqref{eq:existenceProblem1Front}),
\[
\mathcal{L}_0 u_{u,1} = - F(u_\mathrm{up}(x;\phi),u_{\mathrm{up},x}(x;\phi),x).
\]
The form of solutions $u_{u,1}(x;\phi)$ that stay bounded as $x \rightarrow -\infty$ (i.e. solutions $u_u$ that tend to $u_-^\varepsilon$ as $x \rightarrow -\infty$) is given in \eqref{eq:generalSolutionOrderEpsSystem}, \eqref{eq:defABphi}, \eqref{eq:appendixBoundedUnboundedSolution} (see also Appendix~\ref{sec:orderEpsSystem}) and we have thus obtained the leading order approximation,
\begin{equation}
\widetilde{W}^u(\mathcal{M}_-^\varepsilon) := \left\{ \left( u_\mathrm{up}(x;\phi) + \varepsilon u_{u,1}(x;\phi), \frac{d}{dx} u_{\mathrm{up}}(x;\phi) + \varepsilon  \frac{d}{dx} u_{u,1}(x;\phi), x \right): x \in (-\infty, x_\mathrm{max}) \right\},
\label{eq:tildeWuMepsmin}
\end{equation}
where $x_\mathrm{max}$ is the maximal value for which the approximation holds. More precisely, $x_\mathrm{max}$ is not uniquely determined, choosing $x_\mathrm{max}$ also determines the accuracy of the approximation. It is shown in Appendix~\ref{sec:orderEpsSystem} that
\begin{equation}
\label{eq:limxinfty-uu1}
\lim_{x-\phi \rightarrow \infty} u_{u,1}(x;\phi) e^{-\sqrt{2}[x-\phi]}=\lim_{x-\phi \rightarrow \infty} \frac{1}{8} \int_{-\infty}^{\infty} F(u_\mathrm{up}(z;\phi), u_{\mathrm{up}, x}(z;\phi),z)\ \Psi_b(z;\phi)\ dz e^{\pm \sqrt{2}[x-\phi]} = \frac{1}{8} \mathcal{R}(\phi)
\end{equation}
\eqref{eq:FredholmCondition1Front}. Thus, the choice $x_\mathrm{max} - \phi = |\log \varepsilon|/(2\sqrt{2})$ would yield an $\mathcal{O}(\sqrt{\varepsilon})$ accurate approximation, while \eqref{eq:limxinfty-uu1} does not provide useful information if  $x_\mathrm{max} - \phi = \rho |\log \varepsilon|/\sqrt{2}$ with $\rho > 1$ (unless $\mathcal{R}(\phi) = 0$, for which we have recovered the one-front connection between the limiting solutions $u_{\pm}^\varepsilon(x)$ of Sec.~\ref{sec:1fronts}). Following the procedure of \cite{doelman2022slow}, it can be shown that considering the next order terms in the regular expansion \eqref{eq:regexpuux} does not lead to a potential way out: $u_{u,2}(x; \phi)$ grows as exp$(2\sqrt{2}[x-\phi])$ as $x-\phi \to \infty$ so that $\varepsilon^2 u_{u,2}(x; \phi) = \mathcal{O}(\varepsilon^{2-2\rho})$ for $x_\mathrm{max} - \phi = \rho |\log \varepsilon|/\sqrt{2}$. In fact, it follows (in general) that all higher order terms yield $\mathcal{O}(1)$ contributions as $x_\mathrm{max} - \phi = |\log \varepsilon|/\sqrt{2}$. However, by restricting expansion \eqref{eq:regexpuux} to  $x_\mathrm{max} - \phi = |\log \varepsilon|/(2\sqrt{2})$ ($+\mathcal{O}(1)$), an approximation can be obtained for the first `jump' of $W^{u}(\mathcal{M}_-^\varepsilon)$ from $u_{-}^\varepsilon(x)$ into the $\mathcal{O}(\sqrt{\varepsilon})$ region through which it passes $u_{+}^\varepsilon(x)$ (and more accurate approximations within this region can be obtained from expansion \eqref{eq:regexpuux}). Beyond this region, further approximations of $W^u(\mathcal{M}_-^\varepsilon)$ can be obtained by a similar approach based on the `jump back orbit' $u_\mathrm{down}(x; \phi)$, etc., see again the procedure of \cite{doelman2022slow}.
\\
\begin{figure}[ht!]
	\centering
		\begin{subfigure}[t]{0.325 \textwidth}
			\centering
			\includegraphics[width=\textwidth]{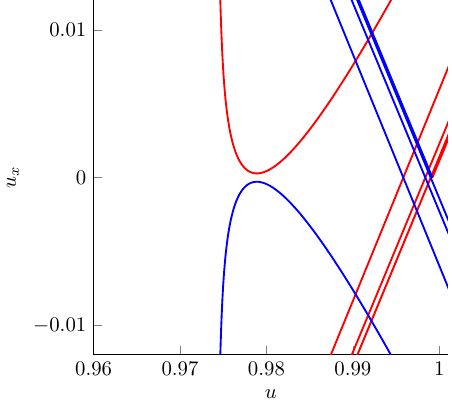}
			\caption{$\alpha_1 = -0.094$}
		\end{subfigure}
		\begin{subfigure}[t]{0.325 \textwidth}
			\centering
			\includegraphics[width=\textwidth]{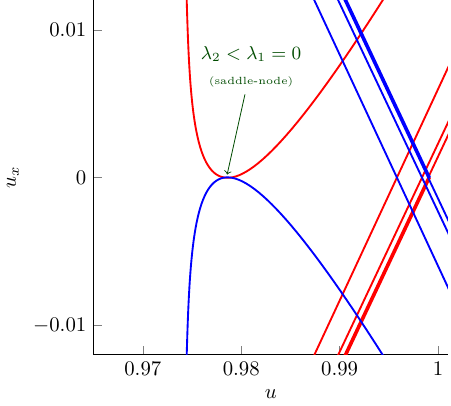}
			\caption{$\alpha_1 = -0.096$}
		\end{subfigure}
		\begin{subfigure}[t]{0.325 \textwidth}
			\centering
			\includegraphics[width=\textwidth]{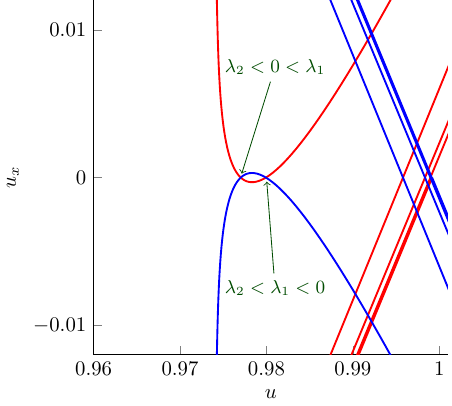}
			\caption{$\alpha_1 = -0.098$}
		\end{subfigure}
		\begin{subfigure}[t]{0.325 \textwidth}
			\centering
			\includegraphics[width=\textwidth]{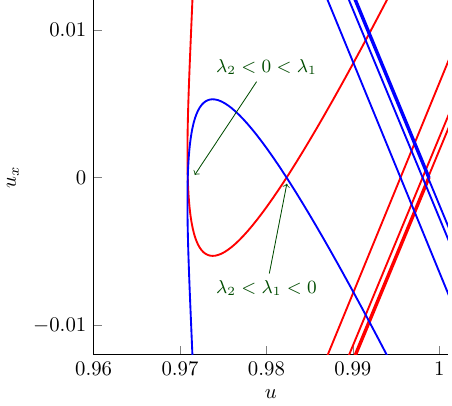}
			\caption{$\alpha_1=-0.131$}
		\end{subfigure}
		\begin{subfigure}[t]{0.325 \textwidth}
			\centering
			\includegraphics[width=\textwidth]{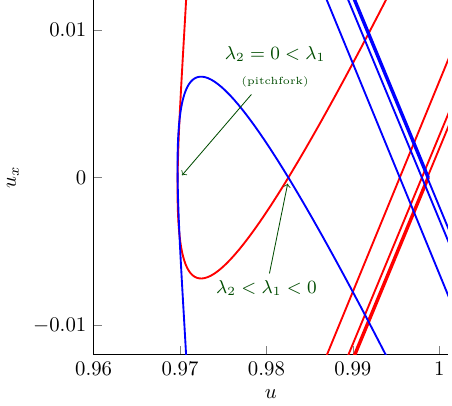}
			\caption{$\alpha_1 = -0.141$}
		\end{subfigure}
		\begin{subfigure}[t]{0.325 \textwidth}
			\centering
			\includegraphics[width=\textwidth]{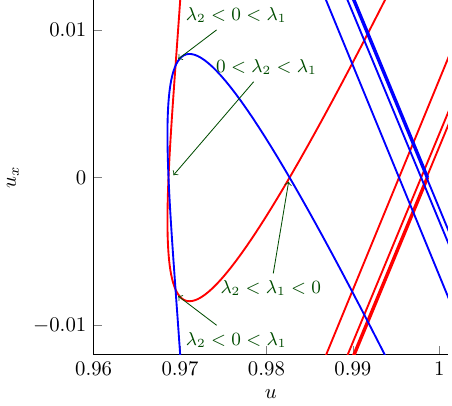}
			\caption{$\alpha_1 = -0.151$}
		\end{subfigure}
		\begin{subfigure}[t]{0.5 \textwidth}
			\centering
			\includegraphics[width=\textwidth]{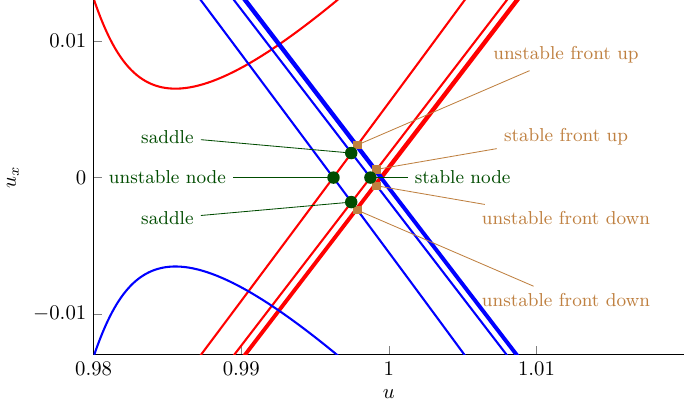}
			\caption{$\alpha_1 = -0.05$}
		\end{subfigure}
	\caption{Analytically obtained approximations of $W_0^{s/u}(\mathcal{M}_-^\varepsilon)$ evaluated at $\varepsilon = 0.1$ with $F(U,V,x) = \alpha_1 \cos(\pi x) U$, for various values for $\alpha_1$. The red lines denote the unstable manifold $W_0^u(\mathcal{M}_-^\varepsilon)$ and the blue lines the stable manifolds $W_0^s(\mathcal{M}_-^\varepsilon)$. (a) No intersection between manifolds, (b) saddle-node bifurcation, (c) two fixed points arise, (d) two fixed points, (e) pitchfork bifurcation, (f) four fixed points, (g) classification of two-front solutions in relation to $1$-front solutions. Results on stability are entirely based on the analysis in Sec.~\ref{sec:ODE2frontper}.}
	\label{fig:manifolds2FrontsSym}
\end{figure}
\\
Here, we do not go further into the details of approximating manifolds $W^{s/u}(\mathcal{M}_\pm^\varepsilon)$ beyond their first `jump' from $\mathcal{M}_\pm^\varepsilon$ to $\mathcal{M}_\mp^\varepsilon$.
(Though we refer the reader to \cite{Kedar1995} for a useful approach for doing this and for finding second-order, third-order, and higher-order intersections of these invariant manifolds.)
Instead, we apply the above procedure to approximate $W^{s/u}(\mathcal{M}_-^\varepsilon)$ up to the $\mathcal{O}(\sqrt{\varepsilon})$ neighborhood of $\mathcal{M}_+^\varepsilon$, which allows us to investigate the first intersections $W^{u}(\mathcal{M}_-^\varepsilon) \cap W^{s}(\mathcal{M}_-^\varepsilon)$, {\it i.e.}, the existence and (dis)appearance of homoclinic orbits of two-front type to $u_-^\varepsilon(x) = - 1 + \mathcal{O}(\varepsilon)$ \eqref{eq:uMinusPlus}. Studying the impact of parameter variations on $W^{u}(\mathcal{M}_\pm^\varepsilon) \cap W^{s}(\mathcal{M}_\pm^\varepsilon)$ thus provides an alternative way to study the bifurcations of (stationary) two-front patterns, compared to that of Sec.~\ref{sec:ODE2frontper}. As usual for these kinds of perturbed integrable systems -- especially for periodically driven systems -- we will consider their restrictions to
specific values of $x$, or equivalently, their intersections with a plane $\{x=x_{\ast}\}$, which we will denote by $W^{s/u}_{x_{\ast}}(\mathcal{M}_\pm^\varepsilon)$ and approximate by $\widetilde{W}^{s/u}_{x_{\ast}}(\mathcal{M}_\pm^\varepsilon)$.
That is,
\[
\widetilde{W}_{x_{\ast}}^u(\mathcal{M}_-^\varepsilon) := \widetilde{W}^u(\mathcal{M}_-^\varepsilon) \cap \left\{ x = x_{\ast} \right\} = \left\{ \left( u_\mathrm{up}(x_{\ast};\phi) + \varepsilon u_{u,1}(x_{\ast};\phi), \frac{d}{dx} u_{\mathrm{up}}(x_{\ast};\phi) + \varepsilon \frac{d}{dx} u_{u,1}(x_{\ast};\phi), x \right) \right\},
\]
and similarly for the approximations of the other stable and unstable manifolds. For functions $F(U,V,x)$ that are periodic in $x$, {\it i.e.}, $F(U,V,x+X) \equiv F(U,V,x)$ for some $X > 0$, we consider the usual Poincar\'e section at, e.g., $\{x=0\}$ and use $W^{u/s}_{0}(\mathcal{M}_\pm^\varepsilon)$ to study the geometry of the system \eqref{eq:stationaryODEsystem}. In Fig. \ref{fig:manifolds2FrontsSym}, we present the results for a function $F$ that possesses the symmetry $F(U,V,x) = - F(-U,-V,x)$ of the topographic case (although the example does not represent a topography $H(x)$). The case of functions $F(U,V,x)$ that possess the symmetry $F(U,V,x) = F(-U,-V,x)$ is similar. Our findings are in agreement with those of Sec.~\ref{sec:ODE2frontper}.
For a changing bifurcation parameter, first a saddle-node bifurcation occurs at which $W_{0}^u(\mathcal{M}_-^\varepsilon)$ and $W_{0}^s(\mathcal{M}_-^\varepsilon)$ are tangent.
This tangency indicates the beginning of the process that leads to a full intersection of the $W_{0}^u(\mathcal{M}_-^\varepsilon)$- and $W_{0}^s(\mathcal{M}_-^\varepsilon)$-lobes. Directly after the saddle-node bifurcation, the `nose' of the $W_{0}^u(\mathcal{M}_-^\varepsilon)$-lobe is inside that of $W_{0}^s(\mathcal{M}_-^\varepsilon)$ (and vice versa), it again punctures through $W_{0}^u(\mathcal{M}_-^\varepsilon)$ by the pitchfork bifurcation (due to the symmetry) -- see again Fig. \ref{fig:manifolds2FrontsSym} and also Fig. \ref{fig:manifolds2FrontsSym2}, that confirms that this crossing through-scenario appears simultaneously at several -- in fact countably many -- locations in the $(u,p)$-plane.
\\
\\
This geometric point of view is especially helpful in understanding the distinction between systems \eqref{eq:stationaryODEsystem} with or without a symmetry. In Fig. \ref{fig:manifolds2FrontsNonSym}, we present the generic results for a function $F(U,V,x)$ that does not possess one of the aforementioned symmetries. In line with the results in Sec.~\ref{sec:ODE2frontper}, it shows that from a generic point of view, having a second saddle node bifurcation is more natural than the pitchfork bifurcation of the symmetric case.
\\
\\
We conclude this (sub)section by the obvious observation that the present geometrical approach based on the intersections $W_{0}^u(\mathcal{M}_\pm^\varepsilon) \cap W_{0}^s(\mathcal{M}_\pm^\varepsilon)$ does not yield any information on the stability of the associated multi-front patterns as solutions of PDE \eqref{eq:mainEquation}. It is natural to expect that there may be a relation between the geometrical nature of a (planar) intersection $W_{0}^u(\mathcal{M}_\pm^\varepsilon) \cap W_{0}^s(\mathcal{M}_\pm^\varepsilon)$ and the stability character of the associated $N$-front pattern in the infinite-dimensional setting. This is a subject of future research.

\subsection{$N$-front dynamics: the existence of nearby attractors}
\label{sec:ODENfrontper}
In this section, we consider general aspects of the dynamics of multi-front patterns in periodically driven Allen-Cahn equations \eqref{eq:mainEquation}. For simplicity we again focus on the topographic case, {\it i.e.} unlike the general case of the (directly) preceding subsections, we consider a driving term $F(U,V,x)$ given by \eqref{eq:Ftopography} with topography $H(x)$ that is periodic with period $X>0$. Thus, we know that $\mathcal{R}_\mathrm{up}(\phi) = \mathcal{R}_\mathrm{down}(\phi)$ ($=\mathcal{R}(\phi)$) \eqref{eq:Rtopoupdown} and that $\mathcal{R}(\phi)$ is also periodic with period $X$ \eqref{eq:Rtopoper}. On the other hand, we consider general periodic topographies, so that $\mathcal{R}(\phi)$ may have more than two zeroes during each full period, unlike the examples considered so far.
\\ \\
As before, we first focus on the existence of stationary $N$-front patterns that will play a central role as attractors of the flow (see Theorem \ref{th:Nfrontsperattr} below).
The situation is different from that of localized topographies considered in Theorem \ref{th:Nfrontslocexp}.
There are countably many stationary $N$-fronts from which (roughly) one out of $2^N$ is stable.
\begin{theorem}
\label{th:Nfrontsper}
Consider the weakly heterogeneous Allen-Cahn equation \eqref{eq:mainEquation} with inhomogeneity $F(U,U_x,x)$ given by the topography $H(x)$ \eqref{eq:Ftopography}. Assume that $H(x)$ is periodic with (minimal) period $X>0$, so that Melnikov function $\mathcal{S}(\psi)$ associated to $H(x)$ -- that is periodic with period $Y = \sqrt{2}X$ (\eqref{eq:Rtopoper}, \eqref{eq:defpsiStau}) -- has an even number $K\geq 0$ of non-degenerate zeroes $0 \leq \psi_{\ast,1} < \psi_{\ast,2} < ... < \psi_{\ast,K} < Y$ on the interval $[0,Y)$ ({\it i.e.}, $\mathcal{S}(\psi_{\ast,k}) = 0$ with $\mathcal{S}'(\psi_{\ast,k}) \neq 0$ and necessarily ${\rm sign} \, \mathcal{S}'(\psi_{\ast,k}) \neq {\rm sign} \, \mathcal{S}'(\psi_{\ast,k+1})$). Denote for $N \geq 2$ a stationary $N$-front pattern -- that connects $u_-^\varepsilon(x) = - 1 + \mathcal{O}(\varepsilon)$ as $x \to - \infty$ to $u_{(-1)^{N+1}}^\varepsilon(x) = (-1)^{N+1} + \mathcal{O}(\varepsilon)$ as $x \to +\infty$ \eqref{eq:uMinusPlus} -- by $(\bar{\psi}_1, \bar{\psi}_2, ..., \bar{\psi}_N)$, where $\bar{\psi}_j$ is the location of the $j$-th front (with $\bar{\psi}_j < \bar{\psi}_{j+1}$, $j=1,...,N$) and let $\Delta \bar{\psi}_j = \bar{\psi}_{j+1} - \bar{\psi}_j$ be the distance between two successive fronts. Let $0 < \varepsilon \ll 1$ sufficiently small.
\\
$\bullet$ If $(\bar{\psi}_1, \bar{\psi}_2, ..., \bar{\psi}_N)$ is a stationary $N$-front pattern with $\Delta \bar{\psi}_j/|\log \varepsilon| = \rho_j$ as $\varepsilon \to 0$, then $\rho_j \geq 1$ for all $j=1,...,N$.
\\
$\bullet$ Consider a sequence $(\tilde{\psi}_{1}, \tilde{\psi}_{2}, ..., \tilde{\psi}_{N}) = (\psi_{\ast,i_1} + n_1(\varepsilon) Y, \psi_{\ast,i_2} + n_2(\varepsilon) Y, ..., \psi_{\ast,i_N} + n_N(\varepsilon) Y)$ of $N$ consecutive zeroes of $\mathcal{S}(\psi)$, {\it i.e.}, $\tilde{\psi}_{j} < \tilde{\psi}_{j+1}$ for $j=1,2,...,N$ (with $n_{j}(\varepsilon) \in \mathbb{Z}$, $i_j = 1,2, ...,K$). If $\Delta \tilde{\psi}_{j} = \tilde{\psi}_{j+1} - \tilde{\psi}_{j}$ increases faster than $|\log \varepsilon|$ as $\varepsilon \to 0$, {\it i.e.}, if $(n_{j+1}(\varepsilon) - n_{j}(\varepsilon))/| \log \varepsilon| > Y$ as $\varepsilon \to 0$ ($j=1,...,N-1$), then there exists a stationary $N$-front pattern $(\bar{\psi}_{1}, \bar{\psi}_{2}, ..., \bar{\psi}_{N})$ asymptotically close to
$(\tilde{\psi}_{1}, \tilde{\psi}_{2}, ..., \tilde{\psi}_{N})$, {\it i.e.}, $\bar{\psi}_{j} = \tilde{\psi}_{j} + r_j(\varepsilon)$ with $r_j(\varepsilon) \to 0$ as $\varepsilon \to 0$ ($j=1,...,N$).
\\
$\bullet$ Associated to any ordered sequence $(n_1(\varepsilon), n_2(\varepsilon), ..., n_N(\varepsilon)) \in \mathbb{Z}^N$ ({\it i.e.}, $n_j(\varepsilon) < n_{j+1}(\varepsilon)$) with $(n_{j+1}(\varepsilon) - n_{j}(\varepsilon))/| \log \varepsilon| > Y$ as $\varepsilon \to 0$ ($j=1,...,N-1$) there are $K^N$ distinct stationary $N$-front solutions denoted by $(\bar{\psi}_{1,i_1}, ..., \bar{\psi}_{N,i_N}) = (\psi_{\ast,i_1} + n_1(\varepsilon) Y + r_{1,i_1}(\varepsilon), ..., \psi_{\ast,i_N} + n_N(\varepsilon) Y + r_{N,i_N}(\varepsilon))$, with $r_{j,i_j}(\varepsilon) \to 0$ as $\varepsilon \to 0$ ($i_j = 1,2, ...,K$, $j=1, ...,N$). The stability of $(\bar{\psi}_{1,i_1}, ..., \bar{\psi}_{N,i_N})$ is determined by $N$ (asymptotically small) eigenvalues $\lambda_{j}(\varepsilon; i_1, ..., i_N)$ that are given to leading order by
\[
\lambda_{j}(\varepsilon; i_1, ..., i_N) = -\varepsilon \mathcal{S}'(\psi_{\ast,i_j}), \; \; i_j = 1,2, ...,K, \; \; j=1, ..., N.
\]
Thus, $(\frac12 K)^N$ of all $K^N$ stationary $N$-front patterns associated to the sequence $(n_1(\varepsilon), n_2(\varepsilon), ..., n_N(\varepsilon))$ are stable, all $(1-\frac{1}{2^N}) K^N$ others are unstable.
\end{theorem}
An essentially identical result holds for stationary $N$-front patterns that connect $u_+^\varepsilon(x) = 1 + \mathcal{O}(\varepsilon)$ at $-\infty$ to $u_{(-1)^{N}}^\varepsilon(x) = (-1)^{N} + \mathcal{O}(\varepsilon)$ at $+\infty$.
\\ \\
The proof of this theorem is completely straightforward. As in the proof of Theorem \ref{th:Nfrontslocexp}, we base our approach on the $N$-front interaction ODE \eqref{eq:dynNfronts}. This is possible by the arguments given at the beginning of this proof ({\it i.e.}, in first paragraph of Appendix \ref{ap:ProofTh}) -- also in the present situation in which we did not provide the analytical details of the proof of the validity of interaction ODE \eqref{eq:dynNfronts} (see the discussion in the second paragraph of Sec.~\ref{sec:interactiondynamics}).
\\ \\
The statement on the non-existence of stationary $N$-front patterns is based on the simple observation that the $\varepsilon \mathcal{S}(\psi_j)$ terms in \eqref{eq:dynNfronts} are not of leading order in the case in which the distances between all successive fronts grow slower than $1 \times |\log \varepsilon|$ as $\varepsilon \to 0$. Thus, the inhomogeneous effects in \eqref{eq:mainEquation} can be neglected which implies that there cannot be stationary $N$-fronts (as in the homogeneous Allen-Cahn equation). On the other hand, if the distances between all successive fronts grow faster than $1 \times |\log \varepsilon|$, then the $\varepsilon \mathcal{S}(\psi_j)$ terms in \eqref{eq:dynNfronts} represent the leading order effects. In other words, the leading order dynamics of each front is only determined by the inhomogeneity: the $j$-th front at $\psi_j(\tau)$ evolves as a solitary front, {\it i.e.} it is governed by $\frac{d \psi_j}{d \tau} = - \varepsilon \mathcal{S}(\psi_j)$ (see \eqref{eq:1frontdyn}, \eqref{eq:defpsiStau}). Both the statement on the existence of stationary $N$-fronts as well as that on their stability in Theorem \ref{th:Nfrontsper} are based on this relatively simple observation.
\\ \\
Transitions between no and countably many stationary $N$-front patterns occur as an $\rho_j$ -- defined in Theorem \ref{th:Nfrontsper} by $\Delta \bar{\psi}_j/|\log \varepsilon| = \rho_j$ as $\varepsilon \to 0$ -- passes through $1$. (Note that it can a priori not be excluded that some critical points may exist for $\rho_j = 1$ in the case that $\mathcal{S}(\psi)$ does not have zeroes ({\it i.e.}, for $K=0$ in Theorem \ref{th:Nfrontsper}); however, $\mathcal{S}(\psi)$ cannot change sign, thus there can be no stationary $N$-front configurations $(\bar{\psi}_1, \bar{\psi}_2, ..., \bar{\psi}_N)$ with both $\mathcal{S}(\bar{\psi}_1) = 16 e^{-(\bar{\psi}_2 - \bar{\psi}_1)}/\varepsilon > 0$ and $\mathcal{S}(\bar{\psi}_N) = -16 e^{-(\bar{\psi}_N - \bar{\psi}_{N-1})}/\varepsilon < 0$.) In Secs.~\ref{sec:ODE2frontper} and \ref{sec:geom-int}, we studied the appearance of critical points for the specific $N=2$ case with a topography $H(x)$ determined by $H'(x) = f_2(x) = \alpha_2 \sin (kx)$. In the general $N$-front setting -- with an $\mathcal{S}(\psi)$ that may have more than two zeroes -- unraveling the transitions becomes technically more involved, but no conceptual novelties are introduced, especially not from the geometrical lobe intersection point of view of Sec.~\ref{sec:geom-int}.
Therefore, we refrain from going further into these details.
\\
\begin{figure}
\centering
	\begin{subfigure}[t]{0.24\textwidth}
		\centering
		\includegraphics[width=\linewidth]{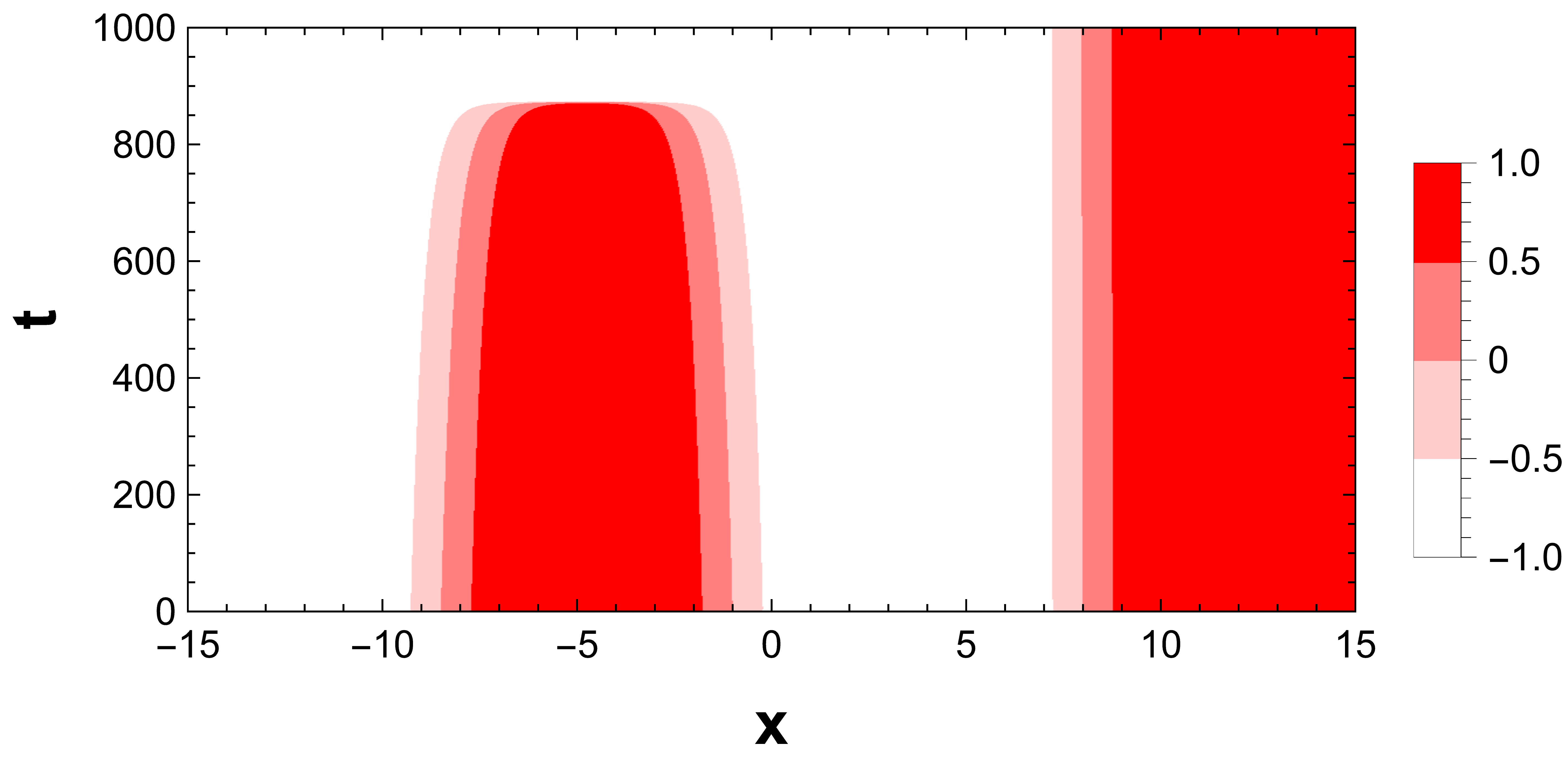}
	\end{subfigure}
	\begin{subfigure}[t]{0.24\textwidth}
		\centering
		\includegraphics[width=\linewidth]{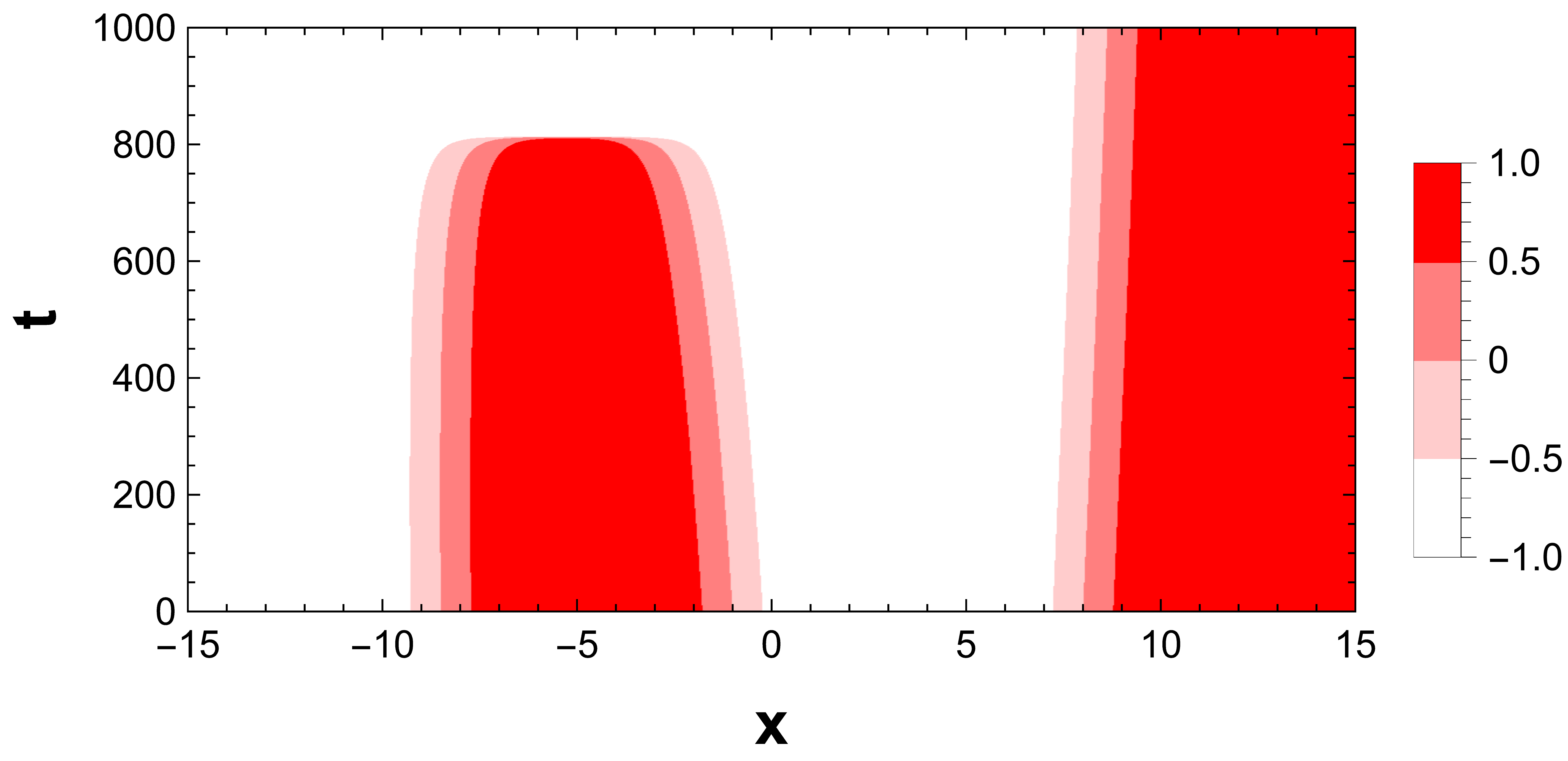}
	\end{subfigure}
	\begin{subfigure}[t]{0.24\textwidth}
		\centering
		\includegraphics[width=\linewidth]{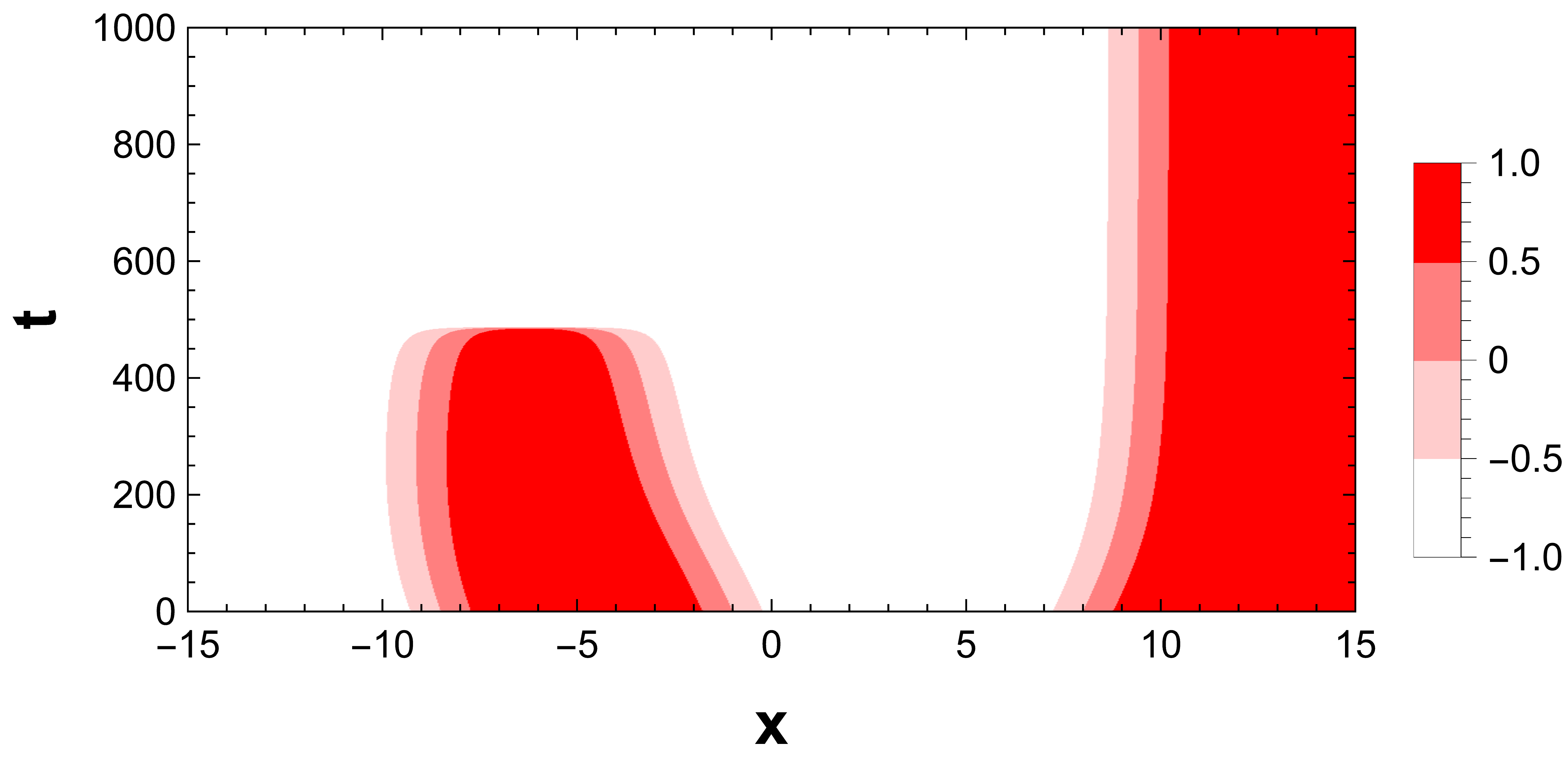}
	\end{subfigure}
	\begin{subfigure}[t]{0.24\textwidth}
		\centering
		\includegraphics[width=\linewidth]{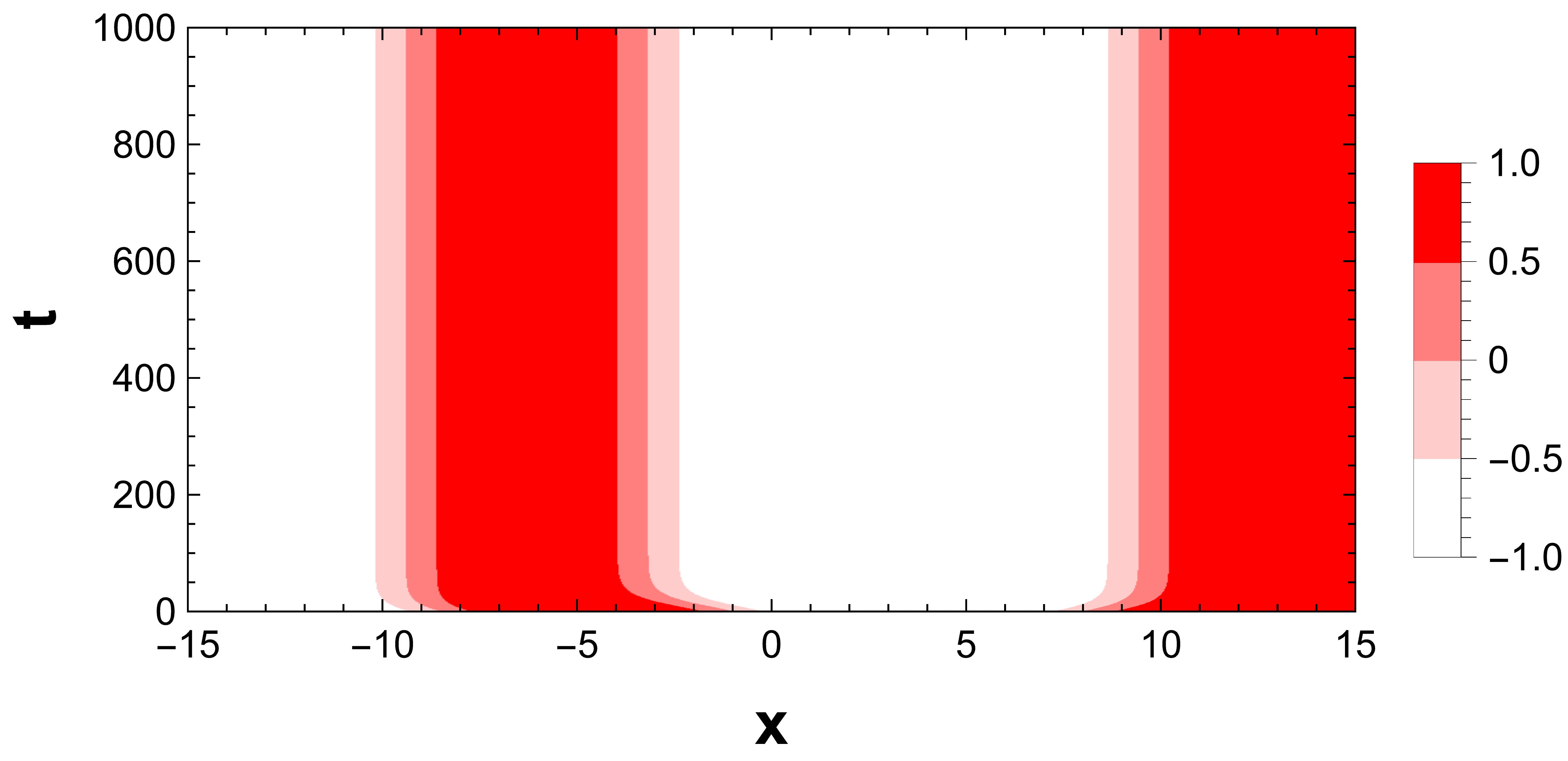}
	\end{subfigure}
	\caption{Simulations of \eqref{eq:mainEquation} on the $x$-interval $(-15,15)$ (with homogeneous Neumann boundary conditions) with spatially periodic inhomogeneity given by \eqref{eq:canonicalExample} with $f_1(x) = \cos x$, $f_2(x) \equiv f_3(x) \equiv 0$ starting from the three-front initial condition $U(x,0) =   \tanh(x + 8.5) -  \tanh(x + 1) + \tanh(x - 8)$ and $\varepsilon = 0$ (a), $0.001$ (b), $0.01$ (c) and $0.1$ (d).}
\label{fig:Numerics3Front}
\end{figure}
\\
The abundance of stable stationary $N$-front patterns established by Theorem \ref{th:Nfrontsper} motivates the upcoming theorem as a natural next step.
It shows that any $N$-front configuration with successive fronts that are sufficiently far apart must evolve toward a nearby attracting stationary $N$-front. This is corroborated by the simulations shown in Fig. \ref{fig:Numerics3Front} of an evolving three-front pattern (similar to the two-front example of Fig. \ref{fig:Numerics2Front}): for $\varepsilon$ sufficiently large, {\it i.e.}, if the separation distance between successive fronts is sufficiently large with respect to the $|\log \varepsilon|$ distance scale, the three-front pattern evolves towards a nearby attracting stationary patterns as described by the upcoming theorem (see Fig. \ref{fig:Numerics3Front}(d)), for smaller $\varepsilon$ (Fig \ref{fig:Numerics3Front}(b),(c)), the three-front dynamics are similar those of the homogeneous Allen-Cahn equation (Fig. \ref{fig:Numerics3Front}(a)) in which two-fronts merge and only a stable attracting one-front pattern remains.
\begin{theorem}
\label{th:Nfrontsperattr}
Consider the hypotheses of Theorem \ref{th:Nfrontsper} and assume $K \geq 0$. Let $(\psi_{0,1}(0), \psi_{0,2}(0), ..., \psi_{0,N}(0))$ (with $\psi_{0,j} < \psi_{0,j+1}$) indicate the initial position of an evolving $N$-front pattern in \eqref{eq:mainEquation} (that connects $u_-^\varepsilon(x) = - 1 + \mathcal{O}(\varepsilon)$ at $x \rightarrow - \infty$ to $u_{(-1)^{N+1}}^\varepsilon(x) = (-1)^{N+1} + \mathcal{O}(\varepsilon)$ at $x \rightarrow +\infty$). If $N \geq 2$ and $\Delta \psi_{0,j} = \psi_{0,j+1} - \psi_{0,j}$ increases faster than $|\log \varepsilon|$ as $\varepsilon \to 0$ (for all $j=1,\ldots,N-1$), then the $N$-front pattern given by $(\psi_{0,1}(\tau), \psi_{0,2}(\tau), ..., \psi_{0,N}(\tau))$ is attracted to a nearby stable stationary $N$-front pattern $(\bar{\psi}_1, \bar{\psi}_2, ..., \bar{\psi}_N)$. More specifically, under the flow of \eqref{eq:mainEquation},
\[
\lim_{\tau \to \infty} (\psi_{0,1}(\tau), \psi_{0,2}(\tau), ..., \psi_{0,N}(\tau)) =
(\bar{\psi}_1, \bar{\psi}_2, ..., \bar{\psi}_N) \; \; {\rm with} \; \; |\psi_{0,j} - \bar{\psi}_j| < Y, j=1,2, ..., N.
\]
If $N=1$, then $\lim_{\tau \to \infty} \psi_{0,1}(\tau) = \bar{\psi}_1$  with $|\psi_{0,j} - \bar{\psi}_j| < Y$ -- without any additional conditions.
\end{theorem}
Note that we refrained from giving a precise description of the $N$-front pattern $(\psi_{0,1}(\tau), \psi_{0,2}(\tau), ..., \psi_{0,N}(\tau))$ and especially of its initial condition $(\psi_{0,1}(0), \psi_{0,2}(0), ..., \psi_{0,N}(0))$.
It needs to be sufficiently close (in the norm of an appropriately chosen function space) to the $N$-dimension manifold to which the dynamics of $N$-front patterns is attracted. In fact, the proof of Theorem~\ref{th:Nfrontsperattr} follows directly from the existence and attractivity of this manifold and the associated validity of $N$-front interaction ODE \eqref{eq:dynNfronts}.
Note also that, in the formulation of this theorem, we excluded (non-generic) situations in which an $N$-front pattern is exactly in the stable manifold of an unstable critical point $(\bar{\psi}_1, \bar{\psi}_2, ..., \bar{\psi}_N)$ of saddle type.

\section{Discussion}
\label{sec:Discussion}
\begin{figure}
\centering
\begin{subfigure}[t]{0.25\textwidth}
\centering
\includegraphics[width=\textwidth]{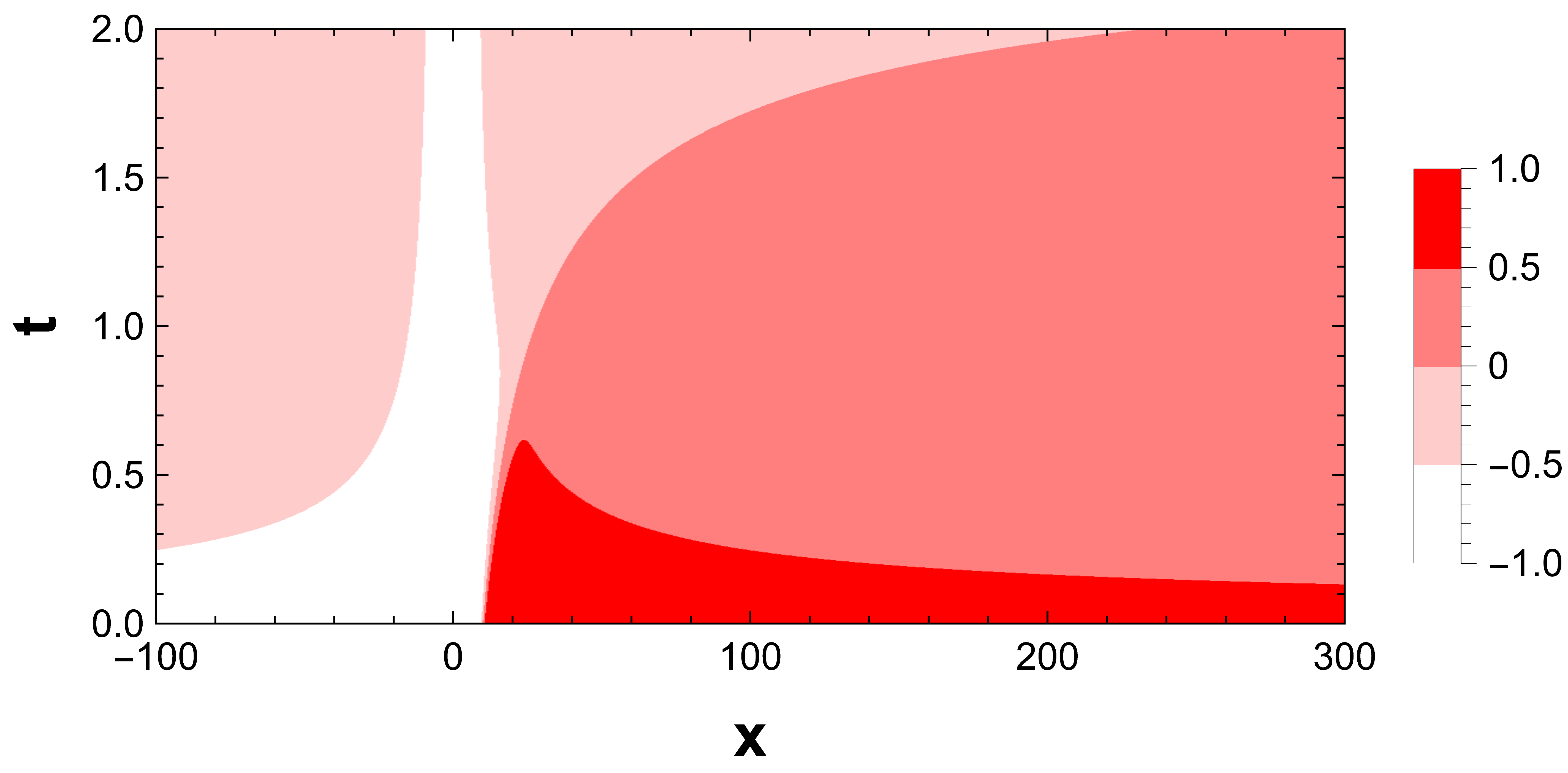}
\end{subfigure}
~
\begin{subfigure}[t]{0.21\textwidth}
\centering
\includegraphics[width = \textwidth]{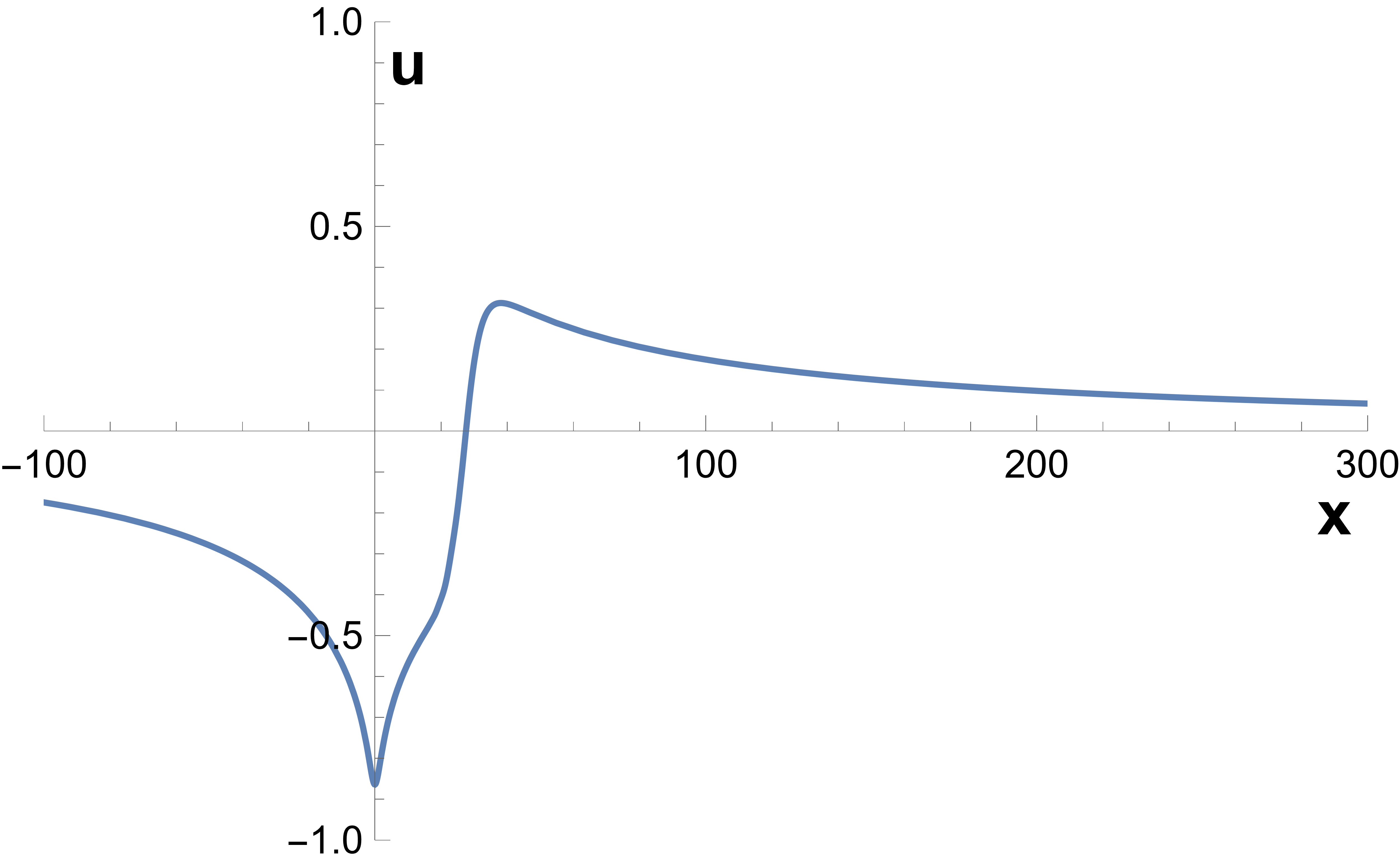}
\end{subfigure}
~
\begin{subfigure}[t]{0.25\textwidth}
\centering
\includegraphics[width = \textwidth]{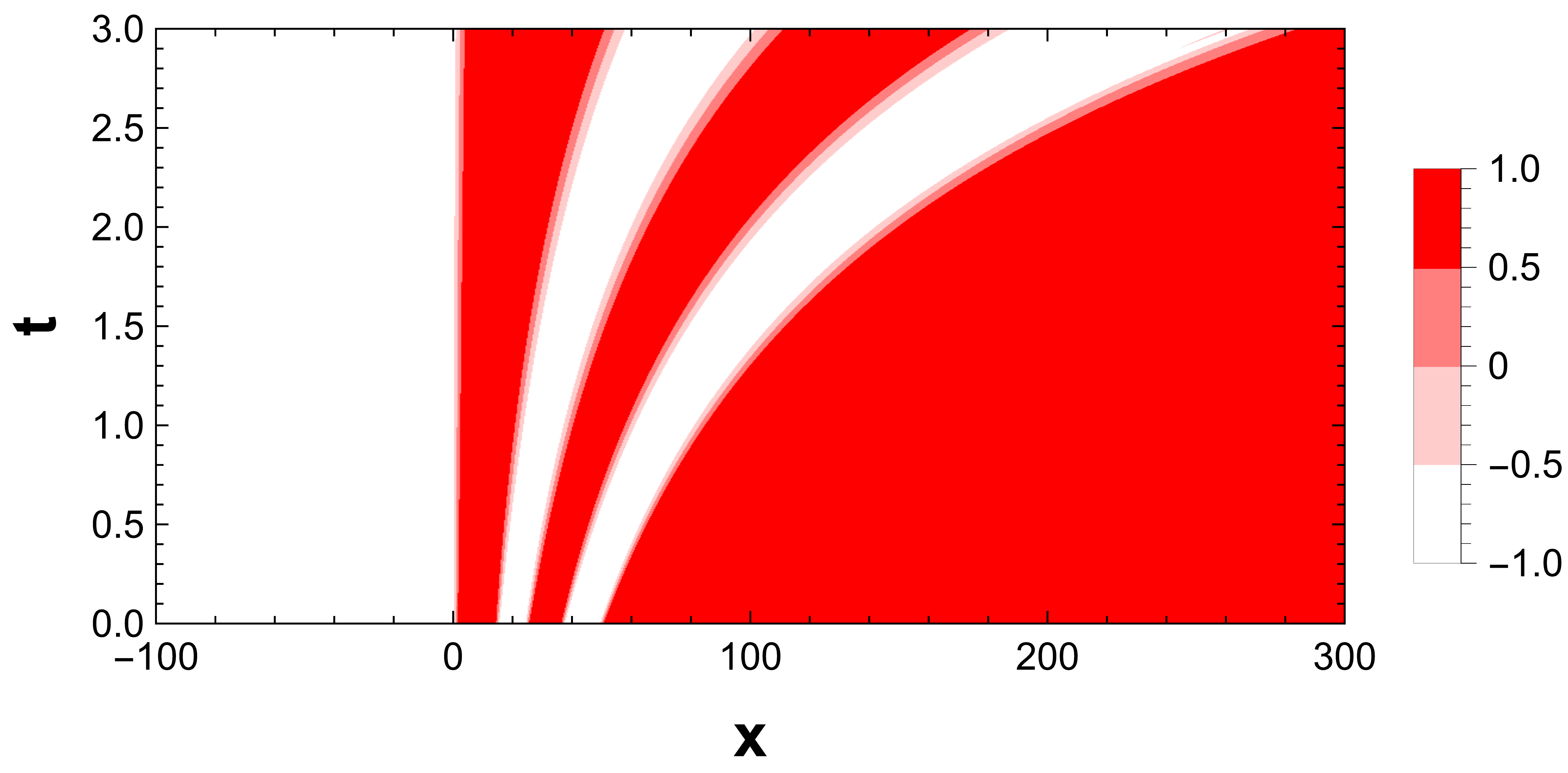}
\end{subfigure}
~
\begin{subfigure}[t]{0.21\textwidth}
\centering
\includegraphics[width = \textwidth]{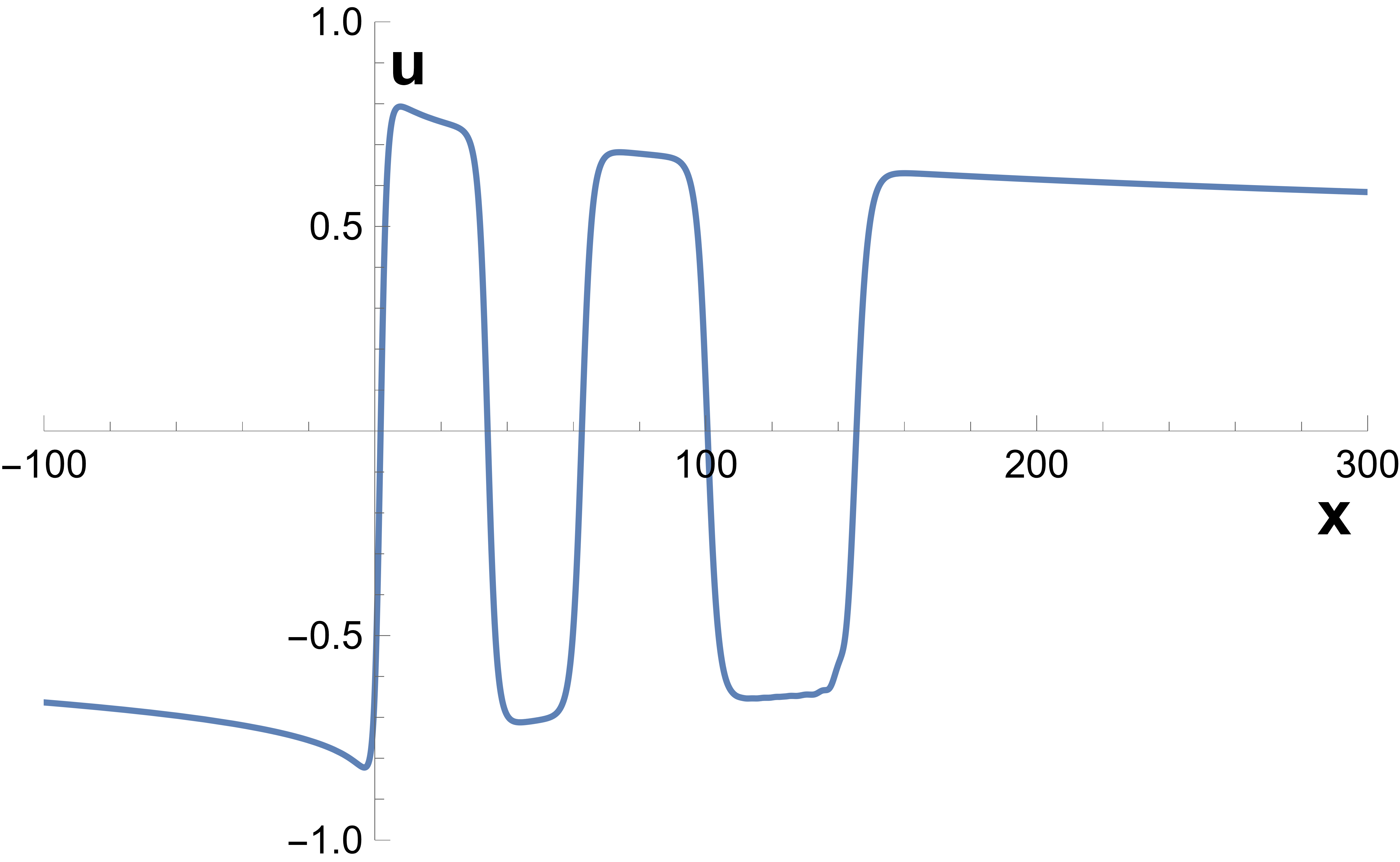}
\end{subfigure}
\caption{Simulations of \eqref{eq:mainEquation} with  topographical inhomogeneity \eqref{eq:Ftopography} $H(x) = -H_{\rm alg}(x; p)$ \eqref{eq:defHuni-alg} for $p<-1$: $H''(x)$ has also become unbounded (with $\varepsilon = 0.1$, $x \in (-100,300)$, homogeneous Neumann boundary conditions). (a)/(b) One-front dynamics for $p=-1.5$ (with initial condition $U(x,0) = \tanh(x - 10)$); $t=1.0$ in (b).  (c)/(d) Five-front dynamics for $p=-1.2$ (with initial condition $U(x,0) = \tanh(x - 1)  - \tanh(x - 15) + \tanh(x - 25) - \tanh(x - 37) + \tanh(x - 50)$); $t=2.0$ in (d). Note that there seems to be a sharp change in $U_x(x,t)$ near $x=0$ in both simulations ({\it i.e.}, in (b) and (d)), however zooming in shows that the solution remains smooth (in the $t$-intervals considered).}
\label{fig:BeyondValidity}
\end{figure}

In this work, we have analyzed several aspects of multi-front dynamics of the Allen-Cahn equation \eqref{eq:mainEquation} driven by general, small-amplitude spatial heterogeneities. Our research is not comprehensive, for simplicity we for instance focused somewhat on `topographical' heterogeneities \eqref{eq:Ftopography}
that are either purely periodic or spatially localized (with uniformly decaying tails that decay with the same rates as $x \to \pm \infty$). Nevertheless, the methods we have employed may be used to investigate the more general -- and typically more technical -- situations. As examples, in this section, we briefly discuss the effect of (localized) heterogeneous driving terms on the process of coarsening and the `mixed' case of a localized inhomogeneity limiting on a spatially periodic background.
In addition, we discuss the apparent validity of our analysis into a realm where it is not a priori expected, {\it i.e.}, for localized topographies with `tails' that no longer decay. Finally, we briefly return to the ecosystem setting that motivated part of our research.
\\
\begin{figure}
\centering
	\begin{subfigure}[t]{0.24\textwidth}
		\centering
		\includegraphics[width=\linewidth]{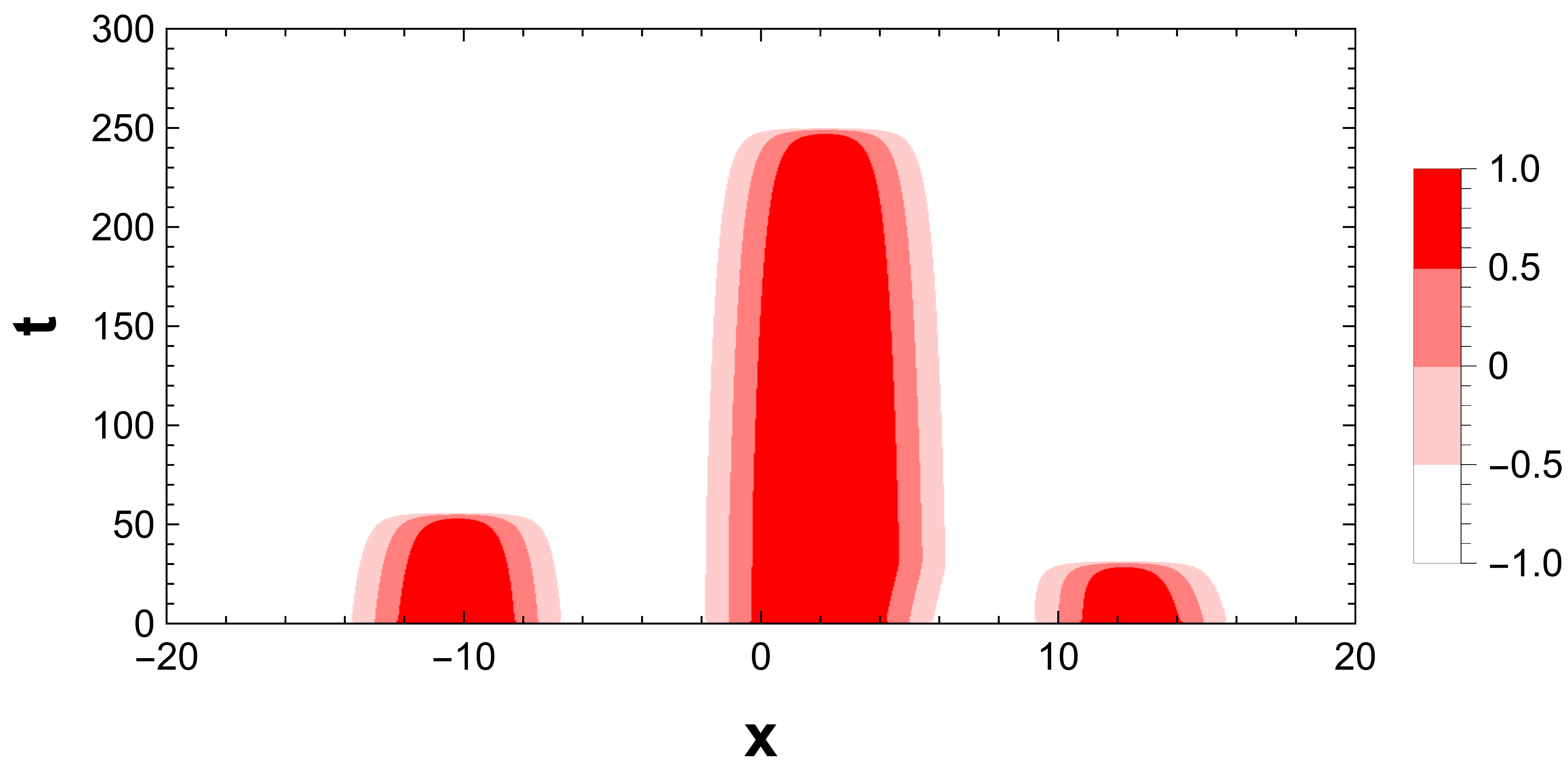}
	\end{subfigure}
	\begin{subfigure}[t]{0.24\textwidth}
		\centering
		\includegraphics[width=\linewidth]{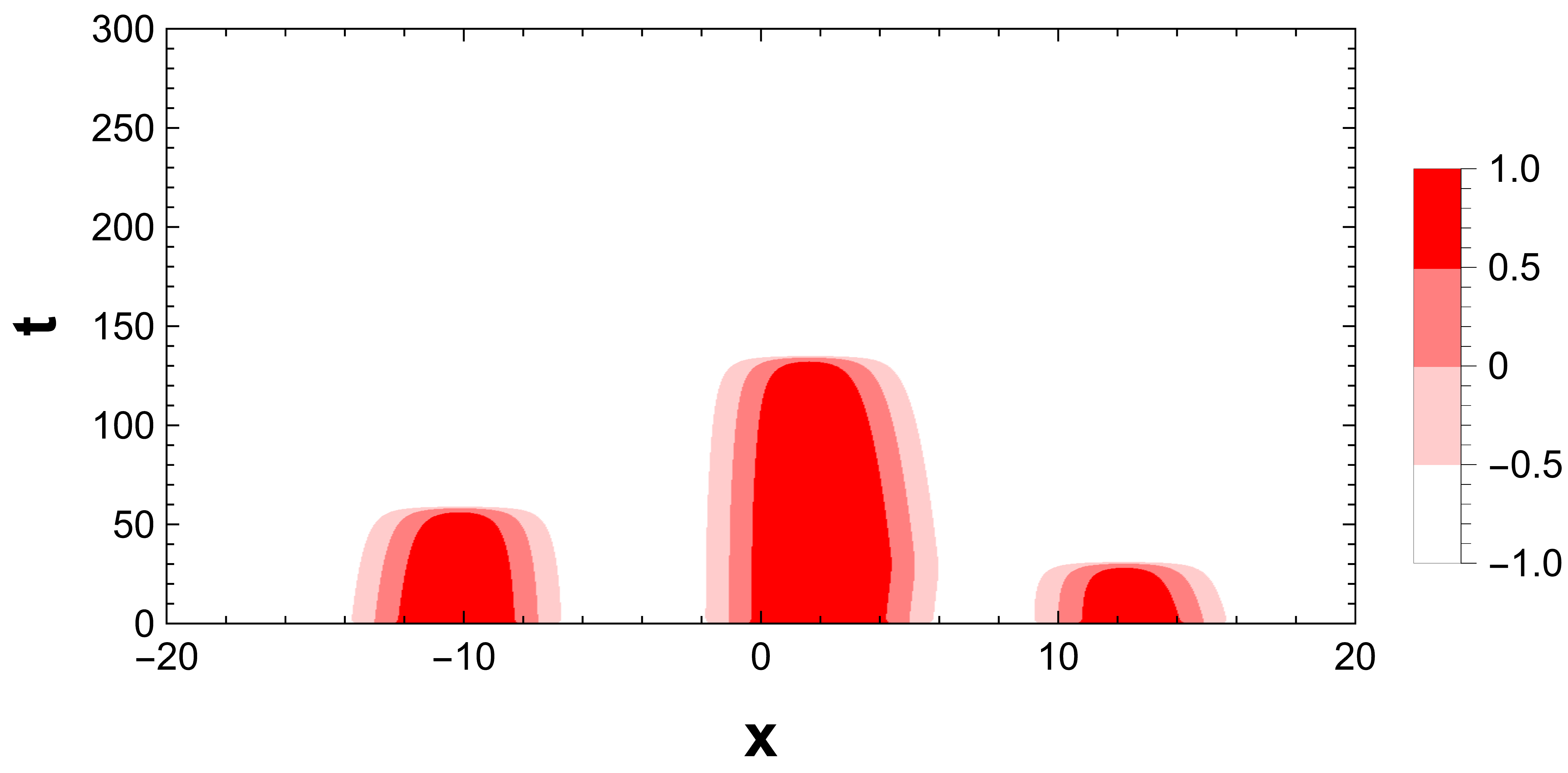}
	\end{subfigure}
	\begin{subfigure}[t]{0.24\textwidth}
		\centering
		\includegraphics[width=\linewidth]{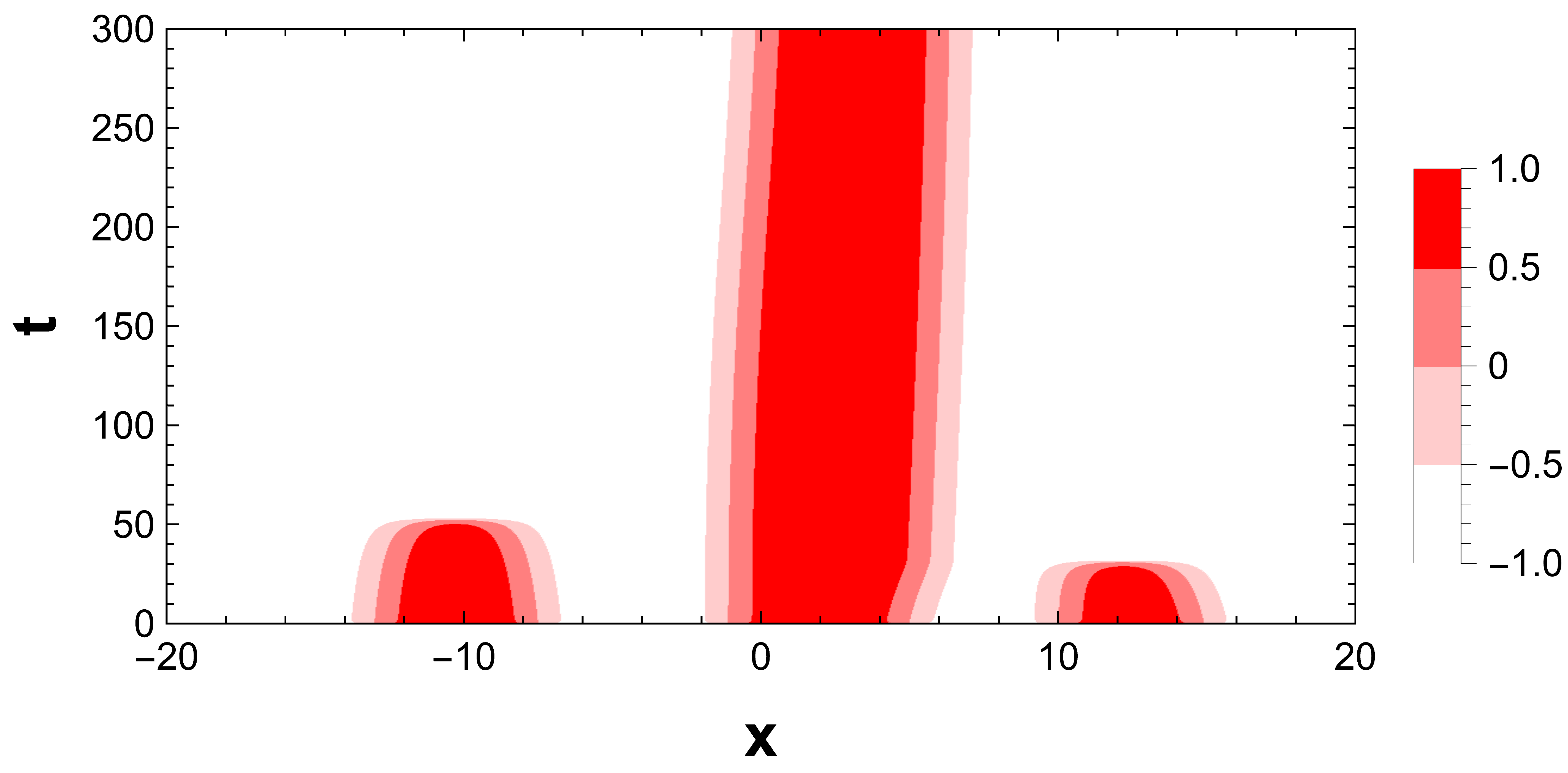}
	\end{subfigure}
	\begin{subfigure}[t]{0.24\textwidth}
		\centering
		\includegraphics[width=\linewidth]{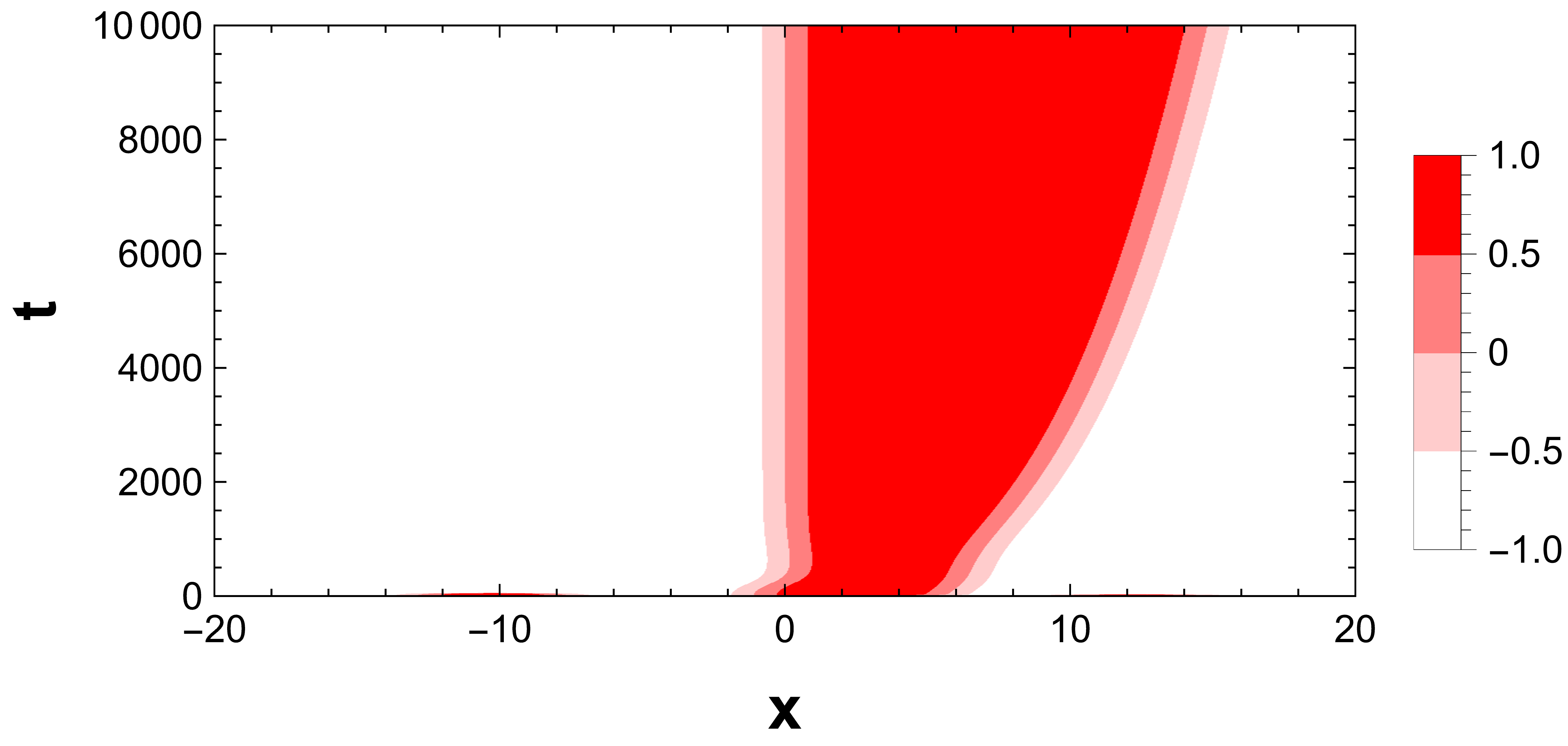}
	\end{subfigure}
	\caption{Simulations of \eqref{eq:mainEquation} with a localized topographical inhomogeneity \eqref{eq:Ftopography} ($\varepsilon = 0.1$) on the interval $(-30,30)$ with homogeneous Neumann boundary conditions (plotted on $(-20,20)$) starting from the six-front initial condition $U(x,0) = -1 + \tanh(x +13.0)  - \tanh(x - 7.5) + \tanh(x + 1.1)  - \tanh(x - 5.0) + \tanh(x - 10.0)  - \tanh(x - 15.9)$. (a) The flat homogeneous Allen-Cahn case: $H(x) \equiv 0$. (b) A weak valley: $H(x) = - H_{\rm alg}(x; 2)$. (c)/(d) A weak hill: $H(x) = H_{\rm alg}(x; 2)$.}
\label{fig:Coarsening}
\end{figure}
\\
{\bf The impact on coarsening.} The Cahn-Hilliard and Allen-Cahn equations are the two prototypical, and thus most studied, models for the process of coarsening -- see \cite{nepomnyashchy2015coarsening, pereira2022steady} and the references therein. As follows from the classical literature \cite{Carr1989,Chen2004,Fusco1989}, any multi-front initial condition of the homogeneous Allen-Cahn equation \eqref{eq:standardAC} evolves towards either one of the stable steady states $U(x,t) \equiv u_{\pm} = \pm 1$ (Sec.~\ref{s:uplusminus}) or towards a stationary one-front solution
$u_{\rm up, down}(x;\phi)$ \eqref{eq:heteroclinicSolution}.
\\
\\
In this article, we have shown that the coarsening process can be altered or even stopped by spatial heterogeneities.
We showed that spatially periodic heterogeneities stop the coarsening process completely in multi-front patterns in which there is a sufficiently large separation distance between the fronts, recall
Theorem \ref{th:Nfrontsperattr}.
(We add that, while Theorem~\ref{th:Nfrontsperattr} only considers topographic heterogeneities, a similar result can be proved for general spatially periodic inhomogeneities $\varepsilon F(U,V,x)$ in \eqref{eq:mainEquation}.)
Also, the simulations of Figs. \ref{fig:5Fronts-Intro}(a) and \ref{fig:Numerics3Front} show that there is an interplay between coarsening and pinning for (spatially) periodically driven multi-front patterns with fronts that are so close to each other that the general result of Theorem \ref{th:Nfrontsperattr} cannot hold.
\\
\\
We also studied spatial heterogeneities which consist of a spatially localized term, and we showed that here the picture is more complex. In Fig.~\ref{fig:Coarsening}, simulations are presented of \eqref{eq:mainEquation} with two different localized topographies, $H(x) = \pm H_{\rm alg}(x; 2)$, that show two very different outcomes. The topography $H(x) = - H_{\rm alg}(x; 2)$ speeds up, and thus strengthens, the coarsening process, as can also be deduced from \eqref{eq:dynNfronts}.
Fronts that are far away from the origin are driven towards the origin by the inhomogeneity, so that the full extent of the multi-front pattern shrinks by a persistent `squeezing process'.
More precisely -- and in general -- the extent of an $N$-front pattern is determined by the distance $\psi_N - \psi_1$ between its two outermost fronts.
In fact, in this case, \eqref{ddtphiNphi1} generalizes to
\[
\frac{d}{dt}(\psi_N - \psi_1) = -\varepsilon \left[\mathcal{S}(\psi_N) - \mathcal{S}(\psi_1) \right] - 16 \left[e^{(\psi_2-\psi_1)} + e^{(\psi_N-\psi_{N-1})} \right],
\]
recall also \eqref{eq:defpsiStau} and \eqref{eq:dynNfronts}.
Thus, it follows from Lemma \ref{lem:algdecayH} for $\psi_1 \ll -1 < 1 \ll \psi_N$ and $H'(x)$ localized with algebraically decaying tails that (to leading order),
\begin{equation}
\label{eq:inhomsqueeze}
\frac{d}{dt}(\psi_N - \psi_1) = -\frac{2^{\frac12(p+3)}}{3} \varepsilon \left[\tilde{h}_+(p) \psi_N^{-p} - \tilde{h}_-(p) |\psi_1|^{-p} \right] - 16 \left[e^{(\psi_2-\psi_1)} + e^{(\psi_N-\psi_{N-1})} \right].
\end{equation}
Since $\tilde{h}_\pm(p) = \pm (1-p)$ for $H_{\rm alg}(x; p)$, it follows that $H(x) = - H_{\rm alg}(x; p)$ with $p > 1$ indeed strengthens the homogeneous squeezing process of \eqref{ddtphiNphi1}. Likewise, the topography $H(x) = + H_{\rm alg}(x; p)$ with $p > 1$ counteracts the homogeneous process. In fact, due to the exponential decay of the homogeneous terms, the inhomogeneous effects of \eqref{eq:inhomsqueeze} dominate when the fronts are sufficiently far apart: the extent of a multi-front pattern will grow for $\tilde{h}_+(p) < 0 < \tilde{h}_-(p)$ -- see the simulation in Fig. \ref{fig:Coarsening}(c),(d) that shows an evolving two-front pattern consisting of a stationary front pinned at $x=0$ and a second front traveling towards $x=+\infty$. Additional examples are presented in Figs. \ref{fig:NumericsTwoThreeFronts-Intro}(b),(e) and \ref{fig:5Fronts-Intro}(b). Each of these examples shows that a spatially localized inhomogeneity may yield an endstate that is significantly different from the endstate of the homogeneous system.
Here, for example, we compare Fig. \ref{fig:NumericsTwoThreeFronts-Intro}(a) to \ref{fig:NumericsTwoThreeFronts-Intro}(b) and Fig. \ref{fig:NumericsTwoThreeFronts-Intro}(d) to \ref{fig:NumericsTwoThreeFronts-Intro}(e);
also, in Fig.\ref{fig:5Fronts-Intro}(b), the homogeneous one-front endstate of a five-front initial condition is replaced by a situation in which eventually $U \to -1$ (to leading order in $\varepsilon$) on the entire domain, while in Fig. \ref{fig:Coarsening}(d) $U \to +1$ (to leading order in $\varepsilon$) on the positive $x$-axis, instead of $U \to -1$ on $\mathbb{R}$ as in Fig. \ref{fig:Coarsening}(a).
Apart from these effects, spatially localized inhomogeneities may also significantly delay the coarsening process (see again Fig. \ref{fig:5Fronts-Intro}(b)).
\\
\begin{figure}
\centering
	\begin{subfigure}[t]{0.24\textwidth}
		\centering
		\includegraphics[width=\linewidth]{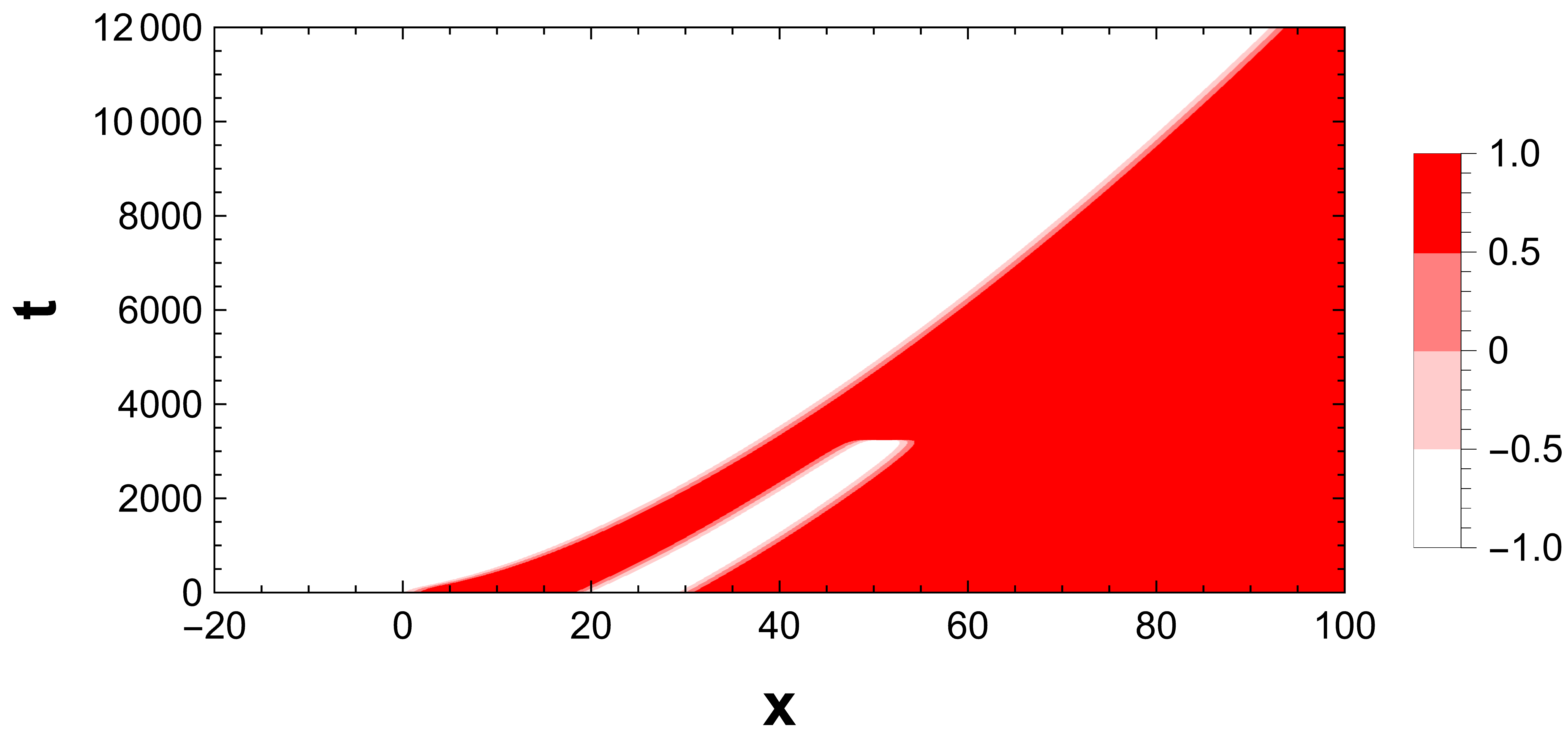}
	\end{subfigure}
	\begin{subfigure}[t]{0.24\textwidth}
		\centering
		\includegraphics[width=\linewidth]{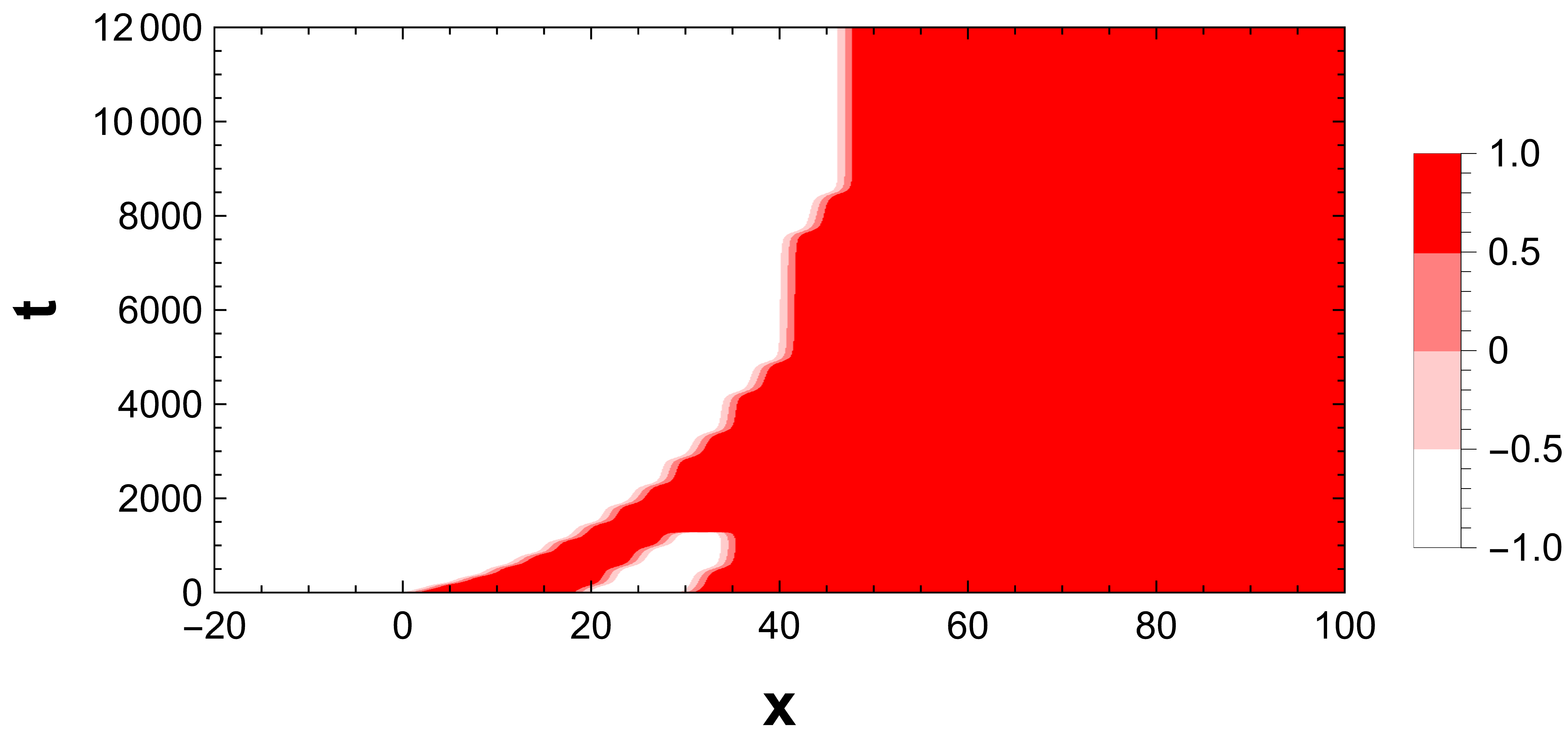}
	\end{subfigure}
	\begin{subfigure}[t]{0.24\textwidth}
		\centering
		\includegraphics[width=\linewidth]{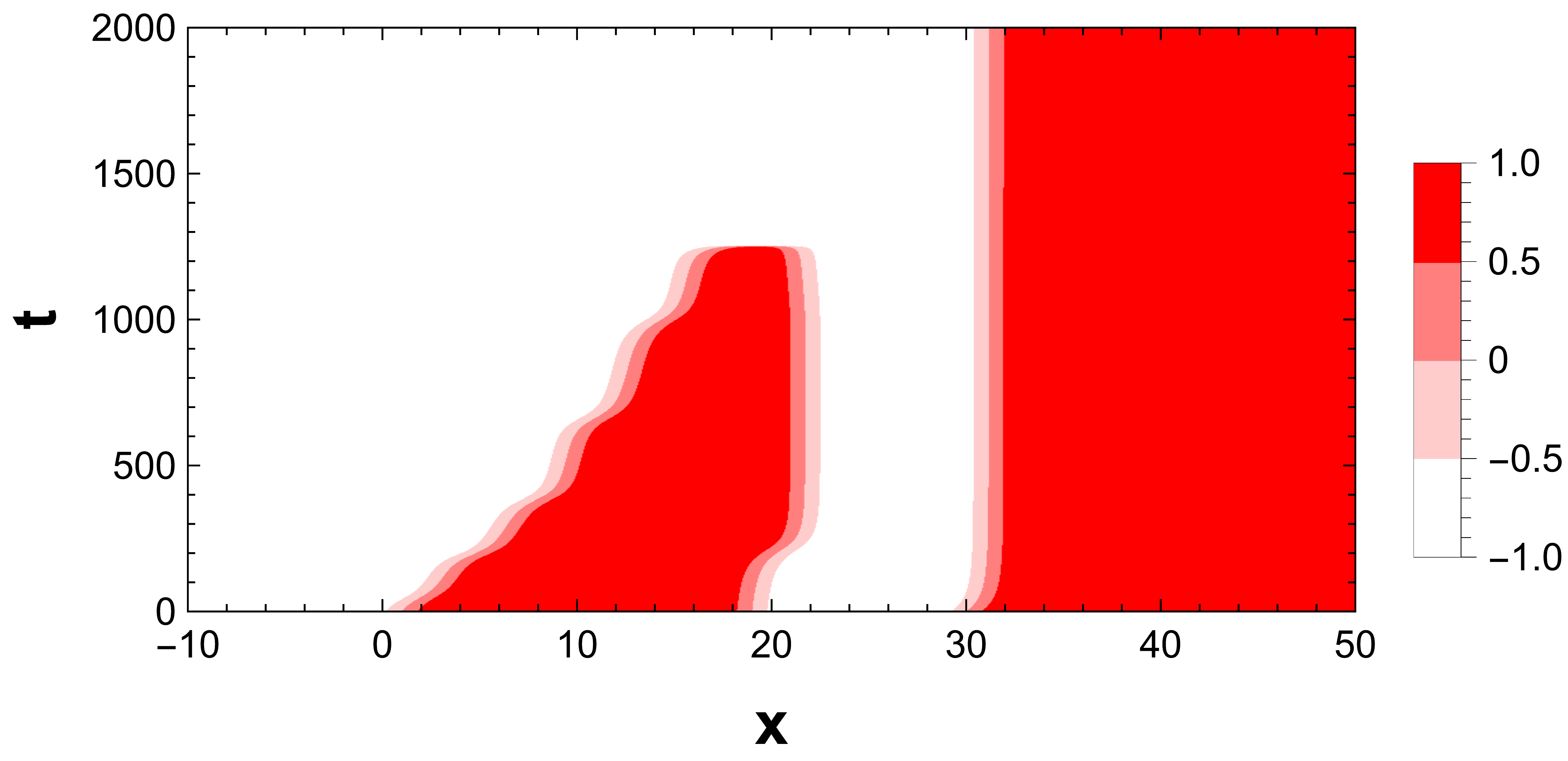}
	\end{subfigure}
	\begin{subfigure}[t]{0.24\textwidth}
		\centering
		\includegraphics[width=\linewidth]{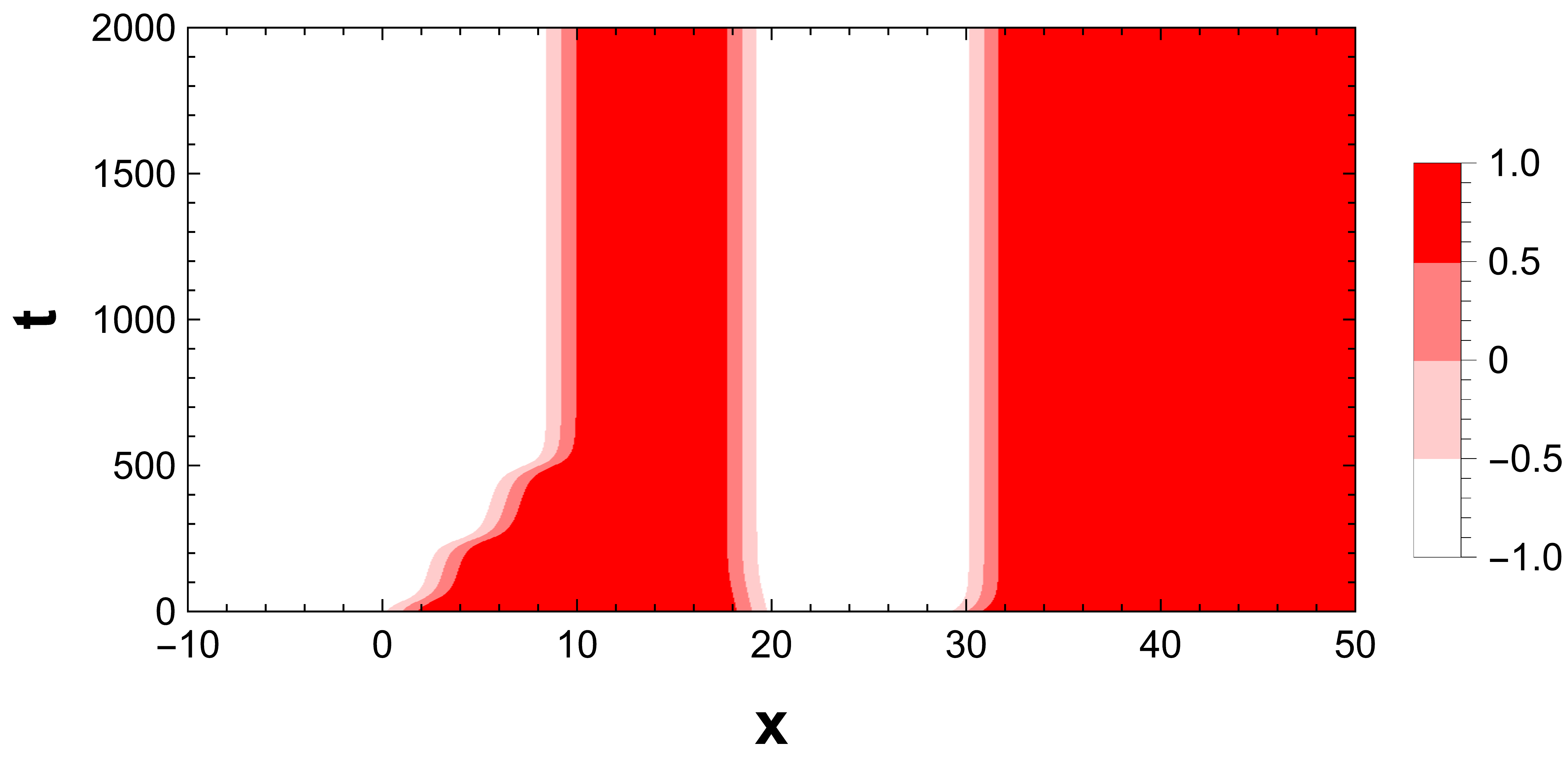}
	\end{subfigure}
	\caption{Simulations of \eqref{eq:mainEquation} with `mixed' topographical inhomogeneity \eqref{eq:Ftopography} given by $H(x;\delta) =-H_{\rm alg}(x; 0.5) + \delta \sin 3x$ \eqref{eq:defHuni-alg}, $\varepsilon = 0.1$ and three-front initial condition $U(x,0) = \tanh(x - 1)  - \tanh(x - 19) + \tanh(x - 30)$; simulations are done for $x \in (-200,200)$ with homogeneous Neumann boundary conditions. (a) $\delta = 0$, (b) $\delta = 0.025$, (c) $\delta = 0.04$ and (d) $\delta = 0.06$.}
\label{fig:3FrontLocPer}
\end{figure}
\\
All these phenomena can be understood by the approach developed here in combination with the classical insights of \cite{Carr1989,Chen2004,Fusco1989}, see also the final paragraph of Sec.~\ref{sec:ODENfrontloc}. This is not the case for the simulations shown in Fig. \ref{fig:5Fronts-Intro}(c) and below in Fig. \ref{fig:5FrontDynVaryingP}(c), which suggest that spatial inhomogeneities may not only stop the coarsening process, they may also drive a spreading process that (continuously) increases the extent of the multi-front pattern without reducing the number of fronts. See the discussion below.
\\
\\
{\bf General topographies.} Other topographies, besides the spatially localized and periodic topographies cases considered in Secs.~5 and 6, involving general inhomogeneous terms $F(U,V,x)$ \eqref{eq:mainEquation} can be studied by our methods.
As a simple example of this, we present in Fig. \ref{fig:3FrontLocPer} simulations of \eqref{eq:mainEquation} with a mixed topographical inhomogeneity consisting of the sum of a localized topography and $\delta$ times a spatially periodic topography.
The associated driving Melnikov function $\mathcal{S}(\psi)$ \eqref{eq:dynNfronts} can be written as $\mathcal{S}(\psi) = \mathcal{S}_{\rm loc}(\psi) + \delta \mathcal{S}_{\rm per}(\psi)$. For $\delta$ sufficiently small, the pinning impact of $\mathcal{S}_{\rm per}(\psi)$ only kicks in for very large $\psi$, {\it i.e.,} when $\log \psi = p |\log \delta|$ to leading order in $\delta$ (recall Lemma \ref{lem:algdecayH} and equation \eqref{eq:defpsiStau}). Thus, while the $\psi_2$-fronts and $\psi_3$-fronts merge and $\psi_1$ travels as a solitary front to $+\infty$ (with decaying speed) in the homogeneous three-front simulation of Fig. \ref{fig:3FrontLocPer}(a) (see Sec.~\ref{sec:ODENfrontloc} and Figs.~\ref{fig:5FrontDynamics}, \ref{fig:5FrontDynVaryingP}(a)), a small spatially periodic term slows down the $\psi_{1,2,3}-$fronts a bit and the remaining $\psi_1$-front eventually gets pinned.
Hence, the configuration now approaches a (pinned) stationary $\psi_1$ one-front pattern (Fig. \ref{fig:3FrontLocPer}(b)).
As $\delta$ increases further, $\psi_3$ gets pinned before it can be reached by $\psi_2$. However, since $\psi_2$ is also slowed down by the periodic effects, it is overtaken by $\psi_1$.
That is, the solution also approaches a pinned stationary one-front pattern, but now originating from $\psi_3$.
Finally, for still larger values of $\delta$, the $\psi_1$-front gets pinned before it reaches (the also pinned) $\psi_2$ front, and the entire configuration approaches a stationary three-front.
\\
\\
{\bf The validity of the $N$-front interaction equations.}
As noted at the beginning of Sec.~\ref{sec:interactiondynamics} and discussed in somewhat more detail at the end of Sec.~\ref{sec:ODEderivation}, both the existence of an attracting approximate $N$-dimensional manifold $\mathcal{M}_N^0$ for the (semi-)flow generated by \eqref{eq:mainEquation} and the validity of  \eqref{eq:NFrontODE}/\eqref{eq:dynNfronts} as leading order descriptions of the flow on $\mathcal{M}_N^0$ can be established by a direct application of existing methods \cite{ei2002pulse, promislow2002renormalization}.
However, manifold $\mathcal{M}_N^0$ has boundaries.
The dynamics near the boundaries associated to fronts getting too close to each other (or to fronts beyond these boundaries) can again be controlled by existing classical methods, see especially \cite{Chen2004}.
Then, for the opposite situations in which the fronts are widely separated, one also needs to verify that the neglected $\mathcal{O}(\varepsilon^2)$ terms $\varepsilon^2 R^{j,N}_{2}(\psi_1,...,\psi_N;\varepsilon)$ ($j=1,...,N$) in \eqref{eq:dynNfronts} remain of higher order.
In general, one expects that the magnitudes of the $R^{j,N}_{2}(\psi_1,...,\psi_N;\varepsilon)$ terms will decrease as $|\psi_i| \gg 1$.
Nevertheless, at this depth of the asymptotic analysis, it could be that the front interactions involve more than just nearest neighbor interactions and depend on other fronts further away -- or even on all other fronts, as is the case for localized structures in semi-strong interaction \cite{doelman2003semistrong, van2010front}.
\\
\\
Relaxation of the hypothesis that $F$ is bounded provides another compelling reason to analyze the validity of the ODE system \eqref{eq:dynNfronts} beyond the classical weak interaction cases considered in \cite{ei2002pulse, promislow2002renormalization}.
Application of these existing methods requires that the inhomogeneous term $\varepsilon F(U,V,x)$ in \eqref{eq:mainEquation} is bounded.
Hence, for topographical inhomogeneities \eqref{eq:Ftopography}, we have required in this article that $H'(x)$ is bounded (for all $x \in \mathbb{R}$). However, as we already noted in \eqref{eq:uMinusPlusTopo} at the very beginning of our analysis, the asymptotic results we derived may extend beyond the case in which $H'(x)$ is bounded.
It may suffice to require only that $H''(x)$ is bounded.
\\
\\
In order to study the impact of this relaxation of the hypothesis on the (coarsening) dynamics of $N$-front configurations, we redo the one-front analysis of Sec.~\ref{sec:ODE2frontloc} and consider the evolution of a solitary front $\psi(\tau)$ that starts out with $\psi(0) = \psi_0 \gg 1$, but now at the (positive) tail of a Melnikov function $\mathcal{S}(\psi)$ associated to a topography with algebraic tails, {\it i.e.}, we use Lemma \ref{lem:algdecayH} instead of Lemma \ref{lem:expdecayH}.
We find that, under the assumptions that $\tilde{h}_+(p) < 0$ and that the higher order corrections to \eqref{eq:dynNfronts} remain of higher order as $\mathcal{S(\psi)}$ becomes small, the algebraic counterparts of \eqref{eq:1frontfareqexp} and \eqref{eq:1frontfarsolexp} are given by
\begin{equation}
\label{eq:1frontfareqsolalg}
\frac{d\psi}{d\tau} = \varepsilon \tilde{H}_+ \psi^{-p},
\; \ \ \ \ \
\tilde{H}_+(p) = \frac{2^{\frac12(p+3)}}{3}|\tilde{h}_+(p)| > 0, \; \ \ \ \ \
\psi(\tau) = \left(\varepsilon \tilde{H}_+ \tau + \psi_0^{p+1}\right)^{\frac{1}{p+1}}, \; \; p \neq -1.
\end{equation}
Hence, the leading order distance between adjacent fronts $\tilde{\psi}(\tau) < \psi(\tau)$, which start out sufficiently far from each other, evolves as
\begin{equation}
\label{eq:deltafrontsalg}
\psi(\tau) - \tilde{\psi}(\tau) =
\frac{1}{(p+1)^{\frac{2p+1}{p+1}}} \frac{\psi_0^{p+1} - \tilde{\psi}_0^{p+1}}{(\varepsilon \tilde{H}_+ \tau)^{\frac{p}{p+1}}} + \;   {\rm h.o.t.} \ \ \ {\rm for} \; \; \tau \gg \frac{\psi_0^{p+1}}{\varepsilon}.
\end{equation}
This is the algebraic counterpart of \eqref{eq:deltafrontsexp}.
We now examine the dynamics for $p>0, p=0, $ and $p<0$, respectively, beginning with the former.
\\
\\
For $p > 0$, {\it i.e.}, for the case in which the validity of \eqref{eq:dynNfronts} can be settled following \cite{ei2002pulse, promislow2002renormalization}, the situation is similar to that of exponentially decaying topographies.
The distance $\psi(\tau) - \tilde{\psi}(\tau)$ decreases (algebraically in time), $\tilde{\psi}(\tau)$ catches up (eventually) with $\psi(\tau)$, and the two fronts merge and disappear.
This scenario is confirmed by the simulations shown in Figs. \ref{fig:5Fronts-Intro}(b) and  \ref{fig:5FrontDynVaryingP}(a), in which the initial data consist of five-front data.
\\ \\
The case $p=0$ is special, since \eqref{eq:deltafrontsalg} predicts that $\psi(\tau) - \tilde{\psi}(\tau)$ will approach a constant value.
Although this is only a leading order result, it  too is corroborated by simulations.
As shown for example in Fig.~\ref{fig:5FrontDynVaryingP}(b), the $\psi_1$-front and $\psi_2$-front seem to `lock' to each other and travel with constant speed to the left maintaining a constant distance, while the $\psi_5$-front travels as a solitary front to the right with the same constant speed. (We add that this all happens after the $\psi_3$-front and $\psi_4$-front have merged, as in Fig. \ref{fig:5FrontDynVaryingP}(a)).)
\\ \\
A significant change occurs as $p$ takes on negative values.
Namely, for $-1 < p < 0$, \eqref{eq:deltafrontsalg} predicts that the distance $\psi(\tau) - \tilde{\psi}(\tau)$ between two successive fronts will increase.
Although the inhomogeneous term $F(U,V,x)$ in \eqref{eq:mainEquation} is no longer bounded (since $H'(x)$ is not bounded), this prediction is again confirmed by the numerics.
See Figs. \ref{fig:5Fronts-Intro}(c) and \ref{fig:5FrontDynVaryingP}(c).
In both simulations, all of the distances between neighboring fronts increase.
As a consequence, no coarsening takes place (if the initial distances between neighboring fronts are not too small).
The end-states of the simulations are formed by dispersing trains of fronts, either moving to $-\infty$ on the left (see the $(\psi_1,\psi_2)$-train in  Fig. \ref{fig:5FrontDynVaryingP}(c)) or to $+\infty$ on the right (see the $(\psi_3,\psi_4,\psi_5)$-train in Fig. \ref{fig:5FrontDynVaryingP}(c) and the $(\psi_1,\psi_2,\psi_3,\psi_4,\psi_5)$-train in Fig. \ref{fig:5Fronts-Intro}(c)).
\\ \\
Thus, the simulations over this entire range of $p$ values strongly suggest that the front interaction system \eqref{eq:dynNfronts} remains valid beyond the natural condition that inhomogeneity $F(U,V,x)$ in \eqref{eq:mainEquation} is bounded (for all $x \in \mathbb{R}$). However, both the above analysis and further simulations indicate that there is a second critical value of $p$ for which this may no longer be the case. As $p$ decreases below $-1$, $H''_{\rm alg}(x; p)$ also is no longer bounded -- and thus also the very first step in our asymptotic approach needs to be reconsidered (see \eqref{eq:uMinusPlusTopo}, the background state $u_\pm^\varepsilon(x)$ is no longer bounded on $\mathbb{R}$ if $H''(x)$ is not).
The simulations of Fig. \ref{fig:BeyondValidity} show that the front dynamics of \eqref{eq:mainEquation} indeed have gone through an essential change in character; namely, the multi-front patterns can no longer be seen as perturbations of combinations of the homogeneous one-front patterns
$u_\textrm{up}(x;\phi)$ and $u_\textrm{down}(x;\phi)$ \eqref{eq:heteroclinicSolution}, which is the backbone of our asymptotic approach. Nevertheless, the patterns exhibited by \eqref{eq:mainEquation} still have a clear multi-front nature, and one might expect that also the case $p < -1$, perhaps combined with the requirement that $|p+1| \ll 1$, could be studied analytically. Such a study would be especially interesting, since the (now quite formal) approximations in \eqref{eq:1frontfareqsolalg} suggest that the position $\psi(\tau)$ of a front may have a finite-time blow-up character, {\it i.e.}, its speed may become unbounded as $\tau$ approaches a critical value $\tau_\ast$. Although simulating \eqref{eq:mainEquation} under these conditions becomes quite subtle, the outcome of these simulations do not contradict this possibility. We refer to Fig. \ref{fig:BeyondValidity} where a simulation is shown in which the fronts travel very fast and keep on speeding up.
\\ \\
Thus, a careful study of the dynamics of $N$-front patterns in \eqref{eq:mainEquation} with inhomogeneity $F(U,V,x)$ -- either or not in the form of \eqref{eq:Ftopography} -- that may become on unbounded (in $x$) would be both interesting and challenging.
\\
\begin{figure}
\centering
\begin{subfigure}[t]{0.32\textwidth}
\centering
\includegraphics[width=\textwidth]{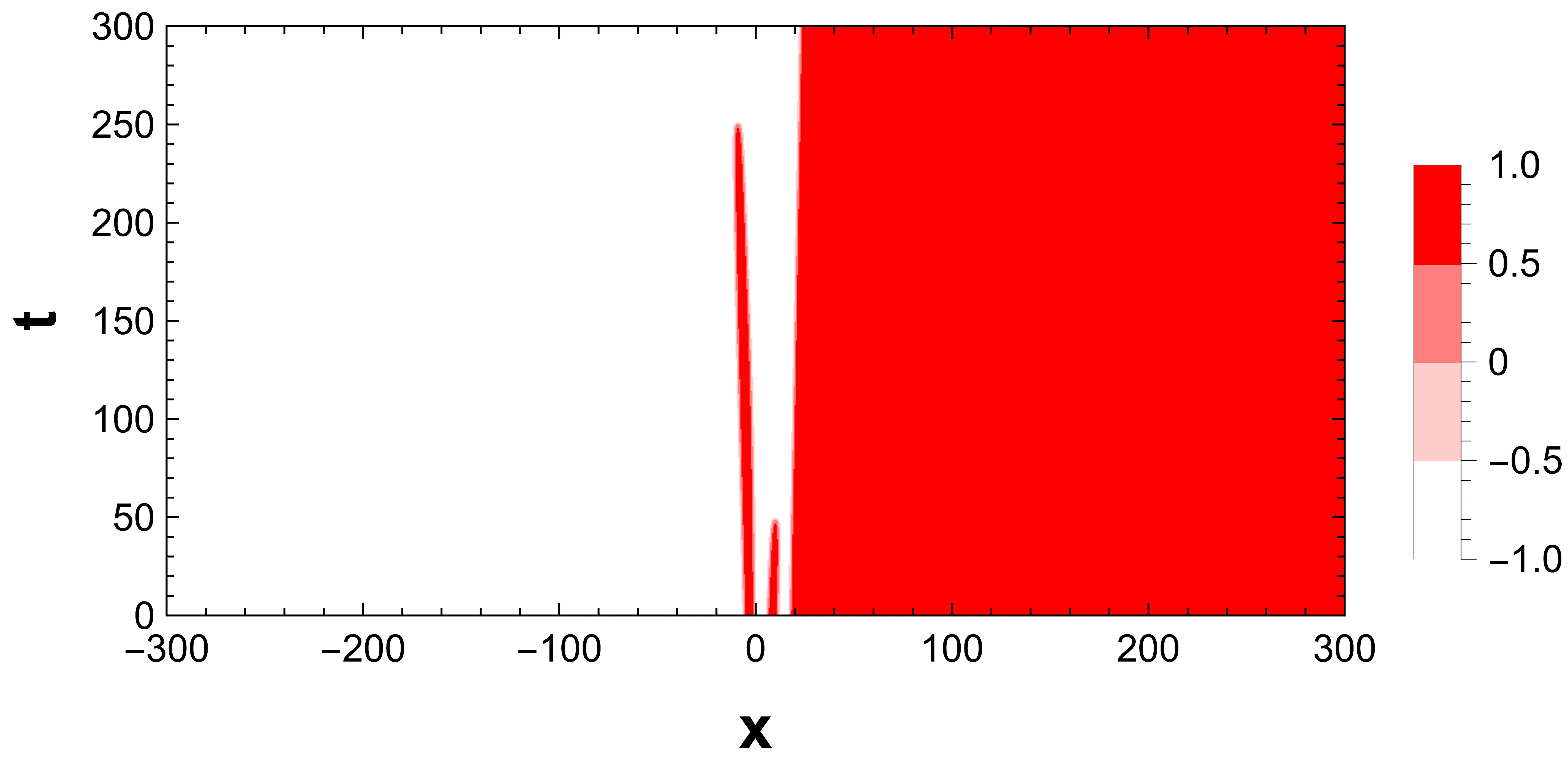}
\end{subfigure}
~
\begin{subfigure}[t]{0.32\textwidth}
\centering
\includegraphics[width=\textwidth]{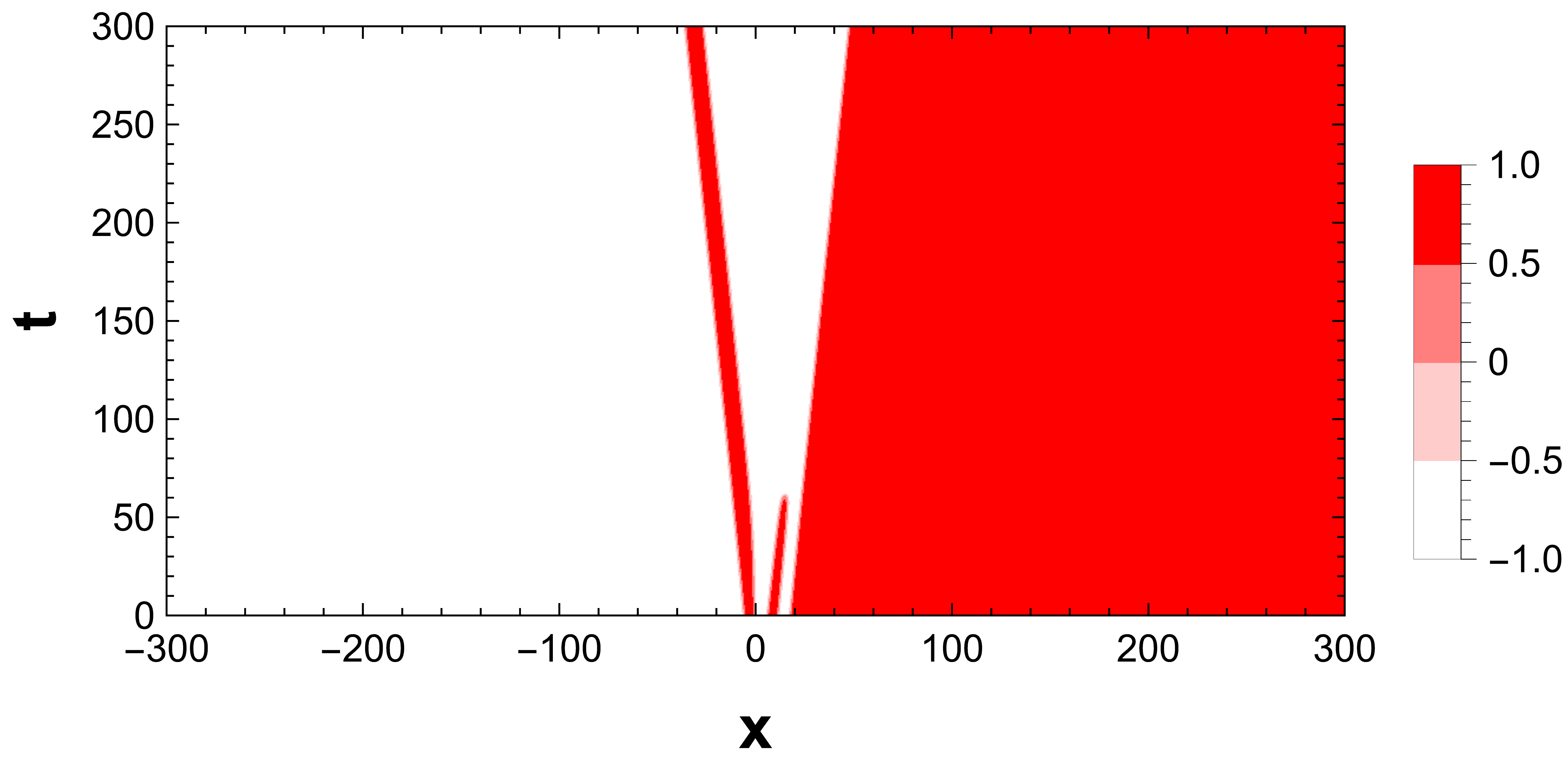}
\end{subfigure}
~
\begin{subfigure}[t]{0.32\textwidth}
\centering
\includegraphics[width = \textwidth]{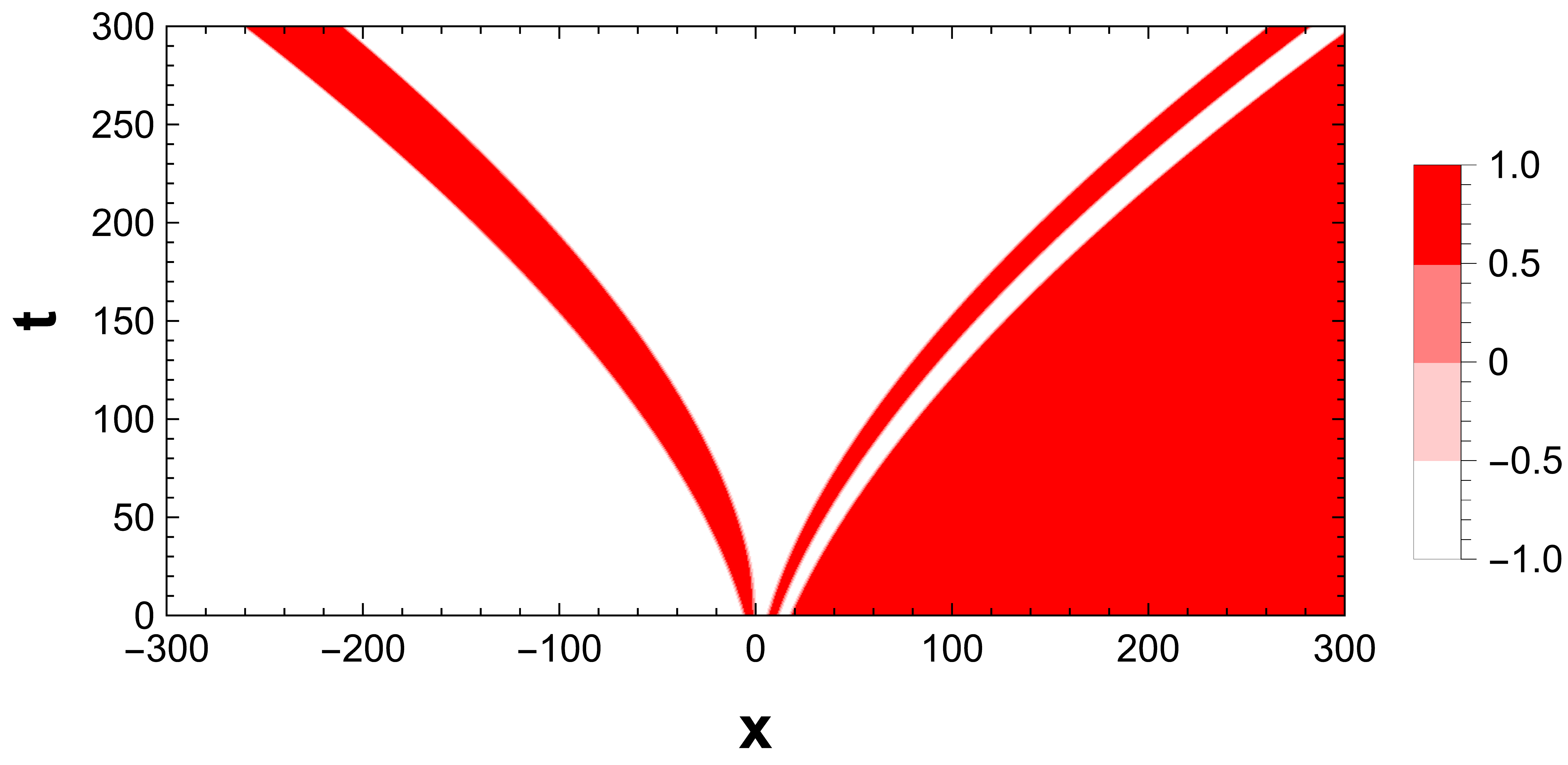}
\end{subfigure}
\caption{Simulations of \eqref{eq:mainEquation} with topographical inhomogeneity \eqref{eq:Ftopography} given by $-H_{\rm alg}(x; p)$ \eqref{eq:defHuni-alg}, $\varepsilon = 0.1$ and five-front initial condition $U(x,0) = \tanh(x + 6)  - \tanh(x + 0.5) + \tanh(x - 6) - \tanh(x - 11.5) + \tanh(x - 17.5)$; $p=0.4$ in (a), $p=0.0$ in (b) and $p=-0.4$ in (c): unlike in Fig. \ref{fig:5FrontDynamics}(c)/(d), the solitary one-front centered around $x=0$ is unstable. All simulations are done for $x \in (-500,500)$ (with homogeneous Neumannn boundary conditions) while only the middle $x \in (-300,300)$ part of the pattern is plotted.}
\label{fig:5FrontDynVaryingP}
\end{figure}
\\
{\bf Implications for ecosystem dynamics.}
The validity of the pulse interaction equations for algebraically `decaying' inhomogeneities with $p < 0$, {\it i.e.}, for topographies with inclinations that grow faster than linear, is a mathematical issue without obvious ecological relevance. However, our investigations into coarsening and pinning of multi-front dynamics within interaction equations \eqref{eq:NFrontODE} and/or its simplification \eqref{eq:dynNfronts} can be expected to be relevant in the context of ecosystems. The most simple ecosystem models are singularly perturbed two-component reaction-diffusion models (see \cite{bastiaansen2019dynamics,jaibi2020existence} and the references therein).
By the methods of \cite{bastiaansen2019dynamics,doelman2007nonlinear,van2010front}, equivalents of \eqref{eq:NFrontODE}/\eqref{eq:dynNfronts} can be derived (and validated) that govern the multi-pulse or multi-front interactions of vegetation patches in ecosystem models (although one should not underestimate the necessary technicalities).
Moreover, these systems will not be identical to \eqref{eq:NFrontODE}/\eqref{eq:dynNfronts}. However, as \eqref{eq:NFrontODE}/\eqref{eq:dynNfronts}, these systems will have an evident division into terms governed by the `external' topography and `internal effects' representing the pulse- or front-interactions that also take place in the homogeneous settings \cite{bastiaansen2019dynamics,van2010front}. In fact, the latter effects will not be exponentially small -- as in the case in \eqref{eq:NFrontODE}/\eqref{eq:dynNfronts} -- due to the semistrong nature of these interactions in singularly perturbed systems \cite{doelman2003semistrong}.
\\ \\
Thus, the semistrong two-component equivalents of \eqref{eq:NFrontODE}/\eqref{eq:dynNfronts} will provide insights into the coarsening and pinning dynamics of vegetation patches in response to the nature of the topography of the domain in which the ecosystem dynamics takes place.
Since desertification of dryland ecosystems can be seen as a coarsening process of vegetation patches in which vegetation is taken over step by step by bare soil states \cite{bastiaansen2020effect} and since pinned vegetation patches persist, understanding the mechanisms that distinguish between coarsening and pinning thus has a clear ecological relevance.
Further, realistic topographies are not necessarily of small amplitude and will vary on various spatial scales; nevertheless, the present work is a step towards building a fundamental understanding of these mechanisms.
\\
\bigskip

\noindent
{\bf Declaration of competing interests}

\noindent
The authors declare that they have no known competing financial interests or personal relationships that could have appeared to influence the work reported in this paper.

\medskip

\noindent
{\bf Acknowledgments}

\noindent
This work was supported by NWO's Mathematics of Planet Earth program. The research of AD is partly supported by the ERC Synergy project RESILIENCE (101071417). The research of TJK was supported in part by NSF DMS-1616064. AD and TJK thank Prof. Tommaso Ruggeri for inviting them to teach at the 2019 Summer School of Mathematical Physics in Ravello, Italy, where part of this research was carried out. The authors are deeply grateful to Margaret Beck, Peter van Heijster, Matt Holzer, Vivi Rottschafer, and Henk Schuttelars for organizing the conference Multiple-Scale Systems: Theory and Application, held at the Lorentz Center in Leiden University, July 8-12, 2024.

\medskip
\noindent{\bf Data availability}

\noindent
The numerical codes used to generate the data shown in the figures are available upon reasonable request from the corresponding author.

\medskip

\appendix

\section{Properties of the system at $\mathcal{O}(\varepsilon)$}
\label{sec:orderEpsSystem}
In the study of the stationary forced ($\varepsilon \neq 0$) system (see especially sections~\ref{sec:frontExistence}, \ref{sec:2frontExistence} and \ref{sec:geom-int}), it is necessary to understand the properties of solutions to the inhomogeneous problem
\begin{equation}
	\mathcal{L}_0 u_1 = - F(u_h,u_h',x),
	\label{eq:inhomogeneousProblem}
\end{equation}
with $\mathcal{L}_0$ as defined in \eqref{eq:defL0}. We focus here on the case $u_h(x;\phi) = u_\mathrm{up}(x;\phi) = + \tanh\left( \frac{\sqrt{2}}{2} \left[x-\phi\right]\right)$ and observe that the situation for $u_h = u_\mathrm{down}$ is equivalent upon replacement of $F(u_\mathrm{up}, u'_\mathrm{up}, x)$ by $F(-u_\mathrm{up},-u'_\mathrm{up},x)$. (We refer to Lemma 2.5 of \cite{doelman2022slow} for comparable more general statements and results.) Without loss of generality, we choose $u_1$ such that $u_1(\phi;\phi)=0$. Foremost, the homogeneous equation $\mathcal{L}_0\Psi = 0$ has a bounded solution $\Psi_b$ and an unbounded solution $\Psi_u$, which are given by
\begin{equation}
\begin{array}{lcl}
\Psi_b(x;\phi) & = & \partial_x u_\mathrm{up}(x;\phi) = \frac{\sqrt{2}}{2} \sech^2\left( \frac{\sqrt{2}}{2} \left[x-\phi\right]\right)
\\
\Psi_u(x;\phi) & = & v(x;\phi) \Psi_b(x;\phi),
\end{array}
\label{eq:appendixBoundedUnboundedSolution}
\end{equation}
where $v$ satisfies
\[
	\partial_x v(x;\phi) =\Psi_b(x;\phi)^{-2} = 2 \cosh^4\left( \frac{\sqrt{2}}{2} \left[x-\phi\right]\right).
\]
By direct integration, $v$ can thus be defined as
\begin{align}
	v(x;\phi)
& = \int_\phi^x \Psi_b(z;\phi)^{-2}\ dz  = \frac{3}{4}\left[x-\phi\right] + \frac{\sqrt{2}}{2} \sinh\left(\sqrt{2}\left[x-\phi\right]\right) + \frac{\sqrt{2}}{16}\sinh\left(2\sqrt{2}\left[x-\phi\right]\right).
\end{align}
The leading order terms for the limit behavior of the functions $v$, $\Psi_b$, and $\Psi_u$ can then be found, confirming the assertion that $\Psi_b$ is bounded and $\Psi_u$ is unbounded in both limits $x \rightarrow \infty$ and $x \rightarrow -\infty$:
\begin{equation}
\begin{array}{cclcc}
v(x;\phi) & \rightarrow & \pm \frac{\sqrt{2}}{32} e^{\pm 2 \sqrt{2} \left[x-\phi\right]} + \; {\rm h.o.t.} & {\rm as} & x-\phi \rightarrow \pm \infty;\\
\Psi_b(x;\phi) & \rightarrow & 2\sqrt{2} e^{\mp \sqrt{2}\left[x-\phi\right]} + \; {\rm h.o.t.} & {\rm as} & x-\phi \rightarrow \pm \infty; \\
\Psi_u(x;\phi) & \rightarrow & \pm \frac{1}{8} e^{\pm \sqrt{2} \left[x-\phi\right]} + \; {\rm h.o.t.} & {\rm as} & x-\phi \rightarrow \pm \infty.
\end{array}
\label{eq:limitbehaviorPsiub}
\end{equation}
The inhomogeneous problem~\eqref{eq:inhomogeneousProblem} can be solved using variation of parameters, which yields the solutions
\[
	u_1(x;\phi) = A(x;\phi) \Psi_b(x;\phi) + B(x;\phi) \Psi_u(x;\phi),
\]
with
\[
\begin{array}{lcl}
	A(x;\phi) & = & \int_\phi^x F(u_h(z;\phi),u_{h,x}(z;\phi),z)\ \Psi_u(z;\phi)\ dz, \\
	B(x;\phi) & = & B_0 - \int_\phi^x F(u_{h}(z;\phi), u_{h,x}(z;\phi),z)\ \Psi_b(z;\phi)\ dz.
\end{array}
\]
Here, $B_0$ is a constant that is determined by the boundary conditions (note $A$ does not carry a constant $A_0$ as $\phi$ has been chosen such that $A(\phi;\phi) = 0$ and thus such that $u_1(\phi;\phi) = 0$). Using the limiting behavior of $\Psi_b$ and $\Psi_u$ \eqref{eq:limitbehaviorPsiub}, we find that the limiting behavior of $A$ and $B$ is as follows,
\[
\begin{array}{rlcl}
	\left| A(x;\phi) \right|  &\leq \left\|F(u_h(\cdot;\phi),u_{h,x}(\cdot;\phi),\cdot)\right\|_\infty\ \frac{\sqrt{2}}{16} e^{\pm \sqrt{2}\left[x-\phi\right]} + \; {\rm  h.o.t.} & {\rm as} &  x-\phi \rightarrow \pm\infty;\\
	\left| B(x;\phi) - B_{\pm}(\phi) \right|  &\leq \left\|F(u_h(\cdot;\phi),u_{h,x}(\cdot;\phi),\cdot)\right\|_\infty\ 2 e^{\mp \sqrt{2}\left[x-\phi\right]} + \; {\rm h.o.t.} & {\rm as} &  x-\phi \rightarrow \pm \infty;
\end{array}
\]
where
\begin{equation}
\label{eq:BpmLimits}
B_{\pm}(\phi) = B_0 - \int_\phi^{\pm \infty} F(u_h(z;\phi), u_{h,x}(z;\phi),z)\ \Psi_b(z;\phi)\ dz.
\end{equation}
Therefore, the limiting behavior of $u_1$ is given by
\begin{equation}
u_1(x;\phi) = \pm \frac{1}{8} B_{\pm} e^{\pm \sqrt{2}[x-\phi]} + \mathcal{O}(1) \; \;  {\rm as} \; \;  x-\phi \rightarrow \pm \infty.
\label{eq:u1limitbehavior}
\end{equation}
Thus, boundedness of $u_1$ for $x \rightarrow -\infty$ (respectively $x \rightarrow +\infty$) is guaranteed when $B_-(\phi) = 0$ (respectively $B_+(\phi) = 0$). That is, $B_0$ needs to be chosen as
\[
B_0 = - \int_{-\infty}^\phi F(u_h(z;\phi), u_{h,x}(z;\phi),z) \Psi_b(z;\phi)\ dz, \; \; {\rm resp.} \; \;
B_0 = \int_{\phi}^\infty F(u_h(z;\phi), u_{h,x}(z;\phi).z) \Psi_b(z;\phi)\ dz
\]
to ensure boundedness as $x \rightarrow -\infty$, resp. $x \rightarrow \infty$. Thus, if we consider an $u_1(x,\phi)$ that remains bounded as $x-\phi \to -\infty$, it follows that its limit behavior for $x-\phi \to \infty$ is indeed given by
\eqref{eq:limxinfty-uu1} (and thus that boundedness in both limits only holds if the solvability condition $\mathcal{R}(\phi) = 0$ \eqref{eq:FredholmCondition1Front} holds).

\section{The proofs of Theorem \ref{th:Nfrontslocexp} and Corollary \ref{cor:Nfrontslocalg}}
\label{ap:ProofTh}
A proof of this theorem and corollary cannot be based a priori on an analysis of the $N$-front interaction ODE \eqref{eq:dynNfronts} -- as was the case in the previous sections for $N = 2,3,4$ -- without a rigorous validation of this reduction (see the brief discussion in the introduction of Sec.~\ref{sec:interactiondynamics}).
On the other hand, the existence and stability of the stationary points of \eqref{eq:dynNfronts} can be established by the approach of Secs.~\ref{sec:1fronts} and \ref{sec:2fronts}.
In fact, there is -- by construction (and as we saw) -- a direct correspondence between the construction and spectral properties of the multi-front patterns in sections \ref{sec:1fronts} and \ref{sec:2fronts} and the critical points of the multi-front interaction ODEs. Therefore, we do consider the stationary $N$-front patterns as critical points of the interaction ODEs in the upcoming proof -- by their direct link to the stationary $N$-front patterns, and thus not necessarily by the validity of these ODEs and the reduction they are based on.
As was already noted in the final paragraph of Sec.~\ref{sec:ODEderivation}, we can also establish the stability of the stationary $N$-front patterns as solutions of  PDE \eqref{eq:mainEquation} by only considering the $N$ eigenvalues of the associated critical points in \eqref{eq:dynNfronts}: apart from the asymptotically small eigenvalues represented by these eigenvalues, the spectrum of a stationary $N$-front pattern in \eqref{eq:mainEquation} is contained in the stable complex half-plane, bounded away from the imaginary axis (see the brief discussion in Sec.~\ref{sec:oneFrontStability}).
\\ \\
Thus, using either of the methods of Sec.~\ref{sec:2fronts} or \ref{sec:interactiondynamics}, we find that  stationary $N$-front patterns correspond to  critical points of the system of $N$ ODEs
\eqref{eq:dynNfronts}
governing the front interactions in $N$-front solutions.
\\ \\
{\bf The proof of Theorem \ref{th:Nfrontslocexp}.} We first focus on a special kind of $N$-front patterns, namely those that are extensions of the two-front patterns of the second kind -- thus we assume that $h_+ < 0$ -- that satisfy
\[
\bar{\psi}_1 = \psi_\ast + \varepsilon^{\nu_1} \tilde{\psi}_1, \;
\bar{\psi}_2 = \nu_2 |\log \varepsilon| + \bar{\ell}_2, \; ...,
\bar{\psi}_N = \nu_N |\log \varepsilon| + \bar{\ell}_N \; \; {\rm with} \; \; \nu_2 < ... < \nu_N,
\]
so that $\bar{\psi}_1 = \mathcal{O}(1) \ll \bar{\psi}_2 \ll ... \ll \bar{\psi}_N$. In other words, we focus on $N$-front patterns that have one front near a zero of $\mathcal{S}(\psi)$ and all other fronts at one side of this zero. The system that determines this $N$-front as critical points of \eqref{eq:dynNfronts} can now be solved inductively, with the two-front of the second kind as starting point. We first consider $N=3$, conclude that $1 \ll \bar{\psi}_2 - \bar{\psi}_1 \ll \bar{\psi}_3 - \bar{\psi}_2$ and that we thus may copy the (leading order) result from Sec.~\ref{sec:ODE2frontloc}: $\nu_1 = \mu/(1-\mu)$, $\nu_2 = 1/(1-\mu)$; moreover, $\bar{\ell}_1$ and $\bar{\ell}_2$ can now be determined -- better: approximated -- uniquely. Then, $\nu_3$ and $\bar{\ell}_3$ are (to leading order) determined by
\[
\varepsilon^{1 + \mu \nu_3} |h_+| w_+ e^{-\mu \bar{\ell}_3} = 16 \varepsilon^{\nu_3 - \nu_2} e^{\bar{\ell}_2 - \bar{\ell}_3}.
\]
Thus, it follows that $\nu_3 = \frac{1}{1-\mu}(1+\nu_2) = \frac{2-\mu}{(1-\mu)^2}$ and that $\bar{\ell}_3$ can be determined uniquely. This process can be repeated for the next $4-,..., N-$front patterns: the pair $(\nu_{j+1},\bar{\ell}_{j+1})$ is determined by
\begin{equation}
\label{eq:iterationnextfront}
\varepsilon^{1 + \mu \nu_{j+1}} |h_+| w_+ e^{-\mu \bar{\ell}_{j+1}} = 16 \varepsilon^{\nu_{j+1} - \nu_j} e^{\bar{\ell}_{j} - \bar{\ell}_{j+1}}.
\end{equation}
so that $\nu_{j+1} = \frac{1}{1-\mu}(1+\nu_j)$ and also $\bar{\ell}_{j+1}$ can be explicitly (and uniquely) determined in terms of $\bar{\ell}_{j}$. Thus, we find
\begin{equation}
\label{eq:nuNcrit}
\nu_N(\mu) = \frac{1}{\mu} \left[ \frac{1}{(1-\mu)^{N-1}} - 1 \right]
\end{equation}
Note that $\lim_{\mu \downarrow 0} \nu_N(\mu) = N-1$ and that $\nu_{j}(\mu) - \nu_{j-1}(\mu) = 1/(1-\mu)^{j-1}$, which implies that the distance between successive fronts increases with $j$ (for $\mu \in (0,1)$). To be certain that this leading order analysis is well-established, we again need to impose that the construction that determines the $j$th front takes place at a magnitude in $\varepsilon$ for which the correction/error terms -- that are {\it prima facie} of $\mathcal{O}(\varepsilon^2)$ and have not been taken into account -- can indeed not give a leading order contribution. It follows from \eqref{eq:iterationnextfront} that we need to impose that $1 + \mu \nu_{j+1} = \nu_{j+1} - \nu_j < 2$ for all $j = 1, ..., N-1$. By \eqref{eq:nuNcrit} we conclude that this is the case if
\begin{equation}
\label{eq:defmuastN}
0 < \mu < \mu_{\ast,1}(N) := 1 - 2^{-\frac{1}{N-1}} \; \; (N \geq 2),
\end{equation}
where the subscript $1$ is an index indicating that we have labeled this N-front solution as the first of in total $\mathcal{N}(N)$ N-front solutions. Thus, to establish the existence of the present type of $N$-front patterns, we need to impose that $0 < \mu < \mu_{\ast,1}(N) < 1$. Clearly, $\mu_{\ast,1}(N)$ is a decreasing function of $N$, for $N \gg1$ $\mu_{\ast,1}(N) = \log 2/(N-1)$  (to leading order) so that $\mu_\ast(N)  \downarrow 0$ as $N \to \infty$: the width of the allowed $\mu$-region decays monotonously to 0 as $N \to \infty$. Nevertheless, we may conclude that the two-front patterns of the second kind can indeed be extended to these $N$-front patterns for any given $N > 2$ as long as $\mu$ is sufficiently close to $0$, {\it i.e.} if the (exponential) decay of the topography $H(x)$ is sufficiently weak (and if $h_+ < 0$). Moreover, we may immediately conclude that these patterns are unstable, since the largest two of the $N$ asymptotically small eigenvalues associated to their spectral stability are (to leading order) given by those of the two-front, {\it i.e.}, by \eqref{eq:eigenv2front2nd}.
\\ \\
In analogy to the situation in Sec.~\ref{sec:ODE2frontloc} (and under the assumption that $h_- > 0$), extending the third kind of critical points in a similar fashion to $N$-front patterns that have one front near a zero of $\mathcal{S}(\psi)$ and all other fronts to the left of this zero ({\it i.e.}, $\bar{\psi}_{1} \ll \bar{\psi}_{2}  \ll ... \bar{\psi}_{N-1} \ll -1$ and $\bar{\psi}_N = \mathcal{O}(1)$) is an identical procedure that leads to the same critical value of $\mu$ as in \eqref{eq:defmuastN}, {\it i.e.},  $\mu_{\ast,2}(N) =  \mu_{\ast,1}(N)$ (and to the same (in)stability result). All other $\mathcal{N}(N) - 2$ $N$-front patterns can be constructed in a similar inductive fashion, with either a two-front of the first, second or third type or the three-front of the fourth type as starting point. For each of these $N$-fronts there will be a critical value $\mu_{\ast,k}(N) > 0$ of $\mu$ so that their existence can be established in the $\mu$-interval $(0,\mu_{\ast,k}(N))$ ({\it i.e.} for these values $\mu$ it can be shown that the $\mathcal{O}(\varepsilon^2)$ corrections to \eqref{eq:defpsiStau} do not (potentially) become of leading order).
\\ \\
Also, the critical value $\mu_\ast(N)$ introduced in the statement of Theorem \ref{th:Nfrontslocexp} is the minimal value of the various $\mu_{\ast,k}(N)$'s ($k=1,2,...,\mathcal{N}(N)$). The fact that the above considered $N$-front patterns have $N-1$ fronts at one side of the first (or last) front at $\psi = \mathcal{O}(1)$ implies that its last (or first) front reaches deeper into the tail of $\mathcal{S}(\psi)$ than that of any of the other $N$-front patterns. In other words, the above considered patterns have fronts at positions for which $\mathcal{S}(\psi)$ is smaller than for any of the other $N$-front patterns. Thus, the necessary assumption that the next order corrections to \eqref{eq:defpsiStau} cannot become of leading order -- {\it i.e.} the condition $\varepsilon |\mathcal{S}(\psi)| \ll \mathcal{O}(\varepsilon^2)$ -- is stronger for these $N$-front patterns than it is for any of the other $N$-front patterns. In other words, $\mu_{\ast,k}(N) \geq \mu_{\ast,1}(N) = \mu_{\ast,2}(N) = \mu_{\ast}(N)$ as given by \eqref{eq:defmuastN} -- which also appears in the statement of the Theorem.
\\ \\
Since the (in)stability of all $N$-fronts patterns is determined by the instability of the underlying patterns of the first, second, third or fourth kind, the statement of Theorem \ref{th:Nfrontslocexp} follows. \hfill $\Box$
\\ \\
{\bf The proof of Corollary \ref{cor:Nfrontslocalg}.} As in the above proof of Theorem \ref{th:Nfrontslocexp}, we again first consider the construction of the $N$-front patterns of the second kind that have $N-1$ fronts to the right of the first front at $\psi = \mathcal{O}(1)$.
\\ \\
For $N=2$, {\it i.e.}, for the two-front of the second kind, we replace \eqref{eq:Ansat2critpts2} by
\begin{equation}
\bar{\psi}_1 = \psi_\ast + \tilde{\psi}_1 \Theta_1(\varepsilon) \; \;
\bar{\psi}_2 = \frac{\tilde{\nu}_2}{\Theta_2(\varepsilon)},
\label{eq:Ansat2critpts2al}
\end{equation}
with $0< \Theta_{1,2}(\varepsilon) \ll 1$ two (asymptotically small) order functions of a priori unknown magnitude in $\varepsilon$ and $\tilde{\psi}_1 = \tilde{\psi}_1(\varepsilon) = \mathcal{O}(1)$, $\tilde{\nu}_2 = \tilde{\nu}_2(\varepsilon) = \mathcal{O}(1)$ two constants to be determined (see \eqref{eq:Ansat2critpts2}). Assuming that $H'(x)$ decays as $1/|x|^p$ \eqref{eq:algdecayHx} it follows by Lemma \ref{lem:algdecayH} (and \eqref{eq:defpsiStau}) that also $\mathcal{S}(\psi)$ decays as $1/|\psi|^p$ (for $\psi| \gg 1$). Thus, under the condition that $\tilde{h}_\pm(p) < 0$ \eqref{eq:algdecayHx}, we find by the second equation of \eqref{eq:intODE2loc} that
\begin{equation}
\label{eq:Theta2}
\varepsilon (\Theta_2(\varepsilon))^p e^{\frac{\tilde{\nu}_2}{\Theta_2(\varepsilon)}}  = C_2 = C_2 (\tilde{\nu}_2),
\end{equation}
for some $\mathcal{O}(1)$ constant $C_2 > 0$. Hence,
\begin{equation}
\label{eq:Theta2-2}
\frac{\tilde{\nu}_2}{\Theta_2(\varepsilon)} = |\log \varepsilon| + p \, |\log \Theta_2(\varepsilon)| + \mathcal{O}(1),
\end{equation}
so that we can determine the order function $0 < \Theta_2(\varepsilon) \ll 1$ and constant $\tilde{\nu}_2(\varepsilon)$ by
\begin{equation}
\label{sol:Theta2}
\Theta_2(\varepsilon) = \frac{1}{|\log \varepsilon| + p \log |\log \varepsilon|} = \frac{1}{|\log \varepsilon|} - p \, \frac{\log |\log \varepsilon|}{|\log \varepsilon|^2}, \; \; \tilde{\nu}_2(\varepsilon) = 1
\end{equation}
to leading order in $\varepsilon$. Similarly, we find from the first equation of \eqref{eq:intODE2loc} that
\[
\varepsilon \Theta_1(\varepsilon) e^{\frac{\tilde{\nu}_2}{\Theta_2(\varepsilon)}}  = C_1 = C_1 (\tilde{\psi}_1) = \mathcal{O}(1).
\]
This yields, by \eqref{eq:Theta2} and \eqref{sol:Theta2}, that
\begin{equation}
\label{sol:Theta1}
\Theta_1(\varepsilon) =
\left(\frac{1}{|\log \varepsilon|} - p \, \frac{\log |\log \varepsilon|}{|\log \varepsilon|^2} \right)^p =
\frac{1}{|\log \varepsilon|^p} - p^2 \, \frac{\log |\log \varepsilon|}{|\log \varepsilon|^{(p+1)}}
\end{equation}
to leading order in $\varepsilon$.
Note that the findings of \eqref{sol:Theta2} and \eqref{sol:Theta1} are natural extensions of the case of a vanishing exponential decay in Theorem \ref{th:Nfrontslocexp}. There, $\nu_2(\mu) = 1/(1-\mu) \to 1$ as $\mu \downarrow 0$ \eqref{eq:Ansat2critpts2}, \eqref{sol:nu12}, which corresponds to $\tilde{\nu}_2(\varepsilon)/\Theta_2(\varepsilon) = |\log \varepsilon|$ \eqref{eq:Ansat2critpts2al} ({\it i.e.}, to leading order, $\tilde{\nu}_2 = \nu_2(0) = 1$ and $1/\Theta_2(\varepsilon) = |\log \varepsilon|$). Moreover, $\nu_1(\mu) = \mu/(1-\mu) \to 0$ as $\mu \downarrow 0$, which would imply that $\bar{\psi}_1$ is no longer asymptotically close to a zero $\psi_\ast$ of $\mathcal{S}(\psi)$ \eqref{eq:Ansat2critpts2}. In the algebraic case this is replaced by
$\bar{\psi}_1$ being $\mathcal{O}(1/|\log \varepsilon|^p)$ close to $\psi_\ast$.
\\ \\
For $N$-front patterns that extend this two-front solution of the second kind to patterns that have one front near a zero of $\mathcal{S}(\psi)$ and all other fronts to the right of this zero, we find in a similar way (to leading order),
\[
\frac{\tilde{\nu}_j}{\Theta_j(\varepsilon)} = (j-1) \left(|\log \varepsilon| + p \, \log |\log \varepsilon|\right), \; j = 2, ..., N
\]
(see \eqref{eq:Theta2-2}). Thus, we can choose all order functions $\Theta_j(\varepsilon)$, $j=3,..,N$, to be equal to $\Theta_2(\varepsilon)$ as given by \eqref{sol:Theta2}. The position of the $j$th front is then determined by (the pre-factor) $\tilde{\nu}_j(\varepsilon) = j-1 $ to leading order, which again agrees with the $\mu \downarrow 0$ limit in the above case of exponential decay: $\nu_j(0) = \tilde{\nu}_j = j-1$ and $1/\Theta_j(\varepsilon) = |\log \varepsilon|$, $j=2,...,N$.
\\ \\
This implies that in the algebraic decay case, the non-determined $\mathcal{O}(\varepsilon^2)$ higher order terms of front interaction system \eqref{eq:dynNfronts} remain of higher order for all stationary $N$-front patterns, since (to leading order)
\[
\varepsilon \mathcal{S}(\bar{\psi}_N) = \varepsilon \frac{\tilde{C}_N}{(\Theta_N(\varepsilon))^p} = \varepsilon \, \tilde{C}_N |\log \varepsilon|^p \gg \varepsilon^2
\]
for some $\mathcal{O}(1)$ constant $\tilde{C}_N=\tilde{C}_N(\varepsilon) > 0$ (where we have again used that $N$-front patterns that have one front near a zero of $\mathcal{S}(\psi)$ and all other fronts to the right (or left) of this zero extend deeper into the tail of $\mathcal{S}(\psi)$ than any of the other $N$-front patterns).
\\ \\
As in the case of exponential decay, the (in)stability of the stationary $N$-front patterns builds on that of the underlying two-front patterns of the first, second and third kind and the three-front pattern of the fourth kind. However, one needs to be careful since the leading order results on the eigenvalues of the underlying basic patterns in the case of exponential decay cannot be automatically extended to the limit $\mu = 0$.
\\ \\
In the case of the two-front patterns of the second kind, the leading order $\mathcal{O}(\varepsilon^{-\mu/(1-\mu)})$ expression in the matrix \eqref{eq:stab2loc2} associated to the stability of the two-front pattern of the second kind merges with the next $\mathcal{O}(1)$ terms in the limit $\mu \downarrow 0$. Thus, one may not take the $\mu \downarrow 0$ limit of the eigenvalues $\lambda_{1,2}(\mu)$ \eqref{eq:lamba21loc} since these are obtained by leading order arguments. Nevertheless, it is straightforward to obtain the counterpart of matrix \eqref{eq:stab2loc2} in the case of algebraic decay,
\begin{equation}
\label{eq:stab2loc2-alg}
E_1
\left(
\begin{array}{cc}
1 - \varepsilon \frac{\mathcal{S}'(\psi_\ast)}{E_1} & - 1
\\
- 1 & 1 - \frac{p}{\bar{\psi}_2}
\end{array}
\right),
\end{equation}
with
\[
E_1 = 16 e^{-(\bar{\psi}_2-\bar{\psi}_1)} = -\varepsilon \mathcal{S}(\bar{\psi}_2) = \mathcal{O}(\varepsilon (\Theta_2(\varepsilon))^p) =
\mathcal{O}\left(\frac{\varepsilon}{|\log \varepsilon|^p}\right).
\]
Thus, similar to the case of exponential decay, the $\varepsilon \mathcal{S}'(\psi_\ast)/E_1$-term in matrix \eqref{eq:stab2loc2-alg} is dominant and we may conclude
\[
\lambda_1 = - \varepsilon \mathcal{S}'(\psi_\ast), \; \;
\lambda_2 = E_1 > 0,
\]
to leading order (see \eqref{eq:eigenv2front2nd}): the two-front of the second kind represented by $(\bar{\psi}_1,\bar{\psi}_2)$ is unstable.
\\ \\
The analysis of the two-front patterns of the third kind is in essence identical to the above analysis.
In the case of two-front patterns of the first kind, the matrix given in \eqref{eq:stab2loc1} also degenerates, in the sense that $\lambda_1(\mu) \uparrow 0$ as $\mu \downarrow 0$ \eqref{eq:lamba21loc}. Although the unstable eigenvalue $\lambda_1(\mu) \downarrow 2 > 0$ -- and one is thus inclined to conclude that the instability persists into the case of algebraic decay -- we need to consider also this case somewhat more carefully. We find that the algebraic counterpart of matrix \eqref{eq:stab2loc1} is very similar to \eqref{eq:stab2loc1},   namely
\begin{equation}
\label{eq:stab2loc1-alg}
E_1
\left(
\begin{array}{cc}
1 - \frac{p}{|\bar{\psi}_1|} & - 1
\\
- 1 &  1 - \frac{p}{\bar{\psi}_2}
\end{array}
\right)
\end{equation}
with $1/|\bar{\psi}_{1,2}| \ll 1$ (in fact, it can be checked that  $\bar{\psi}_{1} = -1/2|\log \varepsilon|$ and $\bar{\psi}_{2} = 1/2|\log \varepsilon|$ to leading order, which again agrees with the $\nu \downarrow 0$ limit of the case of exponential decay). It follows that
\[
\lambda_1(p) = -\frac{p}{2} \left(\frac{1}{|\bar{\psi}_1|} + \frac{1}{\bar{\psi}_2} \right) < 0,
\; \;
\lambda_2(p) = 2 - \frac{p}{2} \left(\frac{1}{|\bar{\psi}_1|} + \frac{1}{\bar{\psi}_2} \right) > 0,
\]
(to leading order). Thus, also in the case of algebraic decay, the two-front pattern of the first kind is unstable.
\\ \\
A similar analysis shows that also the three-front patterns of the fourth kind remain unstable in the case of algebraic decay. Finally, it follows along exactly the same lines as in the exponential case that the eigenvalues of the underlying two-front patterns of the first, second and third kind and those of the three-front patterns of the fourth kind persist as leading order eigenvalues for all $N$-front patterns that are based on it. \hfill $\Box$

\bibliographystyle{plain}
\bibliography{ForcedAllenCahn}

\end{document}